\documentclass[10pt]{article} 
\usepackage{amsthm,amsfonts,amsmath,amscd,amssymb,latexsym,epic,eepic,epsfig,psfrag}

\usepackage[all]{xy}

\author{Dietmar Salamon\\
ETH-Z\"urich
\and
Katrin Wehrheim\thanks
{The second author gratefully acknowledges 
support by the swiss and US national science foundations
and thanks IAS Princeton and FIM Z\"urich for their hospitality, 
during which most of this work was undertaken.}
\\
MIT}

\date{12 August 2007}

\title{Instanton Floer homology with Lagrangian boundary conditions}

\newcommand{\one}{{{\mathchoice \mathrm{ 1\mskip-4mu l} \mathrm{ 1\mskip-4mu l}
\mathrm{ 1\mskip-4.5mu l} \mathrm{ 1\mskip-5mu l}}}}

\newcommand{\dslash}{/\mskip-6mu/}

\def\rd{{\rm d}}
\def\rT{{\rm T}}

\def\rG{{\rm G}}
\def\phi{\varphi}
\def\eps{\varepsilon} 
\def\om{\omega} 

\def\Om{\Omega} 
\def\cg{\mathfrak{g}}

\def\co{\mathfrak{o}}
\def\cA{{\mathcal A}}
\def\cB{{\mathcal B}} 
\def\cC{{\mathcal C}} 
\def\cD{{\mathcal D}} 
\def\cE{{\mathcal E}} 
\def\cF{{\mathcal F}} 
\def\cG{{\mathcal G}} 
\def\cH{{\mathcal H}}

\def\cK{{\mathcal K}}
\def\cL{{\mathcal L}} 
\def\cM{{\mathcal M}}

\def\cR{{\mathcal R}} 
\def\cS{{\mathcal S}}

\def\cU{{\mathcal U}} 
\def\cV{{\mathcal V}}
\def\cW{{\mathcal W}} 
\def\cX{{\mathcal X}}

\def\CS{{\mathcal C}{\mathcal S}}

\def\A{{\mathbb A}}
\def\B{{\mathbb B}}
\def\R{{\mathbb R}} 
\def\Q{{\mathbb Q}} \def\T{{\mathbb T}} 
\def\N{{\mathbb N}} 
\def\C{{\mathbb C}} 
\def\D{{\mathbb D}} 
\def\bF{{\mathbb F}} 
 
\def\Z{{\mathbb Z}}

\def\im{{\rm im}\,}
\def\laplace{\Delta} 
\def\st{\: \big| \:} 
\DeclareMathOperator{\supp}{supp}

\def\dt{{\rm d}t}
\def\ds{{\rm d}s}

\def\dr{{\rm d}r}
\def\pd{\partial}
\def\comp{\circ}

\def\tu{{\tilde{u}}}
\def\tv{{\tilde{v}}}

\def\tA{{\tilde{A}}}
\def\tB{{\tilde{B}}}

\def\tY{{\widetilde{Y}}}

\def\la{\langle\,}
\def\ra{\,\rangle}
\def\dvol{\,\rd{\rm vol}}
\def\CF{{\mathrm{CF}}}
\def\HF{{\mathrm{HF}}}
\def\SU{{\mathrm{SU}}}
\def\su{{\mathfrak s\mathfrak u}}
\def\SO{{\mathrm{SO}}}
\def\Hom{{\mathrm{Hom}}}

\def\tsum{\textstyle\sum}
\def\tint{\textstyle\int}
\def\hol{{\mathrm{hol}}}
\def\Hol{{\mathrm{Hol}}}

\def\Id{{\mathrm{Id}}}
\def\trace{{\mathrm{tr}}}
\def\Crit{{\mathrm{Crit}}}
\def\Met{\mathrm{Met}}
\def\Or{\mathrm{Or}}
\def\spec{\mathrm{spec}}
\def\dom{\mathrm{dom}\,}
\def\flat{{\mathrm{flat}}}
\def\loc{{\mathrm{loc}}}

\def\Abs#1{\left|#1\right|}

\def\Norm#1{\left\|#1\right\|}
\def\Nabla#1{\nabla\kern-.5ex{}_{#1}}
\def\NORM#1{\left|\kern-.2ex{}\|#1|\kern-.2ex{}\right\|}
\newcommand{\p}{\partial}  
\newcommand{\Cinf}{C^\infty}  
\newcommand{\inner}[2]{\langle\, #1, #2\,\rangle}   
\newcommand{\Inner}[2]{\left\langle\, #1, #2\,\right\rangle}  
\newcommand{\winner}[2]{\langle\, #1{\wedge}#2\,\rangle}   
   
\newcommand{\Vect}{\mathrm{Vect}}

\newcommand{\mapdown}[1]{\Big\downarrow
\rlap{$\vcenter{\hbox{$\scriptstyle#1$}}$}}
\newcommand{\mapright}[1]{\smash{
\mathop{\longrightarrow}\limits^{#1}}}

\newtheorem{dfn}{Definition}[section] 
\newtheorem{lem}[dfn]{Lemma} 
\newtheorem{lemma}[dfn]{Lemma} 
\newtheorem{prp}[dfn]{Proposition} 
\newtheorem{proposition}[dfn]{Proposition} 
\newtheorem{thm}[dfn]{Theorem} 
\newtheorem{rmk}[dfn]{Remark} 
\newtheorem{remark}[dfn]{Remark} 
\newtheorem{cor}[dfn]{Corollary} 
\newtheorem{ex}[dfn]{Example}

\begin{document}

\bibliographystyle{plain}

\maketitle



\section*{Contents}

1. Introduction  \hfill \pageref{sec:intro} \linebreak
2. The Chern--Simons functional \hfill \pageref{sec:CS}\linebreak
3. The Hessian \hfill \pageref{sec:HA} \linebreak
4. Operators on the product $S^1\times Y$\hfill\pageref{sec:S1Y} \linebreak
5. Exponential decay \hfill \pageref{sec:exp}\linebreak
6. Moduli spaces and Fredholm theory \hfill \pageref{sec:fredholm}\linebreak
7. Compactness \hfill \pageref{sec:compact}\linebreak
8. Transversality \hfill \pageref{sec:moduli}\linebreak
9. Gluing \hfill \pageref{sec:gluing}\linebreak
10. Coherent orientations \hfill \pageref{sec:orient}\linebreak
11. Floer homology \hfill \pageref{sec:floer}\linebreak
A. Spectral flow \hfill \pageref{app:spec}\linebreak
B. The Gelfand--Robbin quotient\hfill \pageref{app:GR}\linebreak
C. Unique continuation \hfill \pageref{app:uc}\linebreak
D. Holonomy perturbations \hfill \pageref{app:Xf}\linebreak
E. The Lagrangian and its tangent bundle \hfill \pageref{app:Lag}\linebreak


\section{Introduction}\label{sec:intro}

In this paper we define instanton Floer homology groups for a 
pair consisting of a compact oriented $3$-manifold with boundary 
and a Lagrangian submanifold of the moduli space of flat 
$\SU(2)$-connections over the boundary. 
We carry out the construction for a general class of 
{\it irreducible, monotone} boundary conditions.
The main examples of such Lagrangian submanifolds are induced 
from a disjoint union of handle bodies such that the union of 
the $3$-manifold and the handle bodies is an integral homology $3$-sphere.
The motivation for introducing these invariants arises from our program for a 
proof of the Atiyah-Floer conjecture for Heegaard splittings~\cite{A,Sa}. 
We expect that our Floer homology groups are isomorphic to the usual 
Floer homology groups~\cite{F1,Donaldson book} of the closed $3$-manifold
in our main example and thus can be used as a starting point 
for an adiabatic limit argument as in~\cite{DS}. 
On the level of Euler characteristics, the Atiyah-Floer
conjecture was proven by Taubes~\cite{T}.
 
Floer homology groups for $3$-manifolds with boundary were first 
constructed by Fukaya~\cite{Fu} with a different method. 
His setup uses nontrivial $\SO(3)$-bundles 
and thus cannot immediately be used for the proof 
of the Atiyah-Floer conjecture where the bundles 
are necessarily trivial.  Our approach is motivated 
by the construction of a Chern-Simons functional 
on $3$-manifolds with boundary.

Let $Y$ be a compact oriented $3$-manifold with boundary 
and denote
$$
\Sigma:=\p Y,\qquad
\rG:=\SU(2),\qquad 
\cg: = \su(2),\qquad
\inner{\xi}{\eta}:=-\trace(\xi\eta)
$$
for $\xi,\eta\in\cg$.  
While many of the results in this paper carry over to general 
compact Lie groups (and nontrivial bundles),
our construction of Floer homology works in this form 
only for $\rG=\SU(2)$ (where the bundles are necessarily trivial).
The whole story also carries over to nontrivial $\SO(3)$-bundles, 
where the moduli spaces of flat connections are nonsingular and 
monotone, however, in this paper we restrict to the case $\rG=\SU(2)$.

The space $\cA(\Sigma):=\Om^1(\Sigma,\cg)$ 
of connections on $\Sigma$ carries a natural 
symplectic form
\begin{equation}\label{eq:symp}
\omega(\alpha,\beta) 
:= \int_\Sigma \winner{\alpha}{\beta}
\end{equation}
for $\alpha,\beta\in\rT_A\cA(\Sigma)=\Om^1(\Sigma,\cg)$,
the action of the gauge group $\cG(\Sigma):=\cC^\infty(\Sigma,\rG)$
on $\cA(\Sigma)$ is Hamiltonian, and the moment map is the curvature
(see~\cite{AB}). The (singular) symplectic quotient is
the moduli space 
$$
M_\Sigma:= \cA_{\rm flat}(\Sigma)/\cG(\Sigma) 
= \cA(\Sigma)\dslash\cG(\Sigma)
$$
of flat connections.
We assume throughout that $\cL\subset\cA(\Sigma)$ is a gauge invariant, 
monotone, irreducible Lagrangian submanifold in the following sense.

\smallskip\noindent{\bf (L1)}
$\cL$ is a Fr\'echet submanifold of $\cA(\Sigma)$,
each tangent space $\rT_A\cL$ is a Lagrangian 
subspace of $\Om^1(\Sigma,\cg)$, \
$\cL\subset\cA_\flat(\Sigma)$, and $\cL$
is invariant under $\cG(\Sigma)$.

\smallskip\noindent{\bf (L2)}
The quotient of $\cL$ by the based gauge group $\cG_z(\Sigma)$
is compact, connected, simply connected, and
$\pi_2(\cL/\cG_z(\Sigma))=0$. 

\smallskip\noindent{\bf (L3)}
The zero connection is contained in $\cL$ and is nondegenerate
(as a critical point of the Chern-Simons functional).
Moreover, every nontrivial flat connection $A\in\cA(Y)$ 
with $A|_\Sigma\in\cL$ is irreducible. 

\smallskip\noindent
A detailed explanation 
and a finite dimensional characterization 
of these conditions is given in Section~\ref{sec:CS}. 
In particular, the assumptions imply that $\cL$ 
descends to a (singular) Lagrangian submanifold
$L:=\cL/\cG(\Sigma)\subset M_\Sigma$.
If $H$ is a disjoint union of handlebodies with 
$\p H=\bar\Sigma$ then the subset $\cL_H\subset\cA(\Sigma)$
of all flat connections on $\Sigma$ that extend to flat
connections on $Y$ satisfies (L1) and (L2).
It satisfies~(L3) if and only if $Y\cup_\Sigma H$
is an integral homology $3$-sphere. 

The space 
$
\cA(Y,\cL):=\left\{A\in\cA(Y)\,|\,A|_\Sigma\in\cL\right\}
$
of connections on $Y$ with boun\-dary values in $\cL$ 
carries a gauge invariant Chern--Simons functional
$$
\CS_\cL:\cA(Y,\cL)\to\R/4\pi^2\Z,
$$
well defined up to an additive constant, whose differential
is the usual Chern--Simons $1$-form (see Section~\ref{sec:CS}).
The critical points are the flat connections in $\cA(Y,\cL)$.
If we fix a Riemannian metric $g$ on $Y$ then the gradient
flow lines of the Chern--Simons functional with
respect to the $L^2$ inner product are smooth maps
$\R\to\cA(Y): s \mapsto A(s)$ satisfying the 
differential equation
\begin{equation}\label{eq:floerL}
\pd_s A + *F_A = 0, \qquad
A(s)|_\Sigma \in\cL \quad\forall s\in\R.
\end{equation}
As in Floer's original work~\cite{F1} the main idea
is to use the solutions of~(\ref{eq:floerL})
to construct a boundary operator on the chain 
complex generated by the gauge equivalence classes
of the nontrivial flat connections in $\cA_\flat(Y,\cL)$.
This defines the Floer homology groups $\HF(Y,\cL)$.
To make this precise one needs perturbations that 
turn $\CS_\cL$ into a Morse function whose gradient
flowlines satisfy Morse--Smale type transversality conditions.

We shall work with gauge invariant 
holonomy perturbations $h_f:\cA(Y)\to\R$ as 
in~\cite{T,F1,Donaldson book} (see Section~\ref{sec:CS}
and Appendix~\ref{app:Xf}).  The differential 
of $h_f$ has the form 
$
\rd h_f (A) \alpha = \int_Y \winner{X_f(A)}{\alpha}
$
for a suitable map $X_f:\cA(Y)\to\Om^2(Y,\cg)$.
The space of gauge equivalence classes of critical
points of the perturbed Chern--Simons functional 
$\CS_\cL+h_f$ will be denoted by 
$$
\cR_f := \left\{A\in\cA(Y,\cL)\,|\,F_A+X_f(A)=0\right\}/\cG(Y)
$$
and the perturbed gradient flow lines are solutions
of the boundary value problem
\begin{equation}\label{eq:floerLf}
\pd_s A + *\bigl(F_A+X_f(A)\bigr) = 0, \qquad
A(s)|_\Sigma \in\cL \quad\forall s\in\R.
\end{equation}
The space of gauge equivalence classes of solutions
of~(\ref{eq:floerLf}) that are asymptotic to
$[A^\pm]\in\cR_f$ as $s$ tends to $\pm\infty$
will be denoted by $\cM(A^-,A^+;g,f)$. 
In the transverse case with irreducible limits
$[A^\pm]\neq 0$ this moduli space is a manifold whose local dimension
near $[A]\in\cM(A^-,A^+;g,f)$ is given by the Fredholm index
$\delta_f(A)$ of a suitable linearized operator.
A crucial fact is the energy-index relation
$$
\delta_f(A)
= \frac{2}{\pi^2}E_f(A) + \eta_f(A^-) - \eta_f(A^+) 
$$
for the solutions of~(\ref{eq:floerLf}) with energy
$E_f(A)=\int_\R\|\pd_s A\|_{L^2(Y)}^2$, and with  
a function $\eta_f:\cR_f\to\R$.
This is Floer's monotonicity formula; it follows from the 
fact that $\cL/\cG_z(\Sigma)$ is simply connected.
The assumption on $\pi_2$ is only needed for the 
orientability of the moduli spaces. 

Floer's original work corresponds to the case $\p Y=\emptyset$.
The object of the present paper is to show
that all of Floer's ideas carry over to the 
case of nonempty boundary. 
The upshot is that, for a generic perturbation $h_f$,
all critical points of $\CS_\cL+h_f$ are nondegenerate
and so $\cR_f$ is a finite set, and that, for every
pair $[A^\pm]\in\cR_f$ the moduli space 
$\cM^1(A^-,A^+;g,f)$ of index $1$ connecting trajectories
consists of finitely many flow lines up to time shift.
The monotonicity formula plays a central role in this
finiteness theorem. As a result we obtain a Floer chain complex
$$
\CF_*(Y,\cL;f) 
:= \bigoplus_{[A]\in\cR_f\setminus[0]} \Z \, \langle A \rangle
$$
with boundary operator given by
$$
\pd \langle A^-\rangle \,:= \sum_{[A^+]\in \cR_f\setminus[0]} 
\#\bigl(\cM^1(A^-,A^+;g,f)/\R\bigr) \; 
\langle A^+ \rangle .
$$
Here the connecting trajectories are counted with appropriate signs 
determined by coherent orientations of the moduli spaces 
(Section~\ref{sec:orient}). It then follows from gluing 
and compactness theorems (Sections~\ref{sec:compact} 
and~\ref{sec:gluing}) that $\pd^2=0$.
The Floer homology groups are defined by
$$
\HF_*(Y,\cL;f,g):= \ker\pd/\im\pd.
$$
We shall prove that the Floer homology groups are independent 
of the choice of the metric $g$ and the perturbation $f$ 
used to define them (Section~\ref{sec:floer}).

\begin{rmk}\label{rmk:YLH}\rm
In the handle body case we expect the Floer homology groups 
$\HF(Y,\cL_H)$ to be naturally isomorphic to the instanton
Floer homology groups of the homology $3$-sphere 
$Y\cup_\Sigma H$. The proof will be carried out elsewhere.
\end{rmk}

\begin{rmk}\label{rmk:heegaard}\rm
An interesting special case arises from a Heegaard splitting 
$
M={H_0\cup_\Sigma \bar H_1}
$ 
of a homology $3$-sphere into two handle bodies $H_i$ with 
$\pd H_i=\Sigma$. We obtain the Floer homology groups 
$\HF_*([0,1]\times\Sigma,\cL_{H_0}\times \cL_{H_1})$ 
from the following setup: 
The $3$-manifold 
$
Y:=[0,1]\times\Sigma
$
has two boundary components $\pd Y=\bar\Sigma \sqcup \Sigma$,
and attaching the disjoint union of the handle bodies 
$H:=H_0 \sqcup \bar H_1$ yields the homology $3$-sphere 
$
Y\cup_{\bar\Sigma \sqcup \Sigma} H \cong M.
$ 
The Lagrangian submanifold is 
$
\cL_{H_0}\times \cL_{H_1}\cong\cL_H
\subset
\cA(\bar\Sigma\sqcup\Sigma).
$
If this Floer homology is isomorphic to $\HF_*(M)$, as expected,
then the proof of the Atiyah--Floer conjecture for $M$ reduces
to an adiabatic limit argument as in~\cite{DS} which identifies
the symplectic Floer homology group of the pair of Lagrangian
submanifolds $L_{H_0},L_{H_1}$ of the singular symplectic
manifold 
$
M_\Sigma:=\cA_{\mathrm{flat}}(\Sigma)/\cG(\Sigma)
$ 
with the Floer homology groups 
$\HF([0,1]\times\Sigma,\cL_{H_0}\times\cL_{H_1})$
defined in the present paper. Since $M_\Sigma$ is a singular space,
this requires as a preliminary step the very definition of the 
symplectic Floer homology groups of $L_{H_0}$ and~$L_{H_1}$
with $L_{H_i}:=\cL_{H_i}/\cG(\Sigma)$.
\end{rmk}
 
\begin{rmk}\label{rmk:product}\rm
If $H_0, H_1, H_2$ are three handle bodies with boundary $\Sigma$ 
such that the manifold $M_{ij}:=H_i\cup_\Sigma \bar H_j$ 
is a homology $3$-sphere for $i\ne j$, then there is a product 
morphism
$$
\HF_*(Y,\cL_{H_0}\times\cL_{H_1})\times 
\HF_*(Y,\cL_{H_1}\times\cL_{H_2}) \to
\HF_*(Y,\cL_{H_0}\times\cL_{H_2}),
$$
where $Y:=[0,1]\times\Sigma$.
A key ingredient in the definition is the observation 
that~(\ref{eq:floerLf}) is the perturbed anti-self-duality equation 
for a connection on $\R\times Y$ in temporal gauge.
Thus equation~(\ref{eq:floerLf}) can be generalized 
to a $4$-manifold $X$ with a boundary space-time splitting 
and tubular ends (Section~\ref{sec:fredholm}). 
The definition of the product morphism will be based on the 
moduli space for the $4$-manifold $X=\Delta\times\Sigma$,
where $\Delta$ is a triangle (or rather a disc with three cylindrical
ends attached). The details will be carried out elsewhere. 
We expect that our conjectural isomorphisms 
will intertwine the corresponding product structures 
on the symplectic and instanton Floer homologies.  
\end{rmk}

The construction of the Floer homology groups in the 
present paper is based on the foundational analysis 
in~\cite{W Cauchy,W elliptic,W bubb, W lag} for the 
solutions of the boundary value problem~(\ref{eq:floerL}). 
In our exposition we follow the work of Floer~\cite{F1}
and Donaldson~\cite{Donaldson book} 
and explain the details whenever new phenomena arise
from our boundary value problem. 
Recall that the present Lagrangian boundary conditions are a mix
of first order conditions (flatness of the restriction to $\pd Y$) 
and semi-global conditions (pertaining the holonomy on $\pd Y$),
so they cannot be treated by standard nonlinear elliptic methods.

In Section~\ref{sec:CS} we recall the basic properties
of the Chern--Simons functional on a $3$-manifold with boundary
and in Section~\ref{sec:HA} we discuss the Hessian
and establish the basic properties of the linearized operator
on $\R\times Y$. 
Section~\ref{sec:S1Y} examines the spectral flow and the 
determinant line bundle for operators over $S^1\times Y$.
Section~\ref{sec:exp} establishes exponential decay on 
tubular ends.
Section~\ref{sec:fredholm} sets up the Fredholm theory
for general $4$-manifolds with space-time splittings of
the boundary and tubular ends.  In the second half of the 
section we focus on the tube $\R\times Y$, examine the 
spectral flow, and prove monotonicity.
Section~\ref{sec:compact} proves 
the compactness of the moduli spaces, 
based on~\cite{W elliptic, W bubb}.

In Section~\ref{sec:moduli} we establish transversality, 
using holonomy perturbations.  The novel difficulty here 
is that we do not have a geometric description 
of the bubbling effect at the boundary.
So, instead of a gluing theorem converse to bubbling, we 
use monotonicity and work inductively on the energy levels.
The second difficulty is that we need to keep the support of the
perturbations away from the boundary, since the techniques
of~\cite{W bubb} do not extend to the perturbed equation.
As a result we cannot obtain an open and dense set of regular
perturbations but -- still sufficient --
we find a regular perturbation up to index $7$ 
near any given perturbation.
In an appendix to this section we establish the relevant 
unique continuation results. 
In the process we reprove Taubes' 
unique continuation result~\cite{T:ucon} for anti-self-dual
connections that vanish to infinite order at a point.
This is needed to overcome difficulties arising from 
the nonlinear boundary conditions.  After these preparations, 
the construction of the Floer homology follows the standard routine. 
For the gluing results in Section~\ref{sec:gluing} we focus
on the pregluing map and the Banach manifold
setup for the inverse function theorem.
In Section~\ref{sec:orient} we construct
coherent orientations in the Lagrangian setting.
The Floer homology groups are defined
in Section~\ref{sec:floer}.

There are several appendices where we review standard 
techniques and adapt them to our boundary value 
problems. Appendix~\ref{app:spec} deals with
the spectral flow for self-adjoint operator families 
with varying domains. Appendix~\ref{app:GR} discusses
the Gelfand--Robbin quotient, an abstract setting which relates
self-adjoint operators with Lagrangian subspaces.
These results are needed for the index calculations and orientations
in Sections~\ref{sec:S1Y} and~\ref{sec:fredholm}.
Appendix~\ref{app:uc} reviews the Agmon--Nirenberg 
unique continuation technique used in Section~\ref{sec:moduli}.
In Appendix~\ref{app:Xf} we discuss the basic analytic properties 
of the holonomy perturbations and prove a compactness result 
needed in Section~\ref{sec:compact}.
Appendix~\ref{app:Lag} deals with Lagrangian submanifolds 
in the space of connections. We construct an $L^2$-continuous 
trivialization of the tangent bundle $\rT\cL$, 
used in Sections~\ref{sec:HA} and~\ref{sec:fredholm},
and a gauge invariant exponential map for $\cL$,
used in Section~\ref{sec:gluing}.

\smallskip\noindent{\bf Notation.}
We denote the spaces of smooth connections and gauge 
transformations on a manifold $Z$ by $\cA(Z):=\Om^1(Z,\cg)$ 
and $\cG(Z):=\cC^\infty(Z,\rG)$.
The gauge group $\cG(Z)$ acts on $\cA(Z)$ 
by $u^*A:=u^{-1}A u + u^{-1}\rd u$ and the gauge 
equivalence class of $A\in\cA(Z)$ is denoted by $[A]$.
A connection $A\in\cA(Z)$ induces an exterior differential
$\rd_A :\Omega^k(Z,\cg)\to \Omega^{k+1}(Z,\cg)$ via
$
\rd_A \tau := \rd \tau + [A\wedge\tau].
$
Here $[\cdot,\cdot]$ denotes the Lie bracket on~$\cg$.
The curvature of $A$ is the $2$-form 
$
F_A:=\rd A + A\wedge A
$ 
and it satisfies $\rd_A\rd_A\tau=[F_A\wedge\tau]$.
The space of flat connections is denoted by
$
\cA_\flat(Z) := \{ A\in\cA(Z) \st F_A=0 \}.
$
Connections on $X=\R\times Y$ or other $4$-manifolds 
will be denoted by $\A$ or $\Xi$, whereas $A$ denotes 
a connection on a $3$-manifold $Y$ or a $2$-manifold $\Sigma$.
We say that a connection $\A=A+\Phi\ds$ on $\R\times Y$
is in {\bf temporal gauge} on $I\times Y$ if $\Phi|_{I\times Y}\equiv 0$. 


\section{The Chern--Simons functional}\label{sec:CS}

Let $Y$ be a compact oriented $3$-manifold with 
boundary $\p Y=\Sigma$ and $\rG=\SU(2)$.  
The Chern--Simons $1$-form on $\cA(Y)$ is defined by
\begin{equation} \label{CS 1form}
\alpha \;\mapsto\; \int_Y \winner{F_A}{\alpha}
\end{equation}
for $\alpha\in\rT_A\cA(Y)= \Omega^1(Y,\cg)$.
If $Y$ is closed, then (\ref{CS 1form}) is the 
differential of the Chern--Simons functional
$\CS:\cA(Y)\to\R$ given by
$$
\CS(A) := 
\frac 12 \int_Y \Bigl(
\winner{A}{\rd A} + \frac 13 \winner{A}{[A\wedge A]} 
\Bigr).
$$
It changes by 
\begin{equation}\label{eq:CS}
\CS(A) - \CS(u^*A) = 4\pi^2\deg(u)
\end{equation}
under a gauge transformation $u\in\cG(Y)$;  
thus the Chern--Simons functional descends to a circle valued 
function $\cB(Y):=\cA(Y)/\cG(Y)\to \R/4\pi^2\Z$ 
which will still be denoted by $\CS$.
If $Y$ has nonempty boundary $\pd Y=\Sigma$, 
then the differential of (\ref{CS 1form})
is the standard symplectic form~(\ref{eq:symp})
on $\cA(\Sigma)$.  To obtain a closed $1$-form 
we restrict the Chern--Simons $1$-form to a subspace 
of connections satisfying a Lagrangian boundary condition. 

\subsection*{Lagrangian submanifolds}

The relevant Lagrangian submanifolds of $\cA(\Sigma)$ 
were studied in detail in~\cite[Section~4]{W Cauchy}. 
Following~\cite{W Cauchy} we assume that
$
\cL\subset\cA(\Sigma)
$
is a gauge invariant Lagrangian submanifold 
satisfying (L1). This condition can be rephrased 
as follows.

\label{p:L1}
\begin{description}
\item[(L1)]
First, $\cL$ is contained in $\cA_{\mathrm{flat}}(\Sigma)$ 
and is invariant under the action of $\cG(\Sigma)$. 
Second, for some (and hence every) $p>2$ the $L^p$-closure of $\cL$ is a 
Banach submanifold of the space of $L^p$-connections, 
$\cA^{0,p}(\Sigma):=L^p(\Sigma,\rT^*\Sigma\otimes\cg)$.
Third, for every $A\in\cL$ the tangent space 
$\rT_A\cL\subset\Om^1(\Sigma,\cg)$ is Lagrangian, i.e. 
\begin{equation} \label{eq:omega}
\om(\alpha,\beta) = 0\quad
\forall \beta\in \rT_A\cL\qquad
\iff\qquad \alpha\in \rT_A\cL
\end{equation}
for every $\alpha\in\Om^1(\Sigma,\cg)$.
\end{description}

Let $\cL^{0,p}\subset\cA^{0,p}(\Sigma)$ denote the $L^p$-closure 
of $\cL$.  Then $\cL =\cL^{0,p}\cap\cA(\Sigma)$ and the tangent 
space $\rT_A\cL$ of a smooth element $A\in\cL$ - as in (L1) -
is understood as the intersection of the Banach tangent space 
$\rT_A\cL^{0,p}$ with the space of smooth $1$-forms.  
This space is independent of $p>2$ and coincides with the space 
of derivatives of smooth paths in $\cL$ passing 
through $A$.\footnote{
It is not clear whether one could also work with Hilbert 
submanifolds $\cL\subset\cA^{0,2}(\Sigma)$. This is connected 
to subtle questions concerning the gauge action at this 
Sobolev borderline, see~\cite{W lag}.
} 
This follows from a finite dimensional
characterization of the manifold property
which we explain next.

A {\bf base point set} is a finite set $z\subset\Sigma$ which 
intersects each component of $\Sigma$ in precisely one point. 
For every base point set $z$ the based gauge group 
$\cG_z(\Sigma):=\{u\in\cG(\Sigma)\,|\, u(z)\equiv\one\}$
acts freely on $\cA(\Sigma)$.  Let $2g:=\dim\,H_1(\Sigma)$
and pick $2g$ loops in $\Sigma$ that generate $H_1(\Sigma)$
with base points chosen from $z$. The holonomy around these
loops defines a map $\rho_z:\cA_{\rm flat}\to\rG^{2g}$ which is 
invariant under the action  of the based gauge group $\cG_z(\Sigma)$. 
If $\cL$ is a gauge invariant subset of $\cA_\flat(\Sigma)$ 
then $\cL^{0,p}$ is a Banach submanifold of $\cA^{0,p}(\Sigma)$
if and only if the image $\rho_z(\cL)\subset\rG^{2g}$ of the 
holonomy morphism is a smooth submanifold.  
There is however no well defined moment map for the action 
of $\cG_z(\Sigma)$, so the symplectic structure does 
not descend to the quotient.  On the other hand, the quotient 
$L:=\cL/\cG(\Sigma)$ has singularities in general, 
but it intersects the smooth part of the moduli space
$M_\Sigma:=\cA_{\mathrm{flat}}(\Sigma)/\cG(\Sigma)$
in a Lagrangian submanifold.

If $\cL^{0,p}\subset\cA^{0,p}(\Sigma)$ is a Lagrangian submanifold 
then $\cL$ is gauge invariant if and only if 
$\cL\subset\cA_\flat(\Sigma)$; \cite[Sec.~4]{W Cauchy}.
Condition~(L1) implies that $\cL$ is a totally real submanifold
with respect to the Hodge $*$-operator for any metric on $\Sigma$,
i.e.
$$
\Om^1(\Sigma,\cg)=\rT_A\cL\oplus *\rT_A\cL 
\qquad\forall A\in\cL .
$$
The construction of Floer homology groups for the Chern--Simons $1$-form
will require the following additional assumptions on $\cL$.

\begin{description}
\item[(L2)]
The quotient space $\cL/\cG_z(\Sigma)$ is compact, 
connected, simply connected, and $\pi_2(\cL/\cG_z(\Sigma))=0$
for some (and hence every) base point set $z\subset\Sigma$.
\item[(L3)]
The zero connection is contained in $\cL$.
It is nondegenerate in the sense that
$
\rd \alpha = 0\iff\alpha \in \im\rd
$
for every $\alpha\in\rT_0\cA(Y,\cL)$.
Moreover, every flat connection in $\cA(Y,\cL)$ that is not
gauge equivalent to the zero connection is irreducible. 
\end{description}

In (L2) the hypothesis that $\cL/\cG_z(\Sigma)$ is simply 
connected is needed to establish an energy-index relation 
for the Chern-Simons functional.  The hypothesis 
$\pi_2(\cL/\cG_z(\Sigma))=0$ is only used to orient the 
moduli spaces. It can be dropped if one wants to define 
Floer homology with $\Z_2$ coefficients.  These two conditions 
imply that $\pi_1(\cL)$ is isomorphic to 
$\pi_1(\cG_z(\Sigma))\cong\pi_1(\cG(\Sigma))$ and the map 
$\pi_2(\cG_z(\Sigma))\cong\pi_2(\cG(\Sigma))\to \pi_2(\cL)$ 
is surjective. To see this, note that $\cL$ is a fiber bundle 
over the base $\cL/\cG_z(\Sigma)$ (see \cite[Lemma~4.3]{W Cauchy}).
In particular, (L2) implies that $\pi_1(\cL)\cong\Z^{\pi_0(\Sigma)}$
since the fiber $\cG_z(\Sigma)$ has fundamental group $\Z^N$ whenever 
$\Sigma$ has $N$ connected components.
(For a connected component $\Sigma'$ an isomorphism 
$\pi_1(\cG_z(\Sigma'))\cong\Z$ is given by the degree 
of a map $S^1\times\Sigma' \to \SU(2)\cong S^3$.)

The main example of a Lagrangian submanifold 
of $\cA(\Sigma)$ arises from the space of flat connections 
on a disjoint union $H$ of handle bodies\footnote
{
A handle body is an oriented $3$-manifold with 
boundary that is obtained from a $3$-ball by 
attaching $1$-handles. Equivalently, it admits a Morse function
with exactly one minimum, no critical points of index $2$, and
attaining its maximum on the boundary.
}
with boundary $\pd H=\bar\Sigma$. 
Here $\bar\Sigma$ is the same manifold as $\Sigma$ 
but equipped with the opposite orientation. 
Given such a manifold $H$ define
$$
\cL_H \,:=\;  \bigl\{ \tA|_\Sigma \st \tA\in\cA_\flat(H) \bigr\}.
$$

\begin{lemma}\label{le:YH}
Let $H$ be a disjoint union of handle bodies with
$\p Y=\bar\Sigma$.  Then the following holds.

\smallskip\noindent{\bf (i)}
$\cL_H$ is a Lagrangian submanifold of $\cA(\Sigma)$
that satisfies (L1) and (L2) and contains the 
zero connection.

\smallskip\noindent{\bf (ii)}
The zero connection is nondegenerate if and only 
if $Y\cup H$ is a rational homology $3$-sphere 

\smallskip\noindent{\bf (iii)}
Every nontrivial flat connection in $\cA_\flat(Y,\cL_H)$
is irreducible if and only if $Y\cup H$ is an integral 
homology $3$-sphere 
\end{lemma}

\begin{proof}
That $\cL_H$ satisfies (L1) was proved 
in~\cite[Lemma~4.6]{W Cauchy}.  That $\cL_H$ 
contains the zero connection is obvious. 
That it satisfies (L2) follows from the fact 
that the based holonomy map $\rho_z$ induces a 
homeomorphism from $\cL_H/\cG_z(\Sigma)$ to $\rG^g$ 
with $\rG=\SU(2)$ when $\Sigma$ is connected 
and has genus $g$, and that 
$$
\cL_{H_1\sqcup\ldots\sqcup H_m}/
\cG_{\{z_1,\ldots,z_m\}}(\Sigma_1\sqcup\ldots\sqcup\Sigma_m)
\cong \cL_{H_1}/\cG_{z_1}(\Sigma_1)\times\ldots\times 
\cL_{H_m}/\cG_{z_m}(\Sigma_m)
$$
in the case of several connected components.
This proves~(i).

To prove~(ii) we need to consider $\alpha\in\Om^1(Y,\cg)$ with $\rd\alpha=0$.
The linearized Lagrangian boundary condition on $\alpha$ is equivalent to the 
existence of an extension $\tilde\alpha\in\Om^1(Y\cup H,\cg)$
with $\rd\tilde\alpha=0$.
If $H^1(Y\cup H;\R)=0$ (or equivalently $H_1(Y\cup H;\Q)=0$),
then any such $1$-form is exact on $Y\cup H$ and thus on $Y$.
Conversely, if $\tilde\alpha\in\ker\rd$, then nondegeneracy
implies $\tilde\alpha|_Y\in\im\rd$ and hence
$\int_\gamma\tilde\alpha = 0$ for every loop $\gamma\subset Y$.
This implies that $\tilde\alpha$ is also exact on $Y\cup H$ 
since every loop in $Y\cup H$ is homotopic to a loop in $Y$.
This proves~(ii).

We prove~(iii).  Flat connections in $\cA(Y,\cL_H)$ 
can be identified with flat connections in $\cA(Y\cup H)$. 
The gauge equivalence classes of irreducible but 
nontrivial connections are in one-to-one correspondence with nontrivial 
homomorphisms $\pi_1(Y\cup H)\to S^1$. These exist if and only if
$H_1(Y\cup H;\Z)\neq 0$. 
\end{proof}

\subsection*{Lagrangian submanifolds and representations}

We characterize our Lagrangian submanifolds as subsets 
of the representation spaces for Riemann surfaces. 
For simplicity we assume first that $\Sigma$ is connected. 
Fix a base point $z\in\Sigma$ and choose based loops 
$\alpha_1,\dots,\alpha_g,\beta_1,\dots,\beta_g$
representing a standard set of generators\footnote{
The standard generators of $\pi_1(\Sigma,z)$ satisfy the relation
$\prod_{i=1}^g \alpha_i\beta_i\alpha_i^{-1}\beta_i^{-1} = \one $.} 
of the fundamental group.
The based holonomy around the loops $\alpha_i$ and $\beta_i$
gives rise to a map 
$
\rho_z:\cA(\Sigma)\to\rG^{2g}.
$
This map identifies the moduli space $M_\Sigma$ 
of flat connections with the quotient of $f^{-1}(\one)$
by conjugation, where $f:\rG^{2g}\to\rG$ is defined by
\begin{equation}\label{thatsf}
f(x_1,\ldots,x_g,y_1,\ldots y_g):=
{\textstyle\prod_{i=1}^g} x_i y_i x_i^{-1} y_i^{-1} .
\end{equation}
The correspondence between flat connections and representations
is reformulated in~(a) and~(b) below.  Assertions~(c) and~(d) are 
the infinitesimal versions of these observations.  

\begin{remark}\rm 
\begin{description}
\item[(a)]
Let $w=(x_1,\ldots,x_g,y_1,\ldots y_g)\in\rG^{2g}$.
Then there exists a flat connection $A\in\cA_\mathrm{flat}(\Sigma)$
with $\rho_z(A)=w$ if and only if $f(w)=\one$.
\item[(b)]
Let $A,A'\in\cA_\mathrm{flat}(\Sigma)$.  Then $A$ is gauge equivalent
to $A'$ if and only if $\rho_z(A)$ is conjugate to $\rho_z(A')$.
\item[(c)]
Let $A\in\cA_\mathrm{flat}(\Sigma)$, 
$w:=\rho_z(A)$, and $\hat w\in\rT_w\rG^{2g}$.
Then $\rd f(w)\hat w=0$ if and only if there exists an
$\alpha\in\Om^1(\Sigma,\cg)$ such that $\rd_A\alpha=0$ and 
$\rd\rho_z(A)\alpha=\hat w$. 
\item[(d)]
Let $A\in\cA_\mathrm{flat}(\Sigma)$ and $\alpha\in\Om^1(\Sigma,\cg)$.
Denote $w:=\rho_z(A)$ and $\hat w:=\rd\rho_z(A)\alpha$.
Then $\alpha\in\im\rd_A$ if and only if $\hat w$ belongs to
the image of the infinitesimal conjugate action 
$L_w:\cg\to \rT_w\rG^{2g}$ given by $L_w\xi=\xi w - w \xi$.
\end{description}
\end{remark}

While the identity element $\one\in\rG$ is not a regular value
of $f$, it follows from~(c),(d) that the differential
$\rd\rho_z(A):\Om^1(\Sigma,\cg)\to T_w\rG^{2g}$
at a flat connection $A\in\cA_\mathrm{flat}(\Sigma)$
identifies $H^1_A:=\ker\rd_A/\im\rd_A$
(the virtual tangent space of $M_\Sigma$)
with the quotient $\ker\rd f(w)/\im L_w$ at $w=\rho_z(A)$. 
The gauge invariant symplectic form (\ref{eq:symp})
descends to $H^1_A$ and thus induces a symplectic form 
$$
\Om_w:\ker\rd f(w)/\im L_w\times\ker\rd f(w)/\im L_w\to\R
$$
$$
\Om_w(\hat w,\hat w') :=  \int_\Sigma\winner{\alpha}{\alpha'},
$$
where the (infinitesimal) connections 
$A\in\cA_{\rm flat}(\Sigma)$ and $\alpha,\alpha'\in\ker\rd_A$ 
are chosen such that $w=\rho_z(A)$, 
$\hat w = \rd\rho_z(A)\alpha$, and $\hat w' = \rd\rho_z(A)\alpha'$.
An explicit formula for this symplectic form at 
$w=(x_1,\ldots,x_g,y_1,\ldots,y_g)$ 
on the vectors ${\hat w=(\xi_1 x_1,\ldots,\xi_g x_g,\eta_1 y_1,\ldots,\eta_g y_g)}$,
$\hat w'=(\xi'_1 x_1,\ldots,\xi'_g x_g,\eta'_1 y_1,\ldots,\eta'_g y_g)$
is
\begin{align} \label{Omega}
\Omega_w(\hat w,\hat w') &= 
\tsum_{i=1}^g \bigl(
\la ( x_i^{-1}\xi_i x_i + x_i^{-1} \delta_i x_i - \delta_{i-1} ) , 
\eta'_i \ra \\
&\qquad\qquad
- \la ( y_i^{-1} \eta_i y_i + y_i^{-1} \delta_i y_i - \delta_{i-1} ) , 
\xi'_i \ra
\bigr). \nonumber
\end{align}
Here $\delta_j = - (\rd\hol(A)\alpha)\hol(A)^{-1}$ is the infinitesimal holonomy
along the path $\prod_{i=1}^j \alpha_i\beta_i\alpha_i^{-1}\beta_i^{-1}$, i.e.
$$
\delta_j = c_j^{-1} \delta_{c_j} 
+ c_j^{-1}c_{j-1}^{-1} \delta_{c_{j-1}} c_j + \ldots
 + c_j^{-1}\ldots c_{1}^{-1} \delta_{c_{1}} c_2 \ldots c_j ,
$$
$$
c_i := x_i y_i x_i^{-1} y_i^{-1} , \qquad
\delta_{c_i} := x_i y_i \bigl( y_i^{-1} \xi_i y_i - \xi_i + \eta_i
- x_i^{-1} \eta_i x_i  \bigr)  x_i^{-1} y_i^{-1} .
$$
One should compare this with the identities $f(w)=c_1\cdots c_g =\one$ and
$$
\rd f(w)(\hat w) 
= c_1\ldots c_{g-1}\delta_{c_g} 
+ c_1\ldots c_{g-2}\delta_{c_{g-1}} c_{g}
+ \ldots + \delta_{c_1}c_2\ldots c_{g} = 0.
$$
Combining these we see that $\delta_{c_g}=\one$. 
So on the torus $\Sigma=\T^2$ the formula simplifies to 
$\Omega_w(\hat w,\hat w') = 
\la x^{-1} \xi x , \eta' \ra - \la y^{-1} \eta y , \xi' \ra $.
Moreover, if $T\subset\rG$ is any circle and 
$w\in T^{2g}\subset \rG^{2g}$ then the restriction 
of $\Omega_w$ to $\R^{2g}\cong\rT_w T^{2g}\subset\ker\rd f(w)$ 
is the standard symplectic form on Euclidean space.
By construction and assertions~(a-d) above, $\Omega$ descends 
to the symplectic form on the (singular) symplectic quotient 
$f^{-1}(\one)/\rG\cong M_\Sigma=\cA(\Sigma)\dslash\cG(\Sigma)$.
In fact, one can verify directly that $\Omega$ is $\rG$-invariant and that
its kernel at each point is the tangent space to the $\rG$-orbit.
Thus, on the complement of the reducible set, 
$\Omega$ descends to a smooth symplectic structure 
on the $\rG$-quotient.  

In the case of the torus $\Sigma=\T^2$ all points of 
$f^{-1}(\one)$ are reducible; 
in this case $M_{\T^2}$
can be identified with the quotient of the moduli space of flat 
$S^1$-connections by a residual $\Z_2$-action
with four isolated fixed points $\{(\pm\one,\pm\one)\}$
(corresponding to the same four points in $\rG^2$).
For a general surface $\Sigma$, the set of reducibles 
in $f^{-1}(\one)$ is the union $\bigcup_{T\subset\rG} T^{2g}$ 
over all maximal tori $T\subset\rG$. 
For $g>2$ this set has codimension $4g-2>3g$ in $\rG^{2g}$.
So for a half dimensional submanifold $N\subset \rG^{2g}$ the set 
of irreducibles will always be dense in $N$.   In the case of genus $2$ 
the same is true if we require $\Omega|_{\rT N}\equiv 0$,
since the codimension of the set of reducibles is $3g$ but $N$ cannot 
intersect it in an open set since $\Omega$ is nondegenerate on each 
subtorus $T^{2g}$ of the reducibles.

If $\Sigma$ has several connected components we fix a base point
set $z\subset\Sigma$ and obtain the $2$-form $\Omega$ as sum of the
$2$-forms of the connected components.
We then have $M_\Sigma\cong f^{-1}(\one)/\rG^{\pi_0(\Sigma)}$,
where $\rG^{\pi_0(\Sigma)}$ acts by conjugation with a fixed group element on 
each connected component and $f:\rG^{2g}\to\rG^{\pi_0(\Sigma)}$ is the
product of the relations~(\ref{thatsf}) for each connected component.
Now we can reformulate the assumptions (L1-3) on the Lagrangian 
submanifolds $\cL\subset\cA(\Sigma)$ as follows:
{\it $\cL=\rho_z^{-1}(N)\subset\cA(\Sigma)$ is the preimage of a submanifold
$N\subset\rG^{2g}$ satisfying the following conditions.}
\begin{description}
\item[(L1)]
{\it $N \subset f^{-1}(\one)$, $N$ is invariant under $\rG^{\pi_0(\Sigma)}$,
$\dim N=3g$, and $\Omega|_{\rT N}\equiv 0$.}
\item[(L2)]
{\it $N$ is compact, connected, simply connected, and $\pi_2(N)=0$.}
\item[(L3)]
{\it $N\cap \rho_z(\cA_{\rm flat}(Y))$ contains
$(\one,\ldots,\one)$ as isolated point and does not
contain any other reducible points (with respect to the
conjugate action of $\rG^{\pi_0(\Sigma)}$).}
\end{description}
The above discussion of the reducible locus shows that, by condition~(L1),
the quotient $L:=N/G^{\pi_0(\Sigma)}\subset M_\Sigma$ is Lagrangian at 
a dense set of smooth points.

\subsection*{The Chern--Simons functional}

Fix a compact, connected, oriented $3$-manifold $Y$ 
with nonempty boundary ${\pd Y=\Sigma}$ and a
gauge invariant, monotone, irreducible Lagrangian submanifold $\cL\subset\cA(\Sigma)$
satisfying (L1-3) on page~\pageref{p:L1}. Then the restriction of the 
Chern--Simons $1$-form (\ref{CS 1form}) to the submanifold
$$
\cA(Y,\cL):=\{ A\in\cA(Y) \st A|_\Sigma \in \cL\}
$$
is closed. It is the differential of the 
circle valued Chern--Simons functional
$$
\CS_\cL:\cA(Y,\cL)\to \R/4\pi^2\Z
$$ 
given by $\CS_\cL(A):=[\CS(A,B)]$, where
$$
\CS(A,B) := 
\frac 12 \int_Y \Bigl( \winner{A}{\rd A} 
+ \frac 13 \winner{A}{[A\wedge A]}\Bigr)  
- \frac12\int_0^1\int_\Sigma \winner{B(s)}{\pd_s B(s)}\,\ds.
$$
Here $B:[0,1]\to\cL$ is a smooth path satisfying 
$B(0)=A|_\Sigma$ and $B(1)=0$. 

\begin{rmk}\label{rmk:CSAB} \rm
Note that $\CS(A,B)$ is the value of the Chern--Simons functional 
on the connection $\tA$ on $\tY:=Y\cup \bigl([0,1]\times\Sigma\bigr)$ 
given by $A$ on $Y$ and by $B$ on $[0,1]\times\Sigma$. 
Here we glue $\pd Y=\Sigma$ to $\{0\}\times\bar\Sigma$,
and on the new boundary $\pd\tY=\{1\}\times\Sigma$
we have $\tA\equiv 0$.  
\end{rmk}

\begin{lem}\label{le:CS}
{\bf (i)} 
The Chern--Simons functional $\CS(A,B)$ is invariant
under homotopies of $B$ with fixed endpoints.  

\smallskip\noindent{\bf (ii)}
If $u:[0,1]\to\cG(\Sigma)$ satisfies
$u(0)=u(1)=\one$ then 
$$
\CS(A,B) - \CS(A,u^*B) = 4\pi^2\, {\rm deg}\, u .
$$
\smallskip\noindent{\bf (iii)}
If $B_0,B_1:[0,1]\to\cL$ are two paths with 
$B_0(0)=B_1(0)$ and $B_0(1)=B_1(1)=0$
then there is a path $u:[0,1]\to\cG(\Sigma)$ 
with $u(0)=u(1)=\one$ such that $B_1$ is 
homotopic to $u^*B_0$ (with fixed endpoints).

\smallskip\noindent{\bf (iv)}
The circle valued function 
$\CS_\cL:\cA(Y,\cL)\to\R/4\pi^2\Z$ descends to 
the quotient $\cB(Y,\cL):=\cA(Y,\cL)/\cG(Y)$.
\end{lem}

\begin{proof}
The Chern-Simons functional is invariant under homotopies since 
\begin{align*}
- \pd_t \CS(A,B_t) 
&= \frac12\int_0^1\int_\Sigma \Bigl( \winner{\pd_t B_t(s)}{\pd_s B_t(s)}
   + \winner{B_t(s)}{\pd_t\pd_s B_t(s)} \Bigr)\,\ds \\
&= \int_0^1\int_\Sigma \winner{\pd_t B_t(s)}{\pd_s B_t(s)} \,\ds
   + \biggl[ \int_\Sigma \winner{B_t(s)}{\pd_t B_t(s)} \biggr]_0^1.
\end{align*}
for every smooth homotopy $B_t:[0,1]\to\cL$ with fixed
endpoints.  The first term on the right is the symplectic form on 
$\pd_t B_t,\pd_s B_t\in\rT_B\cL$ and the second term vanishes 
since $\pd_t B_t(s)=0$ for $s=0,1$. Hence $\pd_t \CS(A,B_t)=0$. 
This proves~(i).

To prove~(ii), we abbreviate $S^1:=\R/\Z$, define 
$\tu:S^1\times\Sigma\to\SU(2)$ by $\tu(t,z):=u(t)(z)$,
and calculate
\begin{align*}
&2 \bigl( \CS(A,B) - \CS(A,u^*B) \bigr) \\
&= \int_0^1\int_\Sigma 
\Bigl( \winner{u^*B}{\pd_s(u^*B)} - \winner{B}{\pd_s B} \Bigr) \,\ds \\
&= \int_0^1\int_\Sigma 
\Bigl( \winner{B}{\rd_B(\pd_s u \cdot u^{-1}) }
+ \winner{\rd u\cdot u^{-1}}{\bigl(\pd_s B 
+ \rd_B(\pd_s u \cdot u^{-1})\bigr) } 
 \Bigr) \,\ds \\
&= \int_0^1\int_\Sigma 
\Bigl( \winner{B}{\bigl(2\rd(\pd_s u \cdot u^{-1}) 
+ [B , \pd_s u \cdot u^{-1}]\bigr)}
+ \winner{\rd u\cdot u^{-1}}{\rd(\pd_s u \cdot u^{-1})} 
 \Bigr) \,\ds \\
&= 2\int_0^1\int_\Sigma  \winner{F_B}{\pd_s u \cdot u^{-1}}
\; -\frac 13 \int_{S^1\times\Sigma} {\rm tr}\bigl(
\rd \tu\cdot \tu^{-1} \wedge 
\rd \tu \cdot \tu^{-1} \wedge \rd \tu \cdot \tu^{-1} \bigr) \\
&= 8\pi^2\, {\rm deg}\, \tu .
\end{align*}
Here the first equation follows from the definitions,
the second equation uses the formula
$\p_s(u^*B)=u^{-1}(\p_sB+\rd_B(\p_su\cdot u^{-1}))u$,
the third equation uses integration by parts in $s$
and the fact that $\rd u(0)=\rd u(1)=0$,
the fourth equation uses the formula 
$d(\p_su\cdot u^{-1})-\p_s(du\cdot u^{-1})
=[du\cdot u^{-1},\p_su\cdot u^{-1}]$
and integration by parts over $\Sigma$,
and the last equation follows from the fact that
$F_{B(s)}=0$ for every $s$ and that the standard 
volume form on $\SU(2)$ with integral~$1$ is
$
24\pi^2 \tu^*\dvol_{\SU(2)} 
= -\trace\bigl(\rd \tu\cdot \tu^{-1} \wedge \rd \tu \cdot \tu^{-1} \wedge 
\rd \tu \cdot \tu^{-1} \bigr).
$
Thus we have proved~(ii).

To see~(iii) note that the catenation of $-B_0$ and $B_1$ 
is a loop in $\cL$ based at~$0$.  It is contractible in 
the base of the fibre bundle 
$\cG_z(\Sigma)\hookrightarrow \cL \to \cL/\cG_z(\Sigma)$
and hence it is homotopic to a loop $u:[0,1]\to\cG_z(\Sigma)$ in the fibre 
based at ${u(0)=u(1)\equiv\one}$.
Now the catenation of $B_0$, $-B_0$, and $B_1$ is homotopic with fixed
endpoints to $B_1$ on the one hand, and on the other hand to the catenation
of $B_0$ with the loop $u^*0$, which is also homotopic to $u^*B_0$.

It follows from~(i-iii) that the map
$(A,B)\to\CS(A,B)$ induces a circle valued function
$\CS_\cL:\cA(Y,\cL)\to\R/4\pi^2\Z$. 
We prove that this function is invariant under
gauge transformations. 
To see this we can use Remark~\ref{rmk:CSAB}
and extend any given $u\in\cG(Y)$ to 
a gauge transformation $\tilde u\in\cG(\tilde Y)$ 
on $\tY:=Y\cup \bigl([0,1]\times\Sigma\bigr)$ 
with $\tilde u|_{\tilde Y}\equiv\one$.
Such an extension exists because $\cG(\Sigma)$ 
is connected (which in turn follows from the fact that
$\rG=\SU(2)$ is connected, simply connected, and $\pi_2(\rG)=0$).
Hence assertion~(iv) follows from~(\ref{eq:CS}), which directly 
extends to gauge transformations that are trivial over the boundary.
This proves the lemma. 
\end{proof}

\begin{cor}\label{cor:CS}
Let $B_0\in\cL$ and $u:[0,1]\to\cG(\Sigma)$ 
with ${u(0)=u(1)=\one}$.  Then
$$
\int_0^1\int_\Sigma\winner{u(s)^*B_0}{\p_s(u(s)^*B_0)}\,\rd s
= 8\pi^2\deg(u).
$$
\end{cor}

\begin{proof}
The left hand side is twice the difference of the 
Chern-Simons functionals in Lemma~\ref{le:CS}~(ii).
\end{proof}

\subsection*{Perturbations}

We work with holonomy perturbations as 
in~\cite{T,F1,Donaldson book}. Let 
$
\D:=\{z\in\C \st |z|\leq 1\}
$
be the closed unit disc and identify $S^1$ with $\R/\Z$, with 
the real coordinate denoted by $\theta$. Choose embeddings 
$\gamma_i:S^1\times\D\hookrightarrow\mathrm{int}(Y)$ 
for $i=1,\ldots,N$ such that the $\gamma_i$ 
coincide on a neighbourhood of $\{0\}\times\D$.  
We denote by $\rho_i:\D\times\cA(Y) \to \rG$ the map that 
assigns to a pair $(z,A)$ the holonomy of the connection $A$ 
around the loop $[0,1]\to Y : \theta\mapsto \gamma_i(\theta,z)$.
Then the map ${\rho=(\rho_1,\ldots,\rho_N):\D\times\cA(Y) \to \rG^N}$ 
descends to a map between the quotient spaces 
$\D\times\cB(Y) \to \rG^N/\rG$, where the action
of $\rG$ on $\rG^N$ is by simultaneous conjugation
and $\cB(Y):=\cA(Y)/\cG(Y)$.

Now every smooth function $f:\D\times\rG^N\to\R$ that is invariant 
under conjugation and vanishes near the boundary
induces a gauge invariant perturbation 
${h_f:\cA(Y)\to\R}$ given by
$$
h_f (A) := \int_\D f(z,\rho(z,A))  \,\rd^2 z .
$$
The differential $\rd h_f (A) : \rT_A\cA(Y) \to \R$ 
has the form
\begin{equation}\label{Xf}
\rd h_f (A) \alpha = \int_Y \winner{X_f(A)}{\alpha},
\end{equation}
where $X_f:\cA(Y)\to\Om^2(Y,\cg)$ is a smooth map satisfying
\begin{equation}\label{Xf identities}
\rd_A X_f(A)=0,\quad
X_f(u^*A) = u^{-1}X_f(A)u,\quad 
\rd X_f(A) \rd_A\xi = [X_f(A),\xi]
\end{equation}
for $A\in\cA(Y)$, $u\in\cG(Y)$, $\xi\in\Om^0(Y,\cg)$.
This follows from the gauge invariance
of $h_f$ (see Appendix~\ref{app:Xf}).
Since $\rd X_f(A)$ is the Hessian 
of $h_f$ we have
\begin{equation} \label{X symmetry}
\int_Y \winner{\rd X_f(A) \alpha}{\beta} 
= \int_Y \winner{\rd X_f(A) \beta}{\alpha}.
\end{equation}
Moreover, $X_f(A)$ is 
supported in the union of the thickened loops 
${\gamma_i(S^1\times\D)}$
and hence in the interior of $Y$.

\subsection*{Critical points}

The critical points of the perturbed Chern--Simons functional 
$\CS_\cL + h_f$ are the solutions $A\in\cA(Y)$ of the equation
$$
F_A + X_f(A) = 0, \qquad A|_\Sigma\in\cL .
$$
Let $\Crit(\CS_\cL+h_f)$ denote the set of critical points and abbreviate
$$
\cR_f := \Crit(\CS_\cL+h_f)/\cG(Y) .
$$ 
Associated to every critical point $A\in\cA(Y,\cL)$ of $\CS_\cL+h_f$
is a twisted deRham complex
\begin{equation}\label{eq:Acomplex}
\Om^0(Y,\cg)
\;\stackrel{\rd_A}{\longrightarrow}\;
\Om^1_{\rT_A\cL}(Y,\cg)
\;\stackrel{\rd_A+\rd X_f(A)}{\longrightarrow}\;
\Om^2_0(Y,\cg)
\;\stackrel{\rd_A}{\longrightarrow}\;
\Om^3(Y,\cg),
\end{equation}
where
$$
\Om^1_{\rT_A\cL}(Y,\cg):=\left\{\alpha\in\Om^1(Y,\cg)\,|\,
\alpha|_\Sigma\in T_{A|_\Sigma}\cL\right\},
$$
$$
\Om^2_0(Y,\cg):=\left\{\tau\in\Om^2(Y,\cg)\,|\,
\tau|_\Sigma=0\right\}.
$$
The first operator in this complex is the infinitesimal action
of the gauge group, the second corresponds to the Hessian
of the Chern--Simons functional, and the third to the Bianchi
identity.  A critical point $A$ is called {\bf irreducible} if the 
cohomology group $H^0_A$ of (\ref{eq:Acomplex}) 
vanishes, i.e.~the operator
${\rd_A:\Om^0(Y,\cg)\to\Om^1(Y,\cg)}$ is injective. 
It is called {\bf nondegenerate}
if the cohomology group $H^1_{A,f}$ vanishes, 
i.e.~for every $\alpha\in\rT_A\cA(Y,\cL)$ we have
\begin{equation}\label{eq:nondeg}
\rd_A \alpha + \rd X_f (A) \alpha = 0
\quad\iff\quad
\alpha \in \im\rd_A .
\end{equation}
This nondegeneracy means that the Hessian of the Chern-Simons functional 
is nondegenerate on a local slice of the gauge action. 
In Section~\ref{sec:moduli} we will prove that for a 
generic perturbation every critical point is nondegenerate, i.e.\
$\CS_\cL + h_f$ induces a Morse function on the quotient $\cB(Y,\cL)$.

\subsection*{Gradient flow lines}

Fix a metric $g$ on $Y$.  Then a negative gradient flow line 
of the perturbed functional $\CS_\cL+h_f$ 
is a connection $\A\in\cA(\R\times Y)$ in temporal gauge, 
represented by a smooth path $\R\to\cA(Y): s \mapsto A(s)$ 
that satisfies the boundary value problem
\begin{equation}\label{eq:floer}
\pd_s A + *\bigl( F_A + X_f(A) \bigr) = 0, \qquad
A(s)|_\Sigma \in\cL \quad\forall s\in\R .
\end{equation}
The energy of a solution is
$$
E_f(\A) = \frac12 \int_{\R\times Y} 
\Bigl(\left|\pd_s A\right|^2 + \left|F_A + X_f(A)\right|^2\Bigr).
$$
In Section~\ref{sec:exp} we prove that (in the nondegenerate case) a solution 
$\A$ of (\ref{eq:floer}) has finite energy if and only if
there exist critical points $A^+,A^-\in\Crit(\CS_\cL+h_f)$ such that
$A(s)$ converges exponentially to $A^\pm$ as $s$ tends to $\pm\infty$.
Denote the moduli space of connecting trajectories from $[A^-]$ to $[A^+]$ by
$$
\cM(A^-,A^+) :=\left\{
\A\in\cA^{\rm tmp}(\R\times Y) \;\Bigg|
\begin{aligned}
&(\ref{eq:floer}),\; E_f(\A)<\infty, \\
&\lim_{s\to\pm\infty}A(s) \in [A^\pm] 
\end{aligned}
\right\} / {\cG(Y)} ,
$$
where $\cA^{\rm tmp}(\R\times Y)$ denotes the space of connections 
on $\R\times Y$ in temporal gauge.  
The analogue of equation~(\ref{eq:floer}) for connections 
$\A=\Phi\ds + A$ that are not in temporal gauge is 
\begin{equation}\label{eq:floerPhi}
\pd_s A - \rd_A\Phi + *\bigl( F_A + X_f(A) \bigr) = 0, \qquad
A(s)|_\Sigma \in\cL \quad\forall s\in\R .
\end{equation}
This equation can be written in the form
\begin{equation} \label{eq:asd}
F_\A + X_f(\A) + *( F_\A + X_f(\A) ) = 0, \qquad
\A|_{\{s\}\times\pd Y} \in\cL \quad\forall s\in\R ,
\end{equation}
where $X_f(\A)(s,y)=X_f(A(s))(y)$.  In this form it generalizes
to $4$-manifolds with a space time of the boundary
and tubular ends.

The moduli space $\cM(A^-,A^+)$ can also be described as the
quotient of the space of all finite energy solutions 
of~(\ref{eq:floerPhi}) in temporal gauge outside 
of a compact set that converge to $A^\pm$ as $s\to\pm\infty$.
In this case the gauge group consists of gauge transformations
that are independent of $s$ outside of a compact set
and preserve~$A^\pm$ at the ends.
The study of the moduli space is based on the analysis 
of the linearized operator 
for equation~(\ref{eq:floerPhi}). As a first step 
we examine the Hessian of the Chern--Simons functional. 


\section{The Hessian}\label{sec:HA}

In this section we establish the basic analytic properties
of the Hessian of the Chern-Simons functional
and draw some conclusions on the structure of the set of critical
points and the linearized operator of the gradient flow lines.

We continue the notation of Section~\ref{sec:CS}.
The augmented Hessian of the perturbed Chern--Simons functional 
at a connection $A\in\cA(Y,\cL)$ is the operator 
\begin{equation}\label{eq:HA}
\cH_A := \left( \begin{array}{cc} 
*\rd_A + *\rd X_f(A) & - \rd_A \\ -\rd_A^* & 0 
\end{array} \right).
\end{equation}
The additional terms $-\rd_A$ and $-\rd_A^*$ arise from a local 
slice condition. Think of $\cH_A$ as an unbounded operator on the 
Hilbert space
$
{L^2(Y,\rT^*Y\otimes\cg)  \times L^2(Y\otimes\cg)}
$
with dense domain 
$$
\mathrm{dom}\,\cH_A := \{
(\alpha,\phi)\in W^{1,2}(Y,\rT^*Y\otimes\cg) \times W^{1,2}(Y,\cg)|
*\alpha|_{\pd Y} = 0,\alpha|_{\pd Y}\in\rT_A\cL\}.
$$
Here we abbreviate $\rT_A\cL:=\rT_{A|_\Sigma}\cL$ for $A\in\cA(Y,\cL)$.

The operator $\cH_A$ is symmetric: for 
$\alpha,\beta\in\Om^1(Y,\cg)$ and $\phi,\psi\in\Om^0(Y,\cg)$
\begin{equation}\label{H symmetry}
\begin{split}
&
\inner{\cH_A (\alpha,\phi)}{(\beta,\psi)}_{L^2} 
-\inner{(\alpha,\phi)}{\cH_A(\beta,\psi)}_{L^2}  \\
&=
\int_Y \winner{( \rd_A\alpha +\rd X_f(A)\alpha - *\rd_A\phi)}{\beta}
+ \int_Y \winner{(\rd_A*\alpha)}{\psi} \\
&\quad
-\int_Y \winner{\alpha}{( \rd_A\beta + \rd X_f(A)\beta - *\rd_A\psi)}
-\int_Y \winner{\phi}{(\rd_A*\beta)} \\
&= 
\int_{\pd Y} \winner{\alpha}{\beta}
- \int_{\pd Y} \inner{\phi}{*\beta}
+\int_{\pd Y}\inner{*\alpha}{\psi}.
\end{split}
\end{equation}
If both $(\alpha,\phi)$ and $(\beta,\psi)$ belong to the domain of $\cH_A$, 
then the boundary conditions guarantee that the last three integrals vanish. 
In particular, ${\int_{\pd Y} \la \alpha\wedge\beta \ra}$ 
is the symplectic form on $\alpha|_{\pd Y},\beta|_{\pd Y}\in\rT_A\cL$.
An $L^2$-estimate for the Hessian is obtained from the following elementary calculation:
If $(\alpha,\phi)\in\mathrm{dom}\,\cH_A$ then
\begin{align*}
\left\| \cH_A(\alpha,\phi) \right\|_{L^2}^2
&=
\left\|*\rd_A\alpha-\rd_A\phi\right\|_{L^2}^2 
+ \left\|\rd_A^*\alpha\right\|_{L^2}^2 \\
&=
\left\|\rd_A\alpha\right\|_{L^2}^2
+ \left\|\rd_A\phi\right\|_{L^2}^2 
+ \left\|\rd_A^*\alpha\right\|_{L^2}^2
- 2\int_Y\winner{\alpha}{[F_A,\phi]} \\
&\ge 
\delta\left\|(\alpha,\phi)\right\|_{W^{1,2}}^2
- C\left\|(\alpha,\phi)\right\|_{L^2}^2 .
\end{align*}
Here the second equation follows from integration by parts.
The inequality, with suitable constants $\delta>0$ and $C$, 
follows from the Cauchy--Schwarz inequality and~\cite[Theorem~5.1]{W} with $p=2$.
The resulting estimate
$\|(\alpha,\phi)\|_{W^{1,2}}
\le \delta^{-1/2} \| \cH_A(\alpha,\phi) \|_{L^2}
+ (C/\delta)^{1/2}\|(\alpha,\phi)\|_{L^2} $
implies that $\cH_A$ has a finite dimensional kernel and a closed image.  
In Proposition (\ref{prop:HA}) below (which is the main result of this section)
we will identify the cokernel $(\im\cH_A)^\perp)$ with the kernel and thus
prove that the Hessian is a Fredholm operator and self-adjoint.
We moreover establish the estimate for the Hessian in general $W^{k,p}$-Sobolev spaces.
This will be used in the analysis of the linearized operator on $\R\times Y$
and for the exponential decay analysis.

\begin{proposition}\label{prop:HA}
{\bf (i)}  $\cH_A$ is a self-adjoint Fredholm operator.

\smallskip\noindent{\bf (ii)}
For every $A\in\cA(Y,\cL)$ and every integer $k\geq 0$ and every $p>1$
there exists a constant $C$ such that the following holds.
If $(\alpha,\phi)\in\mathrm{dom}\,\cH_A$ and 
$\cH_A(\alpha,\phi)$ is of class $W^{k,p}$,
then $(\alpha,\phi)$ is of class $W^{k+1,p}$ and
$$
\bigl\|(\alpha,\phi)\bigr\|_{W^{k+1,p}(Y)} 
\le C\bigl( \bigl\|\cH_A(\alpha,\phi)\bigr\|_{W^{k,p}(Y)} 
+ \bigl\|(\alpha,\phi)\bigr\|_{L^p(Y)} \bigr).
$$

\smallskip\noindent{\bf (iii)}
If $F_A+X_f(A)=0$ then
$
\ker\,\cH_A = H^1_{A,f} \times H^0_A,
$
where 
\begin{equation}\label{eq:OmA}
\begin{split}
H^0_A&
:=\ker\,\rd_A\subset\Om^0(Y,\cg), \\
H^1_{A,f}&
:=\ker\,(\rd_A+\rd X_f(A))\cap\ker\rd_A^*\subset\Om_A^1(Y,\cg), \\
\Om^1_A(Y,\cg)&
:=\bigl\{\alpha \in \Om^1(Y,\cg) \st *\alpha|_{\pd Y} = 0,\, 
\alpha|_{\pd Y} \in \rT_A\cL \bigr\}.
\end{split}
\end{equation}
\end{proposition}

\begin{dfn}\label{def:nondegenerate}
Let $A\in\cA(Y,\cL)$ be a critical point of the perturbed
Chern--Simons functional, i.e.\ $F_A+X_f(A)=0$.
The connection $A$ is called {\bf nondegenerate} if $H^1_{A,f}=0$;
it is called {\bf irreducible} if $H^0_A=0$. 
\end{dfn}

\begin{rmk}\label{rmk:nondegenerate old} \rm
{\bf (i)}
The vector spaces $H^0_A$ and $H^1_{A,f}$ in 
Proposition~\ref{prop:HA} are isomorphic to the first
two cohomology groups in the complex~(\ref{eq:Acomplex});
they are the spaces of harmonic representatives. Hence
a critical point $A\in\cA(Y,\cL)$ is nondegenerate in the sense of 
Definition~\ref{def:nondegenerate} if and only if 
it satisfies~(\ref{eq:nondeg}). 

\smallskip\noindent{\bf (ii)}
Hypothesis~(L3) says that $A=0$ is nondegenerate for the zero 
perturbation $f=0$.  Since the differential $\rd X_f(A)$ vanishes
at $A=0$ for every $f$ (see Appendix~\ref{app:Xf}) 
it follows that $A=0$ is nondegenerate for any perturbation. 
\end{rmk}

The proof of Proposition~\ref{prop:HA} requires some preparation.
First, we need to introduce norms for the boundary terms in the
upcoming estimates.
Let $p^*$ denote the dual exponent of $p$ 
given by $1/p + 1/p^* = 1$. We define the following norms 
(which strictly speaking depend on $Y$)
for a smooth function $\phi:\Sigma=\p Y\to\cg$
\begin{equation*}
\begin{split}
\Norm{\phi}_{bW^{1-1/p,p}(\Sigma)}
&:= \inf \bigl\{ \|\tilde\phi \|_{W^{1,p}(Y)} \,\big|\,
\tilde\phi|_\Sigma = \phi\bigr\}, \\
\Norm{\phi}_{bW^{-1/p,p}(\Sigma)}
&:= \sup_{0\ne\psi\in\Om^0(\Sigma,\cg)}
\frac{\Abs{\int_\Sigma\inner{\phi}{\psi}\dvol_\Sigma}}
{\Norm{\psi}_{bW^{1-1/p^*,p^*}(\Sigma)}}.
\end{split}
\end{equation*}
For a $2$-form $\tau\in\Om^2(\Sigma,\cg)$ the corresponding 
norms are understood as the norms 
of the function $*\tau\in\Om^0(\Sigma,\cg)$.
The following estimates for these boundary Sobolev norms will be useful.

\begin{lem}\label{le:norm}
For $A\in\cA(Y)$ and $\alpha\in\Om^1(Y,\cg)$ we have
$$
\Norm{\rd_{A|_\Sigma}(\alpha|_\Sigma)}_{bW^{-1/p,p}(\Sigma)}
=\sup_{0\ne\psi\in\Om^0(Y,\cg)} 
\frac{\Abs{\int_{Y} \bigl(\winner{\rd_A \alpha}{\rd_A\psi}
- \winner{\alpha}{[F_A,\psi]}\bigr)}}
{\Norm{\psi}_{W^{1,p^*}(Y)}}.
$$
Moreover, if $A\in\cA(Y,\cL)$ 
is a critical point of $\CS_\cL+h_f$
then
$$
\Norm{\rd_{A|_\Sigma}(\alpha|_\Sigma)}_{bW^{-1/p,p}(\Sigma)}
\le \bigl(1+\Norm{A}_{L^\infty(Y)}\bigr) 
\bigl\|\rd_A \alpha + \rd X_f(A)\alpha\bigr\|_{L^p(Y)}.
$$
\end{lem}

\begin{proof}
By definition we have 
$$
\Norm{\rd_{A|_\Sigma}(\alpha|_\Sigma)}_{bW^{-1/p,p}(\Sigma)}
= \sup_{\psi\ne0}
\frac{\Abs{\int_\Sigma
\inner{\rd_{A|_\Sigma}(\alpha|_\Sigma)}{\psi}}}
{\Norm{\psi}_{bW^{1-1/p^*,p^*}(\Sigma)}} 
= \sup_{\psi\ne0}
\frac{\Abs{\int_Y\rd\winner{\alpha}{\rd_A\psi}}}
{\Norm{\psi}_{W^{1,p^*}(Y)}},
$$
where the supremum runs over all nonzero functions
$\psi\in\Om^0(Y,\cg)$.
Now the first identity follows from
$\rd\winner{\alpha}{\rd_A\psi}=\winner{\rd_A \alpha}{\rd_A\psi}
- \winner{\alpha}{[F_A,\psi]}$.
If ${A\in\cA(Y,\cL)}$ is a critical point of $\CS_\cL+h_f$
then $F_A+X_f(A)=0$ and hence
\begin{align*}
\Norm{\rd_{A|_\Sigma}(\alpha|_\Sigma)}_{bW^{-1/p,p}(\Sigma)} 
&=\sup_{\psi\ne0} 
\frac{\Abs{\int_Y
\winner{(\rd_A \alpha+ \rd X_f(A) \alpha)}{\rd_A\psi}}}
{\Norm{\psi}_{W^{1,p^*}(Y)}} \\
&\le \bigl( 1+\Norm{A}_{L^\infty(Y)}\bigr) 
\bigl\|\rd_A \alpha + \rd X_f(A)\alpha\bigr\|_{L^p(Y)},
\end{align*}
where we have used~(\ref{Xf identities}) and~(\ref{X symmetry}).
This proves the lemma.
\end{proof}

The following lemma provides the basic estimates for Proposition~\ref{prop:HA}.
The first part is a regularity statement which goes a long way towards identifying
the dual domain of $\cH_A$ with its domain (thus establishing self-adjointness).
The second part is an estimate for the Hessian on pairs $(\alpha,\phi)$ that do
not necessarily satisfy the boundary conditions. This degree of generality is 
necessary since the Lagrangian boundary conditions are nonlinear, so 
differences in $\cA(Y,\cL)$ or derivatives of tangent vectors only satisfy
the boundary conditions up to some small curvature term.

\begin{lem} \label{lem:weak reg}
The following holds for every $p>1$ and every $A\in\cA(Y,\cL)$.

\smallskip\noindent{\bf (i)}
If $(\alpha,\phi)\in L^p(Y,\rT^*Y\otimes\cg) \times L^p(Y,\cg)$ 
and there is a constant $c$ such that
\begin{equation} \label{weak H}
\Abs{\int_Y  \inner{\alpha}{(*\rd_A\beta -\rd_A\psi )}
 - \int_Y \inner{\phi}{\rd_A^*\beta}} 
\leq c\Norm{(\beta,\psi)}_{L^{p^*}(Y)}
\end{equation}
for every $(\beta,\psi)\in \Om^1(Y,\cg) \times\Om^0(Y,\cg)$ 
with 
$\beta|_{\pd Y}\in\rd_{A|_{\Sigma}}\Om^0(\Sigma,\cg)$
and $*\beta|_{\pd Y} = 0$,  
then $(\alpha,\phi)\in W^{1,p}(Y,\rT^*Y\otimes\cg) \times W^{1,p}(Y,\cg)$ 
and it satisfies $*\alpha|_{\pd Y}=0$ and
$\rd_{A|_{\pd Y}}(\alpha|_{\pd Y})=0$ in the weak sense.

\smallskip\noindent{\bf (ii)}
There is a constant $C$ such that
\begin{align*}
\bigl\|(\alpha,\phi)\bigr\|_{W^{1,p}(Y)} 
&\leq C \Bigl( \bigl\|*\rd_A\alpha-\rd_A\phi\bigr\|_{L^{p}(Y)} 
+ \bigl\|\rd_A^*\alpha\bigr\|_{L^{p}(Y)} 
+  \bigl\|(\alpha,\phi)\bigr\|_{L^p(Y)} \\
&\qquad\quad
+ \bigl\|*\alpha|_\Sigma\bigr\|_{bW^{1-1/p,p}(\Sigma)} 
+ \Norm{\rd_{A|_\Sigma}(\alpha|_\Sigma)}_{bW^{-1/p,p}(\Sigma)} \Bigr) 
\end{align*}
for all $\alpha\in \Om^1(Y,\cg)$ and $\phi\in \Om^0(Y,\cg)$.
\end{lem} 

Before we prove this lemma let us draw a conclusion that will be useful
for the exponential decay analysis.

\begin{cor} \label{cor:nondegenerate}
Let $p>1$ and $A\in\cA(Y, \cL)$ be a nondegenerate critical point
of $\CS_\cL+h_f$. Then there is a constant $C$ such that
\begin{align*}
\bigl\|\alpha\bigr\|_{W^{1,p}(Y)} 
&\leq 
C \Bigl( \bigl\|\rd_A\alpha + \rd X_f(A)\alpha\bigr\|_{L^p(Y)} 
+ \bigl\|\rd_A^* \alpha\bigr\|_{L^p(Y)} \\
&\qquad\quad
+ \bigl\|*\alpha|_\Sigma\bigr\|_{bW^{1-1/p,p}(\Sigma)} 
+ \Norm{\Pi_A^\perp(\alpha|_\Sigma)}_{L^p(\Sigma)} \Bigr)
\end{align*}
for every $\alpha\in\Om^1(Y,\cg)$, where 
$\Pi_A^\perp:\Om^1(\Sigma,\cg)\to T_A\cL^\perp$ 
denotes the $L^2$ orthogonal projection onto the $L^2$
orthogonal complement of $\rT_A\cL$. 
\end{cor}

\begin{proof}
By Lemma~\ref{lem:weak reg}~(ii) with $\phi=0$ we have
\begin{align*}
\bigl\|\alpha\bigr\|_{W^{1,p}(Y)} 
&\le C \Bigl( \bigl\|\rd_A\alpha\bigr\|_{L^{p}(Y)} 
+ \bigl\|\rd_A^*\alpha\bigr\|_{L^{p}(Y)} 
+  \bigl\|\alpha\bigr\|_{L^p(Y)} \\
&\qquad\quad
+ \bigl\|*\alpha|_\Sigma\bigr\|_{bW^{1-1/p,p}(\Sigma)} 
+ \Norm{\rd_{A|_\Sigma}(\alpha|_\Sigma)}_{bW^{-1/p,p}(\Sigma)} 
\Bigr) \\
&\leq C'\Bigl( \bigl\|\rd_A\alpha+\rd X_f(A) \alpha\bigr\|_{L^{p}(Y)} 
+ \bigl\|\rd_A^*\alpha\bigr\|_{L^{p}(Y)}  \\
&\qquad\quad
+ \bigl\|*\alpha|_\Sigma\bigr\|_{bW^{1-1/p,p}(\Sigma)}  
+ \Norm{\Pi_A^\perp(\alpha|_\Sigma)}_{L^p(\Sigma)}
+ \bigl\|\alpha\bigr\|_{L^p(Y)}\Bigr).
\end{align*}
Here we have used the estimate 
$\Norm{\rd X_f(A)\alpha}_{L^p(Y)}\leq c\Norm{\alpha}_{L^p(Y)}$
of Proposition~\ref{prop:Xf}~(iv) and Lemma~\ref{le:norm}.
We added the term $\Norm{\Pi_A^\perp(\alpha|_\Sigma)}_{L^p(\Sigma)}$ 
on the right since 
$$
\Pi_A^\perp (\alpha|_\Sigma)=0\qquad\iff\qquad 
\alpha|_\Sigma\in\rT_A\cL
$$
and the restriction of the operator $\cH_A$ to the subspace 
$\{(\alpha,0)\}\subset{\rm dom}\,\cH_A$ is injective.  
Hence the operator 
$\alpha\mapsto \bigl(\rd_A\alpha+\rd X_f(A) \alpha,\rd_A^*\alpha,
*\alpha|_\Sigma,\Pi_A^\perp(\alpha|_\Sigma)\bigr)$ is injective
and it follows that the compact term 
$\Norm{\alpha}_{L^p(Y)}$ on the right
can be dropped.  This proves the corollary.
\end{proof}

\begin{proof}[Proof of Lemma~\ref{lem:weak reg}.]
It suffices to prove the lemma in the case $*A|_{\pd Y}=0$.
The general case can be reduced to this by a compact perturbation of the operator
(leaving the boundary conditions fixed).
To prove (i) consider a pair $(\alpha,\phi)\in L^p(Y,\rT^*Y\otimes\cg) \times L^p(Y,\cg)$ 
that satisfies (\ref{weak H}) with a constant $c$.
Let $\zeta \in\Om^0(Y,\cg)$ with $\tfrac{\pd \zeta }{\pd\nu}|_{\pd Y}=0$
and choose $(\beta,\psi)=(\rd_A \zeta ,0)$.
Then $*\beta|_{\pd Y}=0$ and 
$\beta|_{\pd Y}=\rd_{A|_{\Sigma}}(\zeta |_{\Sigma})$
and hence, by~(\ref{weak H}),

\begin{equation} \label{weak 0}
\begin{split}
\Abs{\int_Y \inner{\phi}{\laplace_A \zeta }}
&\le
c \,\Norm{\rd_A \zeta }_{L^{p^*}(Y)} 
+ \Abs{\int_Y \inner{\alpha}{* [F_A, \zeta  ]}}   \\
&\le 
\left(c + c \Norm{A}_{L^\infty(Y)} 
+ \Norm{F_A}_{L^\infty(Y)} \Norm{\alpha}_{L^p(Y)} 
\right)  \Norm{\zeta }_{W^{1,p^*}(Y)} 
\end{split}
\end{equation}
Hence it follows from the regularity theory for the Neumann problem 
(\cite{ADN} or e.g.~\cite[Theorem~2.3']{W}) 
that $\phi\in W^{1,p}(Y,\cg)$ and
\begin{equation}\label{phireg}
\Norm{\phi}_{W^{1,p}(Y)} \le 
C\bigl( c + \Norm{(\alpha,\phi)}_{L^p(Y)}\bigr),
\end{equation}
for a suitable constant $C=C(A)$.

Now fix a vector field $Z\in\Vect(Y)$ 
with $\Norm{Z}_{L^\infty(Y)}\leq 1$ that 
is perpendicular to $\pd Y$. 
Then it follows from~(\ref{weak H}) with
$\beta=0$ and $\psi=\cL_Z \zeta $ that 
\begin{equation} \label{weak 1}
\begin{split}
\Abs{\int_Y \inner{\alpha}{\rd(\cL_Z \zeta )}} 
&\le 
c\Norm{\cL_Z \zeta }_{L^{p^*}(Y)} 
+ \Abs{\int_Y \inner{\alpha}{[A,\cL_Z \zeta ]}} \\
&\le 
\left(c + \Norm{A}_{L^\infty(Y)}\Norm{\alpha}_{L^p(Y)}
\right) \Norm{\zeta }_{W^{1,p^*}(Y)}
\end{split}
\end{equation}
for every $\zeta \in\Om^0(Y,\cg)$.  Choosing $\psi=0$ and 
$\beta= *\bigl(\iota_Z g\wedge\rd \zeta \bigr)$ gives 
\begin{align} \label{weak 2}
&\Abs{\int_Y \inner{\alpha}{\rd^*(\iota_Z g\wedge\rd \zeta )}}
\nonumber \\
&\le 
c\,\Norm{\iota_Z g \wedge\rd \zeta }_{L^{p^*}(Y)} 
+ \Abs{\int_Y  \inner{\phi}{\rd_A (\iota_Z g\wedge\rd \zeta )}}
+ \Abs{\int_Y \inner{\alpha}{*[A\wedge*(\iota_Z g\wedge\rd \zeta )]}} 
\nonumber \\
&\le 
\left(c + C_Z\Norm{\phi}_{L^p(Y)} 
+ \Norm{A}_{L^\infty(Y)} \Norm{(\alpha,\phi)}_{L^p(Y)}
\right) \Norm{\zeta }_{W^{1,p^*}(Y)}  
\end{align}
for every $\zeta \in\Om^0(Y,\cg)$ with $\zeta |_{\pd Y}=0$,
where $C_Z:=\Norm{\rd\iota_Z g}_{L^\infty(Y)}$.
Here we have used~(\ref{weak H}) with
$*\beta|_{\pd Y}=0$ and $\beta|_{\pd Y}=0$.
Combining~(\ref{weak 1}) and~(\ref{weak 2}) we obtain 
the estimate 
$$
\Abs{\int_Y \inner{\alpha(Z)}{\laplace \zeta }} 
\le \left( 2c + C' \Norm{(\alpha,\phi)}_{L^p(Y)}\right) 
\Norm{\zeta }_{W^{1,p^*}(Y)} 
$$
for every $\zeta \in\Om^0(Y,\cg)$ with $\zeta |_{\pd Y}=0$ 
and a suitable constant constant $C'=C'(A,Z)$
(see~\cite[Theorem~5.3~(ii)]{W}).
This implies $\alpha(Z)\in W^{1,p}(Y,\cg)$ and
\begin{equation}\label{aZ reg}
\Norm{\alpha(Z)}_{W^{1,p}(Y)} 
\le C\bigl( c + \Norm{(\alpha,\phi)}_{L^p(Y)}\bigr),
\end{equation}
where the constant $C$ depends on $A$ and the vector field $Z$.
This proves the interior regularity of $\alpha$ as well 
as the regularity of its normal component. 
Moreover, partial integration now shows that,
for every $\zeta \in\Om^0(Y,\cg)$ with $\zeta |_{\pd Y}=0$, we have
$$
\Abs{\int_{\pd Y}\inner{\alpha(Z)}{\tfrac{\pd \zeta }{\pd\nu}}} 
\le \left(
2c + C'\Norm{(\alpha,\phi)}_{L^p(Y)} 
+ \Norm{\rd (\alpha(Z))}_{L^{p}(Y)}
\right) \Norm{\zeta }_{W^{1,p^*}(Y)}.
$$
In particular, we can fix any normal derivative 
$\tfrac{\pd \zeta }{\pd\nu}=g\in\Om^0(\pd Y,\cg)$ and 
find an admissible function $\zeta \in\Om^0(Y,\cg)$ with 
$\zeta |_{\pd Y}=0$ and $\Norm{\zeta }_{W^{1,p^*}(Y)}$ arbitrarily 
small.  Thus we have $\int_{\pd Y} \inner{\alpha(Z)}{g} =0$ 
for all $g\in\Om^0(\pd Y,\cg)$, and hence $\alpha(Z)=0$ 
for normal vector fields $Z$, i.e.~$*\alpha|_{\pd Y}=0$.

To deal with the tangential components near the boundary $\pd Y=\Sigma$ 
we use normal geodesics to identify a neighbourhood of the boundary
with $[0,\eps)\times\Sigma$ with the split metric $\rd t^2 + g_t$, 
where $(g_t)_{t\in[0,\eps)}$ is a smooth family of metrics on~$\Sigma$.
In this splitting we write 
$$
\alpha=\alpha_{\scriptscriptstyle\Sigma} + a\,\dt
$$ 
for ${\alpha_{\scriptscriptstyle\Sigma}
\in L^p([0,\eps)\times\Sigma,\rT^*\Sigma\otimes\cg)}$
and $a\in W^{1,p}([0,\eps)\times\Sigma,\cg)$.  Then 
$$
a|_{t=0}=0, \qquad
\Norm{a}_{W^{1,p}}\leq C 
\bigl( c + \Norm{(\alpha,\phi)}_{L^p(Y)}\bigr)
$$ 
by~(\ref{aZ reg}).   From now on $*$, $\rd$, and $\rd^*$ 
will denote the Hodge operator, the exterior derivative, 
and its adjoint on $\Sigma$. We abbreviate
$I:=[0,\eps)$ and denote by $\cC_0^\infty(I\times\Sigma)$ 
the space of functions with compact support in $(0,\eps)\times\Sigma$.
Then the inequality~(\ref{weak H}) can be rewritten as
\begin{multline*} 
\biggl|
\int_{I\times\Sigma} 
\inner{\alpha_{\scriptscriptstyle\Sigma}}
{\bigl( *\pd_t\beta_{\scriptscriptstyle\Sigma}-*\rd b+\rd\psi\bigr)} \\
\quad\;\;\; - \int_{I\times\Sigma} 
\inner{a}{\bigl(\pd_t\psi-*\rd\beta_{\scriptscriptstyle\Sigma}\bigr)} 
 + \int_{I\times\Sigma} 
\inner{\phi}{\bigl(\pd_t b-\rd^*\beta_{\scriptscriptstyle\Sigma}\bigr)} 
\biggr| 
\le c \bigl\| ( \beta_{\scriptscriptstyle\Sigma},b,\psi ) 
\bigr\|_{L^{p^*}(I\times\Sigma)}
\end{multline*}
for all $\beta_{\scriptscriptstyle\Sigma}
\in\cC_0^\infty(I\times\Sigma,\rT^*\Sigma\otimes\cg)$ 
and $b,\psi\in\cC_0^\infty(I\times\Sigma,\cg)$.
Partial integration in the terms involving $a$ and $\phi$ 
then yields
$$
\Abs{\int_{I\times\Sigma} 
\inner{\alpha_{\scriptscriptstyle\Sigma}}
{\bigl(\pd_t\beta_{\scriptscriptstyle\Sigma}-\rd b-*\rd\psi\bigr)}} 
\le
\bigl(c+\Norm{a}_{W^{1,p}}+\Norm{\phi}_{W^{1,p}}\bigr) 
\Norm{(\beta_{\scriptscriptstyle\Sigma},b,\psi)}_{L^{p^*}}.
$$
Since $\cC^\infty_0(I\times\Sigma)$ is dense in $L^{p^*}(I\times\Sigma)$ 
we obtain $\pd_t\alpha_{\scriptscriptstyle\Sigma}
\in L^{p}(I\times\Sigma,\rT^*\Sigma\otimes\cg)$ and 
$*\rd\alpha_{\scriptscriptstyle\Sigma},\rd^*\alpha_{\scriptscriptstyle\Sigma}
\in L^{p}(I\times\Sigma,\cg)$ with corresponding estimates.
Hence $\nabla_\Sigma\alpha_{\scriptscriptstyle\Sigma}$
is of class $L^p$ (see e.g.~\cite[Lemma~2.9]{W elliptic});
so $\alpha_{\scriptscriptstyle\Sigma}$ is of class $W^{1,p}$
and satisfies the estimate
$$
\Norm{\alpha_{\scriptscriptstyle\Sigma}}_{W^{1,p}}
\leq C \bigl(
c+\Norm{a}_{W^{1,p}}+\Norm{\phi}_{W^{1,p}} 
+ \Norm{\alpha_{\scriptscriptstyle\Sigma}}_{L^{p}}
\bigr) 
$$
with yet another constant $C$.
In combination with (\ref{phireg}) and (\ref{aZ reg}) this
proves the regularity claimed in (i) and the estimate
$$
\Norm{(\alpha,\phi)}_{W^{1,p}(Y)}
\leq C \bigl(  c + \Norm{(\alpha,\phi)}_{L^{p}(Y)} \bigr) .
$$
To prove the second boundary condition on $\alpha|_{\pd Y}$
we use partial integration in~(\ref{weak H}) to obtain
$$
\Abs{\int_\Sigma\winner{\alpha}{\beta}}
\le
\left(c + \Norm{\rd_A\alpha}_{L^p(Y)} 
+ \Norm{\rd_A \phi}_{L^p(Y)}\right)
\Norm{\beta}_{L^{p^*}(Y)}
$$
for every $\beta\in \Om^1(Y,\cg)$ with $*\beta|_\Sigma = 0$ 
and $\beta|_\Sigma\in\rd_{A|_\Sigma}\Om^0(\Sigma,\cg)$.
In particular, we can fix $\beta|_\Sigma=\rd_{A|_\Sigma}\xi$
for any $\xi\in\cC^\infty(\Sigma,\cg)$ and find
admissible $\beta\in\Om^1(Y,\cg)$ with $*\beta|_\Sigma=0$ 
and $\Norm{\beta}_{L^{p^*}(Y)}$ arbitrarily small.
Thus we have ${\int_\Sigma\winner{\alpha}{\rd_{A|_\Sigma}\xi} =0}$ 
for all $\xi\in\Om^0(\Sigma,\cg)$, that is 
$\rd_{A|_\Sigma} (\alpha|_\Sigma) = 0$ in the weak sense.
This proves~(i).

To prove (ii) let $(\alpha,\phi)\in\Om^1(Y,\cg)\times\Om^0(Y,\cg)$ be given
and choose $\gamma\in\Om^1(Y,\cg)$ such that 
$$
*\gamma|_\Sigma=*\alpha|_\Sigma,\qquad
\gamma|_\Sigma=0,\qquad
\Norm{\gamma}_{W^{1,p}(Y)}
\le2\|*\alpha|_\Sigma\|_{bW^{1-1/p,p}(\Sigma)},
$$
and denote 
$
\alpha':=\alpha-\gamma. 
$
There exists a constant $C_0=C_0(A)>0$ such that
$
\Norm{\cH_A(\gamma,0)}_{L^p(Y)} 
\le C_0\|*\alpha|_\Sigma\|_{bW^{1-1/p,p}(\Sigma)}
$
and hence
$$
\Norm{\cH_A(\alpha',\phi)}_{L^p(Y)} 
\le \Norm{\cH_A(\alpha,\phi)}_{L^p(Y)} 
+ C_0 \Norm{*\alpha|_\Sigma}_{bW^{1-1/p,p}(\Sigma)}
=: c.
$$
Then it follows from~(\ref{H symmetry}) that, for every
pair $(\beta,\psi)\in\Om^1(Y,\cg)\times\Om^0(Y,\cg)$
with $*\beta|_\Sigma = 0$, we have 
\begin{equation}\label{eq:Hweak}
\Abs{\inner{(\alpha',\phi)}{\cH_A(\beta,\psi)}} 
\le c\Norm{(\beta,\psi)}_{L^{p^*}(Y)}
+ \Abs{\int_\Sigma \winner{\alpha}{\beta}}.
\end{equation}
Let $\zeta \in\Om^0(Y,\cg)$ with $\frac{\pd \zeta }{\pd\nu}|_{\pd Y} = 0$
and choose $(\beta,\psi)=(\rd_A \zeta ,0)$.  Then, 
by Lemma~\ref{le:norm}, we have 
$$
\Abs{\int_\Sigma \winner{\alpha}{\rd_A \zeta }}
\le \Norm{\rd_{A|_\Sigma}(\alpha|_\Sigma)}_{bW^{-1/p,p}(\Sigma)}
\Norm{\zeta }_{W^{1,p^*}(Y)}
$$
and hence, by~(\ref{eq:Hweak}),
\begin{align*}
&\Abs{\int_Y \inner{\phi}{\laplace_A\zeta }} 
= \Abs{\inner{(\alpha',\phi)}{\cH_A(\rd_A\zeta ,0)}
- \int_Y\inner{\alpha'}{*[F_A,\zeta ]}} \\
&\le
c \Norm{\rd_A \zeta }_{L^{p^*}(Y)} 
+ \Abs{\int_\Sigma \winner{\alpha}{\rd_A \zeta }}  
+ \Norm{\alpha'}_{L^p(Y)}\Norm{F_A}_{L^\infty(Y)}
\Norm{\zeta }_{L^{p^*}(Y)} \\
&\le 
C\Bigl(
\bigl\|\cH_A(\alpha,\phi)\bigr\|_{L^p(Y)} 
+ \bigl\|*\alpha|_\Sigma\bigr\|_{bW^{1-1/p,p}(\Sigma)}  \\
&\qquad\qquad 
+ \bigl\|\rd_{A|_\Sigma}(\alpha|_\Sigma)\bigr\|_{bW^{-1/p,p}(\Sigma)} 
+ \bigl\|\alpha\bigr\|_{L^p(Y)} 
\Bigr)  \bigl\|\zeta\bigr\|_{W^{1,p^*}(Y)}
\end{align*}
for a suitable constant $C=C(A)$. 
(Compare this with~(\ref{weak 0}).)
As in the proof of~(i) this implies 
\begin{align*}
\Norm{\phi}_{W^{1,p}(Y)}
&\le C\Bigl( \bigl\|\cH_A(\alpha,\phi)\bigr\|_{L^p(Y)} 
+ \bigl\|*\alpha|_\Sigma\bigr\|_{bW^{1-1/p,p}(\Sigma)} \\
&\qquad\quad
+ \Norm{\rd_{A|_\Sigma}(\alpha|_\Sigma)}_{bW^{-1/p,p}(\Sigma)} 
+ \bigl\|(\alpha,\phi)\bigr\|_{L^p(Y)}
\Bigr)
\end{align*}
with a possibly larger constant $C$. 
(Compare this with~(\ref{phireg}).) 
To prove the same estimate for $\alpha'$ (and hence for $\alpha$)
one can repeat the argument in the proof of~(i),
because in this part of the argument the 
inequality~(\ref{eq:Hweak}) is only 
needed for $(\beta,\psi)$ with 
$*\beta|_\Sigma=0$ and $\beta|_\Sigma=0$.
This proves~(ii) and the lemma. 
\end{proof}

\begin{proof}[Proof of Proposition \ref{prop:HA}.]
We prove (ii) by induction. Observe that
\begin{equation}\label{Xf bound}
\bigl\| \rd X_f(A) \alpha \bigr\|_{W^{k,p}(Y)} 
\leq C\left\|\alpha\right\|_{W^{k,p}(Y)}  
\end{equation}
for all $\alpha\in\Om^1(Y,\cg)$ and a constant $C=C(A,f)$,
by Proposition~\ref{prop:Xf}~(iv).  
Hence it suffices to prove the estimate with $f=0$. 
For $k=0$ regularity holds by assumption and the estimate follows from 
Lemma~\ref{lem:weak reg}~(ii), using the fact that 
$\rd_{A|_\Sigma}\Om^0(\Sigma,\cg)\subset\rT_A\cL$, so
$\rd_{A|_\Sigma}(\alpha|_\Sigma)=0$.
(For $p=2$ an elementary proof of the estimate was given
at the beginning of the section.)
Thus we have proved~(ii) for $k=0$. 
It follows that $\cH_A$ has a finite dimensional kernel and a closed image.  

Now let $k\ge 1$ and suppose that~(ii) has been established for $k-1$. 
Let $(\alpha,\phi)\in{\rm dom}\,\cH_A$ and assume that
$\cH_A(\alpha,\phi)$ is of class $W^{k,p}$. By the induction 
hypothesis $(\alpha,\phi)$ is of class $W^{k,p}$ and
$$
\|(\alpha,\phi)\|_{W^{k,p}(Y)} 
\le C\bigl( \|\cH_A(\alpha,\phi)\|_{W^{k-1,p}(Y)} 
+ \|(\alpha,\phi)\|_{L^p(Y)} \bigr) .
$$
Let $X_1,\dots,X_k\in\Vect(Y)$.  Then, 
using the symmetry of $\cH_A$ and integration by parts, we obtain 
for every smooth pair 
$(\beta,\psi)\in\Om^1(Y) \times\Om^0(Y)$ with 
comact support in the interior of $Y$, we have
\begin{align*}
&\Abs{\inner{\cL_{X_1}\cdots\cL_{X_k}(\alpha,\phi)}{\cH_A(\beta,\psi)}}  \\
&= 
\Abs{\inner{(\alpha,\phi)}{\cL_{X_k}^*\cdots\cL_{X_1}^*\cH_A(\beta,\psi)}} \\
&\leq \bigl| \la (\alpha,\phi) \,,\, 
\cH_A \cL_{X_k}^*\cdots\cL_{X_1}^*(\beta,\psi) \ra \bigr| 
+ C_1 \|(\alpha,\phi)\|_{W^{k,p}(Y)}  \|(\beta,\psi)\|_{L^{p^*}(Y)}  \\
&= \bigl| \la \cL_{X_1}\cdots\cL_{X_k}
   \cH_A(\alpha,\phi) \,,\, (\beta,\psi) \ra \bigr| 
+ C_1 \|(\alpha,\phi)\|_{W^{k,p}(Y)}  \|(\beta,\psi)\|_{L^{p^*}(Y)}  \\
&\leq  C_2 \bigl(  \|\cH_A(\alpha,\phi)\|_{W^{k,p}(Y)} 
+ \|(\alpha,\phi)\|_{W^{k,p}(Y)} \bigr)  \|(\beta,\psi)\|_{L^{p^*}(Y)} 
\end{align*}
with uniform constants $C_i$. This estimate extends to the 
$W^{1,p^*}$-closure, so it holds for all $(\beta,\psi)$ 
with zero boundary conditions. However, in order 
to apply Lemma~\ref{lem:weak reg}~(i) to the pair 
$\cL_{X_1}\cdots\cL_{X_k}(\alpha,\phi)$ we would have 
to allow for more general test functions $(\beta,\psi)$. 
Unfortunately, this weak equation does
not extend directly, but we can still 
use the arguments of Lemma~\ref{lem:weak reg}.
For that purpose let the vector fields $X_1,\dots,X_k\in\Vect(Y)$ be
tangential to the boundary. Then the boundary condition
$*\alpha|_{\pd Y}=0$ will be preserved,
and the Lie derivatives $\cL_{X_i}$ in the following all have a dual
$\cL_{X_i}^*$ which does not include a boundary term.
To adapt the proof of Lemma~\ref{lem:weak reg}~(i) to
$\cL_{X_1}\cdots\cL_{X_k}(\alpha,\phi)$ instead of $(\alpha,\phi)$
we replace~(\ref{weak 0}) and~(\ref{weak 1}), 
which use test functions with nonzero boundary values.

Instead of~(\ref{weak 0}) we calculate
for all $\zeta \in\Om^0(Y,\cg)$ with $\tfrac{\pd \zeta }{\pd\nu}|_{\pd Y}=0$
and with a $W^{k,p}$-approximation $\Om^0(Y,\cg)\ni\phi_j\to\phi$
\begin{align*} 
& \bigl| \la \cL_{X_1}\dots\cL_{X_k}\phi \,,\, \laplace_A \zeta  \ra \bigr|
= \lim_{j\to\infty}
\bigl| \la \rd_A \cL_{X_1}\dots\cL_{X_k}\phi_j \,,\, \rd_A \zeta  \ra \bigr| \\
&\leq \lim_{j\to\infty}
\Bigl(\bigl| \la \cL_{X_1}\dots\cL_{X_k} \rd_A\phi_j \,,\, \rd_A \zeta  \ra \bigr| 
+ C_1 \|\phi_j \|_{W^{k,p}} \| \zeta  \|_{W^{1,p^*}}\Bigr) \\
&=
\bigl| \la \cL_{X_2}\dots\cL_{X_k} 
\rd_A \phi \,,\, \cL_{X_1}^*\rd_A \zeta  \ra \bigr| 
+ C_1 \|\phi \|_{W^{k,p}} \| \zeta  \|_{W^{1,p^*}} \\
&\leq
\bigl| \la \cL_{X_2}\dots\cL_{X_k} *\rd_A \alpha \,,\, 
\rd_A\cL_{X_1}^* \zeta  \ra \bigr| 
+ \bigl| \la \cL_{X_1}\dots\cL_{X_k} 
(*\rd_A\alpha - \rd_A\phi) \,,\, \rd_A \zeta  \ra \bigr| \\
&\quad
+ C_2 \|\phi \|_{W^{k,p}} \| \zeta  \|_{W^{1,p^*}} \\
&\leq
\bigl| \la *\rd_A \alpha \,,\, 
\rd_A \cL_{X_k}^*\dots\cL_{X_1}^* \zeta  \ra \bigr| 
+ \bigl| \la *\rd_A \alpha \,,\, 
\bigl[\cL_{X_k}^*\dots\cL_{X_2}^*,\rd_A\bigr]\cL_{X_1}^* \zeta  \ra \bigr| \\
&\quad + C_3 \bigl(  \|*\rd_A\alpha-\rd_A\phi \|_{W^{k,p}} 
+ \|\phi \|_{W^{k,p}} \bigr)  
\| \zeta  \|_{W^{1,p^*}} \\
&\leq C_4 \bigl(  \|\cH_A(\alpha,\phi) \|_{W^{k,p}} 
+ \|(\alpha,\phi) \|_{W^{k,p}} \bigr)  
\| \zeta  \|_{W^{1,p^*}}
\end{align*}
with uniform constants $C_i$.  Here the components of 
$\bigl[\cL_{X_k}^*\dots\cL_{X_2}^*,\rd_A\bigr]\cL_{X_1}^* \zeta $
are sums of derivatives of $\zeta $ including at most one normal derivative, 
so all but one derivative can be moved to the left hand side $*\rd_A\alpha$
by partial integration.  Moreover, we have used the fact that 
$\rd_A \cL_{X_k}^*\dots\cL_{X_1}^* \zeta |_{\pd Y}\in\rT_A\cL$ 
to obtain
\begin{equation*}
\begin{split}
\inner{*\rd_A \alpha}{\rd_A \cL_{X_k}^*\dots\cL_{X_1}^* \zeta }
&=
\inner{\alpha}{*[F_A ,\cL_{X_k}^*\dots\cL_{X_1}^*\zeta  }   \\
&=
\inner{\cL_{X_1}\dots\cL_{X_k}*[F_A\wedge\alpha]}{\zeta }
\end{split}
\end{equation*}
The last term can be estimated by 
$\Norm{\alpha}_{W^{k,p}} \Norm{\zeta }_{L^{p^*}}$.

Instead of~(\ref{weak 1}) we pick a 
$W^{k,p}$-approximation $\Om^1(Y,\cg)\ni\alpha_j\to\alpha$
satisfying the boundary condition $*\alpha_j|_{\pd Y}=0$ and hence
$*\cL_{X_1}\dots\cL_{X_k} \alpha_j|_{\pd Y}=0$.
Then we obtain for all $\zeta \in\Om^0(Y,\cg)$
\begin{align*}
&\left| \int_Y 
\la \cL_{X_1}\dots\cL_{X_k} \alpha , \rd(\cL_Z \zeta ) \ra \right| \\
&\leq\lim_{j\to\infty}
\left(
\left| 
\int_Y \la \cL_{X_1}\dots\cL_{X_k}\alpha_j,\rd_A(\cL_Z \zeta )\ra\right|  
+ C_1 \|\alpha_j\|_{W^{k,p}} \| \zeta \|_{W^{1,p^*}} 
\right) \\
&\leq\lim_{j\to\infty}
\left(
\left| \int_Y \la \cL_{X_1}\dots\cL_{X_k}\rd_A^*\alpha_j,\cL_Z \zeta \ra\right|  
+ C_2 \|\alpha_j\|_{W^{k,p}} \| \zeta \|_{W^{1,p^*}} 
\right) \\
&=
\left| \int_Y \la \cL_{X_2}\dots\cL_{X_k}\rd_A^*\alpha,
\cL_{X_1}^*\cL_Z \zeta  \ra \right|  
+ C_2 \|\alpha\|_{W^{k,p}} \| \zeta \|_{W^{1,p^*}} \\
&\leq
\|\cL_{X_1}\dots\cL_{X_k}\rd_A^*\alpha \|_{L^p} \| \cL_Z \zeta  \|_{L^{p^*}}  
+ C_2 \|\alpha\|_{W^{k,p}} \| \zeta \|_{W^{1,p^*}} \\
&\leq 
C_3 \bigl(  
\|\cH_A(\alpha,\phi) \|_{W^{k,p}} + \|(\alpha,\phi) \|_{W^{k,p}} 
\bigr) \|\zeta \|_{W^{1,p^*}(Y)}
\end{align*}
with uniform constants $C_i$.
Now the remaining arguments of Lemma~\ref{lem:weak reg}~(i) 
go through to prove the regularity 
$\cL_{X_1}\dots\cL_{X_k}(\alpha,\phi)\in W^{1,p}$ 
and the estimate
\begin{equation}\label{tangential}
\Norm{\cL_{X_1}\dots\cL_{X_k} (\alpha,\phi)}_{W^{1,p}}
\leq C \left(\Norm{\cH_A(\alpha,\phi)}_{W^{k,p}} 
+ \Norm{(\alpha,\phi)}_{W^{k,p}} \right)
\end{equation}
for the tangential derivatives and in the interior.
To control the normal derivatives near the boundary 
we use the same splitting as in Lemma~\ref{lem:weak reg}~(i). 
If $\cH_A(\alpha,\phi)\in W^{k,p}$
then this argument shows that 
\begin{align*} 
\pd_t\alpha_{\scriptscriptstyle\Sigma} &\in \rd a - *\rd\phi 
+ W^{k,p}(I\times\Sigma,\rT^*\Sigma\otimes\cg), \\
\pd_t a \,\; &\in \rd^*\alpha_{\scriptscriptstyle\Sigma}
+ W^{k,p}(I\times\Sigma,\cg), \\
\pd_t\phi \; &\in *\rd\beta_{\scriptscriptstyle\Sigma} 
+ W^{k,p}(I\times\Sigma,\cg).
\end{align*}
This can be used iteratively to replace the derivatives 
in~(\ref{tangential}) by normal derivatives.  
It then follows from the assumption $\cH_A(\alpha,\phi)\in W^{k,p}$ 
and the induction hypothesis $(\alpha,\phi)\in W^{k,p}$
that $(\alpha,\phi)\in W^{k+1,p}$ and
$$
\Norm{(\alpha,\phi)}_{W^{k+1,p}}
\leq C \left(\Norm{\cH_A(\alpha,\phi)}_{W^{k,p}} 
+ \Norm{(\alpha,\phi)}_{L^p} \right).
$$
This finishes the proof of~(ii).

We prove~(iii).  If $F_A+X_f(A)=0$ 
and $(\alpha,\phi)\in\ker\,\cH_A$, then the pair $(\alpha,\phi)$ 
is smooth by (ii). Integration by parts shows that 
$*\rd_A\alpha+*\rd X_f(A)\alpha$ is orthogonal to $\rd_A\phi$, 
hence both vanish, so the kernel has the required form.

To prove (i) we first show that the cokernel of $\cH_A$ agrees with its kernel.
Let $(\alpha,\phi)\in L^2(Y,\rT^*Y) \times L^2(Y)$ 
be orthogonal to the image of $\cH_A$. 
Denote by $\cH$ the operator of Lemma~\ref{lem:weak reg} 
for the perturbation $f=0$. Then 
\begin{align*}
\Inner{(\alpha,\phi)}{\cH(\beta,\psi)}_{L^2}
\;=\; -\Inner{\alpha}{*\rd X_f(A)\beta}_{L^2} 
\;\le\; c\left\|(\beta,\psi)\right\|_{L^2}
\end{align*}
for some constant $c$ and every pair 
$(\beta,\psi)\in\Om^1(Y,\cg)\times\Om^0(Y,\cg)$
satisfying the boundary conditions $*\beta|_{\pd Y}=0$
and $\beta|_{\pd Y}\in\rT_A\cL$. Hence it follows from 
Lemma~\ref{lem:weak reg}~(i) that $\alpha\in W^{1,2}(Y,\rT^*Y)$
and $\phi\in W^{1,2}(Y)$. 
So by (\ref{H symmetry}) 
\begin{align*}
0&= \int_Y \winner{( \rd_A\alpha +\rd X_f(A)\alpha - *\rd_A\phi)}{\beta}
+ \int_Y \winner{(\rd_A*\alpha)}{\psi} \\
&\quad
- \int_{\pd Y} \winner{\alpha}{\beta}
+ \int_{\pd Y} \inner{\phi}{*\beta}
- \int_{\pd Y}\inner{*\alpha}{\psi}
\end{align*}
for all $\beta\in\Om^1_A(Y,\cg)$ and $\psi\in\Om^0(Y,\cg)$.
(See~(\ref{eq:OmA}) for the definition of $\Om^1_A(Y,\cg)$.)
Taking $*\beta|_{\pd Y}=0$, $\beta|_{\pd Y}=0$, and $\psi|_{\pd Y}=0$
this implies
$$
*\rd_A\alpha +*\rd X_f(A)\alpha - \rd_A\phi=0,\qquad
\rd_A^*\alpha=0.
$$
Taking $(\beta,\psi)\in{\rm dom}\,\cH_A$ we then get  
$$
\int_{\pd Y} \winner{\alpha}{\beta} +\int_{\pd Y}\inner{*\alpha}{\psi} = 0 
$$
for every $\beta\in\Om^1_A(Y,\cg)$ and every $\psi\in\Om^0(Y,\cg)$. 
This (re-)proves $*\alpha|_{\pd Y}=0$ and, since $\beta|_{\pd Y}$ can
take any value in the Lagrangian subspace $\rT_A\cL$, it also
shows that $\alpha|_{\pd Y}\in\rT_A\cL$.  
Thus we have identified the cokernel of $\cH_A$ with its kernel. 
Since the kernel is finite dimensional, this proves that $\cH_A$
is a Fredholm operator.
Furthermore, every symmetric Fredholm operator with this
property is self-adjoint.  
(Let $x\in\dom\cH^*$, i.e.\ $\la x,\cH y \ra = \la z , y \ra$ for all $y\in\dom\cH$
and some $z$ in the target space. By assumption we can write $z=z_0+\cH x_1$ with
$z_0\in(\im\cH)^\perp$ and $x_1\in\dom\cH$.
Then, using symmetry, we have $\la x-x_1 , \cH y \ra = \la z_0,y\ra = 0$
for all $y\in\im\cH\cap\dom\cH$.
The latter is a complement of $\ker\cH\subset\dom\cH$ so we obtain 
$x-x_1\in(\im\cH)^\perp=\ker\cH\subset\dom\cH$ and hence $x\in\dom\cH$.)
This proves the proposition. 
\end{proof}

\subsection*{The set of critical points}

Using the properties of the Hessian we can now show finiteness 
of the set of gauge equivalence classes of critical points 
of the Chern-Simons functional, where the critical points are 
assumed to be nondegenerate.
More generally, we establish a compactness result that will be needed
to achieve nondegeneracy by a transversality construction.

\begin{proposition}\label{prop:finite}
Fix a Lagrangian submanifold $\cL\subset\cA(\Sigma)$
that satisfies (L1) and an integer $k\ge1$.
Let $f^\nu$ be a sequence of perturbations converging to
$f$ in the $\cC^{k+1}$ topology and $A^\nu\in\cA(Y,\cL)$ be 
a sequence of critical points of $\CS_\cL+h_{f^\nu}$.
Then there is a sequence of gauge transformations 
$u^\nu\in\cG(Y)$ such that $(u^\nu)^*A^\nu$ has a 
$\cC^k$ convergent subsequence. 

Moreover, if all the critical points of $\CS_\cL+h_f$ are
nondegenerate, then $\cR_f$ is a finite set.
\end{proposition}

\begin{proof}
Fix a constant $p>4$. 
The critical points of $\CS_\cL+h_{f^\nu}$ are $S^1$-invariant
solutions of the perturbed anti-self-duality equation
on $S^1\times Y$ and, by Proposition~\ref{prop:Xf}~(iii),
they satisfy a uniform $L^\infty$ bound on the curvature. 
Hence, by Uhlenbeck's weak compactness theorem
(see~\cite{U} or~\cite[Theorem~A]{W}), there is a sequence of gauge tranformations 
$u^\nu\in\cG(Y)$ such that $(u^\nu)^*A^\nu$ is bounded in $W^{1,p}$.
Passing to a subsequence, we may assume that 
$(u^\nu)^*A^\nu$ converges strongly in $\cC^0$ and 
weakly in $W^{1,p}$ to a connection $A\in\cA^{1,p}(Y,\cL)$.
The limit connection is a (weak) solution of $F_A+X_f(A)=0$
and hence, by~\cite[Theorem~A]{W elliptic}, is gauge equivalent
to a smooth solution.  Applying a further sequence of 
gauge transformation we may assume that $A$ is smooth and, 
by the local slice theorem (e.g.\ \cite[Theorem~F]{W}), that
\begin{equation}\label{eq:Ugauge}
\rd_A^*((u^\nu)^*A^\nu-A)=0,\qquad *((u^\nu)^*A^\nu-A)|_{\pd Y} = 0.
\end{equation}
It now follows by induction that $(u^\nu)^*A^\nu$ is uniformly
bounded in $W^{k+1,p}$.  Namely, if $(u^\nu)^*A^\nu$
is uniformly bounded in $W^{j,p}$ for any $j\in\{1,\dots,k\}$
then the curvature $F_{(u^\nu)^*A^\nu}=-X_{f^\nu}((u^\nu)^*A^\nu)$
is uniformly bounded in $W^{j,p}$, by Proposition~\ref{prop:Xf}~(iii), 
and hence $(u^\nu)^*A^\nu$ is uniformly bounded in $W^{j+1,p}$
by~\cite[Theorem~2.6]{W elliptic}.  
Since the Sobolev embedding $W^{k+1,p}\hookrightarrow \cC^k$ is compact,
the sequence $(u^\nu)^*A^\nu$ must have a $\cC^k$ convergent subsequence. 

To prove finiteness in the nondegenerate case it remains 
to show that nondegenerate critical points are isolated
in the quotient $\cA(Y,\cL)/\cG(Y)$. Thus let $A$ be a 
nondegenerate critical point and $A^\nu\in\cA(Y,\cL)$ 
be a sequence of critical points converging to $A$ in 
the $W^{1,p}$ topology (for some $p>2$).  
Then, by the local slice theorem, there exists a sequence 
of gauge transformations $u^\nu\in\cG(Y)$, 
converging to $\one$ in the $W^{2,p}$ topology, 
such that $(u^\nu)^*A^\nu$ satisfies~(\ref{eq:Ugauge}).
Since $\cA^{1,p}(Y,\cL)$ is a gauge invariant Banach submanifold
of $\cA^{1,p}(Y)$ it follows that the intersection with a local slice
gives rise to a Banach submanifold
$$
\cX_A := \left\{\alpha\in W^{1,p}(Y,\rT^*Y\otimes\cg)\,\bigg|\,
\begin{array}{c}
*\alpha|_\Sigma=0,\,\rd_A^*\alpha=0,\,
\Norm{\alpha}_{W^{1,p}}<\eps  \\
A+\alpha\in\cA^{1,p}(Y,\cL)
\end{array}
\right\}
$$
for $\eps>0$ sufficiently small. 
The tangent space of $\cX_A$ at $A$ is 
$$
\rT_A\cX_A =  \left\{\alpha\in W^{1,p}(Y,\rT^*Y\otimes\cg)\,\big|\,
*\alpha|_\Sigma=0,\, \alpha|_\Sigma\in\rT_A\cL,\,\rd_A^*\alpha=0
\right\}.
$$
Define the map 
$
\cF_A:\cX_A\times\left\{\phi\in W^{1,p}(Y,\cg)\,|\,\phi\perp\ker\,\rd_A\right\}
\to L^p(Y,\rT^*Y\otimes\cg)
$
by 
$$
\cF_A(\alpha,\phi) := *(F_{A+\alpha}+X_f(A+\alpha)) - \rd_A\phi.
$$
It has a zero at the origin, and we claim that its differential 
$$
\rd\cF_A(0,0)(\hat\alpha,\hat\phi) 
= *(\rd_A\hat\alpha+\rd X_f(A)\hat\alpha) - \rd_A\hat\phi 
$$
is bijective. The injectivity follows from the nondegeneracy
of $A$ and the fact that $\im\rd_A \perp \im *(\rd_A +\rd X_f(A))$.
To check the surjectivity notice that $\rd\cF_A(0,0)$ is the first
factor of the Hessian $\cH_A$. The Hessian is self-adjoint by 
Proposition~\ref{prop:HA} with cokernel 
$(\im\cH_A)^\perp=\ker\cH_A=H^1_{A,f}\times H^0_A$, so the cokernel
of $\rd\cF_A(0,0)$ is $H^1_{A,f}$, which vanishes by the nondegeneracy
assumption.
This proves that $\rd\cF_A(0,0)$ is bijective.
Since $(u^\nu)^*A^\nu-A\in\cX_A$ converges 
to zero in the $W^{1,p}$ norm and $\cF_A((u^\nu)^*A^\nu-A,0)=0$ 
for every $\nu$, it then follows from the inverse function theorem
that $(u^\nu)^*A^\nu=A$ for $\nu$ sufficiently large.  
This proves the proposition. 
\end{proof}

For nondegenerate critical points (i.e.\ $H^1_{A,f}=0$) we
have the following control on the kernel of the Hessian, $H^0_A=\ker\rd_A\subset\Om^0(Y,\cg)$, which measures reducibility.

\begin{rmk}\label{rmk:jump}\rm
The twisted cohomology groups $H^0_A$ form a 
vector bundle over the space of pairs $(f,A)$ 
with $A$ a nondegenerate critical point of $\CS_\cL+h_f$.
In particular, the dimension cannot jump.
This follows from the general fact that the cohomology groups $H^0$
form a vector bundle over the space of all chain complexes with $H^1=0$. 
To see this consider two chain complexes
$$
C^0\stackrel{\rd^0}{\to} C^1\stackrel{\rd^1}{\to}C^2 ,
\qquad\qquad
C^0\stackrel{\rd^0+P^0}{\longrightarrow} C^1\stackrel{\rd^1+P^1}{\longrightarrow}C^2
$$
of operators with closed images (between Hilbert spaces) and
assume that the first homology of the unperturbed complex vanishes, $H^1=\ker \rd^1/\im \rd^0=0$.
(Then the homology of the other complex, $H^1_P=\ker(\rd^1+P^1)/\im(\rd^0+P^0)$ 
also vanishes for sufficiently small perturbation $P$.)
Choose a complement $D^1\subset C^1$ of $\im\rd^0=\ker\rd^1$ 
and let $\Pi:C^1\to C^1/D^1$ be the projection.  
Then $\Pi\circ\rd^0:C^0\to C^1/D^1$ is surjective and 
the restriction $\rd^1|_{D^1}:D^1\to C^2$ is an injective operator with a closed image.
If $P^i:C^i\to C^{i+1}$ are sufficiently small then $\Pi\circ(\rd^0+P^0):C^0\to C^1/D^1$
is still surjective and ${(\rd^1+P^1)|_{D^1}:D^1\to C^2}$ is still injective.
From the latter and the identity $(\rd^1+P^1)\circ(\rd^0+P^0)=0$ 
it follows that $H^0_P=\ker(\rd^0+P^0)$ agrees with the kernel of the surjective map
$\Pi\circ(\rd^0+P^0)$.
Now let $D^0\subset C^0$ be a complement of $H^0=\ker\rd^0$, then 
$\Pi\circ\rd^0|_{D^0}:D^0\to C^1/D^1$ is bijective, and so is
$\Pi\circ(\rd^0+P^0)|_{D^0}:D^0\to C^1/D^1$ for sufficiently small $P^0$.
Its inverse is an injective map $I_P:C^1/D^1\to C^0$ with image $D^0$ that 
depends continuously on $P$ and satisfies $\Pi\circ(\rd^0+P^0) \circ I_P =\Id$.
Now $\pi_P:=I_P\circ \Pi\circ(\rd^0+P^0) : C^0 \to C^0$ is a 
projection, $\pi_P\circ\pi_P=\pi_P$, with 
$\ker\pi_P=\ker(\Pi\circ(\rd^0+P^0))=\im(1-\pi_P)$ and $\im\pi_P=\im I_P=D^0=\ker(1-\pi_P)$.
The opposite projection $1-\pi_P$ then provides an isomorphism 
$H^0=\ker\rd^0 \to \ker(\Pi\circ(\rd^0+P^0))=H^0_P$ that depends continuously on $P$.

\end{rmk}

\subsection*{The linearized operator on $\R\times Y$}

Next, we shall use the above results on the Hessian to establish
some basic properties of the linearized operator for~(\ref{eq:floer}).
Let $I\subset\R$ be an open interval and 
$\A=A+\Phi\ds\in\cA(I\times Y)$ such that
$A(s)|_{\pd Y}\in\cL$ for every $s\in I$.
A $\cg$-valued $1$-form on $I\times Y$ has the form $\alpha+\phi\ds$ with 
$\alpha(s)\in\Om^1(Y,\cg)$ and $\phi(s)\in\Om^0(Y,\cg)$.
Thus we shall identify $\Om^1(I\times Y,\cg)$ 
with the space of pairs $(\alpha,\phi)$
of smooth maps $\alpha:I\to\Om^1(Y,\cg)$ and $\phi:I\to\Om^0(Y,\cg)$.
For any integer $k\geq 1$ and any $p>1$ let 
$W^{k,p}_\A(I\times Y,\rT^*Y\otimes\cg)$ 
denote the space of $W^{k,p}$-regular $1$-forms 
$\alpha:I\times Y \to\rT^*Y\otimes\cg$ 
that satisfy the boundary conditions
\begin{equation}\label{bc}
*\alpha (s)|_{\pd Y}=0,\qquad 
\alpha(s)|_{\pd Y}\in\rT_{A(s)}\cL  
\end{equation}
for all $s\in I$. 
(The first equation arises from a gauge fixing condition.)

\begin{rmk}\label{rmk:hodge}\rm
The boundary conditions (\ref{bc}) are meaningful for every 
$\alpha$ of class $W^{1,p}$ with $p>1$. In this case we have
$\alpha(s)|_{\pd Y}\in L^p(\Sigma,\rT^*\Sigma\otimes\cg)$ 
for almost all $s\in I$, so there is a Hodge decomposition
$$
\alpha(s)|_{\pd Y}
=\alpha_0 + \rd_{A(s)|_\Sigma}\xi + *\rd_{A(s)|_\Sigma}\eta,
$$
and the second condition in (\ref{bc}) means that
$\eta=0$ and $\alpha_0\in\rT_{A(s)|_\Sigma}\cL$.
In other words, $\alpha(s)|_{\pd Y}$ lies in the  
$L^p$-closure of $\rT_{A(s)|_\Sigma}\cL$.
This $L^p$-closure is Lagrangian in the following sense:
If $\alpha\in L^p(\Sigma,\rT^*\Sigma\otimes\cg)$, then $\alpha$ 
lies in the  $L^p$-closure of $\rT_A\cL$ if and only if
$\int_\Sigma \la\alpha \wedge \beta\ra = 0$ for all smooth $\beta\in\rT_A\cL$.
(This extends the Lagrangian condition (\ref{eq:omega}) to nonsmooth tangent vectors.)
\end{rmk}

On a general $4$-manifold $X$, 
the linearized operator $\cD_\A$ for~(\ref{eq:asd})
with a gauge fixing condition has the form
$$
\Om^1(X,\cg)\to\Om^{2,+}(X,\cg)\times \Om^{0}(X,\cg):
\tilde\alpha \mapsto 
\left( (\rd_\A\tilde\alpha +\rd X_f(\A)\tilde\alpha)^+, 
-\rd_\A^*\tilde\alpha \right).
$$
In the case $X=I\times Y$ we identify $\Om^{2,+}(X,\cg)\times\Om^0(X,\cg)$
with the space of pairs of maps $I\to\Om^1(Y,\cg)$ and $I\to\Om^0(Y,\cg)$, 
using the formula 
$$
\tilde\alpha = \tfrac 12 \left(*\alpha(s) - \alpha(s)\wedge\ds \right)
$$
for self-dual $2$-forms on $I\times Y$. 
With this notation the linearized operator
\begin{multline*}
\cD_\A :
W^{k,p}_\A(I\times Y,\rT^*Y\otimes\cg) 
\times W^{k,p}(I\times Y,\cg) \\
 \to W^{k-1,p}(I\times Y,\rT^*Y\otimes\cg) 
\times W^{k-1,p}(I\times Y,\cg) 
\end{multline*}
for $I\times Y$ is given by
$$
\cD_\A := \nabla_s + \cH_{A(s)} ,
$$
where $\nabla_s := \pd_s + [\Phi,\cdot]$;
explicitly,
\begin{equation}\label{eq:DA}
\cD_\A \left(\begin{array}{c} \alpha \\ \phi \end{array} \right)=
\left(\begin{array}{c} \nabla_s \alpha 
+ *\rd_A\alpha +*\rd X_f(A)\alpha - \rd_A\phi \\
\nabla_s \phi - \rd_A^*\alpha \end{array} \right).
\end{equation}
Here we have dropped the argument $s$ in the notation, 
e.g.\ $\rd_A\phi$ stands for the path 
$s\mapsto \rd_{A(s)}\phi(s)$ of $\cg$-valued $1$-forms on $Y$.

\begin{rmk}\label{rmk:time reversal}\rm  
The formal adjoint operator has the form
$$
\cD_\A^* = -\nabla_s + \cH_{A(s)}.
$$
It is isomorphic to an operator of type 
$\nabla_s + \cH_{A}$ via time reversal.
Namely, if ${\sigma:(-I)\times Y\to I\times Y}$ denotes the
reflection in the $s$-coordinate, then
\begin{equation*}
\cD_\A^*(\beta,\psi) \comp\sigma 
= \cD_{\sigma^*\A}(\beta\comp\sigma,\psi\comp\sigma)
\end{equation*}
for every pair of smooth maps $\beta: I \to \Om^1(Y,\cg)$ 
and $\psi: I \to \Om^0(Y,\cg)$.
\end{rmk}

The following theorem provides the basic regularity (i) and estimate (ii) 
for the Fredholm theory of $\cD_\A$ and will also be needed to prove exponential decay.
The $L^p$-regularity has been established in \cite{W elliptic} by techniques 
that do not extend to $p=2$. Here we prove the $L^2$-regularity using the analytic
properties of the Hessian. A fundamental problem is that its domain varies with the connection,
unlike in the closed case. 
The variation will be controlled in step 1 of the proof, using a trivialization 
of the tangent bundle of $\cL$ in Appendix~\ref{app:Lag}.
This control then allows to apply the general theory of Appendix~\ref{app:spec}

\begin{thm}\label{thm:4reg}
For every integer $k\geq 0$, every $p>1$, 
and every compact subinterval $J\subset I$
there is a constant $C$ such that the following holds.

\smallskip\noindent{\bf(i)}
Assume $k=0$ and define $p^*:=p/(p-1)$.
Let
$$
(\alpha,\phi) \in L^p(I\times Y,\rT^*Y\otimes\cg) \times L^p(I\times Y,\cg)
$$
and suppose that there is a constant $c$ such that
\begin{equation}\label{eq:weak DA}
\left| \int_{I\times Y} 
\la \cD^*_\A (\beta,\psi) \,,\, (\alpha,\phi) \ra \right| 
\leq c \|(\beta,\psi)\|_{L^{p^*}(I\times Y)}
\end{equation}
for every compactly supported smooth map
$(\beta,\psi):I \to \Om^1(Y,\cg)\times\Om^0(Y,\cg)$ satisfying~(\ref{bc}).
Then $(\alpha,\phi)|_{J\times Y}$ is of class $W^{1,p}$ and
satisfies the boundary condition~(\ref{bc}) and the estimate
$$
\Norm{(\alpha,\phi)}_{W^{1,p}(J\times Y)} \leq
C \bigl( \Norm{\cD_\A(\alpha,\phi)}_{L^p(I\times Y)} 
+ \Norm{(\alpha,\phi)}_{L^p(I\times Y)} \bigr).
$$

\smallskip\noindent{\bf (ii)}
Assume $k\geq 1$. 
If $(\alpha,\phi) \in W^{1,p}(I\times Y,\rT^*Y\otimes\cg) 
\times W^{1,p}(I\times Y,\cg)$ satisfies~(\ref{bc}) 
and $\cD_\A(\alpha,\phi)$ is of class $W^{k,p}$, then 
$(\alpha,\phi)|_{J\times Y}$ is of class $W^{k+1,p}$ and
$$
\|(\alpha,\phi)\|_{W^{k+1,p}(J\times Y)} \leq
C \bigl( \|\cD_\A(\alpha,\phi)\|_{W^{k,p}(I\times Y)} +
            \|(\alpha,\phi)\|_{L^p(I\times Y)} \bigr).
$$
\end{thm}

\begin{proof}
Using the estimates on the perturbation $\rd X_f(A)$ 
in Proposition~\ref{prop:Xf}~(iv) we may assume without loss of generality
that $f=0$. Fix $s_0\in J$.  We prove the result for a neighbourhood 
of $s_0$ in four steps. 

\medskip\noindent{\bf Step~1.}
{\it 
After shrinking $I$, there exists a family of bijective linear operators 
$$
Q(s):\Om^1(Y,\cg)\times\Om^0(Y,\cg)\to\Om^1(Y,\cg)\times\Om^0(Y,\cg),
$$
parametrized by $s\in I$, such that the following holds.
\begin{description}
\item[(a)]
For every $s\in I$ and every $(\alpha,\phi)\in\Om^1(Y,\cg)\times\Om^0(Y,\cg)$ 
$$
(\alpha,\phi)\in{\rm dom}\,\cH_{A(s_0)}\qquad\iff\qquad
Q(s)(\alpha,\phi)\in{\rm dom}\,\cH_{A(s)}.
$$
\item[(b)]
For every integer $k\ge 0$ and every $p>1$
the operator family $Q$ induces a continuous linear operator from 
$W^{k,p}_{\rm loc}(I\times Y,\rT^*Y\otimes\cg)\times 
 W^{k,p}_{\rm loc}(I\times Y,\cg)$ to itself.
\end{description}}

\medskip\noindent
Let $\cU\subset\cA(Y,\cL)$ be a neighbourhood of $A(s_0)$
that is open in the $\cC^0$-topology and $\{Q_A\}_{A\in\cU}$ 
be an operator family which satisfies the requirements 
of Theorem~\ref{thm:Q}.  Shrink $I$ so that $A(s)\in\cU$ 
for every $s\in I$. Then the operators $Q(s):=Q_{A(s)}\times\Id$ 
satisfy the requirements of Step~1. 

\medskip\noindent{\bf Step~2.}
{\it We prove~(i) for $p=2$.}

\medskip\noindent
Abbreviate
$$
H:=L^2(Y,\rT^*Y\otimes\cg) \times L^2(Y,\cg)
$$
and let $W(s)\subset H$ be the subspace of
${(\alpha,\phi)\in W^{1,2}(Y,\rT^*Y\otimes\cg) \times W^{1,2}(Y,\cg)}$
that satisfy the boundary conditions
$$
*\alpha|_{\pd Y}=0,\qquad \alpha|_{\pd Y}\in\rT_{A(s)}\cL.
$$
Let $Q$ be as in Step~1, so each $Q(s)$ induces an operator on $H$ 
that descends to a Hilbert space isomorphism from $W(s_0)$ to $W(s)$. 
Then, by Proposition~\ref{prop:HA} with $p=2$, 
the operator family $\cH_{A(s)}:W(s)\to H$ satisfies the
conditions (W1-2) and (A1-2) in Appendix~\ref{app:spec}
for every compact subinterval of $I$.
Hence the estimate in~(i) with $p=2$ follows from Lemma~\ref{lem:spec est}
and a cutoff function argument,
and the regularity statement follows from Theorem~\ref{thm:spec reg}.

\medskip\noindent{\bf Step~3.}
{\it We prove~(i) for $p\ne 2$.}

\medskip\noindent
The result follows from~\cite[Theorem~C]{W elliptic}.
The intervals $I$ and $J$ can be replaced by $S^1$ by using cutoff functions,
and one can interchange $\cD_\A^*$ and $\cD_{\sigma^*\A}$ in (\ref{eq:weak DA})
by reversing time as in Remark~\ref{rmk:time reversal}.
Then \cite[Theorem~C~(iii)]{W elliptic} implies that 
$(\alpha,\phi)\comp\sigma$ is of class $W^{1,p}$ 
(with corresponding estimate). The same holds for 
$(\alpha,\phi)$, and partial integration 
as in~(\ref{H symmetry}) implies that
\begin{align*}
C\Norm{(\beta,\psi)}_{L^{p^*}} 
&\geq 
\left| \int_{I\times Y} \la \cD^*_\A (\beta,\psi) \,,\, (\alpha,\phi) \ra 
- \la (\beta,\psi) \,,\, \cD_\A (\alpha,\phi) \ra \right| \\
&=
\left| \int_{I\times\pd Y} \la \alpha \wedge \beta \ra 
      + \la *\alpha \,,\,\psi \ra \right| .
\end{align*}
Here we can choose any compactly supported 
$\beta|_{I\times\pd Y}:I\to\rT_A\cL\subset\Om^1(\pd Y,\cg)$
and $\psi|_{I\times\pd Y}:I\to\Om^0(\pd Y,\cg)$ and extend 
them to $I\times Y$ with $\|(\beta,\psi)\|_{L^{p^*}}$ arbitrarily small.
Thus the above estimate implies that $\alpha$ satisfies the boundary conditions
$\alpha(s)|_{\pd Y}\in\rT_{A(s)}\cL$ and $*\alpha(s)|_{\pd Y}=0$.

\medskip\noindent{\bf Step~4.}
{\it We prove~(ii).}

\medskip\noindent
The assertion of~(ii) continues to be meaningful for $k=0$;
we prove it by induction on $k$. For $k=0$ the regularity
statement holds by assumption and the estimate follows from~(i).
Fix an integer $k\ge 1$ and assume, by induction, that~(ii) has
been established with $k$ replaced by $k-1$.  Let 
$$
(\alpha,\phi) \in 
W^{1,p}(I\times Y,\rT^*Y\otimes\cg) \times W^{1,p}(I\times Y,\cg)
$$
such that~(\ref{bc}) holds and 
$$
(\beta,\psi):=\cD_\A(\alpha,\phi)\in
W^{k,p}(I\times Y,\rT^*Y\otimes\cg) \times W^{k,p}(I\times Y,\cg).
$$
Denote 
$$
(\alpha',\phi') := Q\pd_s(Q^{-1}(\alpha,\phi)),
$$
and
$$
(\beta',\psi'):=Q\Bigl(\pd_s(Q^{-1}(\beta,\psi))
-\bigl(\pd_s(Q^{-1}\cD_\A Q)\bigr) Q^{-1}(\alpha,\phi)\Bigr).
$$
Then  $(\alpha',\phi')$ satisfies the hypotheses of~(i) and hence
is of class $W^{1,p}$ and satisfies the boundary conditions~(\ref{bc}). 
Thus
$$
\cD_\A(\alpha',\phi') = (\beta',\psi')
$$
is of class $W^{k-1,p}$.  Hence, by the induction hypothesis,
$(\alpha',\phi')$ is of class $W^{k,p}$ and
\begin{align*}
\|(\alpha',\phi')\|_{W^{k,p}(J\times Y)}
&\le C_1\bigl(\|(\beta',\psi')\|_{W^{k-1,p}(I\times Y)}
+\|(\alpha',\phi')\|_{L^{p}(I\times Y)}\bigr) \\
&\le C_2\bigl(\|(\beta,\psi)\|_{W^{k,p}(I\times Y)}
+\|(\alpha,\phi)\|_{W^{k,p}(I\times Y)}\bigr).
\end{align*}
Since $(\alpha',\phi')=(\pd_s\alpha,\pd_s\phi)-(\pd_s Q)Q^{-1}(\alpha,\phi)$, 
this implies that $(\pd_s\alpha,\pd_s\phi)$ is of class $W^{k,p}$ and 
\begin{align*}
\|(\pd_s\alpha,\pd_s\phi)\|_{W^{k,p}(J\times Y)}
&\le C_3\bigl(\|\cD_\A(\alpha,\phi)\|_{W^{k,p}(I\times Y)}
+\|(\alpha,\phi)\|_{W^{k,p}(I\times Y)}\bigr) \\
&\le C_4\bigl(\|\cD_\A(\alpha,\phi)\|_{W^{k,p}(I\times Y)}
+\|(\alpha,\phi)\|_{L^p(I\times Y)}\bigr).
\end{align*}
It remains to establish regularity and estimates for $(\alpha,\phi)$
in $L^p(J,W^{k+1,p}(Y))$. To see it note that
$\cH_A(\alpha,\phi)=\cD_\A(\alpha,\phi)-\nabla_s(\alpha,\phi)$
is of class $L^p(J,W^{k,p}(Y))$. By Proposition~\ref{prop:HA},
$(\alpha(s),\phi(s))\in W^{k+1,p}(Y)$ for almost every $s\in J$ 
and 
\begin{align*}
&\left\|(\alpha,\phi)\right\|_{L^p(J,W^{k+1,p}(Y))}^p\\
&=\int_J\left\|(\alpha(s),\phi(s))\right\|_{W^{k+1,p}(Y)}^p\,\ds    \\
&\le C_4\int_J\left(
\left\|\cH_{A(s)}(\alpha(s),\phi(s))\right\|_{W^{k,p}(Y)}^p
+ \left\|(\alpha(s),\phi(s))\right\|_{L^{p}(Y)}^p
\right)\,\ds    \\
&\le C_5\left(\left\|\cD_\A(\alpha,\phi)\right\|_{W^{k,p}(I\times Y)}^p
+\left\|(\alpha,\phi)\right\|_{L^{p}(I\times Y)}^p\right).
\end{align*}
This completes the proof.
\end{proof}

\begin{rmk} {\rm \label{rmk:4reg}
The proof of Theorem~\ref{thm:4reg} carries over word for word
to the case where the metric and perturbation on $Y$ depend
smoothly on $s\in I$.
}\end{rmk}

We finish this section with a complete description of the linearized
operator for the trivial gradient flow line at an irreducible, nondegenerate
critical point.

\begin{thm}\label{thm:iso}
Let $A\in\cA(Y,\cL)$ be a critical point of the perturbed 
Chern--Simons functional $\CS_\cL+h_f$ such that
$H^0_A=0$ and $H^1_{A,f}=0$.  Then the operator
$$
\cD_A := \frac{\pd}{\pd s} +\cH_A
$$
on $L^p(\R\times Y,\rT^*Y\otimes\cg)\times L^p(\R\times Y,\cg)$
with domain 
\begin{multline*}
\mathrm{dom}\,\cD_A:=
\Bigl\{
(\alpha,\phi)\in W^{1,p}
(\R\times Y,\rT^*Y\otimes\cg)\times W^{1,p}(\R\times Y,\cg)\,\Big| \\
*\alpha(s)|_{\pd Y}=0,\,\alpha(s)|_{\pd Y}\in\rT_{A}\cL
\;\forall s\in\R
\Bigr\}
\end{multline*}
is a Banach space isomorphism for every $p>1$. 
\end{thm}

\begin{proof}
For $p=2$ it follows from~\cite[Theorem~A]{RS} 
and Proposition~\ref{prop:HA}  that $\cD_A$ is a Fredholm
operator of index zero; that it is bijective follows from
the inequality~(8) in~\cite{RS}.
Another argument is given in~\cite[Proposition~3.4]{Donaldson book};
it is based on the fact that $\cH_A$ is a bijective self-adjoint Fredholm 
operator, and on the local $L^2$-regularity (Theorem~\ref{thm:4reg}).
The case $p\ne 2$ can be reduced to the case $p=2$ by 
Donaldson's argument in~\cite[Proposition~3.21]{Donaldson book};
it uses in addition the local $L^p$-regularity in Theorem~\ref{thm:4reg}.
(For an adaptation of Donaldson's argument to the symplectic
case see~\cite[Lemma~2.4]{Sal}.)
\end{proof}


\section{Operators on the product $S^1\times Y$}\label{sec:S1Y}

In this section we study the anti-self-duality operator on $\SU(2)$-bundles
over the product $S^1\times Y$ with Lagrangian boundary conditions.
Our goal is, first, to establish a formula for the Fredholm index and,
second, to prove that the relevant determinant line bundle is orientable. 
Both results are proved with the same technique.  The problem 
can be reduced to the case of a suitable closed $3$-manifold 
$Y\cup_\Sigma Y'$ by means of an abstract argument involving 
the Gelfand--Robbin quotient.   

Throughout we fix a compact connected oriented $3$-manifold $Y$ 
with non\-empty boundary $\p Y=\Sigma$ and a 
gauge invariant, monotone Lagrangian submanifold 
$\cL\subset\cA(\Sigma)$ satisfying (L1-2) on page~\pageref{p:L1}. 
We identify $S^1\cong\R/\Z$.
Every gauge transformation $v:Y\to\rG=\SU(2)$ termines 
a principal $\SU(2)$-bundle
$P_v\to S^1\times Y$ defined by 
$$
P_v := \frac{\R\times Y\times\rG}{\Z},\qquad
[s,y,u]\equiv[s+1,y,v(y)u].
$$
A connection on $P_v$ with Lagrangian boundary conditions
is a pair of smooth maps $A:\R\to\cA(Y,\cL)$ and $\Phi:\R\to\Om^0(Y,\cg)$
satisfying 
\begin{equation}\label{eq:Av}
A(s+1) = v^*A(s),\qquad \Phi(s+1) = v^{-1}\Phi(s)v.
\end{equation}
The space of such connections will be denoted by $\cA(P_v,\cL)$
and we write $\A=\Phi\rd s+A$ or $(A,\Phi)$ for the elements 
of $\cA(P_v,\cL)$.  The space 
$$
\cA(S^1\times Y,\cL) := \bigl\{(v,\A)\,\big|\, \A\in\cA(P_v,\cL)\bigr\}
$$
is a groupoid.  
We will see that it has several connected components, corresponding
to $\pi_1(\cL/\cG_z(\Sigma))$ respectively the degree of $v:(Y,\pd Y)\to(\rG,\one)$.
A morphism from $(v_0,\A_0)$ 
to $(v_1,\A_1)$ is a smooth gauge transformation 
${u:\R\to\cG(Y)}$ on $\R\times Y$ satisfying 
\begin{equation}\label{eq:uv}
\begin{split}
v_1 &= u(s)^{-1}v_0u(s+1),\\
\A_1 &= u^*\A_0 .
\end{split}
\end{equation}
We  abbreviate~(\ref{eq:uv}) by $(v_1,\A_1)=:u^*(v_0,\A_0)$.
In the case $v_0=v_1=v$ a map $u$ that satisfies 
the first equation in (\ref{eq:uv}) 
is a gauge transformation on $P_v$.   
Since the gauge group $\cG(Y)$ is connected
there is, for every pair $v_0,v_1\in\cG(Y)$, a gauge transformation 
$u:\R\to\cG(Y)$ that satisfies the first equation in~(\ref{eq:uv}). 

Fix a perturbation $X_f$. 
Then every pair $(v,\A)=(v,A,\Phi)\in\cA(S^1\times Y,\cL)$
determines Sobolev spaces
\begin{equation*}
\begin{split}
W_v^{k,p}(S^1\times Y,\cg) 
&:= \bigl\{\phi\in W_\loc^{k,p}(\R\times Y,\cg)\,\big|\,
\phi(s+1)=v^{-1}\phi(s)v\bigr\},\\
W_v^{k,p}(S^1\times Y,\rT^*Y\otimes\cg) 
&:= \bigl\{\alpha\in W_\loc^{k,p}(\R\times Y,\cg)\,\big|\,
\alpha(s+1)=v^{-1}\alpha(s)v\bigr\},\\
W_{v,\A}^{k,p}(S^1\times Y,\rT^*Y\otimes\cg) 
&:= \bigl\{\alpha\in W_v^{k,p}(S^1\times Y,\rT^*Y\otimes\cg)\,\big|\,
(\ref{bc})\bigr\}
\end{split}
\end{equation*}
and an anti-self-duality operator
\begin{multline*}
\cD_{v,\A} :
W^{k,p}_{v,\A}(S^1\times Y,\rT^*Y\otimes\cg) 
\times W_v^{k,p}(S^1\times Y,\cg) \\
 \to W_v^{k-1,p}(S^1\times Y,\rT^*Y\otimes\cg) 
\times W_v^{k-1,p}(S^1\times Y,\cg) 
\end{multline*}
given by
$
\cD_{v,\A} := \nabla_s + \cH_{A(s)}
$
respectively by~(\ref{eq:DA}) as in Section~\ref{sec:HA}.

\begin{dfn}\label{def:deg}
The {\bf degree} of a pair $(v,\A)=(v,A,\Phi)\in\cA(S^1\times Y,\cL)$
is the integer 
\begin{equation*}
\begin{split}
\deg(v,\A) 
&:= - \frac{1}{4\pi^2}\int_0^1\int_Y\winner{F_A}{\p_sA}\,\rd s  
\end{split}
\end{equation*}
\end{dfn}

\begin{rmk}\label{rmk:degree}\rm
{\bf (i)}
The degree is an integer because it is the difference of the
Chern-Simons functionals.  Explicitly,
\begin{align*}
\deg(v,\A) &= 
- \frac{1}{8\pi^2} \int_0^1\int_\Sigma \winner{A}{\pd_s A} \ds  \\
&\quad
- \frac{1}{8\pi^2} \biggl[ \int_Y \Bigl( \winner{A}{\rd A} 
+ \frac 13 \winner{A}{[A\wedge A]}\Bigr) \biggr]^{s=1}_{s=0} \\
&= \frac{1}{4\pi^2} \Bigl( \CS(A(0),A|_{\Sigma}\# B) - \CS(A(1),B) \Bigr) \\
&\equiv \frac{1}{4\pi^2} \Bigl( \CS_\cL([A(0)]) - \CS_\cL([A(0)]) \Bigr) \;=\;0 \;\in\;\R/\Z .
\end{align*}
Here $B:[0,1]\to\cL$ is a smooth path from $B(0)=A(1)|_\Sigma$ to $B(1)=0$
and $A|_{\Sigma}\# B$ is the catenation of $A|_\Sigma:[0,1]\to\cL$
with $B$.

\smallskip\noindent{\bf (ii)}
If $v|_\Sigma\equiv\one$ then $A(s+1)|_\Sigma=A(s)|_\Sigma$
and, by (\ref{eq:CS}),
$$
\deg(v,\A) = \deg(v) 
- \frac{1}{8\pi^2}\int_0^1\int_\Sigma\winner{A}{\p_sA}\,\rd s.
$$
The last term is the symplectic action 
of the loop $\R/\Z\to\cL:s\mapsto A(s)|_\Sigma$,
multiplied by the factor $1/4\pi^2$. 

\smallskip\noindent{\bf (iii)}
If $v\equiv\one$ and $A(s)|_\Sigma=u(s)^*A(0)|_\Sigma$ with
$u(s+1)=u(s)\in\cG(\Sigma)$ then $\deg(v,\A)$
is minus the degree of the map $u:S^1\times\Sigma\to\rG$,
see Corollary~\ref{cor:CS}.
\end{rmk}

\begin{thm}\label{thm:S1Y}
Fix $p>1$ and an integer $k\ge 1$, then the following holds.

\smallskip\noindent{\bf (i)}
Two pairs $(v,\A),(v',\A')\in\cA(S^1\times Y,\cL)$ belong to 
the same component of $\cA(S^1\times Y,\cL)$ if and only
if they have the same degree. 

\smallskip\noindent{\bf (ii)}
For every pair $(v,\A)\in\cA(S^1\times Y,\cL)$ the operator
$\cD_{v,\A}$ is Fredholm and
$$
\mathrm{index}(\cD_{v,\A}) = 8\deg(v,\A).
$$

\smallskip\noindent{\bf (iii)}
The determinant line bundle $\det\to\cA(S^1\times Y,\cL)$ with fibers
$\det(\cD_{v,\A})$ is orientable. 

\smallskip\noindent{\bf (iv)}
Let $u:\R\to\cG(Y)$ be a morphism from 
$(v,\A)$ to $(v',\A')=(u^*v,u^*\A)$. 
Then $(v,\A)$ and $(u^*v,u^*\A)$ have the same 
degree and the induced isomorphism
$$
u^*:\det(\cD_{v,\A})\to\det(\cD_{(u^*v,u^*\A)})
$$
is orientation preserving (i.e.\ the map on orientations agrees
with the one induced by a homotopy).
\end{thm}

The proof of (ii) will be based on an identification of the index with
the spectral flow of the Hessian. 
Both the index and orientation results in (ii)-(iv) require a description
of the space of self-adjoint boundary conditions for the Hessian on a
pair of domains with matching boundary. We will use it to 
homotop from Lagrangian boundary conditions to the diagonal 
(representing the closed case). More precisely, 
we will use the abstract setting of Appendix~\ref{app:GR}.

We think of the div-grad-curl operator on $Y$ as an unbounded operator 
$$
D:=\left( \begin{array}{cc} 
*\rd & - \rd \\ -\rd^* & 0 
\end{array} \right) \;:\; W_0 \to H 
$$ 
on the Hilbert space
$$
H:= L^2(Y,\rT^*Y\otimes\cg)\oplus L^2(Y,\cg)
$$
with the dense domain 
$$
\mathrm{dom}\,D := W_0 := W^{1,2}_0(Y,\rT^*Y\otimes\cg)\oplus W^{1,2}_0(Y,\cg) .
$$
With this domain $D$ is symmetric and injective and has a closed image, 
see Lemma~\ref{le:HA} below.
Hence $D$ satisfies the assumptions of Appendix~\ref{app:GR} and thus
defines a symplectic Hilbert space, the Gelfand--Robbin quotient
$$
V:=\dom D^*/\dom D =  W / W_0 ,
\qquad 
\omega(\xi,\eta):=\inner{D^*\xi}{\eta}-\inner{\xi}{D^*\eta},
$$
where $W:=\dom D^*$ is the domain of the adjoint operator $D^*$.
The crucial property of the Gelfand--Robbin quotient is the fact that
self-adjoint extensions of $D$ are in one-to-one
correspondence with Lagrangian subspaces of $V$.

If $A\in\cA(Y)$ is a smooth connection on $Y$ 
then the restricted (unperturbed) Hessian $\cH_A|_{W_0}:W_0\to H$ 
is an unbounded operator on $H$ with domain~$W_0$.
It is a compact perturbation of the div-grad-curl
operator $D$.  The next lemma shows how these operators fit
into the setting of Appendix~\ref{app:GR}.

\begin{lemma}\label{le:HA}
\smallskip\noindent{\bf (i)}
For every smooth connection $A\in\cA(Y)$ 
on $Y$ the operator $\cH_A|_{W_0}:W_0\to H$ is symmetric,
injective, and has a closed image. Its domain $W_0$
is dense in $H$, the graph norm of $\cH_A$ on $W_0$ 
is equivalent to the $W^{1,2}$-norm, 
and the inclusion $W_0\to H$ is compact.

\smallskip\noindent{\bf (ii)}
For every ${A\in\cA(Y)}$ the domain of the dual operator 
$(\cH_A|_{W_0})^*$ is equal to~$W$ and
the symplectic form on the quotient $W/W_0$ 
is given by 
$$
\om(\xi,\eta) = 
\int_{\pd Y} \winner{\alpha}{\beta}
- \int_{\pd Y} \inner{\phi}{*\beta}
+\int_{\pd Y}\inner{*\alpha}{\psi}.
$$
for smooth elements $\xi=(\alpha,\phi)$ 
and $\eta=(\beta,\psi)$ in $W$.

\smallskip\noindent{\bf (iii)}
The kernel of $(\cH_A|_{W_0})^*$ determines a Lagrangian subspace 
$$
\Lambda_0(A) := \frac{\ker(\cH_A|_{W_0})^*+W_0}{W_0}\subset V.
$$
If two connections $A,A'\in\cA(Y)$ coincide in a neighbourhood 
of the boundary $\pd Y$ then $\Lambda_0(A')$ 
is a compact perturbation of $\Lambda_0(A)$.
\end{lemma}

\begin{proof}
The operator $\cH_A|_{W_0}$ is symmetric by~(\ref{H symmetry})
and it has a closed image by Lemma~\ref{lem:weak reg}~(ii).
To prove that it is injective let $(\alpha,\phi)\in\ker\cH_A\cap W_0$.
Extend $A$ to an $S^1$-invariant connection $\Xi$ on $S^1\times Y$
and $(\alpha,\phi)$ to an $S^1$-invariant $1$-form 
$\xi=\alpha+\phi\ds$ on $S^1\times Y$.
Then $\rd_\Xi^+\xi=0$, $\rd_\Xi^*\xi=0$, and $\xi$ vanishes 
on the (nonempty) boundary.  Near the boundary we choose 
coordinates $(s,t,z)\in S^1\times[0,\eps)\times\Sigma$
so that $(t,z)$ are normal geodesic coordinates on $Y$. 
Interchanging $s$ and $t$ we can first bring $\Xi$ into
temporal gauge with respect to $t$ and then use 
Lemma~\ref{le:ucon-loc}~(ii) to deduce that $\xi$ 
vanishes near the boundary. 
Since $Y$ is connected it follows from an open and closed 
argument that $\xi$ vanishes identically. 
The graph norm of $\cH_A$ on $W_0$ is given 
by~\eqref{eq:graphnorm} below.
The boundary term vanishes on $W_0$ and hence
this norm is equivalent to the $W^{1,2}$ norm.
The compactness of the inclusion $W_0\to H$ 
follows from Rellich's theorem. This proves~(i).

The domain of the dual operator and the symplectic form are
independent of $A$ because the difference 
$\cH_A|_{W_0}-D=(\cH_A-\cH_0)|_{W_0} :W_0\to H$ 
extends to a bounded self-adjoint operator from $H$ to itself.
The formula for the symplectic form follows 
from~\eqref{H symmetry}.

Assertion~(iii) follows from Lemma~\ref{le:DPL}.
This uses the fact that the difference operator 
$\Delta:=(\cH_A|_{W_0})^*-(\cH_{A'}|_{W_0})^* : W \to H$
is compact since it coincides with $\Delta\circ \iota \circ\Psi$. 
Here $\Psi:W\to W_0$ is a bounded map, given 
by multiplication with a cutoff function $\psi\in\cC^\infty_0(Y,[0,1])$, 
$\psi|_{\supp(A-A')}\equiv 1$, the inclusion $\iota:W_0\to H$ 
is compact by~(i),
and $\Delta:H\to H$ is bounded. 
This proves the lemma.
\end{proof}

\begin{rmk} \label{hurrah} \rm
\smallskip\noindent{\bf (i)}
The symplectic Hilbert space $(V,\omega)$ can be viewed 
as a space of boundary data for the Hessian, containing the space
$$
\frac{\Omega^1(Y,\cg)\times 
\Omega^0(Y,\cg)}{W_0\cap(\Omega^1(Y,\cg)\times \Omega^0(Y,\cg))}
\cong\Omega^1(\Sigma,\cg)\times \Omega^0(\Sigma,\cg)
\times \Omega^0(\Sigma,\cg)
$$
of smooth boundary data as a dense subspace; 
see Lemma~\ref{hurray} below.
The isomorphism is by $[(\alpha,\phi)]\mapsto 
(\alpha|_{\partial Y}, \phi|_{\partial Y},
*_{\scriptscriptstyle\Sigma}(*\alpha|_{\partial Y}) )$.
In this notation, an explicit formula for the symplectic 
form is given in Lemma~\ref{le:HA}~(ii).

\smallskip\noindent{\bf (ii)}
The space $\Omega^1(Y,\cg)\times\Omega^0(Y,\cg)$ 
of smooth pairs $(\alpha,\phi)$ is contained in 
the domain of the dual operator, and the restriction 
of $(\cH_A|_{W_0})^*$ to this subspace agrees with $\cH_A$.
The graph norm on 
$\Omega^1(Y,\cg)\times\Omega^0(Y,\cg)\subset\dom(\cH_A|_{W_0})^*$ is
\begin{equation} \label{eq:graphnorm} 
\begin{split}
\|(\alpha,\phi)\|^2_{(\cH_A|_{W_0})^*} &= 
\|(\alpha,\phi)\|^2_{L^2(Y)} + \|\rd_A\alpha\|^2_{L^2(Y)} 
+ \|\rd_A^*\alpha\|^2_{L^2(Y)} \\
&\quad 
+ \|\rd_A\phi\|^2_{L^2(Y)} 
+ 2\int_Y \inner{\phi}{[F_A,\alpha]} 
- 2 \int_{\pd Y} \inner{\phi}{\rd_A\alpha} .
\end{split}
\end{equation}
The dual domain $W=\dom(\cH_A|_{W_0})^*$ is the completion of
$\Omega^1(Y,\cg)\times\Omega^0(Y,\cg)$ with respect to this norm.
It is bounded by the $W^{1,2}$-norm and hence
$$
W^{1,2}_Y:= W^{1,2}(Y,\rT^*Y\otimes\cg)\oplus W^{1,2}(Y,\cg)
\subset W.
$$
Moreover, it follows from interior elliptic regularity that every 
element of the dual domain $W$ is of class $W^{1,2}$ on every compact 
subset of the interior of $Y$.  However, $W$ is not contained 
in $W^{1,2}_Y$, see Lemma~\ref{hurray} below.
\end{rmk}

The next Lemma gives a precise description for the spaces $W$ 
and $V$, including some parts of weak regularity. 
However, our theory does not depend 
on the explicit description of these spaces.
In our applications we only use the fact that the 
Gelfand--Robbin quotient is independent of the connection,
see Lemma~\ref{le:HA}. In the following we slightly abuse
notation and identify the Gelfand--Robbin quotient $V=W/W_0$ 
with the orthogonal complement of $W_0$ in $W$ in the graph norm 
of $D^*$.  Remark~\ref{rmk:V}~(ii) shows that it is given by
$$
V = \left\{\xi\in\dom D^*\,|\,
D^*\xi\in\dom D^*,\,D^*D^*\xi+\xi=0\right\}.
$$

\begin{lem}\label{hurray}
\smallskip\noindent{\bf (i)}
The space $V$ admits an orthogonal Lagrangian splitting
$$
V = \Lambda_0\oplus\Lambda_1,\qquad
\Lambda_0 := D^*\Lambda_1,\qquad
\Lambda_1 := V\cap\im D,
$$
where $\Lambda_0$ is the orthogonal projection of 
the kernel of $D^*$ onto $V$. 

\smallskip\noindent{\bf (ii)}
The space $W$ admits an orthogonal splitting
$W = W_0\oplus\Lambda_0\oplus\Lambda_1$,
where $W_0$ and $\Lambda_1$ are closed subspaces 
of $W^{1,2}_Y$ and $\Lambda_0$ is a closed subspace 
of~$H=:L^2_Y$. 

\smallskip\noindent{\bf (iii)}
The spaces of smooth elements are dense 
in $\Lambda_0$, $\Lambda_1$, $V$, and $W$
(with respect to the graph norm of $D^*$). 
The restriction map 
\begin{equation}\label{restriction}
\xi=(\alpha,\phi)\mapsto
\xi|_\Sigma := 
(\alpha|_\Sigma,\phi|_\Sigma, 
*_{\scriptscriptstyle\Sigma}(*\alpha|_{\Sigma})  )
\end{equation}
on the smooth elements extends continuously 
to $\Lambda_0$ and $\Lambda_1$.  This gives rise to 
injective operators
$$
\Lambda_0\to W^{-1/2,2}_\Sigma,\qquad
\Lambda_1\to W^{1/2,2}_\Sigma
$$
with closed images. Here we denote
$W^{-1/2,2}_\Sigma:=\bigl(W^{1/2,2}_\Sigma\bigr)^*$
and
$$
W^{1/2,2}_\Sigma:= 
W^{1/2,2}(\Sigma,\rT^*\Sigma\otimes\cg)\oplus W^{1/2,2}(\Sigma,\cg)
\oplus W^{1/2,2}(\Sigma,\cg).
$$
\end{lem}

\begin{proof}
The splitting in~(i) is the one in Remark~\ref{rmk:V}~(iii)
with $\Lambda_1=\Lambda_0^\perp$. To prove~(ii) we
examine the operator $D^*D$ of Lemma~\ref{le:D*D}. 
On smooth elements this is the Laplace-Beltrami operator. 
Hence its domain is
$$
\dom(D^*D)=
\left\{ \xi\in W_0 \,\Big|\,
\sup_{\eta\in W_0}\frac{\inner{D\xi}{D\eta}_{L^2}}
{\Norm{\eta}_{L^2}}<\infty
\right\}
=
W_0\cap W^{2,2}_Y
$$
by elliptic regularity. This implies that
$\dom D^*\cap\im D={D(W_0\cap W^{2,2}_Y)}$ 
is a closed subspace of $W^{1,2}_Y$. 
One can also think of $D^*$ as a bounded linear 
operator from $L^2_Y$ to $W^{-1,2}_Y:=(W_0)^*$, 
see the proof of Lemma~\ref{le:D*D}. 
Then the operator
\begin{equation}\label{eq:T}
W^{1,2}_Y\to W^{-1,2}_Y\times W^{1/2,2}_\Sigma :
\xi \mapsto (D^*D^*\xi + \xi,\xi|_\Sigma)
\end{equation}
is bijective, by elliptic regularity and the 
Sobolev trace theorem, and $V\cap W^{1,2}_Y$ 
is the preimage of $\{0\}\times W^{1/2,2}_\Sigma$
under this operator. Hence $V\cap W^{1,2}_Y$
is also a closed subspace of $W^{1,2}_Y$ and so is
the space
$$
\Lambda_1=(V\cap W^{1,2}_Y)\cap (\dom D^*\cap\im D).
$$
Next, the kernel of $D^*$ is a closed subspace of $L^2_Y$ 
and hence, so is the space 
$$
\Lambda_0=\bigl\{\xi-\left(\one+D^*D\right)^{-1}\xi\,\big|\,
\xi\in\ker D^*\bigr\}.
$$
See Remark~\ref{rmk:V}~(ii) for the projection $W\to V$;
the formula simplifies for $\xi\in\ker D^*$. This proves~(ii). 

We prove that the spaces of smooth elements are dense 
in $\Lambda_0$, $\Lambda_1$, $V$, and $W$.
Any element in $\Lambda_1$ can be approximated by a 
smooth sequence in $\Lambda_1$: The $W^{1,2}$-approximation 
by any smooth sequence converges in the graph norm of $D^*$ 
and projects under the map $\Pi_0$ in Remark~\ref{rm:llv} 
to a convergent smooth sequence in $\Lambda_1$.
Since $\Lambda_0=D^*\Lambda_1$, this shows that
the smooth elements are dense in $\Lambda_0$
as well as in $W=W_0\oplus\Lambda_0\oplus\Lambda_1$.

That the the restriction map~\eqref{restriction} extends
to an injective bounded linear operator from 
$\Lambda_1$ onto a closed subspace of~$W^{1/2,2}_\Sigma$
follows by restricting the isomorphism~\eqref{eq:T} to
the closed subspace $\Lambda_1$ of $V\cap W^{1,2}_Y$.
Next we prove that the map~\eqref{restriction} sends $\Lambda_0$ 
to a closed subspace of $W^{-1/2,2}(\Sigma)$.
For this it is convenient to use the following norms
for $\xi\in W$:
$$
\Norm{\xi|_\Sigma}_{W^{-1/2,2}_\Sigma}
:= \sup_{\eta\in W^{1,2}_Y}
\frac{\om(\xi,\eta)}{\Norm{\eta}_{W^{1,2}_Y}},\qquad
\Norm{\xi}_{D^*} 
:= \sqrt{\Norm{\xi}_{L^2_Y}^2+\Norm{D^*\xi}_{L^2_Y}^2}
$$
By definition there is a constant $c>0$ such that 
$$
\Norm{\xi|_\Sigma}_{W^{-1/2,2}_\Sigma}\le c\Norm{\xi}_{D^*}
$$
for every $\xi\in W$. Thus~\eqref{restriction}
is a bounded linear operator from $W$ to $W^{-1/2,2}_\Sigma$.
Moreover, $\Lambda_1$ is complete both with respect to the 
graph norm of $D^*$ and the $W^{1,2}$-norm,
and the former is bounded above by the latter.
Hence, by the open mapping theorem, there is a constant 
$\delta>0$ such that 
$$
\Norm{\eta}_{D^*}\ge\delta \Norm{\eta}_{W^{1,2}_Y}
\qquad \forall\eta\in\Lambda_1.
$$ 
Now let $\xi\in\Lambda_0$ be given.  
Then $D^*\xi\in\Lambda_1\subset W^{1,2}_Y$
and hence
$$
\Norm{\xi|_\Sigma}_{W^{-1/2,2}_\Sigma}
\ge
\delta\sup_{\eta\in W^{1,2}_Y}
\frac{\om(\xi,\eta)}{\Norm{\eta}_{D^*}}
\ge
\delta
\frac{\om(\xi,D^*\xi)}{\Norm{D^*\xi}_{D^*}} 
=
\delta\Norm{\xi}_{D^*}.
$$
Since $\Lambda_0$ is a closed subspace of $W$, the operator 
$\Lambda_0\to W^{-1/2,2}_\Sigma:\xi\mapsto\xi|_\Sigma$
is injective and has a closed image.  This proves the lemma.
\end{proof}

\begin{remark}\label{rmk:harmonic}\rm
The dual domain $W$ admits another orthogonal splitting
$$
W = (\dom D^*\cap\im D)\oplus \ker D^*
$$
where $\dom D^*\cap\im D$ is a closed
subspace of $W^{1,2}_Y$ and the kernel of $D^*$
is a closed subspace of $L^2_Y$.  It can be described  
as the image under $D^*$ of the space of harmonic 
pairs $\xi=(\alpha,\beta)\in W^{1,2}_Y$:
$$
\ker D^* = \left\{(*\rd\alpha-\rd\phi,-\rd^*\alpha)\,|\,
(\alpha,\phi)\in W^{1,2}_Y,\,
\rd^*\rd\alpha+\rd\rd^*\alpha=0,\,\rd^*\rd\phi=0\right\}.
$$
This can also be used to prove that the restriction 
map~\eqref{restriction} maps the kernel of $D^*$
to $W^{-1/2,2}_\Sigma$: If $\phi$ is a $W^{1,2}$
harmonic function on $Y$ then its restriction to the boundary
is of class $W^{1/2,2}$ and its normal derivative on the 
boundary is of class $W^{-1/2,2}$.   

Yet another splitting of $W$ can be obtained from 
eigenspace decompositions along the lines of 
Atiyah--Patodi--Singer~\cite{APS}.  The operator 
$D$ has the form $J(\p_t+B)$ near the boundary,
where $J^2=-\one$ and $B$ is a self-adjoint first order 
Fredholm operator over $\Sigma$.  The decomposition 
involves the eigenspaces of $B$~\cite{BossZhu}.
\end{remark}

\begin{proof}[Proof of Theorem~\ref{thm:S1Y}.]
It suffices to prove the theorem for 
$X_f=0$ because any two perturbations are homotopic
and result in compact perturbations of the operators
$\cD_{v,\A}$ and hence in isomorphic determinant
line bundles.

We prove~(i).
By Lemma~\ref{le:CS} the degree depends only 
on the homotopy class of $(v,\A)$.
Given such a pair, there is a smooth 
path $[0,1]\to\cG(Y):\tau\mapsto v^\tau$
with $v^0=v$ and $v^1=\one$, because $\cG(Y)$ is connected.
Let ${u^\tau:\R\to\cG(Y)}$ be the smooth path of gauge
transformations constructed in Lemma~\ref{le:Xuv} below
with $X=\mathrm{pt}$ and define 
$$
\A^\tau:=(u^\tau)^*\A.
$$  
Then $\tau\mapsto(v^\tau,\A^\tau)$
is a smooth path in $\cA(S^1\times Y,\cL)$
connecting $(v^0,\A^0)=(v,\A)$ to a pair 
of the form $(\one,A^1)$.  Hence we may assume 
without loss of generality that
$v=v'=\one$ and $\A,\A'\in\cA(P,\cL)$ where
$P=P_\one=S^1\times Y\times\rG$. Now the map 
$$
\cA(P,\cL)\to\Cinf(S^1,\cL) \,:\; \A\mapsto A|_{S^1\times\Sigma}
$$
is a homotopy equivalence. Hence~(i) follows
from the fact that, by (L2), every loop in $\cL$
is homotopic to a loop of the form 
$\R/\Z\to\cL:s\mapsto u(s)^*A_0$ with 
${u(s+1)=u(s)\in\cG(\Sigma)}$, and that the homotopy
class of such a loop is characterized by the 
degree of the map $u:S^1\times\Sigma\to\rG$.

We prove~(ii).  That the operator $\cD_{v,\A}$
has a finite dimensional kernel and a closed image
follows immediately from the estimate 
in Theorem~\ref{thm:4reg}~(ii) and Rellich's theorem
(see~\cite[Lemma~A.1.1]{MS}).  That it 
has a finite dimensional cokernel follows
from the regularity results in Theorem~\ref{thm:4reg}
and Remark~\ref{rmk:time reversal}. (The dual operator
has a finite dimensional kernel.)
Thus we have proved that $\cD_{v,\A}$ is a Fredholm
operator for every pair $(v,\A)\in\cA(S^1\times Y,\cL)$. 
The regularity theory in Theorem~\ref{thm:4reg}
also shows that its kernel and cokernel, 
and hence also the Fredholm index, are independent 
of $k$ and $p$. Moreover, the Fredholm index 
depends only on the homotopy class of $(v,\A)$;
to see this one can use the argument in the 
proof of Step~1 in Theorem~\ref{thm:4reg}
to reduce the problem to small deformations with
constant domain and then use the stability 
properties of the Fredholm index. 
So by (i) it suffices to consider one pair $(v,\A)$ 
in each degree. Hence we can assume
$$
v|_N = \one,\qquad \Phi=0,\qquad A(s)|_N =0
$$
for all $s$ and an open neighbourhood $N\subset Y$ of $\p Y$.
Then $\deg(v,\A)=\deg(v)$.  
Choose a handle body $Y'$ with $\pd Y'=\bar\Sigma$
and extend $A(s)$ smoothly by the trivial connection 
on $Y'$ to obtain a smooth connection $\tA(s)$ on the 
closed $3$-manifold 
$$
\tY:=Y\cup_\Sigma Y'
$$ 
for every $s$.  Note that $\tA(s+1)=\tv^*\tA(s)$, 
where $\tv\in\cG(\tY)$ agrees with $v$ on $Y$ and 
is equal to $\one$ on $Y'$. 
Let $\cH_0'$ denote the Hessian on $Y'$ (at the trivial connection) 
and $\cH_{A(s)}$ the Hessian on $Y$, both with the 
same boundary Lagrangian $T_0\cL$.
These are self-adjoint Fredholm operators, by 
Proposition~\ref{prop:HA}. The Hessian (\ref{eq:HA}) over the 
closed manifold $\tY$ will be denoted by $\tilde{\cH}_{\tA(s)}$. 
Choose $\eps>0$ such that the operators 
$\cH'_0 + \eps\Id$, $\cH_{A(0)}+ \eps\Id$, 
and $\tilde{\cH}_{\tA(0)}+ \eps\Id$ are all bijective.  
We shall introduce the spectral flow $\mu_\mathrm{spec}$ 
(as defined in Appendix~\ref{app:spec}) and prove that
\begin{equation}\label{eq:index8}
\begin{split}
\mathrm{index}(\cD_{v,\A}) 
&= 
\mu_\mathrm{spec}\bigl(\bigl\{
\bigl(\cH_{A(s)}+ \eps\Id\bigr)\oplus
\bigl(\cH'_0+ \eps\Id\bigr)
\bigr\}_{s\in[0,1]}\bigr) \\ 
&= 
\mu_\mathrm{spec}
\bigl(\bigl\{\tilde{\cH}_{\tA(s)}+ \eps\Id\bigr\}_{s\in[0,1]}\bigr) \\
&= 
\mathrm{index}(\cD_{\tv,\widetilde{\A}}) 
= 8\deg(\tv) 
= 8\deg(v,\A).  
\end{split}
\end{equation}
Here $\cD_{\tv,\widetilde{\A}}=\nabla_s + \tilde{\cH}_{\tA(s)}$ denotes
the anti-self-duality operator on the twisted bundle $P_{\tilde{v}}$ 
over $S^1\times\tilde{Y}$.

To prove~(\ref{eq:index8}) we may assume $k=1$ and $p=2$.
In this case the first and third equations
follow from Theorem~\ref{thm:indexS1}, the fourth 
equation follows from the Atiyah--Singer index theorem
(the second Chern class of the principal bundle 
$P_\tv\to S^1\times\tY$ is the degree of $\tv$),
and the last equation is obvious from the definitions. 
To prove the second equation in~(\ref{eq:index8}) 
consider the operator family
$$
D(s):= \bigl(\cH_{A(s)}+ \eps\Id \bigr) 
\oplus \bigl(\cH'_0  + \eps\Id \bigr)
$$
on the Hilbert space
$$
H:= L^2(Y,\rT^*Y\otimes\cg)\oplus L^2(Y,\cg)
\oplus L^2(Y',\rT^*Y'\otimes\cg)\oplus L^2(Y',\cg) 
$$
with the constant dense domain $\mathrm{dom}\,D(s) = W_0$, where
$$
W_0:= W^{1,2}_0(Y,\rT^*Y\otimes\cg)\oplus W^{1,2}_0(Y,\cg)
\oplus W^{1,2}_0(Y',\rT^*Y'\otimes\cg)\oplus W^{1,2}_0(Y',\cg).
$$
As in Remark~\ref{hurrah}, this choice of domain makes $D(s)$ 
closed, symmetric, and injective. Moreover, the Gelfand--Robbin quotient
and its symplectic structure 
$$
V:=\mathrm{dom}\,D(s)^*/\mathrm{dom}\,D(s) =  W / W_0
$$
are independent of $s$.
Now, by Appendix~\ref{app:GR}, self-adjoint extensions 
of $D(s)$ are in one-to-one
correspondence with Lagrangian subspaces of $V$.
The operators in the first row of~(\ref{eq:index8})
all correspond to the Lagrangian subspace 
$$
\Lambda_1 := \left\{(\alpha,\phi,\alpha',\phi')\in W^{1,2} 
\,\biggl|\,
\begin{array}{c}
*\alpha|_{\pd Y},*\alpha'|_{\pd Y'}=0, \\
\alpha|_{\pd Y},\alpha'|_{\pd Y'}\in \rT_0\cL
\end{array}
\right\} / W_0 \subset V ,
$$
where $W^{1,2}:=W^{1,2}(Y,\rT^*Y\otimes\cg\oplus\cg)
\times W^{1,2}(Y',\rT^*Y'\otimes\cg\oplus\cg)\subset W$. 
The operators in the second row of~(\ref{eq:index8})
all correspond to the `diagonal'
$$
\Lambda_2 :=
\left\{(\alpha,\phi,\alpha',\phi')\in W^{1,2} \,\left|\,
\begin{array}{c}
\phi|_{\pd Y} = \phi'|_{\pd Y'}, \\
\alpha|_{\pd Y} = \alpha'|_{\pd Y'},\\
*\alpha|_{\pd Y}=*\alpha'|_{\pd Y'}
\end{array}
\right.\right\} / W_0 \subset V .
$$
For $i=1,2$ and $s\in\R$ let $D(s)_{\Lambda_i}:\dom D(s)_{\Lambda_i}\to H$ 
denote the restriction of $D(s)^*$ to the preimage of $\Lambda_i$
under the projection $W\to W/W_0$. Then $D(s)_{\Lambda_i}$ is 
self-adjoint.
Moreover, we have $D(s+1)=Q^{-1}D(s) Q$, where $Q:H\to H$ is given by 
conjugation with the gauge transformation $v$ and satisfies
$\xi-Q\xi\in W_0$ for all $\xi\in W$ since $v\equiv\one$ near $\pd Y$.
This implies that
$$
\Lambda_0:=(\ker D(0)^* \oplus W_0)/W_0
= (\ker D(1)^* \oplus W_0)/W_0 .
$$ 
Then, by the choice of $\eps$, the 
Lagrangian subspaces $\Lambda_1$ and $\Lambda_2$ are transverse to $\Lambda_0$. 
Moreover, they are compact perturbations of $\Lambda_0^\perp$ 
by Lemma~\ref{le:DLcpct}, since the graph norm on 
$\dom D(s)_{\Lambda_i}$ is equivalent to the $W^{1,2}$-norm, see (\ref{eq:graphnorm}).
The second identity in (\ref{eq:index8}) follows from Remark~\ref{rmk:loopmaslov}, 
which asserts that the spectral flow of 
$\{D(s)_\Lambda\}_{s\in[0,1]}$ is independent of the Lagrangian 
subspace $\Lambda\subset V$ that is transverse to $\Lambda_0$ and a compact
perturbation of $\Lambda_0^\perp$. 
This proves~(\ref{eq:index8}) and thus~(ii).

We prove~(iii) and~(iv).  That two isomorphic pairs
$(v_0,\A_0)$ and $(v_1,\A_1)=u^*(v_0,\A_0)$ have the 
same degree follows from~(ii) and the fact that 
conjugation by~$u$ identifies kernel and co\-kernel
of the operator $\cD_{v_0,\A_0}$ with kernel and
co\-kernel of $\cD_{v_1,\A_1}$.  For every 
$(v,\A)\in\cA(S^1\times Y,\cL)$ denote by $\Or(\cD_{v,\A})$ 
the two element set of orientations of $\det(\cD_{v,\A})$.
Then the remaining assertions in~(iii) and~(iv) can be 
rephrased as follows.

\smallbreak

\medskip\noindent{\bf Claim:}
{\it Let $\left\{(v_\lambda,\A_\lambda)\right\}_{0\le\lambda\le 1}$
be a smooth path in $\cA(S^1\times Y,\cL)$ 
and $u:\R\to\cG(Y)$ be a morphism from $(v_0,\A_0)$
to $(v_1,\A_1)$.  Then the isomorphism
$$
u^*:\Or(\cD_{v_0,\A_0})\to\Or(\cD_{v_1,\A_1})
$$
agrees with the isomorphism induced by the path 
$\lambda\mapsto(v_\lambda,\A_\lambda)$.}

\medskip\noindent
When $u\equiv\one$, the claim asserts that the 
automorphism of $\det(\cD_{v_0,\A_0})$ induced by a loop
in $\cA(S^1\times Y,\cL)$ is orientation preserving and 
hence the determinant bundle over $\cA(S^1\times Y,\cL)$
is orientable. Throughout we write 
$\A_\lambda=\Phi_\lambda(s)\,\rd s+A_\lambda(s)$
We prove the claim in five steps.

\smallskip\noindent{\bf Step~1.}
{\it It suffices to assume that $v_\lambda=\one$ 
for every $\lambda$.}

\smallskip\noindent
Since $\cG(Y)$ is connected, there exists a smooth homotopy
$[0,1]\times[0,1]\to\cG(Y):(\tau,\lambda)\mapsto v^\tau_\lambda$ 
from $v^0_\lambda=v_\lambda$ to $v^1_\lambda=\one$.
By Lemma~\ref{le:Xuv} below with $X=[0,1]$, there exists 
a smooth map 
$
[0,1]\times[0,1]\times\R\to\cG(Y):
(\tau,\lambda,s)\mapsto u^\tau_\lambda(s)
$
such that 
$$
v^\tau_\lambda 
= u^\tau_\lambda(s)^{-1}v_\lambda u^\tau_\lambda(s+1),\qquad
u^0_\lambda(s) = \one.
$$
Define 
$$
\A^\tau_\lambda:=(u^\tau_\lambda)^*\A_\lambda,\qquad
u^\tau := (u^\tau_0)^{-1}uu^\tau_1.
$$
Then $(v^\tau_\lambda)^*A^\tau_\lambda(s) = A^\tau_\lambda(s+1)$, 
$\A^\tau_1 = (u^\tau)^*\A^\tau_0$, 
and
$
v^\tau_1=u^\tau(s)^{-1}v^\tau_0u^\tau(s+1).
$ 
Hence $(v^\tau_\lambda,\A^\tau_\lambda)\in\cA(S^1\times Y,\cL)$ 
for all $\tau$ and $\lambda$, and $u^\tau$ is a morphism from 
$(v^\tau_0,\A^\tau_0)$ to $(v^\tau_1,\A^\tau_1)$ for every $\tau$.  
By continuity, the claim holds for $\tau=0$ 
if and only if it holds for $\tau=1$. 
Since $v^1_\lambda=\one$ for every $\lambda$,
this proves Step~1. 

\medskip\noindent{\bf Step~2.}
{\it It suffices to assume that $v_\lambda=\one$ 
and $u|_{S^1\times\Sigma}=\one$.}

\smallskip\noindent
By Step~1 we can assume $v_\lambda=\one$.    
The restriction of the map ${u:S^1\times Y\to\rG}$ 
to the boundary has degree zero (see e.g.\ \cite[\S 5,Lemma~1]{MILNOR}). 
Hence there exists a smooth path
$[0,1]\to\cG(P):\tau\mapsto u^\tau$ such that 
$u^0=u$ and $u^1|_{S^1\times\Sigma}=\one$. 
Composing the paths $\{\A_\lambda\}_{0\le\lambda\le 1}$
and $\{(u^{\lambda\tau})^*\A_0\}_{0\le\lambda\le 1}$
we obtain a homotopy of homotopies
$\tau\mapsto\{\A^\tau_\lambda\}_{0\le\lambda\le 1}$
with $\A^0_\lambda=\A_\lambda$ and $\A^\tau_1=(u^\tau)^*\A^\tau_0$.
Hence Step~2 follows as in Step~1 by continuity.

\medskip\noindent{\bf Step~3.}
{\it 
Using (L2) we see that
it suffices to assume that $v_\lambda=\one$,
$u|_{S^1\times\Sigma}=\one$, and there exists a smooth 
map $[0,1]\times S^1\to\cG_z(\Sigma):(\lambda,s)\mapsto w_\lambda(s)$
satisfying 
$
A_\lambda(s)|_\Sigma = w_\lambda(s)^{-1}\rd w_\lambda(s)
$
and $w_\lambda(s+1)=w_\lambda(s)$, $w_0(s)=w_1(s)$,
$w_0(0)=\one$.}

\medskip\noindent
By Step~2 we can assume $v_\lambda=\one$ and $u|_{S^1\times\p Y}=\one$.
Then $A_\lambda(s+1)=A_\lambda(s)$ and $A_0(s)=A_1(s)$
for all $s$ and $\lambda$.  Since $\cL/\cG_z(\Sigma)$ is connected 
and simply connected, the loops 
${[0,1]\to\cL:\lambda\mapsto A_\lambda(0)|_\Sigma}$ 
and ${S^1\to\cL:s\mapsto A_0(s)|_\Sigma}$ 
are homotopic to loops in the based gauge equivalence 
class of the zero connection in $\cL$. 
This implies that there is a smooth homotopy
$[0,1]^2\times S^1\to\cL:(\tau,\lambda,s)\mapsto B^\tau_\lambda(s)$
of homotopies of loops, satisfying
$$
B^\tau_\lambda(s+1)=B^\tau_\lambda(s),\qquad 
B^\tau_0(s)=B^\tau_1(s),
$$
starting at $B^0_\lambda(s)=A_\lambda(s)|_\Sigma$ and 
ending at a homotopy of loops satisfying
$$
B^1_\lambda(0),B^1_0(s)
\in\left\{w^{-1}\rd w\,\big|\,w\in\cG_z(\Sigma)\right\}.
$$
The composition of the map 
$[0,1]^2\to\cL:(\lambda,s)\mapsto B^1_\lambda(s)$
with the projection $\cL\to\cL/\cG_z(\Sigma)$ maps the boundary 
to a point.  Since ${\pi_2(\cL/\cG_z(\Sigma))=0}$ the homotopy
$\tau\mapsto B^\tau$ can be extended to the interval $0\le\tau\le 2$ 
so that 
$
B^2_\lambda(s) = w_\lambda(s)^{-1}\rd w_\lambda(s) .
$
This determines the map 
$[0,1]\times\R\to\cG_z(\Sigma):(\lambda,s)\mapsto w_\lambda(s)$
uniquly, hence $w$ satisfies the requirements of Step~3. 
Since the restriction map $\cA(Y,\cL)\to\cL$ 
is a homotopy equivalence, there exists a smooth homotopy
$[0,2]\times [0,1]\to\cA(P_\one,\cL):(\tau,\lambda)\mapsto\A^\tau_\lambda$ 
with $\A^\tau_1=u^*\A^\tau_0$ from $\A^0_\lambda=\A_\lambda$ to
$\A^2_\lambda$ satisfying $A^2_\lambda(s)|_\Sigma = B^2_\lambda(s)$. 
Step~3 follows since, by continuity, the claim holds for $\tau=0$ 
if and only if it holds for $\tau=2$. 

\medskip\noindent{\bf Step~4.}
{\it It suffices to assume that $v_\lambda=v$ is independent of $\lambda$
and there exists a neighbourhood $N\subset Y$ of $\p Y$ such that 
$v|_N=\one$, $A_\lambda(s)|_N=0$, 
$\Phi_\lambda(s)|_N=0$, and $u(s)|_N=\one$.}

\medskip\noindent
By Step~3 we can assume $v_\lambda=\one$, $u|_{S^1\times\Sigma}=\one$, 
and $A_\lambda(s)|_\Sigma = w_\lambda(s)^{-1}\rd w_\lambda(s)$
for a smooth map $w:[0,1]\times S^1\to\cG_z(\Sigma)$.
By a further homotopy argument we may assume that $w$ is transversally 
constant near the edges of the square, 
$\p_\lambda w_\lambda(s)=0$ for $\lambda\simeq 0$ and $\lambda\simeq 1$, 
and $\p_sw_\lambda(s)=0$ for $s\simeq 0$ and $s\simeq 1$. 
Since every gauge transformation on $\Sigma$ extends 
to a gauge transformation on $Y$ and the same holds for
families parametrized by contractible domains, there is a 
smooth map $[0,1]^2\to\cG(Y):(\lambda,s)\mapsto u_\lambda(s)$ 
such that
$$
u_\lambda(s)|_\Sigma = w_\lambda(s)^{-1}.
$$
This map can be chosen such that
$\p_\lambda u_\lambda(s)=0$ for $\lambda\simeq 0$ and $\lambda\simeq 1$,
and $\p_su_\lambda(s)=0$ for $s\simeq 0$ and $s\simeq 1$.
Moreover, we can achieve $\lambda$-independence of  
$v'_\lambda := u_\lambda(0)^{-1}u_\lambda(1)$.
To see this, note that $v'_\lambda|_\Sigma=\one$ and
there is a $\delta>0$ such that 
$\pd_\lambda v'_\lambda=0$ for $\lambda\not\in (\delta,1-\delta)$.
Let $\beta:[0,1]\to[0,1]$ be a smooth monotone cutoff function 
such that $\beta(\lambda)=\lambda$ for $\lambda\in[\delta,1-\delta]$,
$\beta\equiv 0$ for $\lambda\simeq 0$, 
and $\beta\equiv 1$ for $\lambda\simeq 1$.
Now we can replace $u_\lambda(s)$ by 
$u_\lambda(s) (v'_{\beta(s)\beta(\lambda)})^{-1}$.
The resulting map $(\lambda,s)\mapsto u_\lambda(s)$ satisfies
$u_\lambda(1)=u_\lambda(0)v'$ with $v'$ independent of $\lambda$,
as claimed.
Hence it extends to $[0,1]\times\R$ such that
$
v'=u_\lambda(s)^{-1}u_\lambda(s+1)
$
for all $\lambda$ and $s$.  Define
$$
\A_\lambda' := u_\lambda^*\A_\lambda\in\cA(P_{v'},\cL),\qquad
u' := u_0^{-1}uu_1.
$$
Then $v'|_\Sigma=\one$, $A'_\lambda|_\Sigma\equiv0$, $u'|_\Sigma\equiv\one$,
and ${u'}^*(v',\A_0')=(v',\A_1')$.
Moreover $u_\lambda$ is a morphism from 
$(\one,\A_\lambda)$ to $(v',\A_\lambda')$ for every $\lambda$.
This gives a commuting diagram
$$
\begin{array}{ccc}
\det(\cD_{\one,\A_0}) & \mapright{u^*} &  \det(\cD_{\one,\A_1}) \\
\mapdown{u_0^*}  & & \mapdown{u_1^*}  \\
\det(\cD_{v',\A'_0}) & \mapright{{u'}^*} &\det(\cD_{v',\A'_1}).
\end{array}
$$
There is a second diagram where the horizontal arrows
are induced by the paths $\lambda\mapsto(\one,\A_\lambda)$
and $\lambda\mapsto(v',\A_\lambda')=u_\lambda^*(\one,\A_\lambda)$.
That this second diagram commutes as well follows from a homotopy 
argument; namely the space of smooth maps
$[0,1]^2\to\cG(Y):(s,\lambda)\mapsto u_\lambda(s)$ is connected
and the diagram obviously commutes when $u_\lambda(s)\equiv\one$.
This shows that the claim holds for $(u,\one,\A_\lambda)$ if and only 
if it holds for $(u',v',\A'_\lambda)$.   Hence Step~4 follows from 
a further homotopy argument (to achieve the relevant 
boundary conditions and vanishing of $\Phi$ in a neighbourhood of $\p Y$).

\medskip\noindent{\bf Step~5.}
{\it We prove the claim.}

\medskip\noindent
By Step~4, we may assume that $v_\lambda=v$ and 
there exists a neighbourhood $N\subset Y$ of $\p Y$ such that 
$v|_N=\one$, $A_\lambda(s)|_N=0$, 
$\Phi_\lambda(s)|_N=0$, and $u(s)|_N=\one$.
We shall argue as in the proof of~(ii), namely 
choose a handle body $Y'$ with $\pd Y'=\bar\Sigma$ 
and transfer the problem to the closed 
$3$-manifold $\tilde Y:=Y\cup_\Sigma Y'$.

Since the map on orientations induced by the 
path $\lambda\mapsto\A_\lambda$ is invariant under homotopy
we may assume that the path is the straight line
$$
\A_\lambda=(1-\lambda)\A + \lambda u^*\A,
$$
where $\A\in\cA(P_v)$ vanishes near the boundary and 
$u\in\cG(P_v)$ is equal to the identity near the boundary. 
Since $v\in\cG(Y)$ is the identity near the boundary
we can extend it to a gauge transformation
$\tv\in\cG(\tY)$ via $\tv|_{Y'}:=v':=\one$. 
Then $u\in\cG(P_v)$ extends to a gauge transformation 
$\tu\in\cG(P_\tv)$ via $\tu(s)|_{Y'}:=\one$
and $\A$ extends to a connection $\tilde\A\in\cA(P_\tv)$
via $\tilde\A|_{S^1\times Y'}:=\A'=0$.
As in the proof of~(ii) we have three
Fredholm operators $\cD_{v,\A}$ on $S^1\times Y$,
$\cD_{v',\A'}$ on $S^1\times Y'$ (both with boundary conditions
$*\alpha|_{\p Y}=0$ and $\alpha|_{\p Y}\in\rT_0\cL$), and $\cD_{\tv,\tilde\A}$
on $S^1\times\tY$ (without boundary conditions).  We must prove that the isomorphism
$$
u^*:\Or(\cD_{v,\A})\to\Or(\cD_{v,u^*\A})
$$
agrees with the isomorphism determined by the homotopy. 
Since both the gauge transformation and the homotopy
act trivially on $\det(\cD_{v',\A'})$ this means that
the isomorphism
\begin{equation}\label{eq:u*1}
u^*\otimes\Id:
\Or(\cD_{v,\A}\times\cD_{v',\A'})
\to \Or(\cD_{v,u^*\A}\times\cD_{v',\A'})
\end{equation}
agrees with the homotopy isomorphism.  
As in the proof of~(ii) we choose a family of Lagrangian 
subspaces connecting $\Lambda_1$ to $\Lambda_2$
to obtain two continuous families of isomorphisms 
(see Lemma~\ref{lem:kernel}; we use the fact
that the Lagrangian subspaces can be chosen as 
compact perturbations of $\Lambda_0^\perp$).
For $\Lambda_1$ the gauge transformation 
induces the isomorphism~(\ref{eq:u*1}) 
and for $\Lambda_2$ the isomorphism
\begin{equation}\label{eq:u*2}
\tu^*:\Or(\cD_{\tv,\tilde\A})
\to \Or(\cD_{\tv,\tu^*\tilde\A})
\end{equation}
and similarly for the homotopy induced isomorphisms.
For $\Lambda_2$ both isomorphisms agree by the standard
theory for self-duality operators on closed $4$-manifolds
(see~\cite{Donaldson orient}).  Hence they agree for 
$\Lambda_1$. This proves the claim and the theorem. 
\end{proof}

\begin{lemma}\label{le:Xuv}
Let $X$ be a manifold and 
$
[0,1]\times X\to\cG(Y):(\tau,x)\mapsto v^\tau_x
$
be a smooth map.  Then there is a smooth map 
$$
{[0,1]\times X\times\R\to\cG(Y):(\tau,x,s)\mapsto u^\tau_x(s)}
$$
such that
\begin{equation}\label{eq:Xuv}
v^\tau_x = u^\tau_x(s)^{-1}v^0_xu^\tau_x(s+1),\qquad
u^\tau_x(0)=\one.
\end{equation}
\end{lemma}

\begin{proof}
Choose a cutoff function $\beta:[0,1]\to[0,1]$
such that $\beta(s)=0$ for $s\simeq 0$ and $\beta(s)=1$
for $s\simeq 1$. Define 
$$
u^\tau_x(s):=(v^0_x)^{-1}v^{\beta(s)\tau}_x,\qquad
0\le s\le 1. 
$$
Then $u^\tau_x(s)=\one$ for $s\simeq 0$
and $u^\tau_x(s)=(v^0_x)^{-1}v^\tau_x$ for $s\simeq 1$.
Hence $u^\tau_x$ extends uniquely to a smooth map
from $\R$ to $\cG(Y)$ that satisfies~(\ref{eq:Xuv});
the extension to $(1,\infty)$ is given by 
$u_x^\tau(s+1):=(v_x^0)^{-1} u^\tau_x(s) v_x^\tau$
and the extension to $(-\infty,0)$ by $u_x^\tau(s-1):=v_x^0 u^\tau_x(s) (v_x^\tau)^{-1}$,
in both cases for $s> 0$.  Moreover, the resulting map 
$[0,1]\times X\times\R\to\cG(Y)$ is smooth in all variables.
\end{proof}


\section{Exponential decay}\label{sec:exp}

Let $Y$ be a compact oriented $3$-manifold with boundary $\pd Y=\Sigma$
and let ${\cL\subset\cA(\Sigma)}$ be a gauge invariant, monotone 
Lagrangian submanifold satisfying (L1-2) on page~\pageref{p:L1}. 
(Actually this section only requires the compactness of $\cL/\cG_z(\Sigma)$ from (L2).)
We fix a perturbation ${X_f:\cA(Y)\to\Om^2(Y,\cg)}$ 
as in Section~\ref{sec:CS}. The purpose of this section is 
to establish the exponential decay for finite energy 
solutions in the following two Theorems. 
The unperturbed Yang-Mills energy 
of a connection $\A\in\cA(\R\times Y)$ 
is $\frac 12 \int |F_\A|^2$. 
In the presence of a holonomy perturbation 
the gauge invariant energy of $\A=A + \Phi\ds$ is
$$
E_f(\A) = \frac 12 \int_{\R\times Y} \bigl| F_\A + X_f(\A) \bigr|^2 
= \frac 12 \int_{\R\times Y} \left( \bigl| \pd_s A -\rd_A\Phi \bigr|^2 
+ \bigl| F_A + X_f(A) \bigr|^2 \right) .
$$
An anti-self-dual connection in temporal gauge
satisfies $\pd_s A + *\bigl( F_A + X_f(A) \bigr) = 0$ and $\Phi=0$ 
and the energy simplifies to $E_f(\A) = \int_{\R\times Y} |\pd_s A|^2$.

\begin{thm}\label{thm:decay}
Suppose that every critical point of the perturbed 
Chern--Simons functional $\CS_\cL+h_f$ is nondegenerate.  
Then there is a constant $\delta>0$ such that the following holds.
If $A:[0,\infty)\to\cA(Y)$ is a smooth solution of 
\begin{equation}\label{eq:bvp}
\pd_sA+*(F_A+X_f(A)) =0,\qquad A(s)|_{\pd Y}\in \cL,
\end{equation}
satisfying
$$
\int_0^\infty\int_Y
\left|\pd_sA\right|^p\dvol_Y\,\ds < \infty,\qquad
p\ge 2,
$$
then there is a connection 
$A_\infty\in\cA(Y, \cL)$ such that
$F_{A_\infty}+X_f(A_\infty)=0$ and
$A(s)$ converges to $A$ as $s\to\infty$.
Moreover, there are constants $C_0,C_1,C_2,\dots$ 
such that
$$
\left\| A - A_\infty \right\|_{\cC^k([s-1,s+1]\times Y)}
\le C_ke^{-\delta s}
$$
for every $s\ge1$ and every integer $k\ge0$.
\end{thm}

\begin{remark}\label{rmk:norms}\rm
Let $X$ be a compact Riemannian manifold with boundary.
We shall need gauge invariant Sobolev norms on the spaces
$\Om^\ell(X,\cg)$ depending on a connection $\A\in\cA(X)$. 
For $p\ge1$ and an integer $k\ge0$ we define
$$
\left\|\alpha\right\|_{W^{k,p},\A}
:= \biggl(
\sum_{j=0}^k\int_X\bigl|\nabla_\A^j\alpha\bigr|^p\biggr)^{1/p}
$$ 
for $\alpha\in\Om^\ell(Y,\cg)$, where
$\nabla_\A^j\alpha$ denotes
the $j$th covariant derivative of $\alpha$ twisted by $\A$. 
For $p=\infty$ we define
$$
\left\|\alpha\right\|_{W^{k,\infty},\A}
:= \left\|\alpha\right\|_{\cC^k,\A}
:= \max_{0\leq j\leq k}\, \sup_X\, \bigl|\nabla_\A^j\alpha\bigr| .
$$ 
These norms are gauge invariant in the sense that
$$
\left\|u^{-1}\alpha u\right\|_{W^{k,p},u^*\A}
= \left\|\alpha\right\|_{W^{k,p},\A}
$$
for every gauge transformation $u\in\cG(X)$.
In particular, for $k=0$ the $L^p$-norms are gauge
invariant and do not depend on the connection $\A$.
\end{remark}

\begin{thm}\label{thm:long}
Suppose that every critical point of the perturbed 
Chern--Simons functional is nondegenerate. 
Then, for every $p>1$, there are positive 
constants $\eps$, $\delta$, $C_0,C_1,\dots$ 
such that the following holds for every $T\ge1$.
If ${A:[-T,T]\to\cA(Y)}$ is a smooth solution of~(\ref{eq:bvp}) 
satisfying
\begin{equation}\label{eq:longenergy}
\int_{-T}^T\int_Y\left|\pd_sA\right|^2\dvol_Y\,\ds < \eps ,
\end{equation}
then, for every $s\in[0,T-1]$ and every integer $k\ge0$,
\begin{equation} \label{eqn:dsA}
\left\|\pd_sA\right\|_{\cC^k([-s,s]\times Y),\A}
\le C_ke^{-\delta(T-s)}
\left\|\pd_sA\right\|_{L^2(([-T,1-T]\cup[T-1,T])\times Y)},
\end{equation}
where $\A\in\cA([-T,T]\times Y)$
is the connection associated to the path $A$.
Moreover, there is a connection $A_0\in\cA(Y,\cL)$ 
with $F_{A_0}+X_f(A_0)=0$ such that
\begin{multline} \label{eqn:a-a0}
\left\|A-A_0\right\|_{\cC^0([-s,s]\times Y)}
+ \left\|A-A_0\right\|_{W^{1,p}([-s,s]\times Y),A_0} \\
\quad
\le C_0 e^{-\delta(T-s)}
\left\|\pd_sA\right\|_{L^2(([-T,1-T]\cup[T-1,T])\times Y)}
\end{multline}
for every $s\in[0,T-1]$.
\end{thm}

The proofs of these results will be given below. 
Theorem~\ref{thm:decay} guarantees the existence of a limit
for each finite energy solution of~(\ref{eq:bvp}), however, the 
constants in the exponential decay estimate depend on the solution. 
With the help of Theorem~\ref{thm:long} one can show
that these constants can be chosen independent of the 
solution of~(\ref{eq:bvp}) and depend only on the 
limit $A_\infty$. 
This will be important for the gluing analysis.

\begin{cor}\label{cor:decay}
Let $A_\infty$ be a nondegenerate critical point of the perturbed 
Chern--Simons functional $\CS_\cL+h_f$.  
Then there are positive constant $\delta$, $\eps$, 
$C_0,C_1,\dots$  such that the following holds.
If $A:[0,\infty)\to\cA(Y)$ is a smooth solution of~(\ref{eq:bvp})
satisfying
$$
\int_0^\infty\int_Y
\left|\pd_sA\right|^2\dvol_Y\,\ds < \eps,\qquad
\lim_{s\to\infty}A(s) = A_\infty, 
$$
then 
$$
\Norm{A - A_\infty}_{\cC^k([s,\infty)\times Y)}
\le C_ke^{-\delta s}\Norm{\pd_sA}_{L^2([0,\infty)\times Y)}
$$
for every $s\ge1$ and every integer $k\ge0$.
\end{cor}

\begin{proof}
Let $\delta$, $\eps$, $C_k'$ be the constants
of Theorem~\ref{thm:long}.  Then
$$
\Norm{\p_sA}_{\cC^k([s,\infty)\times Y),\A}
\le C_k'e^{-\delta s}\Norm{\pd_sA}_{L^2([0,\infty)\times Y)}
$$
for $k=0,1,2,\dots$ and $s\ge1$. 
For $k=0$ the desired estimate follows by integrating from
$s$ to $\infty$ because the $\cC^0$-norm is independent 
of the reference connection $\A$.
Now argue by induction.  If the result 
has been established for any $k$ then there is a constant
$c_k$, depending on $C_k$, such that 
$$
\Norm{\alpha}_{\cC^{k+1}([s,\infty)\times Y)}
\le c_k\Norm{\alpha}_{\cC^{k+1}([s,\infty)\times Y),\A}
$$
for every $\alpha:[1,\infty)\to\Om^1(Y,\cg)$.
Applying this to $\alpha=\p_sA$ we obtain
$$
\Norm{\p_sA}_{\cC^{k+1}([s,\infty)\times Y)}
\le c_kC_{k+1}'e^{-\delta s}\Norm{\pd_sA}_{L^2([0,\infty)\times Y)}
$$
and the required $\cC^{k+1}$-estimate follows again by integrating
from $s$ to $\infty$. This proves the corollary. 
\end{proof}

The proof of Theorems~\ref{thm:decay} and \ref{thm:long} is based on
the following three lemmas concerning solutions on a long cylinder with
little energy.
We show that such solutions are uniformly close to a critical point
and establish uniform estimates for the Hessian and the linearized operator.

\begin{lemma}\label{le:decay1}
For every $\kappa>0$, $\rho>0$, and $p>1$
there is an $\eps>0$ such that the following holds.  
If $A:[-\rho,\rho]\to\cA(Y)$ is a solution
of~(\ref{eq:bvp}) that satisfies
$$
\int_{-\rho}^{\rho}\int_Y\left|\pd_sA\right|^2\dvol_Y\,\ds < \eps
$$
then there is a connection $A_\infty\in\cA(Y,\cL)$ with
$F_{A_\infty}+X_f(A_\infty)=0$ such that
\begin{align}\label{eq:ass}
\left\|A(0)-A_\infty\right\|_{W^{1,p}(Y),A_\infty} 
+ \left\|A(0)-A_\infty\right\|_{L^\infty(Y)} 
+ \left\|\pd_sA(0)\right\|_{L^\infty(Y)} < \kappa.
\end{align}
\end{lemma}

\begin{proof}
Assume by contradiction that this is wrong. 
Then there exist constants $\kappa>0$, $\rho>0$, and $p>1$
and a sequence $A_\nu:[-\rho,\rho]\to\cA(Y)$ 
of solutions of~(\ref{eq:bvp}) such that
\begin{equation}\label{eq:limlong}
\lim_{\nu\to\infty}
\int_{-\rho}^{\rho}\int_Y
\left|\pd_sA_\nu\right|^2\dvol_Y\,\ds =0
\end{equation}
but~(\ref{eq:ass}) fails.
Let $\A_\nu\in\cA([-\rho,\rho]\times Y)$ denote the connection
in temporal gauge associated to the path $A_\nu$.
Then $F_{\A_\nu}+X_f(\A_\nu)$ converges to zero 
in the $L^2$-norm, by~(\ref{eq:limlong}) and~(\ref{eq:bvp}).
Now it follows from the energy quantization 
in~\cite[Theorems~1.2,~2.1]{W bubb} 
(for general Lagrangians see~\cite{W lag}, 
and for the perturbed version 
see Theorem~\ref{thm:compactness})
that $\A_\nu$ satisfies an $L^\infty$-bound on the curvature.
Hence, by~\cite[Theorem~B]{W elliptic} and 
Theorem~\ref{thm:compactness}, there is a subsequence 
(still denoted by $\A_\nu$) and a sequence of gauge transformations 
${u_\nu\in\cG([-\rho/2,\rho/2]\times Y)}$ such that $u_\nu^*\A_\nu$ 
converges to $\A_\infty=A_\infty(s)+\Phi_\infty(s)\,\ds
\in\cA([-\rho/2,\rho/2]\times Y)$ in the 
$\cC^\infty$-topology. By~(\ref{eq:bvp}) and~(\ref{eq:limlong}) 
the limit connection satisfies
$$
\pd_sA_\infty(s) - \rd_{A_\infty(s)}\Phi_\infty(s)=0,\quad
F_{A_\infty(s)}+X_f(A_\infty(s))=0,\quad
A_\infty(s)|_\Sigma\in\cL
$$
for every $s\in[-\rho/2,\rho/2]$. 
After modifying the gauge transformations $u_\nu$ we may
assume in addition that $\Phi_\infty(s)=0$ 
and $A_\infty(s)=A_\infty$ is independent of~$s$.
It then follows that $u_\nu^{-1}\pd_s u_\nu$ 
converges to zero in the $\cC^\infty$-topology.
So after a further modification we can assume 
that the $u_\nu(s)=u_\nu$ is independent of $s$, 
and so the convergent connections $u_\nu^*\A_\nu$ 
are in temporal gauge, given by the paths
$[-\rho/2,\rho/2] \to \cA(Y): s \mapsto u_\nu^*A_\nu(s) $.
Hence
$$
\lim_{\nu\to\infty} 
\bigl\|A_\nu(0)-(u_\nu^{-1})^*A_\infty
\bigr\|_{W^{1,p}(Y),u_\nu^{-1\;*}A_\infty}
=\lim_{\nu\to\infty} 
\bigl\|(u_\nu^*\A_\nu - \A_\infty)(0)\bigr\|_{W^{1,p},A_\infty}
=0,
$$
$$
\lim_{\nu\to\infty} 
\bigl\|A_\nu(0)-(u_\nu^{-1})^*A_\infty\bigr\|_{L^\infty(Y)}
=\lim_{\nu\to\infty} 
\bigl\|(u_\nu^*\A_\nu - \A_\infty)(0)\bigr\|_{L^\infty(Y)}
=0 ,
$$
$$
\lim_{\nu\to\infty} \left\|\pd_sA_\nu(0)\right\|_{L^\infty(Y)}
=\lim_{\nu\to\infty} \left\|\pd_s(u_\nu^*\A_\nu)(0)\right\|_{L^\infty(Y)}
=\left\|\pd_s A_\infty\right\|_{L^\infty(Y)}
=0.
$$
This contradicts the assumption that (\ref{eq:ass}) fails,
and thus proves the lemma. 
\end{proof}

\begin{lemma}\label{le:decay est}
Suppose that every critical point of the perturbed 
Chern--Simons functional $\CS_\cL+h_f$ is nondegenerate.  
Then, for every $\rho>0$, there are positive constants 
$c_0$ and $\eps$ with the following significance.  
If $A:[-\rho,\rho]\to\cA(Y)$ is a solution
of~(\ref{eq:bvp}) such that 
$$
\int_{-\rho}^{\rho}\int_Y\left|\pd_sA\right|^2\dvol_Y\,\ds < \eps ,
$$
then for every $\alpha\in\Om^1_{A(0)}(Y,\cg)$
$$
\|\alpha\|_{L^6(Y)} + \|\alpha\|_{L^4(\pd Y)}
\le c_0\bigl(
\bigl\|\rd_{A(0)}\alpha+ \rd X_f(A(0))\alpha\bigr\|_{L^2(Y)} 
+ \bigl\|\rd_{A(0)}^* \alpha\bigr\|_{L^2(Y)} \bigr) .
$$
\end{lemma}

\begin{proof}
Assume by contradiction that this is wrong. 
Then there is a constant $\rho>0$, a sequence 
$A_\nu:[-\rho,\rho]\to\cA(Y)$ of solutions of~(\ref{eq:bvp})
with~(\ref{eq:limlong}), 
and a sequence $\alpha_\nu\in\Om^1_{A_\nu(0)}(Y,\cg)$ 
such that 
\begin{align} \label{contrass}
\frac {\|\alpha_\nu\|_{L^6(Y)}
+ \|\alpha_\nu\|_{L^4(\pd Y)} }
{ \bigl\|\rd_{A_\nu(0)}\alpha_\nu
+ \rd X_f(A_\nu(0))\alpha_\nu\bigr\|_{L^2(Y)} 
+ \bigl\|\rd_{A_\nu(0)}^* \alpha_\nu\bigr\|_{L^2(Y)} } 
\;\; \underset{\nu\to\infty}{\longrightarrow}\infty .
\end{align}
Arguing as in the proof of Lemma~\ref{le:decay1} we find 
a subsequence, still denoted by $A_\nu$, and a sequence of
gauge transformations $u_\nu\in\cG(Y)$ such that 
$u_\nu^*A_\nu(0)$
converges in the $\cC^\infty$-topology to a connection
$A_\infty\in\cA(Y,\cL)$ that satisfies 
${F_{A_\infty}+X_f(A_\infty)=0}$.
By assumption $A_\infty$ is 
nondegenerate, so by Corollary~\ref{cor:nondegenerate}
there is a constant $C$ such that
\begin{equation} \label{eq:Ainfty}
\bigl\|(\alpha,0)\bigr\|_{W^{1,2}(Y)} 
\leq C \bigl\|\cH_{A_\infty}(\alpha,0)\bigr\|_{L^2(Y)}
\end{equation}
for every $(\alpha,0)\in{\rm dom}\,\cH_{A_\infty}$.
By Theorem~\ref{thm:Q} this estimate is stable under 
$\cC^1$-small perturbations of $A_\infty$, and by gauge
invariance it continues to hold with $A_\infty$ replaced by 
$A_\nu(0)$. Precisely, let $\cU\subset\cA(Y,\cL)$ be a neighbourhood
of $A_\infty$ and $\{Q_A\}_{A\in\cU}$ be an operator family 
that satisfies the requirements of Theorem~\ref{thm:Q}.
Then $u_\nu^*A_\nu(0)\in\cU$ for large $\nu$ 
adnd the isomorphisms ${Q_\nu:=Q_{u_\nu^*A_\nu(0)}\times\Id}$ 
from ${\rm dom}\,\cH_{A_\infty}$ to ${\rm dom}\,\cH_{u_\nu^*A_\nu(0)}$
converge to 
$Q_{A_\infty}\times\Id=\Id$ in both $\cL(W^{1,2})$ and $\cL(L^2)$;
so the sequence $Q_\nu^{-1}\cH_{u_\nu^*A_\nu(0)}Q_\nu$
has the constant domain ${\rm dom}\,\cH_{A_\infty}$, and it
converges to $\cH_{A_\infty}$ in the operator norm on $\cL(W^{1,2},L^2)$.
Hence, for large~$\nu$,  we can replace $\cH_{A_\infty}$ 
by $Q_\nu^{-1}\cH_{u_\nu^*A_\nu(0)} Q_\nu$ 
in~(\ref{eq:Ainfty}) to obtain estimates with a uniform constant $C$.
Since $Q_\nu$ converges to the identity in the relevant operator 
norms we obtain the following estimate with uniform constants $C_i$ 
but varying domain:
\begin{align*}
 \|\alpha\|_{L^6(Y)} 
+ \|\alpha\|_{L^4(\pd Y)}
&\leq C_1 \bigl\|(\alpha,0)\bigr\|_{W^{1,2}(Y)} 
\leq C_2 \bigl\|\cH_{u_\nu^*A_\nu(0)}(\alpha,0)\bigr\|_{L^2(Y)}
\end{align*}
for every $(\alpha,0)\in{\rm dom}\,\cH_{u_\nu^*A_\nu(0)}$.
Here we have used the Sobolev embedding
$W^{1,2}(Y)\hookrightarrow L^6(Y)$ and the trace theorem
$W^{1,2}(Y)\hookrightarrow L^4(\pd Y)$.
Since $\rT_{u^*A}\cL = u^{-1}(\rT_A\cL) u$
we can apply the last estimate to 
$\bigl( u_\nu^{-1}\alpha_\nu u_\nu , 0\bigr)
\in{\rm dom}\,\cH_{u_\nu^*A_\nu(0)}$.
Since the norms on the left and right hand side are all 
gauge invariant the resulting inequality contradicts~(\ref{contrass}).
This proves the lemma. 
\end{proof}

\begin{lemma}\label{le:decay est D}
Suppose that every critical point of the perturbed 
Chern--Simons functional $\CS_\cL+h_f$ is nondegenerate.  
Then, for every $\rho>\rho'>0$, there are positive 
constants $c_0,c_1,\dots$ and $\eps$ 
with the following significance.  
If $A:[-\rho,\rho]\to\cA(Y)$ is a solution of~(\ref{eq:bvp}) 
such that 
$$
\int_{-\rho}^{\rho}\int_Y\left|\pd_sA\right|^2\dvol_Y\,\ds < \eps ,
$$ 
then, for every smooth path
$[-\rho,\rho]\to\Om^1(Y,\cg)\times\Om^0(Y,\cg):
s\mapsto(\alpha(s),\phi(s))$ 
satisfying $\alpha(s)\in\Om^1_{A(s)}(Y,\cg)$
and every integer $k\geq 0$, we have
\begin{multline*}
\bigl\|(\alpha,\phi)\bigr\|_{\cC^k([-\rho',\rho']\times Y),\A} \\ 
\le c_k\Bigl(
\bigl\|\cD_\A (\alpha,\phi) \bigr\|_{W^{k+2,2}([-\rho,\rho]\times Y),\A} 
+ \bigl\|(\alpha,\phi)\bigr\|_{L^2([-\rho,\rho]\times Y)}
\Bigr).
\end{multline*}
\end{lemma}

\begin{proof}
If this is wrong, then there exist constants 
$k\geq 0$, $\rho>\rho'>0$ and a sequence 
$A_\nu:[-\rho,\rho]\to\cA(Y)$ of solutions of~(\ref{eq:bvp})
with (\ref{eq:limlong}), for which the constant in the estimate blows up.
As in the proof of Lemma~\ref{le:decay1} we find 
a subsequence of the connections on $[-\rho,\rho]\times Y$, 
still denoted by $\A_\nu$, and 
gauge transformations $u_\nu\in\cG(Y)$ such that $u_\nu^*\A_\nu$
converges in the $\cC^\infty$-topology on $[-\rho/2,\rho/2]\times Y$
to a constant connection $\A_\infty=A_\infty\in\cA(Y,\cL)$.
Now by Theorem~\ref{thm:4reg} and the Sobolev embedding theorem, and 
with the norms of Remark~\ref{rmk:norms}, there is a constant $C$ such that
for every $(\alpha,\phi)$ satisfying $\alpha(s)\in\Om^1_{A_\infty}(Y,\cg)$
\begin{multline} 
\bigl\|(\alpha,\phi)\bigr\|_{\cC^k([-\rho',\rho']\times Y),\A_\infty} \\
\leq C \bigl(
\bigl\|\cD_{\A_\infty} (\alpha,\phi) 
\bigr\|_{W^{k+2,2}([-\rho,\rho]\times Y),\A_\infty} 
+ \bigl\|(\alpha,\phi)\bigr\|_{L^2([-\rho,\rho]\times Y),\A_\infty}
\bigr) .\label{eq:AAinfty}
\end{multline}
The same argument as in the proof of Lemma~\ref{le:decay est}
(with the sequence of operators $Q_\nu(s):=Q_{u_\nu^*A_\nu(s)}\times\Id$)
shows that this estimate continues to hold with 
$\A_\infty$ replaced by $u_\nu^*\A_\nu$.
Note that
$\cD_{u_\nu^*\A_\nu}u_\nu^{-1}(\alpha_\nu,\phi_\nu) u_\nu
= u_\nu^{-1}\bigl(\cD_{\A_\nu}(\alpha_\nu,\phi_\nu)\bigr) u_\nu$.
So since the norms are gauge invariant, the above estimate also holds
with $\A_\infty$ replaced by $\A_\nu$, which contradicts the choice of 
$A_\nu$ and thus proves the lemma. 
\end{proof}

\begin{proof}[Proof of Theorem~\ref{thm:decay}.]
The proof has three steps.

\medskip\noindent{\bf Step~1.}
{\it 
There is a uniform constant $\delta>0$ 
(independent of the solution $A$)
and a constant $C$ (which depends on $A$) such that
$$
\bigl\|\pd_sA(s)\bigr\|_{L^2(Y)} \le C e^{-\delta s}
\qquad \text{for}\;s\ge 0.
$$
}
\noindent
Define
$$
g(s):= \frac12 \int_Y \bigl| \pd_s A \bigr|^2
=  \frac12 \int_Y \bigl| F_A + X_f(A) \bigr|^2.
$$
Then
$$
g'(s) =  \int_Y 
\la \bigl(\rd_A\p_s A+\rd X_f(A)\p_sA\bigr) 
\wedge *\bigl(F_A + X_f(A)\bigr) \ra,
$$
and hence
\begin{align*}
g''(s) 
&= \int_Y 
\bigl| \rd_A \pd_s A + \rd X_f(A) \pd_s A  \bigr|^2 
-  \int_Y 
\winner{\bigl([\pd_sA\wedge\pd_s A]+\rd_A\pd_s^2A\bigr)}
{\pd_s A} \\
&\quad
- \int_Y 
\winner{\bigl(\rd^2 X_f(A)(\pd_s A,\pd_s A)+\rd X_f(A)\pd_s^2 A\bigr)}
{\pd_s A} \\
&= \int_Y
\bigl|\rd_A\pd_sA+\rd X_f(A)\pd_sA\bigr|^2 
- \int_Y
\winner{\pd_s^2A}{\bigl(\rd_A\pd_sA+\rd X_f(A)\pd_sA\bigr)} \\
&\quad
- \int_Y 
\winner{\bigl([\pd_sA\wedge\pd_sA]+\rd^2X_f(A)(\pd_sA,\pd_sA)\bigr)}
{\pd_s A}  
- \int_\Sigma \winner{\pd_s^2A}{\pd_sA} \\
&\geq 
2 \bigl\|\rd_A\pd_sA+\rd X_f(A)\pd_sA\bigr\|_{L^2(Y)}^2 \\
&\quad - c_1\bigl\|\pd_sA\bigr\|_{L^\infty(Y)} \bigl\|\pd_sA\bigr\|_{L^2(Y)}^2
- c_1\bigl\|\pd_sA\bigr\|_{L^3(\pd Y)}^3 \\
&\geq  
\left(4\delta^2-c_2\bigl\|\pd_sA\bigr\|_{L^\infty(Y)}\right)
\left(\bigl\|\pd_sA\bigr\|_{L^2(Y)}^2
+ \bigl\|\pd_sA\bigr\|_{L^3(\pd Y)}^2\right) \\
&\geq  
2\delta^2 \bigl\|\pd_sA\bigr\|_{L^2(Y)}^2
\end{align*}
for uniform constants $c_i$ and $\delta>0$ and $s$ sufficiently large.
Here we used~(\ref{eq:bvp}).
In the first inequality the term $\int_\Sigma\winner{\pd_s^2A}{\pd_sA}$ 
is controlled by $\|\pd_s A\|_{L^3(\pd Y)}^3$, see~\cite[Lemma~2.3]{W bubb}
and \cite{W lag} for general Lagrangian submanifolds. 
The first inequality also uses the estimate on $\rd^2 X_f(A)$
from Proposition~\ref{prop:Xf}~(v). 
For the second inequality note that every solution of~(\ref{eq:bvp}) 
satisfies $\pd_s A(s)\in\rT_{A(s)}\cL$ and 
\begin{equation} \label{eq:dsA} 
\begin{split}
*\pd_s A|_{\pd Y} &= -(F_A + X_f(A))|_{\pd Y}=0, \\
\rd_A^*\pd_s A &= *\rd_A(F_A + X_f(A))=0. 
\end{split}
\end{equation}
These identities use~(\ref{Xf identities}) and the Bianchi identity as well as
the facts that the perturbation vanishes near $\pd Y$ and that the 
Lagrangian submanifold $\cL$ is contained in the flat connections on $\pd Y$. 
Now we can apply Lemma~\ref{le:decay est} to the paths
$[-1,1]\to\cA(Y):\sigma\mapsto A(s+\sigma)$ 
(whose derivative is $L^2$-small due to the finite $L^p$-energy of the path)
and to the $1$-forms $\alpha=\pd_s A(s)\in\Om^1_{A(s)}(Y,\cg)$, 
for sufficiently large $s\geq 0$ to obtain
$$
\bigl\| \pd_s A \bigr\|_{L^2(Y)}^2 
+ \bigl\|\pd_s A \bigr\|_{L^3(\pd Y)}^2 
\le 
(2\delta^2)^{-1} \bigl\| \rd_A \pd_s A + \rd X_f(A)\pd_s A  \bigr\|_{L^2(Y)}^2.
$$
Here we have chosen $(2\delta^2)^{-1}= (c_0 c)^2$ with the constant 
$c_0$ from Lemma~\ref{le:decay est} and a further 
Sobolev constant $c$, so $\delta>0$ is independent of the 
solution $A$. The last inequality in the estimate of $g''$ is due to 
$\|\pd_s A(s)\|_{L^\infty(Y)}\leq 2\delta^2 c_2^{-1}$
for $s$ sufficiently large. 
This follows from Lemma~\ref{le:decay1} 
applied to the paths ${[-1,1]\to\cA(Y): \sigma \mapsto A(s+\sigma)}$.
So we have 
$
g''(s)\geq 4\delta^2 g(s)
$
for $s$ sufficiently large.  This implies the assertion of Step~1, 
i.e.~$g(s)\leq C^2 e^{-2\delta s}$, by a standard argument 
(see e.g.~the proof of~\cite[Lemma~2.11]{Sal}).

\medskip\noindent{\bf Step~2.}
{\it Let $\delta>0$ be the constant of Step~1
and $\A\in\cA([0,\infty)\times Y)$ be the connection
associated to the path $A$.
For every integer $k\ge 0$ there is a constant 
$C_k$ such that for every $s\ge1$}
$$
\bigl\|\pd_sA\bigr\|_{\cC^k([s-1,s+1]\times Y),\A}
\le C_ke^{-\delta s}.
$$

\medskip\noindent
Fix $k\geq 0$ and consider the connections 
$\A_\sigma\in\cA([-2,2]\times Y)$ given by
the paths $A_\sigma(s):=A(\sigma+s)$.
Due to the finite $L^p$-energy of $A$ on $[0,\infty)$ for some $p\geq 2$ 
these paths on $[-2,2]$ satisfy 
${\|\pd_s A_\sigma\|_{L^2([-2,2]\times Y)}\to 0}$ as $\sigma\to\infty$.
So by Lemma~\ref{le:decay est D} there is a constant $c_k$ such that 
for all sufficiently large $\sigma$
\begin{multline*}
\bigl\|(\alpha,\phi)\bigr\|_{\cC^k([-1,1]\times Y),\A_\sigma} \\
\le c_k\left(
\bigl\|\cD_{\A_\sigma}(\alpha,\phi)
\bigr\|_{W^{k+2,2}([-2,2]\times Y),\A_\sigma}
+\bigl\|(\alpha,\phi)\bigr\|_{L^2([-2,2]\times Y)}
\right)
\end{multline*}
for every smooth $\phi:[-2,2]\to\Om^0(Y,\cg)$
and $\alpha:[-2,2]\to\Om^1(Y,\cg)$
satisfying $\alpha(s)\in\Om^1_{A_\sigma(s)}(Y,\cg)$.
Now apply the estimate to the pair 
$$
\alpha(s):=\pd_s A(\sigma + s),\qquad \phi(s):=0.
$$  
Differentiate~(\ref{eq:bvp}) and recall~(\ref{eq:dsA}) 
to see that $(\alpha,\phi)\in\ker\cD_{\A_\sigma}$
and hence 
$$
\left\|\pd_s A\right\|_{\cC^k([\sigma-1,\sigma+1]\times Y),\A}
\le c_k\left\|\pd_s A\right\|_{L^2([\sigma-2,\sigma+2]\times Y)}
\le c_k C (2\delta)^{-\frac 12} e^{2\delta} e^{-\delta \sigma}.
$$
The last inequality follows from Step~1 and proves Step~2. 

\medskip\noindent{\bf Step~3.}
{\it Let $\delta>0$ be the constant of Step~1. 
Then there is a connection ${A_\infty\in\cA(Y,\cL)}$ such that
$
F_{A_\infty}+X_f(A_\infty)=0
$
and a sequence of constants $C_0,C_1,C_2,\dots$ such that
\begin{equation}\label{eq:Cklim}
\bigl\|A-A_\infty\bigr\|_{\cC^k([s-1,s+1]\times Y)}
\le C_ke^{-\delta s}
\end{equation}
for every integer $k\ge 0$ and every $s\ge1$.}

\medskip\noindent
By Step~2 we have 
$
\left\|\pd_sA(s)\right\|_{L^\infty(Y)}\le C_0 e^{-\delta s}
$
for every $s\ge 0$.  Hence the integral
$$
A_\infty \,:=\; A(0) + \int_0^\infty\pd_sA(s)\,\ds
\;=\; \lim_{s\to\infty} A(s)
$$
converges in $L^\infty(Y,\rT^*Y\otimes\cg)$ and defines
a $\cC^0$-connection on $Y$.  
This directly implies $A_\infty|_\Sigma\in\cL$.
Moroever, (\ref{eq:Cklim}) holds with $k=0$.
We prove by induction on $k$ that $A_\infty$ is a $\cC^k$
connection that satisfies~(\ref{eq:Cklim}). 
For $k=0$ this is what we have just proved. 
Fix an integer $k\ge 1$ and suppose that $A_\infty$
is a $\cC^{k-1}$ connection that satisfies~(\ref{eq:Cklim})
with  $k$ replaced by $k-1$. Then $\A$ is bounded in $\cC^{k-1}$
and so there is a constant $C$ such that
\begin{equation} \label{norm equiv}
\left\|\alpha\right\|_{\cC^\ell([s-1,s+1]\times Y)}
\le C\left\|\alpha\right\|_{\cC^\ell([s-1,s+1]\times Y),\A}
\end{equation}
for every $\ell\leq k$,  $s\ge 1$, and every 
$\alpha\in\Om^1([s-1,s+1]\times Y,\rT^*Y\otimes\cg)$.
So it follows from Step~2 that 
$$
\left\|\pd_s A \right\|_{\cC^k([s-1,s+1]\times Y)}
\le C C_ke^{-\delta s}.
$$
Hence for $s_1\ge s_0\ge 0$
$$
\left\|A(s_0)-A(s_1)\right\|_{\cC^k(Y)}
\le \int_{s_0}^{s_1}\left\|\pd_s A\right\|_{\cC^k(Y)}\,\ds
\le \frac{C C_k}{\delta}e^{-\delta s_0} .
$$
This shows that $A_\infty$ is a $\cC^k$ connection with
$$
\left\|A(s)-A_\infty\right\|_{\cC^k(Y)}
\le \frac{C C_k}{\delta}e^{-\delta s}.
$$
The exponential decay of ${\pd_s^\ell (A(s)-A_\infty)=\pd_s^\ell A(s)}$ 
in $\cC^{k-\ell}(Y)$ for $\ell=1,\dots,k$ follows from Step~2 and
(\ref{norm equiv}), so this implies~(\ref{eq:Cklim}). 
Moreover,
$$
F_{A_\infty}+X_f(A_\infty)
=\lim_{s\to\infty} \bigl( F_{A(s)}+X_f(A(s))\bigr)
=-\lim_{s\to\infty}*\pd_s A(s)=0.
$$
This proves Step~3 and the lemma.
\end{proof}

\begin{proof}[Proof of Theorem~\ref{thm:long}.]
Let $\delta>0$ be the constant of Step~1 in the proof 
of Theorem~\ref{thm:decay}.
We prove that there are constants $C$ and $\eps>0$
such that the following holds for every $T\ge1$.
If $A:[-T,T]\to\cA(Y)$ is a solution
of~(\ref{eq:bvp}) that satisfies~(\ref{eq:longenergy}),
then it also satisfies
\begin{equation}\label{eq:L2exp}
\left\|\pd_sA(s)\right\|_{L^2(Y)}
\le C e^{-\delta(T-|s|)}\left\|\pd_sA
\right\|_{L^2(([-T,1-T]\cup[T-1,T])\times Y)}
\end{equation}
for $|s|\le T-1/2$.  
Let $\eps>0$ be the constant of Lemma~\ref{le:decay est} with $\rho=\frac 14$
and assume that~(\ref{eq:longenergy}) holds with this 
constant $\eps$.  Define $f:[-T,T]\to\R$ by
$$
f(s):=\tfrac12\|\pd_sA(s)\|_{L^2(Y)}^2.
$$
Then the same argument as in Step~1 in the proof 
of Theorem~\ref{thm:decay} shows that there is 
a constants $c_2$, independent of $A$, such that for $|s|\le T-1/4$  
$$
f''(s)\ge 2 \left(4\delta^2-c_2\left\|\pd_sA(s)\right\|_{L^\infty(Y)}\right)
\left( f(s) + \|\pd_s A(s)\|_{L^3(\pd Y)}^2 \right).
$$
Shrinking $\eps$ if necessary we may assume,
by Lemma~\ref{le:decay1} with $\rho=1/4$, 
that $\left\|\pd_sA(s)\right\|_{L^\infty(Y)}\le2\delta^2/c_2$ and hence 
$$
f''(s)\ge 4\delta^2 f(s)\qquad \text{for}\; |s|\le T-1/4.
$$
Now~(\ref{eq:L2exp}) follows from Lemma~\ref{le:long} below with $\rho=1/4$, 
$\delta$ replaced by $2\delta$, and $T$ replaced by $T-1/4$.

Integration of~(\ref{eq:L2exp}) yields
$$
\left\|\pd_sA\right\|_{L^2([\sigma-3/2,\sigma+3/2]\times Y)}
\le C'' e^{-\delta(T-|\sigma|)}
\left\|\pd_sA\right\|_{L^2(([-T,1-T]\cup[T-1,T])\times Y)}
$$
for every $\sigma\in[-T+2,T-2]$ with $C''=C e^{3\delta/2}\delta^{-1/2}$.
Now, shrinking $\eps$ if necessary, we can apply 
Lemma~\ref{le:decay est D} with $\rho=3/2$ and $\rho'=1$ to the 
paths shifted by $\sigma$.
Since $(\pd_s A,0)\in\ker\cD_\A$
(as in Step~2 of the proof of Theorem~\ref{thm:decay})
we obtain constants $C_k$ and $C_k'$ for every $k\geq 0$ such that
\begin{align*}
\left\|\pd_sA\right\|_{\cC^k([\sigma-1,\sigma+1]\times Y),\A}
&\leq  C_k' \left\|\pd_sA\right\|_{L^2([\sigma-3/2,\sigma+3/2]\times Y)} \\
&\le C_k e^{-\delta(T-|\sigma|)}
\left\|\pd_sA\right\|_{L^2(([-T,1-T]\cup[T-1,T])\times Y)} .
\end{align*}
for every $\sigma\in[-T+2,T-2]$.
Taking the supremum over $\sigma\in[-s+1,s-1]$ then proves the
assertion \eqref{eqn:dsA} on $\pd_s A$.

To prove \eqref{eqn:a-a0} it remains to estimate the derivatives tangent to $Y$.
We fix any two constants $\kappa>0$ and $p>1$ 
and then, by Lemma~\ref{le:decay1}, 
find a connection ${A_0\in\cA(Y,\cL)}$ such that
$F_{A_0}+X_f(A_0)=0$ and 
$$
\Norm{A(0)-A_0}_{W^{1,p}(Y),A_0}
+ \Norm{A(0)-A_0}_{L^\infty(Y)}\leq\kappa.
$$
After a gauge transformation on $A_0$ we can assume that $A(0)$ 
lies in the local slice $S_{A(0)}$ of $A(0)$, that is
$\rd_{A_0}^*(A(0)-A_0)=0$ and $*(A(0)-A_0)|_{\pd Y}=0$.
Since all critical points are nondegenerate,
Corollary~\ref{cor:nondegenerate} provides a universal 
constant $c_0$ depending on $q>\max\{3,p\}$
such that 
for all $\alpha\in\Om^1(Y,\cg)$ with $*\alpha|_{\pd Y}=0$
\begin{multline*}
\Norm{\alpha}_{L^\infty(Y)} + \Norm{\alpha}_{W^{1,p}(Y),A_0} \\
\le c_0 \left( 
\Norm{\rd_{A_0}\alpha + \rd X_f(A_0)\alpha}_{L^q(Y)}
+ \Norm{\rd_{A_0}^*\alpha}_{L^q(Y)} 
+ \Norm{\Pi_{T_A\cL}^\perp (\alpha|_{\pd Y})}_{L^q(\pd Y)} 
\right) .
\end{multline*}
When applying this to $\alpha=A(0)-A_0$ we can use the estimate
$$
\Norm{\Pi_{T_A\cL}^\perp (\alpha|_{\pd Y})}_{L^q(\pd Y)} 
\leq c_1 \Norm{(A(0)-A_0)|_{\pd Y}}_{L^q(\pd Y)}^2
$$
with a uniform constant $c_1$ since $A(0)|_{\pd Y}$ and $A_0|_{\pd Y}$ 
both lie in the submanifold $\cL\subset\cA^{0,q}(\Sigma)$.
More precisely, we abbreviate $A'_0:=A_0|_{\pd Y}$, 
then we can use the exponential map in Lemma~\ref{lem:lagg}
to write 
$$
A(0)|_{\pd Y}=\Theta_{A_0'}(\beta)
=A'_0 + \beta 
+ \int_0^1 \bigl( D\Theta_{A'_0}(\tau\beta) - D\Theta_{A'_0}(0) \bigr) \beta \; \rd\tau
$$
for some $\beta\in\rT_{A_0}\cL$, using the identities
$\Theta_A(0)=A$ and $D\Theta_{A}={\rm Id}$.
The map $\Theta$ is smooth and gauge invariant, and $\cL/\cG(\Sigma)$ is compact, so 
by the choice of $\kappa>0$ we obtain arbitrarily small bounds on $\|\beta\|_{L^q(\Sigma)}$
and a uniform linear bound $\|D\Theta_{A'_0}(\tau\beta) - D\Theta_{A'_0}(0)\|\leq c'\|\beta\|_{L^q(\Sigma)}$.
This implies the uniform estimate
\begin{align*}
\bigl\| \bigl(A(0)|_{\pd Y} - A'_0 \bigr) - \beta \bigr\|_{L^q(\Sigma)}
\leq c'\|\beta\|_{L^q(\Sigma)}^2
\leq c_1 \|(A(0)-A_0)|_{\pd Y}\|_{L^q(\pd Y)}^2 .
\end{align*}
We also use the identity
$\rd_{A_0}\alpha = F_{A(0)} - F_{A_0} - \frac 12 [\alpha\wedge\alpha]$
to obtain
\begin{align*}
&\|A(0)-A_0\|_{L^\infty(Y)} + \|A(0)-A_0\|_{W^{1,p}(Y),A_0} \\
&\leq 
c_0 \Bigl( \|F_{A(0)} + X_f(A(0)) \|_{L^q(Y)}
+ \|\tfrac 12 [\alpha\wedge\alpha]\|_{L^q(Y)} \\
&\qquad\quad
+ \|X_f(A_0+\alpha) - X_f(A_0) - \rd X_f(A_0)\alpha\|_{L^q(Y)} 
+ c_1 \|\alpha|_{\pd Y}\|_{L^q(\pd Y)}^2 \Bigr) \\
&\leq 
c_0 \|\pd_s A(0) \|_{L^q(Y)} + c_2 \kappa \|A(0)-A_0\|_{L^\infty(Y)} .
\end{align*}
Here $c_2$ is another uniform constant and we have used
Proposition~\ref{prop:Xf}~(v) for the perturbation term.
If we choose $\kappa=(2c_2)^{-1}$ and the corresponding
$\eps>0$ from Lemma~\ref{le:decay1}, then this proves
$$
\|A(0)-A_0\|_{L^\infty(Y)}+ \|A(0)-A_0\|_{W^{1,p}(Y),A_0}
\leq 2 c_0 \|\pd_s A(0) \|_{L^q(Y)}.
$$
Now \eqref{eqn:a-a0} follows by integrating over the 
estimate \eqref{eqn:dsA} for $\pd_s A$.
\end{proof}

\begin{lemma}\label{le:long}
For every $\delta>0$ and every $\rho>0$ there exists a constant $C$
such that the following holds.  If $T\ge\rho$ and 
$f:[-T,T]\to\R$ is a $\cC^2$-function satisfying
\begin{equation}\label{eq:f}
f''(s)\ge \delta^2 f(s),\qquad f(s)\ge 0
\end{equation}
for all $s\in[-T,T]$, then
\begin{equation}\label{eq:f1}
f(s) \le C e^{-\delta(T-|s|)}E_\rho(f)
\end{equation}
for all $|s|\le T-\rho$, where 
$$
E_\rho(f):=\int_{-T}^{\rho-T}f(s)\,\ds + \int_{T-\rho}^{T}f(s)\,\ds.
$$
\end{lemma}

\begin{proof}
We claim that there is a constant $C_0=C_0(\delta,\rho)>0$
such that every $\cC^2$-function $f:[-T,T]\to\R$ with $T\ge\rho$ 
that satisfies~(\ref{eq:f}) also satisfies
\begin{equation}\label{eq:long}
f'(s)-\delta f(s) \ge -C_0 e^{-\delta T}E_\rho(f)
\end{equation}
for all $0\leq s\leq T$.
To see this note that, for every $s\in[-T,T]$, we have
$$
\frac{\rd}{\rd s} e^{\delta s}\bigl(f'(s)-\delta f(s)\bigr)
= e^{\delta s}\bigl(f''(s)-\delta^2f(s)\bigr)\ge 0.
$$
Hence 
$$
f'(s)-\delta f(s) \ge e^{\delta(r-s)}\bigl(f'(r)-\delta f(r)\bigr)
$$
for all $-T\le r \le s \le T$.
Integrating this over the interval 
$t\le r\le t+\rho/2$ for 
$-T\le t\le -\rho/2$ 
$-T\le t\le \rho/2 - T$ 
and $s\geq 0$ gives
\begin{align*}
f'(s)-\delta f(s)
&\ge  \frac{2e^{-\delta s}}{\rho}
\int_{t}^{t+\rho/2}e^{\delta r}\bigl(f'(r)-\delta f(r)\bigr)\,\dr \\
&= \frac{2e^{-\delta s}}{\rho}
\int_{t}^{t+\rho/2}
\left(\frac{\rd}{\rd r}
(e^{\delta r}f(r))-2\delta e^{\delta r} f(r)\right)\,\dr \\
&\ge 
- \frac{2}{\rho}e^{\delta t} f(t) 
- \frac{4\delta e^{\delta\rho/2}}{\rho}e^{\delta t}E_\rho(f) .
\end{align*}
Integration over the interval 
$-T\le t\le\rho/2-T$ yields~(\ref{eq:long}) with 
$C_0:=12\rho^{-2}e^{\delta\rho}$.
By~(\ref{eq:long}), we have
\begin{align*}
\frac{\rd}{\rd s} e^{-\delta s}f(s) 
&= e^{-\delta s}\bigl(f'(s)-\delta f(s)\bigr) 
\ge -C_0e^{-\delta(s+T)}E_\rho(f)
\end{align*}
for $0\le s\le T$ and hence
$$
e^{-\delta t}f(t)-e^{-\delta s}f(s) \ge - C_1e^{-\delta T}E_\rho(f)
$$
for $0\le s\le t\le T$, where $C_1:=C_0/\delta$. 
For $s\le T-\rho\le t\le T$ this implies
\begin{align*}
f(s)
&\le e^{\delta(s-t)}f(t) + C_1e^{\delta(s-T)}E_\rho(f) 
\le e^{\delta(s-T)}\bigl(e^{\delta\rho}f(t)+C_1E_\rho(f)\bigr).
\end{align*}
Integrating this inequality over the interval $T-\rho\le t\le T$ 
gives~(\ref{eq:f1}) for $0\le s\le T-\rho$
with $C:=C_1+\rho^{-1}e^{\delta\rho}$.
To prove the estimate for $-T+\rho \le s \le 0$ 
replace $f$ by the function $s\mapsto f(-s)$.
\end{proof}

We close this section with a useful exponential 
estimate for the solutions of the linearized equation.

\begin{thm}\label{thm:lin-exp}
Let $A:[0,\infty)\to\cA(Y,\cL)$ be a finite energy
solution of~(\ref{eq:bvp}) that converges to a nondegenerate 
critical point $A^+\in\cA(Y,\cL)$ of $\CS_\cL+h_f$. 
Then there exists a constant $\delta>0$ with the 
following significance.
If $\alpha:[0,\infty)\to\Om^1(Y,\cg)$
is a smooth solution of the equation
$$
\p_s\alpha(s) = *\bigl(
\rd_{A(s)}\alpha+\rd X_f(A(s))\alpha(s)\bigr),\qquad
\rd_{A(s)}^*\alpha(s)=0
$$
satisfying the boundary conditions
$\alpha(s)|_\Sigma\in\rT_{A(s)}\cL$
and $*\alpha(s)|_\Sigma=0$, and
$$
\int_0^\infty e^{-\delta s}\Norm{\alpha(s)}^2_{L^2(Y)}\,\ds 
< \infty ,
$$
then there are constants $C_k$ such that, 
for every $s\ge1$ and every integer $k\ge 0$,
$$
\Norm{\alpha}_{\cC^k([s-1,s+1]\times Y)}
\le C_ke^{-\delta s}.
$$
\end{thm}

\begin{proof}
We prove first that
\begin{equation}\label{eq:expL2}
\Norm{\alpha(s)}^2_{L^2(Y)}
\le Ce^{-\delta s}.
\end{equation}
Since the limit connection is nondegenerate,
Corollary~\ref{cor:nondegenerate} provides an estimate
$$
\Norm{\alpha(s)}_{W^{1,2}(Y)}
\le c\Norm{\rd_{A(s)}\alpha+\rd X_f(A(s))\alpha(s)}_{L^2(Y)}
$$
for $s$ sufficiently large.  
This implies that the function 
$$
g(s) := \tfrac12\Norm{\alpha(s)}_{L^2(Y)}^2
$$
satisfies 
\begin{align} \label{eqn:f''}
g''(s) &= \Norm{\pd_s\alpha}^2_{L^2}
+ \inner{(\rd_{A}+\rd X_f(A))\pd_s\alpha}{\alpha}
+ \inner{\bigl([\pd_s A,\alpha] 
+ \rd^2 X_f(\pd_s A,\alpha)\bigr)}{\alpha} \nonumber\\
&\geq 2 \Norm{\rd_{A}\alpha + \rd X_f(A)\alpha}_{L^2(Y)}^2 
+ \int_{\pd Y}\winner{\pd_s\alpha|_{\pd Y}}{\alpha|_{\pd Y}}
- C \|\pd_s A\|_\infty \Norm{\alpha}^2_{L^2(Y)} \nonumber\\
&\geq 2 c^{-2}\Norm{\alpha}_{W^{1,2}(Y)}^2 
- C \|\pd_s A\|_\infty \Norm{\alpha|_{\pd Y}}^2_{L^2(\pd Y)}
- C \|\pd_s A\|_\infty \Norm{\alpha}^2_{L^2(Y)} \nonumber\\
&\geq \delta^2 g(s)
\end{align}
for some $\delta>0$ and all $s\ge s_0$.
Here we used Proposition~\ref{prop:Xf}~(v) to estimate
$\|\rd^2 X_f(\pd_s A,\alpha)\|_{L^2(Y)}$
and Theorem~\ref{thm:P} to write $\alpha(s)|_{\pd Y}=P_{A(s)|_{\pd Y}}\beta(s)$
for tangent vectors $\beta(s)\in\rT_{A_0}\cL$ at the limit connection
$A_0:=\lim_{s\to\infty}A(s)|_{\pd Y}$.
This gives the estimate
$$
\int_{\pd Y}\winner{\pd_s\alpha|_{\pd Y}}{\alpha|_{\pd Y}}
=\int_{\pd Y}\winner{\bigl(\pd_s P_{A(s)|_{\pd Y}}\bigr)\beta}{\alpha|_{\pd Y}}
\leq C \|\pd_s A\|_\infty \|\alpha|_{\pd Y}\|_{L^2(\pd Y)}^2.
$$
The final inequality in (\ref{eqn:f''}) follows from the exponential decay
of $\pd_s A$ (see Theorem \ref{thm:decay}) with any $\delta<2 c^{-1}$ and 
sufficiently large $s_0$.
This shows that the function 
$h(s):=e^{-\delta s}(g'(s)+\delta g(s))$
is monotonically increasing for $s\ge s_0$. 
We claim that $h(s)\le 0$ for all $s\ge s_0$.  
Suppose otherwise that there is an $s_1\ge s_0$ 
such that $c_1:=h(s_1)>0$.  Then $h(s)\ge c_1$
for all $s\ge s_1$, hence
$$
\frac{d}{ds}(e^{\delta s}g(s)) 
= e^{2\delta s} h(s) \ge e^{2\delta s} c_1,\qquad s\ge s_1,
$$
and hence, by integration, 
$$
e^{\delta s}g(s) \ge 
\frac{c_1}{2\delta} e^{2\delta s}
- \left( \frac{c_1}{2\delta} e^{2\delta s_1} 
- e^{\delta s_1}g(s_1) \right) .
$$
But this means that the function 
$s\mapsto e^{-\delta s}g(s)$ is not integrable,
in contradiction to our assumption.
Thus we have proved that $h(s)\le 0$
and hence $g'(s)\le -\delta g(s)$ for every $s\ge s_0$.
Hence either $g$ vanishes identically for all 
sufficiently large $s$ or $g>0$ for all $s\ge s_0$ 
and $(\log g)'\le -\delta$. This proves~(\ref{eq:expL2}).

To obtain bounds on the derivatives of $\alpha$ we use 
Theorem~\ref{thm:4reg}~(ii) with $\cD_\A$ replaced by the
adjoint $-\cD_\A^*=\nabla_s - \cH_A$. 
Since $A(s)$ converges in the $\Cinf$ topology for $s\to\infty$
we obtain
$\|\alpha\|_{W^{k+1,2}([s-1,s+1]\times Y)} 
\le C_k \|\alpha\|_{L^2([s-2,s+2]\times Y)}$
with a uniform constant $C_k$ for each integer $k$ and all $s\ge2$.
The result then follows from the Sobolev embeddings 
$W^{k+3,2}([-1,1]\times Y)\hookrightarrow \cC^k([-1,1]\times Y)$.
\end{proof}


\section{Moduli spaces and Fredhom theory}\label{sec:fredholm}

In this section we set up the Fredholm theory
for the boundary value problem~(\ref{eq:asd}).
For the purpose of this paper we could restrict the
discussion to the case of a tube $\R\times Y$ as base manifold.
In view of a future definition of product structures however,
we take some time to introduce a more general class of base manifolds
and develop the basic Fredholm theory for these.
For the index computations we then restrict to the case of a tube.
We begin by introducing the basic setup followed by a discussion
of the relevant moduli spaces. The main part of this section then
discusses the properties of the linearized operators.

\subsection*{Instanton data}

\begin{dfn}\label{def:4}
A {\bf $\mathbf 4$-manifold with boundary space-time 
splitting and tubular ends} is a triple 
$(X,\tau,\iota)$ consisting of 
\begin{description}
\item[$\bullet$]
an oriented smooth $4$-manifold $X$ 
with boundary, 
\item[$\bullet$]
a tuple $\tau = (\tau_1,\dots,\tau_m)$
of orientation preserving embeddings
$$
\tau_i:\cS_i\times\Sigma_i\to \pd X, \qquad i=1,\dots m,
$$
where each $\Sigma_i$ is a compact oriented
$2$-manifold and each $\cS_i$ is either $\R$ or 
${S^1\cong\R/\Z}$, 
\item[$\bullet$]
a tuple $\iota=(\iota_1,\dots,\iota_n)$
of orientation preserving embeddings
$$
\iota_j:(0,\infty)\times Y_j\to X,\qquad j=1,\dots,n,
$$
where $Y_j$ is a compact oriented $3$-manifold with boundary,
\end{description}
satisfying the following conditions.
\begin{description}
\item[(i)]
The images of the embeddings $\tau_1,\dots,\tau_m$ have disjoint
closures and
$$
\pd X = \bigcup_{i=1}^m \tau_i(\cS_i\times\Sigma_i).
$$
\item[(ii)]
For $j=1,\dots,n$ the image $U_j:=\iota_j((0,\infty)\times Y_j)$ 
of $\iota_j$ is an open subset of $X$, the closures of the sets
$U_j$ are pairwise disjoint, and the set $X\setminus\bigcup_{j=1}^nU_j$
is compact. 
\item[(iii)]
For every $j\in\{1,\dots,n\}$ there is a subset
$I_j\subset\{1,\dots,m\}$ and a map
$\eps_j:I_j\to\{\pm1\}$ such that
$$
\pd Y_j = \bigsqcup_{i\in I_j}\Sigma_i,\qquad
\iota_j(s,z)=\tau_i(\eps_j(i)(s+1),z)
$$
for $s>0$, $i\in I_j$, and $z\in\Sigma_i$.
The orientation of $\Sigma_i$ coincides with
the boundary orientation of $Y_j$ iff $\eps_j(i)=-1$.  
\end{description}
\end{dfn}

\begin{dfn}\label{def:4g}
Let $(X,\tau,\iota)$ be a $4$-manifold with 
boundary space-time splitting and tubular ends.
A Riemannian metric $g$ on $X$ is called 
{\bf compatible with the boundary space-time splitting
and the tubular ends} if 
\begin{description}
\item[(i)]
on each tubular end the metric is of split form 
$$
\iota_j^*g=\ds^2 + g_j,
$$ 
where $g_j$ is a metric on $Y_j$ independent 
of $s\in(0,\infty)$,
\item[(ii)]
each $\tau_i$ can be extended to an embedding 
$\bar\tau_i:\cS_i\times[0,\eps_i)\times\Sigma_i \to X$
for some $\eps_i>0$ such that 
$$
\bar\tau_i^*g=\ds^2+\dt^2+g_{i,s,t},
$$
where $g_{i,s,t}$ is a smooth family of 
metrics on $\Sigma_i$.
\end{description}
A quadruple $(X,\tau,\iota,g)$ with these properties 
is called a {\bf Riemannian $\mathbf 4$-man\-i\-fold with boundary
space-time splitting and tubular ends}.
\end{dfn}

\begin{rmk}\rm
{\bf (i)}
On the tubular ends condition~(ii) in Definition~\ref{def:4g} 
follows from~(i).
Indeed, on $U_j$ the extension $\bar\tau_i$
for $i\in I_j$ is obtained by composing $\iota_j$ with 
the embedding $[0,\eps)\times\Sigma_i\to Y_j$
associated to geodesic normal coordinates.

\smallskip\noindent{\bf (ii)}
Let $(X,\tau,\iota,g)$ be a Riemannian $4$-manifold with boundary
space-time splitting and tubular ends.  Then $X$ can be exhausted 
by compact deformation retracts. Hence the triple $(X,\tau,g)$
is a Riemannian $4$-manifold with a boundary space-time splitting 
in the sense of \cite[Definition~1.2]{W elliptic}.
\end{rmk}

\begin{ex}\label{ex:RY}\rm
Let $Y$ be a compact oriented $3$-manifold with nonempty boundary
${\pd Y=\Sigma}$. Then 
$X:=\R\times Y$
satisfies the requirements of Definition~\ref{def:4}
with the obvious inclusion $\tau:\R\times\Sigma\to\pd X$,
$Y_1:=Y$, $Y_2:=\bar Y$ (which has the reversed orientation), 
$\iota_1(s,y):=(s+1,y)$, $\iota_2(s,y):=(-s-1,y)$. 
For any metric $g_Y$ on $Y$ the metric $\ds^2+g_Y$ 
on $\R\times Y$ satisfies the conditions of 
Definition~\ref{def:4g}. If $g_\pm$ are two metrics
on $Y$ then, by~\cite[Example~1.4]{W elliptic}, there 
is a metric $g$ on $\R\times Y$ that satisfies
the conditions of Definition~\ref{def:4g}
and has the form $g=\ds^2+g_\pm$ for $\pm s\ge 1$. 
\end{ex}

The following result will be needed in the proof of
independence of the Floer homology from the choice of a metric.

\begin{lemma}\label{le:metric}
Let $(X,\tau,\iota)$ be a $4$-manifold with boundary
space-time splitting and tubular ends and, 
for $j=1,\dots,n$, let $g_j$ be a metric on $Y_j$. 
Then there is a metric $g$ on $X$, compatible 
with the boundary space-time splitting and the 
tubular ends, such that~(i) in Definition~\ref{def:4g}
holds with the given metrics $g_j$. 

Moreover, the space of such metrics $g$ is contractible
if we restrict the consideration to those metrics 
with $\eps_i\geq\eps$ in (ii) for any fixed $\eps>0$.
\end{lemma}

\begin{proof}
The construction of a metric with given ends works as in
\cite[Example~1.4]{W elliptic}.
Denote by $\Met(X,\tau,\iota)$ the set of 
metrics on $X$ that satisfy~(i) in Defintion~\ref{def:4g}
and $\tau_i^*g=\ds^2+g_{i,s}$ for $i=1,\dots,m$ and some families
of metrics $(g_{i,s})_{s\in\cS_i}$ on $\Sigma_i$. 
Then $\Met(X,\tau,\iota)$ is convex and hence contractible. 
Fix $\eps>0$ and let $\Met_\eps(X,\tau,\iota)\subset\Met(X,\tau,\iota)$
denote the subset of all metrics that are compatible
with the boundary space-time splitting and the tubular ends
as in Definition~\ref{def:4g} with $\eps_i\geq\eps$ in (ii).  
To prove that $\Met_\eps(X,\tau,\iota)$ is contractible
it suffices to construct a continuous left inverse of the inclusion
$\Met_\eps(X,\tau,\iota)\hookrightarrow\Met(X,\tau,\iota)$.

Every metric $g\in\Met(X,\tau,\iota)$ determines
embeddings 
$$
\bar\tau_{g,i}:\cS_i\times[0,\delta)\times\Sigma_i\to X
$$ 
defined by 
$$
\bar\tau_{g,i}(s,t,z) := \exp_{\tau_i(s,z)}(t\,\nu_i(s,z)) ,
$$
where $\nu_i:\cS_i\times\Sigma_i\to\tau_i^*\rT X$ 
denotes the inward unit normal.
The constant $\delta>0$ for which the $\bar\tau_{g,i}$ are embeddings
can be chosen uniform on a $\cC^1$-neighbourhood of the metric.
Taking a locally finite refinement of the cover 
of $\Met(X,\tau,\iota)$ by these neighbourhoods 
and using a partition of unity one can construct a function
$$
\bar\delta:\Met(X,\tau,\iota)\to(0,\eps],
$$
continuous with respect to the $\cC^\infty$-topology, 
such that the maps $\bar\tau_{g,i}$ are embeddings for
$0<\delta\le\bar\delta(g)$. 

For $g\in\Met(X,\tau,\iota)$ and $i=1,\dots,m$ define 
the metrics $h_{g,i}$ on the strips 
${\cS_i\times[0,\bar\delta(g))\times\Sigma_i}$
by 
$$
h_{g,i}:=\ds^2+\dt^2+g_{i,s,t}, 
$$
where the metric $g_{i,s,t}$ on $\Sigma_i$ is the pullback
of the metric on $X$ under the embedding 
$z\mapsto\bar\tau_{g,i}(s,t,z)$. 
We fix a smooth cutoff function $\lambda:[0,1]\to[0,1]$
such that $\lambda(t)=0$ for $t$ near $0$ and $\lambda(t)=1$
for $t$ near~$1$. 
Then for $\delta>0$ we define 
$\lambda_\delta:\cS_i\times[0,\delta)\times\Sigma_i\to[0,1]$
by 
$$
\lambda_\delta(s,t,z):=\lambda(t/\delta).
$$
Now we can define the map 
$\Met(X,\iota,\tau)\to\Met_\eps(X,\iota,\tau):g\mapsto \tilde g$
by 
$$
\tilde g :=(\bar\tau_{g,i})_*\bigl(
\lambda_{\bar\delta(g)}\bar\tau_{g,i}^*g
+\bigl(1-\lambda_{\bar\delta(g)}\bigr)h_{g,i}
\bigr)
$$
on the image of $\bar\tau_{g,i}$ for $i=1,\dots,m$
and by $\tilde g:=g$ on the complement. This map is 
the identity on $\Met_\eps(X,\tau,\iota)$ since 
$\eps_i\geq\eps\geq\bar\delta(g)$.
So we have constructed the required left inverse of the inclusion 
$\Met_\eps(X,\tau,\iota)\hookrightarrow\Met(X,\tau,\iota)$.
\end{proof}

\begin{dfn}\label{def:data}
Let $(X,\tau,\iota)$ be a $4$-manifold with boundary space-time 
splitting and tubular ends.  {\bf Instanton data} on $X$  are given 
by a triple $(g,\cL,f)$ with the following properties.
\begin{description}
\item[$\bullet$] 
$g$ is a Riemannian metric on $X$ compatible with the 
boundary space-time splitting and the tubular ends.
\item[$\bullet$]
$\cL=(\cL_1,\dots,\cL_m)$ is an $m$-tuple of gauge invariant, monotone
Lagrangian submanifolds $\cL_i\subset\cA(\Sigma_i)$, satisfying (L1-2)
on page~\pageref{p:L1}. 
\item[$\bullet$]
$X_f:\cA(X)\to\Om^2(X,\cg)$ is a holonomy perturbation as in the introduction
such that, on every tubular end and for every $\A\in\cA(X)$, the $2$-form 
$\iota_j^*X_f(\A)\in\Om^2((0,\infty)\times Y_j,\cg)$ is induced by the path
$s\mapsto X_{f_j}(A_j(s))$, where $\iota_j^*\A=:A_j(s)+\Phi_j(s)\ds$.
Here $X_{f_j}:\cA(Y_j)\to\Om^2(Y_j,\cg)$ is as in~(\ref{Xf}).
The perturbation $f$ involves a choice of thickened loops, 
i.e.\ embeddings $\gamma_i:S^1\times\Om\to{\rm int}(X)$,
where $\Om\subset\R^3$ is a contractible open set.
\end{description}
\end{dfn}

\subsection*{The moduli space}

Let $(X,\tau,\iota)$ be a $4$-manifold with boundary space-time 
splitting and tubular ends and let $(g,\cL,f)$ be instanton data on $X$.
The perturbed anti-self-duality equation with Lagrangian
boundary conditions has the form
\begin{equation}\label{eq:asdX}
F_\A + X_f(\A) + *\bigl(F_\A+X_f(\A)\bigr)=0,\qquad 
\tau_{i,s}^*\A\in\cL_i \quad \forall s\in\cS_i.
\end{equation}
Here the embedding $\tau_{i,s}:\Sigma_i\to X$ is defined by 
$\tau_{i,s}(z):=\tau_i(s,z)$. 
The energy of a solution is
$$
E_f(\A):= \frac 12 \int_X \bigl| F_\A + X_f(\A) \bigr|^2 .
$$
By Theorem~\ref{thm:decay} every finite energy solution of 
(\ref{eq:asdX}) that is in temporal gauge on the tubular ends
converges to critical points $A_j$ of the 
perturbed Chern--Simons functionals, i.e.
\begin{equation}\label{eq:limX}
\lim_{s\to\infty}\left\|\iota_j^*\A-A_j\right\|_{\cC^k([s-1,s+1]\times Y_j)}=0
\end{equation}
for every $j\in\{1,\dots,n\}$ and every integer $k\ge 0$.  
This equation is understood as follows.
We denote by $\cA(X,\cL)$ the set of smooth connections $\A\in\cA(X)$
that satisfy the Lagrangian boundary conditions 
$\tau_{i,s}^*\A\in\cL_i$ for all $i\in\{1,\ldots,m\}$ and $s\in\cS_i$.  
On a tubular end, any such connection decomposes as
$$
\iota_j^*\A = B_j + \Phi_j\ds
$$
with $\Phi_j:(0,\infty)\to\Om^0(Y_j,\cg)$ and 
$B_j:(0,\infty)\to\cA(Y_j,\cL)$.  
Here $\cA(Y_j,\cL)$ denotes the set of smooth connections $B\in\cA(Y_j)$ 
that satisfy the Lagrangian boundary conditions $B|_{\Sigma_i}\in\cL_i$ 
for all $i\in I_j$.
The temporal gauge condition means that $\Phi_j\equiv0$.
For $j=1,\dots,n$ the connection $A_j\in\cA(Y_j,\cL)$ in~(\ref{eq:limX}) is
a critical point of the perturbed Chern--Simons functional for $Y_j$, i.e.
$$
F_{A_j}+X_{f_j}(A_j)=0.
$$
The space of solutions of (\ref{eq:asdX}) and (\ref{eq:limX}) 
that are in temporal gauge on the tubular ends
will be denoted by 
$$
\widetilde\cM(A_1,\dots,A_n;X_f)\subset\cA(X,\cL).
$$
Let us denote by $\cG_{A_j}\subset\cG(Y_j)$ the isotropy
subgroup of $A_j$.  Then the group
$\cG(A_1,\dots,A_n)$ of all gauge transformations $u\in\cG(X)$ 
that satisfy $u\circ\iota_j\equiv u_j\in\cG_{A_j}$ for $j=1,\dots,n$, 
acts on the space $\widetilde\cM(A_1,\dots,A_n;X_f)$.
The quotient will be denoted by 
\begin{equation}\label{eq:GAA}
\cM(A_1,\dots,A_n;X_f):=
\widetilde\cM(A_1,\dots,A_n;X_f)/
\cG(A_1,\dots,A_n).
\end{equation}
In the case of the tube $X=\R\times Y$, this moduli space
can easily be identified with the one that is mentioned
in the introduction.
Similarly, the moduli space $\cM(A'_1,\dots,A'_n;X_f)$
for gauge equivalent limits $A'_i\in[A_i]$ can be
identified with $\cM(A_1,\dots,A_n;X_f)$.

\subsection*{The linearized operator}

Fix critical points $A_j\in\cA(Y_j,\cL)$, $j=1,\dots,n$,
of the perturbed Chern--Simons functionals and let 
$\A\in\cA(X,\cL)$ be a connection satisfying (\ref{eq:limX}).
Denote by $\Om^1_\A(X,\cg)$ the space of smooth 
$1$-forms that satisfy the boundary conditions
\begin{equation}\label{eq:alphabc}
*\alpha|_{\pd X}=0, \qquad
\tau_{i,s}^*\alpha\in\rT_{\tau_{i,s}^*\A}\cL_i
\end{equation}
for $i\in\{1,\dots,m\}$ and $s\in\cS_i$.
Then $\A$ determines  a differential operator
$$
\cD_\A:\Om^1_\A(X,\cg)\to \Om^{2,+}(X,\cg)\times \Om^0(X,\cg),
$$ 
\begin{equation}\label{eq:DAX}
\cD_\A\alpha := \left( (\rd_\A \alpha 
+ \rd X_f(\A)\alpha )^+ ,-\rd_\A^*\alpha \right),
\end{equation}
where 
$
\om^+:=\frac12(\om+*\om)
$
denotes the self-dual part of a $2$-form $\om\in\Om^2(X,\cg)$. 
This is a generalization of the linearized operator on 
$\R\times Y$ in (\ref{eq:DA}).  
The formal adjoint operator 
$$
\cD_\A^*:\Om^{2,+}_\A(X,\cg)\times \Om^0(X,\cg)\to \Om^1(X,\cg) 
$$
is given by 
$$
\cD_\A^*(\om,\phi) =  \rd_\A^*\om + \rd X_f(\A)^*\om -\rd_\A\phi.
$$
Here $\Om^{2,+}_\A(X,\cg)$ denotes the space of 
self-dual $2$-forms $\om$ on $X$ that satisfy the 
boundary condition
\begin{equation}\label{eq:bc-dual}
\tau_{i,s}^*\om=0, \qquad
\iota(\pd/\pd s)\tau_i^*\om|_{\{s\}\times\Sigma}\in\rT_{\tau_{i,s}^*\A}\cL_i
\end{equation}
for $i\in\{1,\dots,m\}$ and $s\in\cS_i$.

To obtain a Fredholm operator we must impose decay conditions 
on $\alpha$ at the tubular ends and extend the operator to 
suitable Sobolev completions. 
For any integer $k\geq 1$ and any $p>1$ denote 
by $W^{k,p}_\A(X,\rT^*X\otimes\cg)$ the space of $1$-forms 
on $X$ of class $W^{k,p}$ with values in  $\cg$ that satisfy the 
boundary conditions~(\ref{eq:alphabc})\footnote
{
Note that the subscript $\A$ in $W^{k,p}_\A$ indicates 
boundary conditions for the $1$-forms in this space. 
This is not to be confused with the norms $\|\cdot\|_{W^{k,p},\A}$
in Remark~\ref{rmk:norms}, where the subscript indicates that the covariant
derivatives are twisted by $\A$.
}
and by $W^{k,p}_\A(X,\Lambda^{2,+}\rT^*X\otimes\cg)$ the space 
of self-dual $2$-forms on $X$ of class $W^{k,p}$ with values in $\cg$ 
that satisfy the boundary conditions~(\ref{eq:bc-dual}).
The following theorem summarizes the Fredholm properties of $\cD_\A$
and $\cD_\A^*$. The regularity results (ii) and (iii) are steps towards 
the proof of (i).

\begin{thm}\label{thm:fredholm}
Suppose the limit connections $A_j$ are nondegenerate
and irreducible, i.e.~$H^0_{A_j}=0$ and $H^1_{A_j,f_j}=0$
for $j=1,\dots,n$. 
Then the following holds for every connection
$\A\in\cA(X,\cL)$ that satisfies~(\ref{eq:limX}).

\smallskip\noindent{\bf(i)}
The operators
$$
\cD_\A:W^{k,p}_\A(X,\rT^*X\otimes\cg)
\to W^{k-1,p}(X,\Lambda^{2,+}\rT^*X\otimes\cg)
\times W^{k-1,p}(X,\cg),
$$
$$
\cD_\A^*: W^{k,p}_\A(X,\Lambda^{2,+}\rT^*X\otimes\cg)
\times W^{k,p}(X,\cg)
\to W^{k-1,p}(X,\rT^*X\otimes\cg)
$$
are Fredholm for every integer $k\ge 1$ and every $p>1$.
Their Fredholm indices
$$
\delta_f(\A) := \mathrm{index}\,\cD_\A = - \mathrm{index}\,\cD_\A^*
$$
are independent of $k$ and $p$ and depend only 
on the homotopy class of $\A$ subject to~(\ref{eq:limX}). 

\smallskip\noindent{\bf(ii)}
If $\alpha\in L^p(X,\rT^*X\otimes\cg)$, 
$\om\in W^{k-1,p}(X,\Lambda^{2,+}\rT^*X\otimes\cg)$,
$\phi\in W^{k-1,p}(X,\cg)$ satisfy the equation
\begin{equation}\label{eq:asd-weak}
\int_X\inner{\cD_\A^*(\om',\phi')}{\alpha} 
= \int_X\Bigl(\inner{\om'}{\om}+\inner{\phi'}{\phi}\Bigr)
\end{equation}
for every compactly supported smooth
$(\om',\phi')\in\Om^{2,+}_\A(X,\cg)\times \Om^0(X,\cg)$,
then ${\alpha\in W^{k,p}_\A(X,\rT^*X\otimes\cg)}$ and
$\cD_\A\alpha=(\om,\phi)$.

\smallskip\noindent{\bf(iii)}
If $\om\in L^p(X,\Lambda^{2,+}\rT^*X\otimes\cg)$, $\phi\in L^p(X,\cg)$, 
$\alpha\in W^{k-1,p}(X,\rT^*X\otimes\cg)$ satisfy the equation
\begin{equation}\label{eq:asd-dual-weak}
\int_X\inner{(\om,\phi)}{\cD_\A\alpha'}
= \int_X\inner{\alpha}{\alpha'}
\end{equation}
for every compactly supported smooth $1$-form 
$\alpha'\in\Om^1_\A(X,\cg)$, then we have
${\om\in W^{k,p}_\A(X,\Lambda^{2,+}\rT^*X\otimes\cg)}$, 
$\phi\in W^{k,p}(X,\cg)$, and $\cD_\A^*(\om,\phi)=\alpha$.
\end{thm}

\begin{proof}
Assertions~(ii) and (iii) follow from Theorem~\ref{thm:4reg}
and Remark~\ref{rmk:time reversal}. 
(To obtain global $W^{k,p}$-regularity one sums up
estimates on compact domains 
-- with and without boundary -- exhausting $X$.)
To prove~(i)  we combine Theorems~\ref{thm:4reg} and~\ref{thm:iso}
with a cutoff function argument to obtain the estimate 
\begin{equation}\label{eq:prefred}
\Norm{\alpha}_{W^{k,p}(X)}
\le c\bigl(\Norm{\cD_\A\alpha}_{W^{k-1,p}(X)} 
+ \Norm{\alpha}_{W^{k-1,p}(K)}\bigr)
\end{equation}
for a sufficiently large compact subset $K\subset X$. 
(See~\cite[p.50]{Donaldson book}, or~\cite{RS} for the case 
$X=\R\times Y$, $k=0$, and $p=2$).  This estimate shows 
that $\cD_\A$ has a finite dimensional kernel and a closed image.
(See for example~\cite[Lemma~A.1.1]{MS}.)
By~(iii) the cokernel of $\cD_\A$ agrees with the kernel
of $\cD_\A^*$.  Since $\cD_\A^*$ satisfies a similar estimate
as~(\ref{eq:prefred}), it follows that the cokernel is finite dimensional
as well.   Hence $\cD_\A$ and $\cD_\A^*$ are Fredholm operators.
By~(ii) and~(iii), their Fredholm indices add up to zero and
are independent of $k$ and $p$.  That they depend only on the homotopy 
class of $\A$ follows from the stability properties of the Fredholm index.  
\end{proof}

In the case $\pd X=\emptyset$ the space of connections 
satisfying~(\ref{eq:limX}) is convex and so the index of 
$\cD_\A$ depends only on the limit connections $A_j$. 
The change of the index under gauge transformations 
on $Y_j$ depends on the degrees of the gauge transformations.   
By contrast, in the case $\p X\ne\emptyset$ and $\p Y_j\ne\emptyset$
the space of gauge transformations on $Y_j$ is connected, 
but the Lagrangian submanifolds $\cL_i$ have nontrivial 
fundamental groups.  So the index of $\cD_\A$ also depends 
on the homotopy classes of the paths in $\cL_i$ 
that are given by~$\A|_{\pd X}$.

\subsection*{Weighted theory}

In order to deal with reducible critical points we set up a refined
Fredholm theory on weighted Sobolev spaces.
Fix small 
nonzero real numbers $\delta_1,\dots,\delta_n$ 
and choose a smooth function $w:X\to(0,\infty)$ such that 
on all tubular ends
$$
w(\iota_j(s,y)) = e^{\delta_j s} \qquad\text{for}\;s\ge1 ,
$$ 
$w$ is independent of $y\in Y_j$ for $s\in[0,1]$, 
and $w\equiv 1$ on the complement. 
We introduce the weighted spaces
$$
W^{k,p}_{\A,\delta}(X,\rT^*X\otimes\cg)
:=\bigl\{ \alpha:X\to\rT^*X\otimes\cg 
\st w\alpha\in W^{k,p}_\A(X,\rT^*X\otimes\cg) \bigr\},
$$
and similarly for $W^{k,p}_\delta(X,\cg)$ and 
$W^{k,p}_\delta(X,\Lambda^{2,+}\rT^*X\otimes\cg)$.
The function $w$ does not appear in the notation because the 
spaces only depend on the choice of the~$\delta_j$.
The weighted inner product on $L^2_\delta(X,\rT^*X\otimes\cg)$ is 
$$
\la \alpha \,,\, \beta \ra_{L^2_\delta} 
:= \int_X w^2 \la \alpha\wedge *\beta \ra ,
$$
and similarly for $L^2_\delta(X,\cg)$.
The adjoint operator of $\rd_\A$ with respect to these 
two inner products 
is given by
$$
\rd_\A^{*,\delta}:=w^{-2}\rd_\A^*w^2: 
W^{k,p}_{\A,\delta}(X,\rT^*X\otimes\cg) 
\to W^{k-1,p}_\delta(X,\cg).
$$
It has the form ${(\alpha,\phi)\mapsto 
\rd_A^*\alpha-\nabla_s\phi - 2\delta_j\phi}$
on the tubular ends.
We will be using the following generalized Hodge decomposition.

\begin{lem}\label{lem:weighted Hodge}
Let $k$ be a positive integer and $p>1$ and 
suppose $\A\in\cA(X,\cL)$ satisfies (\ref{eq:limX}).
Then the operator 
$$
\rd_\A^{*,\delta}\rd_\A:
W^{k+1,p}_{\A,\delta}(X,\cg)\to W^{k-1,p}_\delta(X,\cg)
$$
with domain 
$
W^{k+1,p}_{\A,\delta}(X,\cg)
:=\bigl\{ \xi\in W^{k+1,p}_\delta(X,\cg) 
\st *\rd_\A\xi|_{\pd X}=0 \bigr\}
$
is bijective and there is a Hodge decomposition
$$
W^{k,p}_{\A,\delta} (X,\rT^*X\otimes\cg)
= \ker \rd_\A^{*,\delta} \oplus 
\rd_\A W^{k+1,p}_{\A,\delta}(X,\cg),
$$
\end{lem}

\begin{proof}
This Hodge decomposition is standard 
(see e.g.~\cite[Section 4.3]{Donaldson book})
except for the boundary conditions.
The two subspaces do not intersect since
$$
\la \rd_\A\xi , \alpha \ra_{L^2_\delta} 
- \la \xi , \rd_\A^{*,\delta}\alpha \ra_{L^2_\delta}
\;=\; \int_{\pd X} w^2 \inner{\xi}{* \alpha}
\;=\;0 
$$
for all $\alpha\in W^{k,p}_{\A,\delta}(X,\rT^*X\otimes\cg)$.
Assuming the operator $\rd_\A^{*,\delta}\rd_\A$ is bijective
we obtain the Hodge decomposition of 
$\beta\in W^{k,p}_{\A,\delta}(X,\rT^*X\otimes\cg)$
by solving the Neumann problem
$$
\rd_\A^{*,\delta}\rd_\A\xi = \rd_\A^{*,\delta}\beta,\qquad
*\rd_\A\xi|_{\pd X} = 0
$$
for $\xi\in W^{k+1,p}_\delta(X,\cg)$.
Since $\rd_\A\xi$ satisfies the Lagrangian boundary condition 
we have $\alpha:=\beta-\rd_\A\xi\in W^{k,p}_{\A,\delta}
(X,\rT^*X\otimes\cg)$.

To prove that the operator $\rd_\A^{*,\delta}\rd_\A$
is bijective we work with the weight function
$w=e^V:X\to(0,\infty)$ given by $V(s)=\delta_js$
on the tubular ends. Since $w$ has normal derivative
zero the function $\xi':=w\xi\in W^{k+1,p}(X,\cg)$
satisfies the boundary condition $*\rd_\A\xi'|_{\p X}=0$
whenever $\xi$ does. On the tubular ends we have 
$$
w\,\rd_\A^{*,\delta}\rd_\A w^{-1}
=\rd_{A_j}^*\rd_{A_j} - \Nabla{s}\Nabla{s} + \delta_j^2.
$$
This operator is bijective on $W^{k+1,p}_\A(\R\times Y_j,\cg)$ 
since it is Fredholm, symmetric, and positive definite.
So, as in the proof of Theorem~\ref{thm:fredholm},
one can use a cutoff function argument to show 
that $\rd_\A^{*,\delta}\rd_\A$ is a Fredholm operator.
Partial integration then shows that its kernel and cokernel are equal 
to the kernel of $\rd_\A$. To prove that the kernel is zero
let $\xi\in W^{k+1,p}_{\A,\delta}(X,\cg)$ with $\rd_\A\xi=0$
and assume w.l.o.g.~that $\A$ is in temporal gauge on the 
tubular ends. Then on each tubular end we have
$\p_s\xi_j\equiv0$, hence $\xi_j\equiv 0$ by the decay
condition, and hence $\xi\equiv 0$.
This proves the lemma. 
\end{proof}

Every connection $\A\in\cA(X,\cL)$ that satisfies~(\ref{eq:limX})
determines a differential operator
$$
\cD_{\A,\delta}:W^{k,p}_{\A,\delta}(X,\rT^*X\otimes\cg)
\to W^{k-1,p}_\delta(X,\Lambda^{2,+}\rT^*X\otimes\cg)
\times W^{k-1,p}_\delta(X,\cg)
$$
given by 
$$
\cD_{\A,\delta}\alpha 
:= \bigl( (\rd_\A\alpha + \rd X_f(\A)\alpha )^+,
-\rd_\A^{*,\delta}\alpha \bigr).
$$
Different choices of $w$ with the same $\delta_j$ give rise
to compact perturbations of~$\cD_{\A,\delta}$.

\begin{thm}\label{thm:weighted fredholm}
For $j=1,\dots,n$ let $A_j\in\cA(Y_j,\cL)$ 
and $\A\in\cA(X,\cL)$ such that ${F_{A_j}+X_{f_j}(A_j)=0}$ 
and $\A\in\cA(X,\cL)$ satisfies~(\ref{eq:limX}).
Then the following holds.

\smallskip\noindent{\bf (i)}
The operator $\cD_{\A,\delta}$ is Fredholm for every 
integer $k\ge 1$, every $p>1$, and every
$n$-tuple of sufficiently small nonzero
real numbers $\delta_1,\dots,\delta_n$.

\smallskip\noindent{\bf (ii)}
The Fredholm index of $\cD_{\A,\delta}$ 
is independent of $k$ and $p$; it depends only 
on the signs of the $\delta_j$ and on the homotopy class 
of $\A$ subject to~(\ref{eq:limX}).

\smallskip\noindent{\bf (iii)}
If the limit connections $A_j$ are all nondegenerate 
and irreducible, then ${\rm index}\,\cD_{\A,\delta}={\rm index}\,\cD_\A$.

\smallskip\noindent{\bf (iv)}
If the limit connections $A_j$ are all nondegenerate
and $\A$ satisfies~(\ref{eq:asdX}) then the cokernel 
of $\cD_{\A,\delta}$ is independent of the 
weight function (up to natural isomorphisms)
as long as the $|\delta_j|$ are sufficiently small. 
\end{thm}

\begin{proof}
The operator $w\cD_{\A,\delta}w^{-1}$ differs from 
$\cD_\A$ by a zeroth order perturbation which makes the 
operators on the tubular ends invertible.  
Hence assertions~(i-iii) follow by adapting the proof 
of Theorem~\ref{thm:fredholm} to the present case. 
To prove~(iv) we observe that
the restriction of the second component
$\rd_\A^{*,\delta}$ of $\cD_{\A,\delta}$
to the image of $\rd_\A$ is surjective
and, when $\A$ satisfies~(\ref{eq:asdX}), 
the image of $\rd_\A$ is contained in the kernel of
the first component ${(\rd_\A +\rd X_f(\A))^+}$ of 
$\cD_{\A,\delta}$.  Hence every element in the cokernel 
of $\cD_{\A,\delta}$ has the form $(\eta,0)$.
Moreover, $(\eta,0)$ belongs to the kernel 
of the adjoint operator $\cD^*_{\A,\delta}$
(with respect to the $L^2$-inner product determined
by $w$) if and only if $\eta=w^{-2}\zeta$, where
\begin{equation}\label{eq:D*eta}
w^{-1}\zeta\in W^{k-1,p}_\A(X,\Lambda^{2,+}\rT^*X\otimes\cg),\qquad
(\rd_\A+\rd X_f(\A))^*\zeta = 0.
\end{equation}
The subscript in $W^{k-1,p}_\A$ indicates the
dual boundary condition. It follows from linear 
exponential decay in Theorem~\ref{thm:lin-exp}
that every solution $\zeta$ of~(\ref{eq:D*eta}) decays 
exponentially. Hence the space of solutions of~(\ref{eq:D*eta})
is independent of the choice of the weight function $w$
as long as the $\Abs{\delta_j}$ are sufficiently small. 
This proves the theorem.
\end{proof}

\begin{rmk}\label{rmk:halfonto}\rm 
{\bf(i)}
The linearized operator is gauge equivariant in the sense that
$
\cD_{u^*\A,\delta} (u^{-1}\alpha u) 
= u^{-1} ( \cD_{\A,\delta} \alpha ) u
$
for all $\alpha\in W^{k,p}_{\A,\delta}(X,\rT^*X\otimes\cg)$
and all gauge transfomations $u\in\cG(X)$ that satisfy
$u\comp\iota_j\equiv u_j\in\cG(Y_j)$.

\smallskip\noindent{\bf(ii)}
In contrast to Theorem~\ref{thm:weighted fredholm}~(iv), 
the kernel of $\cD_{\A,\delta}$ is not
independent of the sign of the $\delta_j$
unless the $A_j$ are also irreducible. 

\smallskip\noindent{\bf(iii)}
On a tube $X=\R\times Y$ we will 
use weight functions of the form 
\begin{equation}\label{eq:wV}
w(s,y)=\exp(V(s))
\end{equation}
with $V\in\cC^\infty(\R)$ such that $V(s)=\pm \delta s$ for 
$\pm s\ge1$ (i.e.~$\delta_1=\delta_2=:\delta>0$).
Then $\cD_{\A,\delta}$ can -- as in Section~\ref{sec:HA} -- 
be identified with the operator
\begin{multline}\label{eq:DAd}
\cD_{\A,\delta}:W^{k,p}_{\A,\delta}(\R\times Y,\rT^*Y\otimes\cg)
\times W^{k,p}_\delta(\R\times Y,\cg) \\
\to W^{k-1,p}_\delta(\R\times Y,\rT^*Y\otimes\cg)
\times W^{k-1,p}_\delta(\R\times Y,\cg) 
\end{multline}
given by
$$
\cD_{\A,\delta} := \nabla_s + \cH_{A(s)}  
+ \left(
\begin{array}{cc}
0 & 0 \\ 0 & 2\lambda
\end{array}\right),\qquad 
\lambda:=\p_sV.
$$
The formal $L^2_\delta$-adjoint operator of $\cD_{\A,\delta}$ has the form
$$
\cD_{\A,\delta}^*(\alpha,\phi) 
:= - \nabla_s + \cH_{A(s)} 
- \left(
\begin{array}{cc}
2\lambda & 0 \\ 0 & 0
\end{array}\right).
$$

\smallskip\noindent{\bf(iv)}
The operator~(\ref{eq:DAd}) is conjugate
to the operator
\begin{equation}\label{eq:Dw}
w \cD_{\A,\delta} w^{-1} =
\Nabla{s} + \cH_{A(s)} - I_{\lambda(s)},\qquad 
I_\lambda := \left(
\begin{array}{cc}
\lambda & 0 \\ 0 & -\lambda
\end{array}\right),
\end{equation}
on the unweighted Sobolev spaces.
By Theorem~\ref{thm:weighted fredholm}~(iv) 
and its proof, this operator is surjective 
if and only if the operator
$
\Nabla{s} +\cH_{A(s)} - I_\delta
$
is surjective, provided $\delta\in\R\setminus\{0\}$ is 
sufficiently small and 
$\A\in\widetilde{\cM}(A^-,A^+;X_f)$ is a Floer 
connecting trajectory with nondegenerate ends.
\end{rmk}

\subsection*{The nonlinear setup}

In the remainder of this section we fix the constants 
$\delta_1=\cdots=\delta_n=\delta>0$.  Then the operators 
$\cD_{\A,\delta}$ have the following significance
for the study of the moduli space $\cM(A_1,\dots,A_n;X_f)$.
Let $\A\in\widetilde\cM(A_1,\dots,A_n;X_f)$ and suppose that
$\cD_{\A,\delta}$ is surjective.
If the $A_j$ are all nondegenerate and irreducible
and $\delta=0$, then $\cM(A_1,\dots,A_n;X_f)$ is a smooth manifold
near $[\A]$ whose tangent space is the kernel of $\cD_\A=\cD_{\A,\delta}$.
In general, the kernel of $\cD_{\A,\delta}$ is the tangent space of the
quotient
$$
\cM_0(A_1,\dots,A_n;X_f):=\widetilde\cM(A_1,\dots,A_n;X_f) / \cG_0(X),
$$
where $\cG_0(X)$ denotes the group of gauge transformations
$u\in\cG(X)$ that satisfy $u\comp\iota_j\equiv\one$ for every $j$.
Hence the dimension of $\cM(A_1,\dots,A_n;X_f)$ is equal to
\begin{equation}\label{eq:deltaf}
\delta_f(\A) := {\rm index}\,\cD_{\A,\delta} - \sum_{j=1}^n\dim H^0_{A_j}.
\end{equation}
(This agrees with the notation in Theorem~\ref{thm:fredholm}.)
To prove these assertions one can set up the nonlinear theory 
as follows. Fix an integer $k\ge 1$ and a real number $p>2$.
Associated to a tuple $A_j\in\cA(Y_j,\cL)$, $j=1,\dots,n$,
of critical points of the perturbed Chern--Simons functionals 
is a Banach manifold 
\begin{equation}\label{eq:Akpdelta}
\cA^{k,p}_\delta(X,\cL;A_1,\dots,A_n)
:= \left\{\A=\A_0+\alpha\,\Bigg|\,
\begin{array}{l}
\alpha\in W^{k,p}_\delta(X,\rT^*X\otimes\cg) \\
\tau_{i,s}^*\A\in\cL_i\;\forall i\,\forall s\in\cS_i
\end{array}\right\}
\end{equation}
where $\A_0\in\cA(X,\cL)$ is a reference connection 
satisfying $\iota_j^*\A_0\equiv A_j$ for all~$j$.
The tangent space of $\cA^{k,p}_\delta(X,\cL;A_1,\dots,A_n)$
is 
$$
\rT_\A\cA^{k,p}_\delta(X,\cL;A_1,\dots,A_n)
= \bigl\{\alpha\in W^{k,p}_\delta(X,\rT^*X\otimes\cg)\,\big|\,
\tau_{i,s}^*\alpha\in\rT_{\tau_{i,s}^*\A}\cL_i\bigr\}.
$$
Banach submanifold charts for
${\cA^{k,p}_\delta(X,\cL;A_1\dots A_n)
\subset \A_0 + W^{k,p}_\delta(X,\rT^*X\otimes\cg)}$
can be constructed with the help of the Banach submanifold
coordinates for $\cL_i\subset\cA^{0,p}(\Sigma_i)$ 
in~\cite[Lemma~4.3]{W Cauchy} (see Appendix~\ref{app:Lag}).
The gauge group 
\begin{equation}\label{eq:Gkpdelta}
\cG^{k+1,p}_\delta(X):=\Bigl\{
u:X\to\rG\,\Big|\,u^{-1}du\in W^{k,p}_\delta(X,\cg),\,
\lim_{s\to\infty}u\circ\iota_j = \one\Bigr\}
\end{equation}
acts freely on $\cA^{k,p}_\delta(X,\cL;A_1,\dots,A_n)$.
Its Lie algebra is the Banach space $W^{k+1,p}_\delta(X,\cg)$
and the quotient 
$\cA^{k,p}_\delta(X,\cL;A_1,\dots,A_n)/\cG^{k+1,p}_\delta(X)$
is a Banach manifold.  There is a gauge equivariant smooth map
$$
\cA^{k,p}_\delta(X,\cL;A_1,\dots,A_n)\to
W^{k,p}_\delta(X,\Lambda^{2,+}\rT^*X\otimes\cg):
\A\mapsto (F_\A+X_f(\A))^+
$$
and the moduli space $\cM_0(A_1,\dots,A_n;X_f)$
can be identified with the quotient of the zero set of this map
by the action of $\cG^{k+1,p}_\delta(X)$.  The operator 
$\cD_{\A,\delta}$ arises from linearizing this setup
in a local slice of the gauge group action
and hence, if this operator is surjective, it follows
from the implicit function theorem that $\cM_0(A_1,\dots,A_n;X_f)$
is a smooth manifold near $\A$, whose tangent space can
be identified with the kernel of $\cD_{\A,\delta}$.
The isotropy group $\cG_{A_1}\times\cdots\times\cG_{A_n}$
still acts on $\cM_0(A_1,\dots,A_n;X_f)$ and the quotient
by this action is the moduli space $\cM(A_1,\dots,A_n;X_f)$. 
If all limit connections $A_j$ are irreducible then the action
is free, so the moduli space is smooth.

\subsection*{The spectral flow}

We now specialize to the case $X:=\R\times Y$ and establish 
index identities for the linearized operator.
The main results are Theorem~\ref{thm:index} and Corollary~\ref{cor:index} below.
They will be proven by identifying the index with a spectral flow.

We fix a gauge invariant, monontone Lagrangian submanifold $\cL\subset\cA(\pd Y)$ 
satisfying (L1-2) on page~\pageref{p:L1} such that the zero connection
is contained in $\cL$ and is nondegenerate.
Choose a perturbation $h_f:\cA(Y)\to\R$ as in the introduction 
with a conjugation invariant function $f:\D\times\rG^N\to\R$.
Then the zero connection is a (nondegenerate) critical point 
of the perturbed Chern--Simons functional.  
For $A\in\mathrm{Crit}(\CS_\cL+h_f)$ and a path
$B:[0,1]\to\cL$ from $B(0)=A|_\Sigma$ to $B(1)=0$
we define an integer $\mu_f(A,B)$
as follows.  Choose a smooth path $A:[0,1]\to\cA(Y,\cL)$ 
such that $A(0) = A$, $A(1)=0$, and $A(s)|_\Sigma=B(s)$.
Define 
$$
\mu_f(A,B) := \mu_\mathrm{spec}
\left(\left\{\cH_{A(s)}+I_\eps\right\}_{s\in[0,1]}\right),\qquad
I_\eps := \left(\begin{array}{cc}\eps & 0 \\ 0 & -\eps\end{array}\right),
$$
where $\mu_\mathrm{spec}$ denotes the upward spectral flow
(see e.g.\ \cite{RS} and Appendix~\ref{app:spec}) 
and $\eps>0$ is sufficiently small. 
This integer is independent of the choice of the path $A$ 
and the constant $\eps$ used to define it. 
(The space of paths $A$ with fixed endpoints 
and boundary values is in fact convex.
Moreover, the kernel $\ker\,\cH_A=H^1_{A,f}\times H^0_A$ 
splits at the endpoints $A=A(0),A(1)$ by Proposition~\ref{prop:HA}.)

The significance of the following theorem is that the index resp.\ 
local dimension of the moduli space $\cM(A^-,A^+)$ is determined
modulo $8$ by the limit connections $A^-,A^+$.

\begin{thm} \label{thm:index}
{\bf (i)}
Let $A^\pm\in\cA(Y,\cL)$ be critical points of $\CS_\cL+h_f$ 
and ${\A\in\cA(\R\times Y)}$ be the connection 
associated to a smooth path $A:\R\to\cA(Y,\cL)$ with limits 
\begin{equation}\label{eq:Alimit}
\lim_{s\to\pm\infty}\left\|A-A^\pm\right\|_{\cC^1([s-1,s+1]\times Y)}=0 .
\end{equation}
Choose paths $B^\pm:[0,1]\to\cL$ 
from $B^\pm(0)=A^\pm|_\Sigma$ to $B^\pm(1)=0$
such that $B^-$ is homotopic to the catenation of 
the path $\R\to\cL:s\mapsto A(s)|_\Sigma$ with $B^+$. 
Then 
\begin{align*}
{\rm index}\,\cD_{\A,\delta}
&= \mu_\mathrm{spec}\bigl(\left\{\cH_{A(s)}
- I_{\lambda(s)}
\right\}_{s\in\R}\bigr)
\end{align*}
and
\begin{equation} \label{eq:deltamu}
\begin{split}
\delta_f(\A)
&:= {\rm index}\,\cD_{\A,\delta} 
- \dim H^0_{A^-} - \dim H^0_{A^+}  \nonumber \\
&= \mu_f(A^-,B^-)-\mu_f(A^+,B^+) 
- \dim\,H^0_{A^-} - \dim\,H^1_{A^+,f}.
\end{split}
\end{equation}

\smallskip\noindent{\bf (ii)}
If $A\in\cA(Y,\cL)$ is a critical point of $\CS_\cL+h_f$ and 
$B:[0,1]\to\cL$ is a path from $B(0)=A|_\Sigma$ to $B(1)=0$, then
for every loop $u:[0,1]\to\cG(\Sigma)$ with $u(0)=u(1)=\one$
$$
\mu_f(A,B) -\mu_f(A,u^*B) = 8\, {\rm deg}\, u .
$$
\end{thm}

\begin{proof}
Multiplication by $w$ defines an 
isomorphism $W^{k,p}_\delta\to W^{k,p}$, 
so $\cD_{\A,\delta}$ has the same index as 
the operator $w \cD_{\A,\delta} w^{-1}$ on 
$W^{k,p}_\A(\R\times Y,\rT^*(\R\times Y)\otimes \cg)$.
Hence, by~(\ref{eq:Dw}) and Theorem~\ref{thm:Fredholm}, 
the index of the operator $\cD_{\A,\delta}$ is given by 
\begin{align*}
\mathrm{index}(\cD_{\A,\delta}) 
&= \mu_\mathrm{spec}\bigl(
\left\{\cH_{A(s)} - I_{\lambda(s)}\right\}_{s\in\R}
\bigr) \\
&= \mu_\mathrm{spec}\bigl(
\left\{\cH_{A(s)}+I_\delta\right\}_{s\in\R}
\bigr) 
 - \dim\,H^1_{A^+,f}  + \dim\,H^0_{A^+}.
\end{align*}
Here $\lambda:=\pd_s V:\R\to\R$ satisfies 
$\lambda(s)=-\delta$ for $s\le -1$ 
and $\lambda(s)=\delta$ for $s\ge1$.
The second equation follows from a homotopy argument.
Namely, the path $\cH_{A(s)} - I_{\lambda(s)}$
is homotopic to the catenation of the path 
$\cH_{A(s)} + I_\delta$ with $\cH_{A^+} - I_{\lambda(s)}$.
Now the catenation of the path $\cH_{A(s)} + I_\eps$
with the path in the definition of $\mu_f(A^+,B^+)$ 
yields a path homotopic to the one in the definition of $\mu_f(A^-,B^-)$.
(By assumption the paths are homotopic over the boundary $\pd Y$, 
and this homotopy can be extended to the interior.)
Hence
\begin{align*}
\mu_f(A^-,B^-)
&= \mu_\mathrm{spec}
\bigl(\left\{\cH_{A(s)}+I_\eps\right\}_{s\in\R}\bigr) 
+ \mu_f(A^+,B^+).
\end{align*}
For $\delta>0$ sufficiently small we can choose $\eps=\delta$
and obtain
\begin{equation*}
\mu_f(A^-,B^-)- \mu_f(A^+,B^+)
= \mathrm{index}(\cD_{\A,\delta}) 
- \dim\,H^0_{A^+} + \dim\,H^1_{A^+,f}.
\end{equation*}
This proves~(i).

To prove~(ii) choose a path $A(s):[0,1]\to\cA(Y,\cL)$ with
$A(0)=A$, $A(1)=0$, and $B(s)=A(s)|_\Sigma$.
By homotopy invariance we may assume that 
$A(s)=0$ for $s\ge 1/2$.  Now let $u:[0,1]\to\cG(\Sigma)$ be a loop
with $u(0)=u(1)=\one$ and choose a path $A':[0,1]\to\cA(Y,\cL)$
such that $A'(0)=A$, $A'(1)=0$ and $A'(s)|_\Sigma=u(s)^*B(s)$.
Assume w.l.o.g.~that $u(s)=\one$ and $A'(s)=A(s)$ for $s\le 1/2$.
Then the spectral flow of the path $\cH_{A'(s)}+ I_\eps$
on the interval $0\le s\le 1/2$ is equal to $\mu_f(A,B)$. 
On the other hand, by Theorem~\ref{thm:indexS1} 
and a homotopy from $\cH_{A'}+I_\eps$ to $\cH_{A'}$,
the spectral flow on the interval $1/2\le s\le 1$
is equal to $\mathrm{index}(\cD_{\one,\A})$ for a connection
$\A=\tilde u^{-1}\rd\tilde u\in\cA(S^1\times Y,\cL)$ on the bundle $P_\one$
in the notation of Section~\ref{sec:S1Y}.
Here $\tilde u\in\cG(S^1\times Y)$ is homotopic to $u$ on $[1/2,1]\times Y$
and identically $\one$ on the complement. Hence 
\begin{equation*}
\begin{split}
\mu_f(A,B) - \mu_f(A,u^*B)
&= - \mu_\mathrm{spec}\bigl(
\left\{\cH_{A'(s)}+ I_\eps\right\}_{1/2\le s\le1}
\bigr) \\
&= - \mathrm{index}(\cD_{\one,\A}) 
\;=\;  8\deg(\one,\A) \;=\;  8\deg(u).
\end{split}
\end{equation*}
Here the third identity follows from Theorem~\ref{thm:S1Y}~(ii)
and the last from Remark~\ref{rmk:degree}~(iii). 
This proves the theorem.
\end{proof}

For every critical point $A\in\cA(Y,\cL)$ of the perturbed
Chern--Simons functional we define the real number $\eta_f(A)$
by 
$$
\eta_f(A) := \mu_f(A,B) - \frac{2}{\pi^2}
\Bigl(\CS(A,B) + h_f(A)\Bigr),
$$
where $B:[0,1]\to\cL$ is a path from $B(0)=A|_\Sigma$ to $B(1)=0$,
and $\CS(A,B)$ denotes the value of the Chern-Simons functional 
for the connection given by $A$ and $B$.  

\begin{cor}\label{cor:eta mu}
{\bf (i)}
The spectral flow $(A,B)\mapsto\mu_f(A,B)$ descends to a
circle valued function $\mu_f:\cR_f\to\Z/8\Z$. 

\smallskip\noindent{\bf (ii)}
The function $\eta_f:\Crit(\CS_\cL+h_f)\to\R$ is well defined and 
descends to a real valued function on $\cR_f$.
\end{cor}

\begin{proof}
Lemma~\ref{le:CS}~(iii), the homotopy invariance of the spectral flow,
and Theorem~\ref{thm:index}~(ii) imply that $\mu_f(A,B)\in\Z/8\Z$ 
is independent of the choice of $B$. 
Given a gauge transformation $u\in\cG(Y)$ we can connect it to 
the identity by a smooth path 
$\tu:[0,1]\to\cG(Y)$ from $\tu(0)=u$ to $\tu(1)=\one$.
Let $A:[0,1]\to\cA(Y,\cL)$ be the path in the definition of $\mu_f(A,B)$,
then $\mu_f(u^*A,(\tu|_{\pd Y})^*B)$ is defined as the spectral flow along 
the path $s\mapsto \tu(s)^*A(s)$ and hence, by the gauge equivariance of the 
Hessian,
$$
\mu_f(A,B) = \mu_f(u^*A,(\tu|_{\pd Y})^*B) .
$$
This proves (i).
That $\eta_f$ is well defined (i.e.\ independent of the choice of $B$)
follows from Lemma~\ref{le:CS} and Theorem~\ref{thm:index}~(ii).
To see that $\eta_f$ is gauge invariant it remains to check that
$$
\CS(A,B) = \CS(u^*A,(\tu|_{\pd Y})^*B) .
$$
This follows from the same argument as Lemma~\ref{le:CS}~(iv).
Namely, $\CS(A,B)$ is the Chern-Simons functional on 
$\tilde Y=Y\cup([0,1]\times\Sigma)$ of a connection $\tA$ given by $A$ and $B$.
The connection given by $u^*A$ and $(\tu|_{\pd Y})^*B$ is 
$\hat{u}^*\tA$, where the gauge transformation $\hat{u}\in\cG(\tilde Y)$
is given by $u$ and $\tu|_{\pd Y}$. It satisfies $\hat{u}|_{\pd Y}=\one$ 
and has degree zero since a homotopy to $\one$ is given by combining
$\tu(\sigma)$ on $Y$ with $s\mapsto\tu(s + (1-s)\sigma)|_{\pd Y}$
on $[0,1]\times\Sigma$.
Hence the equality of the Chern-Simons functionals follows from
the analogon of (\ref{eq:CS}) for manifolds with boundary and
gauge transformations that are trivial on the boundary.
\end{proof}

\begin{rmk}\label{rmk:etacont}\rm
The function $(f,A)\mapsto\eta_f(A)$ is continuous
on the space of nondegenerate pairs $(f,A)$. 
To see this note that the dimension of $H^0_A$
cannot jump, by Remark~\ref{rmk:jump}, and hence
one can locally work with the same constant
$\eps>0$ for the definition of $\mu_f$
in a neighbourhood of a pair $(f,A)$.
\end{rmk}

We can now state further index identities.
The monotonicity formula in (i) below -- a linear relationship
between index and energy -- will be central for excluding bubbling effects.

\begin{cor}\label{cor:index}
{\bf (i)}
Let $\A\in\cA(\R\times Y)$ be the connection associated to a smooth 
solution $A:\R\to\cA(Y,\cL)$ of (\ref{eq:floer}).
Suppose that it satisfies~(\ref{eq:Alimit}) with the critical points
$A^\pm\in\cA(Y,\cL)$ of $\CS_\cL+h_f$. Then 
$$
\delta_f(\A)
= \frac{2}{\pi^2}E_f(\A) + \eta_f(A^-) - \eta_f(A^+) 
- \dim\,H^0_{A^-} - \dim\,H^1_{A^+,f}.
$$
{\bf (ii)}
If $A,A':\R\to\cA(Y,\cL)$ are paths connecting $A^-$ to $B$,
respectively $B$ to $A^+$, then the index of their catenation
is given by
$$
\delta_f (\A\#\A') =
\delta_f(\A) + \delta_f(\A') 
+ \dim\,H^0_{B} + \dim\,H^1_{B,f}.
$$
{\bf (iii)}
If $A:\R\to\cA(Y,\cL)$ is a self--connecting path with limits 
${A^-=A^+=:A_0}$ and $s\mapsto A(s)|_\Sigma$  is homotopic
to $s\mapsto u(s)^*A_0|_\Sigma$ for $u:\R\to\cG(\Sigma)$
with ${u(\pm\infty)=\one}$, then
$$
\delta_f(\A)
= 8 \deg(u) -\dim\,H^0_{A_0} - \dim\,H^1_{A_0,f}.
$$
\end{cor}

\begin{proof}
Assertions (ii) and (iii) follow immediately from 
Theorem~\ref{thm:index}.  Assertion (i) follows from 
the definition of $\eta_f$, Theorem~\ref{thm:index},
and the following energy identity.  For a path 
$A:\R\to\cA(Y,\cL)$ satisfying 
$$
\pd_s A = - * ( F_A + X_f(A) )
$$
choose paths $B^\pm:[0,1]\to\cL$ from $B^\pm(0)=A^\pm|_\Sigma$ 
to $B^\pm(1)=0$ such that $B^-$ is homotopic to the catenation 
of $A(s)|_\Sigma$ with $B^+$.  Then
\begin{align*}
- E_f(\A)
&= \int_\R \int_Y \la \pd_s A \wedge \bigl( F_A + X_f(A) \bigr) \ra \,\ds \\
&= \int_\R  \biggl( 
 \frac 12 \frac\pd{\pd s}\int_Y \Bigl( \la A\wedge \rd A \ra 
 + \frac 13 \la A \wedge [A\wedge A] \ra \Bigr) \\
&\qquad\qquad\qquad\qquad
+ \frac 12 \int_\Sigma \la A \wedge \pd_s A \ra    
\;+\; \frac{\pd}{\pd s}h_f(A) 
 \biggl)\,\ds \\
&= \CS(A^+,B^+) +  h_f(A^+) - \CS(A^-,B^-) - h_f(A^-) .
\end{align*}
Here the second equation follows from~(\ref{Xf}) and the
fact that
\begin{align*}
2 \int_Y \la F_A \wedge \pd_s A \ra
&=  \int_Y \frac\pd{\pd s} 
      \Bigl( \la A\wedge \rd A \ra 
+ \frac 13 \la A \wedge [A\wedge A] \ra \Bigr)
  + \int_{\pd Y} \la A\wedge \pd_s A \ra .
\end{align*}
The last identity follows from the $\cC^1$-convergence 
of $A$ for $s\to\pm\infty$.
Since $B^-$ is homotopic (with fixed endpoints)
to the catenation of $A|_\Sigma$ with $B^+$, 
we have
$$
\int_\R\int_\Sigma \la A \wedge \pd_s A \ra \,\ds  
= \int_0^1\int_\Sigma \la B^- \wedge \pd_s B^- \ra \,\ds  
\;- \int_0^1\int_\Sigma \la B^+ \wedge \pd_s B^+ \ra \,\ds  .
$$
(See the proof of Lemma~\ref{le:CS} above for the invariance 
of this integral under homotopy.) This proves the corollary.
\end{proof}

\begin{remark}\label{rmk:index}\rm
Our notation for the indices is motivated by the following
finite dimensional model.  Let $M$ be a  
Riemannian $n$-manifold, $\rG$ be a compact 
Lie group that acts on $M$ by isometries, and 
$f:M\to\R$ be a $\rG$-invariant Morse--Bott function.
Associated to every critical point $x\in M$ is a chain
complex 
$$
0 \longrightarrow
\cg \stackrel{L_x}{\longrightarrow}
\rT_xM \stackrel{\nabla^2f(x)}{\longrightarrow}
\rT_xM \stackrel{L_x^*}{\longrightarrow}
\cg \longrightarrow 0,
$$
where $L_x$ is the infinitesimal action of $\cg$
and $\nabla^2f(x)$ is the Hessian of $f$ (see~(\ref{eq:Acomplex})).
We denote
$$
\nu_0(x) := \dim\ker L_x,\qquad
\nu_1(x) := \dim\frac{\ker\nabla^2f(x)}{\im L_x},\qquad
\mu(x) := \mathrm{ind}_f(x),
$$
that is $\mu(x)$ is the number of negative eigenvalues of 
the Hessian and $\nu_0(x)$ is the dimension of the isotropy subgroup. 
Now the kernel of the Hessian has dimension ${\nu_1(x)+\dim\rG-\nu_0(x)}$, the unstable manifold
$W^u(x)$ of the orbit $\rG x$ has dimension ${\mu(x)+\dim\rG-\nu_0(x)}$,
the stable manifold $W^s(x)$ of $\rG x$ has dimension ${n-\mu(x)-\nu_1(x)}$,
and, in the tranverse case, the moduli space 
$$
\cM(x^-,x^+):=W^u(x^-)\cap W^s(x^+)/\rG
$$
of connecting trajectories has dimension (compare with~(\ref{eq:deltamu}))
$$
\delta(x^-,x^+):=\dim\cM(x^-,x^+)
= \mu(x^-)-\mu(x^+)-\nu_0(x^-)-\nu_1(x^+) .
$$
\end{remark}


\section{Compactness}\label{sec:compact}

Let $Y$ be a compact oriented Riemannian $3$-manifold 
with boundary $\pd Y=\Sigma$ and $\cL\subset\cA(\Sigma)$ 
be a gauge invariant, monotone, irreducible Lagrangian submanifold 
satisfying (L1-3) on page~\pageref{p:L1}. 
Fix a collection of embeddings $\gamma_i:S^1\times\D\to{\rm int}(Y)$,
$i=1,\dots,m$, as in Section~\ref{sec:CS}. We use the notation
$$
\widetilde\cM(A^-,A^+;X_f) 
:= \left\{ \A\in\cA^\tau(\R\times Y) \left|
\begin{array}{l}
\pd_s A - \rd_A\Phi + * (F_A + X_f(A)) = 0, \\
A(s)|_\Sigma\in\cL \quad\forall s\in\R, \\
E_f(\A)<\infty, \; \lim_{s\to\pm\infty}A(s)= A^\pm
\end{array}\right.\right\}
$$
for the space of Floer connecting trajectories 
associated to a perturbation 
$
f\in\cC^\infty(\D\times\rG^m)^\rG
$
and two critical points $A^\pm\in\cA(Y,\cL)$ of $\CS_\cL+h_f$.
Here $\cA^\tau(\R\times Y)$ denotes the space of connections 
$\Xi=\Phi\ds+A$ on $\R\times Y$ that are in temporal gauge
outside of $[-1,1]\times Y$, i.e.~$\Phi(s)=0$ for $|s|\ge 1$.
The corresponding gauge group $\cG(A^-,A^+)$ 
consists of all gauge transformations
$u:\R\to\cG(Y)$ that satisfy~$u(s)=u^\pm\in\cG_{A^\pm}$ 
for $\pm s\geq 1$ and the quotient 
space will be denoted by 
$$
\cM(A^-,A^+;X_f):=\widetilde\cM(A^-,A^+;X_f)/\cG(A^-,A^+)
$$
The goal of this section is to establish compactness theorems
for these moduli spaces. 
The proofs will be heavily based on the basic compactness
results in \cite{W elliptic,W bubb}.
We start with a summary of the compactness 
for uniformly bounded curvature.

\begin{prp}\label{prop:compact}
Let $f^\nu\in\cC^\infty(\D\times\rG^m)^\rG$ 
be a sequence that converges to 
$f^\infty\in\cC^\infty(\D\times\rG^m)^\rG$ 
in the $C^{k+1}$-topology for some $k\geq 1$.
Let $I^\nu\subset\R$ be a sequence of open intervals 
such that $I^\nu\subset I^{\nu+1}$ for all $\nu$
and denote $I:=\bigcup_\nu I^\nu$. Let 
$
\Xi^\nu=\Phi^\nu\rd s+A^\nu\in\cA(I^\nu\times Y)
$ 
be a sequence of solutions of the Floer equation 
\begin{equation}\label{eq:bvp nu}
\pd_s A^\nu - \rd_{A^\nu}\Phi^\nu 
+*\bigl( F_{A^\nu} + X_{f^\nu}(A^\nu) \bigr) =0, 
\qquad
A^\nu(s)|_\Sigma\in\cL,
\end{equation}
such that the curvature $|F_{\Xi^\nu}|$
is locally uniformly bounded. Then the following holds.

\smallskip\noindent{\bf (i)}
There exists a subsequence, still denoted by $\Xi^\nu$,
and a sequence of gauge transformations
$u^\nu\in\cG(I^\nu\times Y)$ such that $(u^\nu)^*\Xi^\nu$ 
converges in the $\cC^k$ topology 
on every compact subset of $I\times Y$.

\smallskip\noindent{\bf (ii)}
There exists a subsequence, still denoted by $\Xi^\nu$,
and a sequence of gauge transformations
$u^\nu\in\cG(I^\nu\times Y)$ such that $(u^\nu)^*\Xi^\nu$ 
is in temporal gauge and converges in the $\cC^{k-1}$ topology 
on every compact subset of $I\times Y$.

\smallskip\noindent{\bf (iii)}
In both cases, the limit $\Xi^\infty\in\cA(I\times Y)$ 
of the subsequence can be chosen smooth and it
satisfies~(\ref{eq:bvp nu}) with $f^\nu$ replaced by $f^\infty$.
\end{prp}

\begin{proof}
In a neighbourhood of the boundary $I\times\pd Y$,
where the perturbations vanish, compactness for 
anti-self-dual connections with Lagrangian boundary 
conditions was established in~\cite[Theorem~B]{W elliptic}.
The interior compactness follows from 
standard techniques (e.g.~\cite{DK}, \cite{W})
and Remark~\ref{rmk:Xf}. The crucial point 
in the bootstrapping argument is that 
a $W^{k,p}$-bound on $(u^\nu)^*\Xi^\nu$ implies a 
$W^{k,p}$-bound on $X_{f^\nu}((u^\nu)^*\Xi^\nu)$ 
and hence on $F_{(u^\nu)^*\Xi^\nu}^+$.
(The constant in the $W^{k,p}$-estimate 
of Proposition~\ref{prop:Xf}~(iii) depends 
continuously on $f\in\cC^{k+1}$.)
Combining these two compactness results via a general 
patching procedure as in~\cite[Lemma~4.4.5]{DK} 
or~\cite[Proposition 7.6]{W} we deduce that, 
for a suitable subsequence and choice of $u^\nu$,
the sequence $(u^\nu)^*\Xi^\nu$ is bounded in 
$W^{k+1,p}(K)$ for every compact subset $K\subset I\times Y$
and a fixed $p>4$, and hence has a $\cC^k$ convergent subsequence. 
A diagonal argument then proves~(i).  

To prove~(ii) we write 
$
\tilde\Xi^\nu:=(u^\nu)^*\Xi^\nu=:\tilde\Phi^\nu\ds+\tA^\nu
$
where $u^\nu$ is as in~(i).  Then $\tilde\Xi^\nu$
is bounded in $W^{k+1,p}$ on every compact subset of $I\times Y$. 
Define $v^\nu:I^\nu\times Y\to\rG$ as the unique 
solution of the differential equation
$$
\p_sv^\nu+\tilde\Phi^\nu v^\nu=0,\qquad v^\nu(0)=\one.
$$
Then $v^\nu$ is bounded in $W^{k+1,p}$ on every compact subset of $I\times Y$.  
(To check this use the identity $\pd_s(v^{-1}\rd v)=-v^{-1}\Phi v$.)
Hence $(v^\nu)^*\tilde\Xi^\nu=(u^\nu v^\nu)^*\Xi^\nu$ 
is in temporal gauge and is bounded in $W^{k,p}$ 
on every compact subset of $I\times Y$. The compact embeddings 
$W^{k,p}(K)\hookrightarrow\cC^{k-1}(K)$ together with a diagonal
argument then prove~(ii).

The regularity of the limit $\Xi^\infty$ can be achieved 
by a further gauge transformation.  That $\Xi^\infty$ 
solves~(\ref{eq:bvp nu}) follows from the fact that these equations 
are gauge invariant and preserved under weak $W^{k,p}$ convergence.
\end{proof}

The following is the most general compactness result for 
bounded energy.

\begin{thm}\label{thm:compact1}
Let ${f\in\cC^\infty(\D\times\rG^m)^\rG}$ be a perturbation
such that every critical point of ${\CS_\cL+h_f}$ is nondegenerate.
Let $f^\nu\in\cC^\infty(\D\times\rG^m)^\rG$ be a sequence
that converges to $f$ in the $C^{k+1}$-topology
and let $\Xi^\nu=\Phi^\nu\ds+A^\nu
\in\widetilde{\cM}(A^\nu_-,A^\nu_+;X_{f^\nu})$
be a sequence of Floer connecting trajectories 
with bounded energy
$$
\sup_\nu E_{f^\nu}(\Xi^\nu) 
= \sup_\nu \int_{\R\times Y}
\Abs{\pd_s A^\nu - \rd_{A^\nu}\Phi^\nu }^2 < \infty.
$$
Fix $p>1$ and suppose that $A^\nu_\pm$ converges to 
$A^\pm\in\Crit(\CS_\cL+h_f)$ in the $\cC^k$ topology. 
Then there is a subsequence, still denoted by $\Xi^\nu$, 
critical points $B_0,\dots,B_\ell\in\Crit(\CS_\cL+h_f)$
with $B_0=A^-$, $B_\ell=A^+$,
and Floer connecting trajectories
$
\Xi_i\in\widetilde{\cM}(B_{i-1},B_i;X_f)
$
for $i=1,\dots,\ell$,
such that $\Xi^\nu$ converges to the broken trajectory 
$(\Xi_1,\dots,\Xi_\ell)$ in the following sense.

For every $i\in\{1,\dots,\ell\}$ there is a sequence $s^\nu_i\in\R$
and a sequence of gauge transformations $u^\nu_i\in\cG(\R\times Y)$
such that the sequence $s\mapsto((u^\nu_i)^*\Xi^\nu)(s+s^\nu_i)$ 
converges to $\Xi_i$ in the $W^{1,p}$-norm on every compact 
subset of $\R\times Y\setminus Z_i$.
Here $Z_i\subset\R\times Y$ is the bubbling locus 
consisting of finitely many interior points
and finitely many boundary slices; it is nonempty 
whenever $\Xi_i$ has zero energy. 

The broken 
trajectory $(\Xi_1,\dots,\Xi_\ell)$ has energy and index
\begin{equation}\label{eq:Edelta}
\begin{split}
&\sum_{i=1}^\ell E_f(\Xi_i)  
\le \lim_{\nu\to\infty}E_{f^\nu}(\Xi^\nu),\\
&\sum_{i=1}^\ell \delta_f(\Xi_i)  + \sum_{i=1}^{\ell-1} \dim H^0_{B_i}
\le \lim_{\nu\to\infty}\delta_{f^\nu}(\Xi^\nu).
\end{split}
\end{equation}
If $\sup_\nu \|F_{\Xi^\nu}\|_{L^\infty}<\infty$ then 
there is no bubbling (i.e.~$Z_i=\emptyset$ 
for all $i$), equality holds in~(\ref{eq:Edelta}), 
and $(u^\nu)^*\Xi^\nu$ converges in the $\cC^k$ 
topology on every compact set. 
If $\sup_\nu \|F_{\Xi^\nu}\|_{L^\infty}=\infty$ then
there is bubbling (i.e.~$Z_i\ne\emptyset$ for some~$i$)
and
\begin{equation} \label{eq:bubb}
\begin{split}
&\sum_{i=1}^\ell E_f(\Xi_i)  
\le \lim_{\nu\to\infty}E_{f^\nu}(\Xi^\nu) - 4\pi^2,\\
&\sum_{i=1}^\ell \delta_f(\Xi_i)  
+ \sum_{i=1}^{\ell-1} \dim H^0_{B_i}
\le \lim_{\nu\to\infty}\delta_{f^\nu}(\Xi^\nu) - 8. 
\end{split}
\end{equation}
\end{thm}

\begin{rmk}\rm
The assumption that $A^\nu_\pm$ converges in the $\cC^k$ 
topology always holds for a subsequence in a suitable gauge,
by Proposition~\ref{prop:finite}.  
\end{rmk}

\begin{proof}[Proof of Theorem~\ref{thm:compact1}.]
Replacing the uniform bound on the curvature 
in Proposition~\ref{prop:compact} by an energy bound 
on $\Xi^\nu$ allows for bubbling. For the (unperturbed) 
anti-self-duality equation with Lagrangian boundary
conditions this was dealt with in \cite[Theorems~1.2,1.5]{W bubb},
\cite{W lag}, and \cite[Section~3]{W survey};
for the perturbed equation in the interior the (well known) 
result is Theorem~\ref{thm:compactness}.
Combining these one essentially obtains the same 
basic compactness theorem as for anti-self-dual 
connections (see~\cite[Proposition~2.1]{Donaldson book}).
A minor difference is that -- due to the holonomy 
perturbations -- we obtain convergence in the $W^{1,p}$-norm 
for any $p>1$ rather than in the $\cC^\infty$-topology;
so~\cite[Proposition~2.1~(1)]{Donaldson book} is 
replaced by $W^{1,p}$-convergence.  The crucial difference 
is in the knowledge about the bubbling phenomenon.
First, the finite set $\{x_1,\dots,x_\ell\}\subset\R\times Y$ 
of bubbling points is replaced by a more general bubbling locus 
$Z\subset\R\times Y$ consisting of finitely many interior points 
and finitely many boundary slices $\{s\}\times\pd Y$.
On the complement of $Z$, one has local $L^p$-bounds 
on the curvature.  Second, we do not have a geometric 
description of the bubbles (after rescaling) or the 
precise quantum $4\pi^2$ for the energy concentration.
There is however a universal constant $\hbar>0$ that is 
a lower bound for the energy 
concentration at each component of the bubbling locus $Z$;
so~\cite[Proposition~2.1~(2)]{Donaldson book} 
is replaced by $\int_U |F_A+X_f(A)|^2  \leq 
\limsup_{\alpha'\to\infty}\int_U 
|F_{A_{\alpha'}}+X_{f_{\alpha'}}(A_{\alpha'})|^2 - \ell\hbar $,
where $\ell$ is the number of points and boundary slices in $Z$.

The second source of noncompactness, the splitting of trajectories,
is the same as for the usual Floer theories.
With the exponential decay results of Section~\ref{sec:exp} 
and the modified basic compactness above, 
one can adapt the discussion in~\cite[Chapter~5.1]{Donaldson book} 
to prove the convergence to a broken trajectory. 
In particular, exponential decay holds for sufficiently 
$\cC^2$-close perturbations with uniform constants 
(see Theorem~\ref{thm:fcrit} for the nondegeneracy
and Proposition~\ref{prop:Xf}~(v) for the constants).
More precisely we argue as follows.

Throughout we denote the perturbed Yang-Mills energy 
of $\Xi^\nu$ on $I\times Y$ by 
$$
E_{f^\nu}(\Xi^\nu;I) 
:= \int_I\int_Y\Abs{\p_sA^\nu-\rd_{A^\nu}\Phi^\nu}^2.
$$
Passing to a subsequence we may assume that bubbling 
occurs only for finitely many sequences $t^\nu_j$,
$j=1,\dots,m$, with
$$
\hbar_j := \lim_{\delta\to 0}
\lim_{\nu\to\infty}
E_{f^\nu}(\Xi^\nu;[t^\nu_j-\delta,t^\nu_j+\delta])
\ge\hbar.
$$
In particular, the limits exist.
The sequences are chosen such that 
$t^\nu_{j+1}-t^\nu_j>0$ and that these
differences converge either to a positive number 
or to infinity. We may also assume that the curvature
of $\Xi^\nu$ is uniformly bounded on the complement 
of the sets $[t^\nu_j-\delta,t^\nu_j+\delta]\times Y$ 
for every $\delta>0$ and that the following limits exist:
\begin{equation*}
\begin{split}
\eps_0 &:= \lim_{\delta\to 0}\lim_{\nu\to\infty}
E_{f^\nu}(\Xi^\nu;(-\infty,t^\nu_1-\delta]),\\
\eps_j &:= \lim_{\delta\to 0}\lim_{\nu\to\infty}
E_{f^\nu}(\Xi^\nu;[t^\nu_j+\delta,t^\nu_{j+1}-\delta]),\qquad
j=1,\dots,m-1, \\
\eps_m &:= \lim_{\delta\to 0}\lim_{\nu\to\infty}
E_{f^\nu}(\Xi^\nu;[t^\nu_m+\delta,\infty)).
\end{split}
\end{equation*}
Then
$$
\lim_{\nu\to\infty} E_{f^\nu}(\Xi^\nu)
= \eps_0+\hbar_1+\eps_1+\cdots+\hbar_m+\eps_m.
$$
Next we choose a constant $\eps>0$
smaller than the constant in Theorem~\ref{thm:long}
and smaller than $\hbar$. Following~\cite[5.1]{Donaldson book}
we choose the $s^\nu_i\in\R$ inductively such that 
$$
E_{f^\nu}(\Xi^\nu;(-\infty,s^\nu_1]) = \frac{\eps}{2},\qquad
E_{f^\nu}(\Xi^\nu;[s^\nu_i,s^\nu_{i+1}]) 
= E_f(\Xi_i)+\sum_{j\in J_i}\hbar_j,
$$
where $\Xi_i$ is the limit of the sequence 
$\Xi^\nu(s^\nu_i+\cdot)$ modulo gauge and bubbling
and $J_i\subset\{1,\dots,m\}$ denotes the set of all $j$
such that the sequence $t^\nu_j-s^\nu_i$ is bounded. 
This choice guarantees that $s^\nu_{i+1}-s^\nu_i\to\infty$
for all $i$, that $\{1,\dots,m\}$ is the disjoint 
union of the $J_i$, and that $J_i\ne\emptyset$ 
whenever $\Xi_i$ has zero energy.  By Theorem~\ref{thm:long} 
(applied to a temporal gauge of the $\Xi^\nu$
on intervals $[s^\nu_i+T,s^\nu_{i+1}]$ with energy less than $\eps$)
the positive end of $\Xi_i$ is gauge equivalent 
(and hence w.l.o.g.~equal to) the negative end of $\Xi_{i+1}$,
the negative end of $\Xi_1$ is $A^-$, and 
the positive end of $\Xi_\ell$ is $A^+$. 
The total energy of the broken trajectory is  
\begin{equation}\label{eq:enelim}
\sum_{i=1}^\ell E_f(\Xi_i)
= \sum_{j=0}^m\eps_j 
= \lim_{\nu\to\infty}E_{f^\nu}(\Xi^\nu)-\sum_{j=1}^m\hbar_j.
\end{equation}
If the curvature is bounded then $m=0$
and all bubbling loci $Z_i$ are empty.
In this case the energy identity is~(\ref{eq:enelim})
and the index identity follows from the 
monotonicity formula in Corollary~\ref{cor:index}~(i).
If the curvature blows up then $m\ge 1$,
hence $Z_i\ne\emptyset$ for some $i$, and
we obtain the strict inequality 
\begin{equation*}
\begin{split}
\sum_{i=1}^\ell \delta_f(\Xi_i) 
+ \sum_{i=1}^\ell \dim H^0_{B_{i-1}}
&= 
\sum_{i=1}^\ell
\biggl(
\frac{2}{\pi^2}E_f(\Xi_i) 
+ \eta_f(B_{i-1})-\eta_f(B_i)
\biggr) \\
& < \lim_{\nu\to\infty} 
\biggl(\frac{2}{\pi^2}E_{f^\nu}(\Xi^\nu) 
+ \eta_{f^\nu}(A^\nu_-) - \eta_{f^\nu}(A^\nu_+)\biggr) \\
&=
\lim_{\nu\to\infty}\delta_{f^\nu}(\Xi^\nu) 
+ \dim H^0_{A^-}.
\end{split}
\end{equation*}
Here the first step follows from 
Corollary~\ref{cor:index}~(i), 
the second step uses~(\ref{eq:enelim}) and the continuity of the 
function $(f,A)\mapsto\eta_f(A)$ (see Remark~\ref{rmk:etacont}), 
and the last step uses Corollary~\ref{cor:index}~(ii) 
and $\dim H^0_{A^\nu_-}=\dim H^0_{A^-}$ 
for $\nu$ sufficiently large (see Remark~\ref{rmk:jump}).
Each side of our inequality has the form 
$\delta_f(\Xi)+\dim H^0_{A^-}$ for a suitable 
path $\Xi$ running from $A^-$ to $A^+$.  
For the left hand side, by Corollary~\ref{cor:index}~(ii),
$\Xi$ can be chosen as the catenation of the $\Xi_i$ and 
for the right hand side as a small deformation of $\Xi^\nu$ 
for $\nu$ sufficiently large. Since the inequality is strict 
it follows from Theorem~\ref{thm:index}~(i)
and Corollary~\ref{cor:eta mu} that the defect is at least $8$. 
Using monotonicity again we obtain an energy gap of
at least $4\pi^2$.  This proves the theorem. 
\end{proof}

A first consequence of the compactness and index identities is that
we can exclude bubbling in certain moduli spaces by transversality.

\begin{cor}\label{cor:compact} 
Suppose that the sequence of solutions in Theorem~\ref{thm:compact1}
has index 
$$
\delta_{f^\nu}(\Xi^\nu)\leq 7.
$$
Suppose that either bubbling occurs or one of the limit trajectories $\Xi_i$
is a self-connecting trajectory of $[B_{i-1}]=[B_i]=[0]$.
Then one of the limit trajectories $\Xi_j$ must have 
negative index $\delta_{f}(\Xi_j)< 0$ and
at least one of its endpoints $B_{j-1}$ or $B_j$ is not
gauge equivalent to the trivial connection.
\end{cor}

\begin{proof}
Every nontrivial self--connecting trajectory $\Xi_i$ of $[0]$ has
index $\delta_f(\Xi_i)\geq 5$ by Corollary~\ref{cor:index}  
with $E_f(\Xi_i)=4\pi^2\deg(u)>0$.
It also adds $\dim H^0_{[0]}=3$ to the sum of indices.
So to achieve a sum $\leq 7$, one of the other indices must
be negative. 
A trivial self--connecting trajectory of $[0]$ has index $-3$ but 
also adds another $\dim H^0_{[0]}=3$ to the sum of indices.
Hence there must be a trajectory with negative index and at least one
nontrivial end. The same holds in the bubbling case by~(\ref{eq:bubb}).
\end{proof}

We will refine the compactness theorem in two special cases.
First we consider the case of no breaking and no bubbling in which
we obtain actual compactness of moduli spaces.

\begin{thm}\label{thm:compact3}
Fix a constant $p>1$. Let $f,f^\nu$ be as in 
Theorem~\ref{thm:compact1} and $A^\pm\in\cA(Y,\cL)$ 
such that
$
{F_{A^\pm}+X_{f^\nu}(A^\pm)=0}
$
for all $\nu$.  Then there is a $\delta>0$ such that the following holds.
If $\Xi\in\widetilde\cM(A^-,A^+;X_f)$ and, for each $\nu$,
$\Xi^\nu$ is a solution of~(\ref{eq:bvp nu}) that is
gauge equivalent to an element of $\widetilde\cM(A^-,A^+;X_{f^\nu})$
such that $\Xi^\nu$ converges to $\Xi$
in the $\cC^k$ topology on compact sets and
$$
E_f(\Xi) = \lim_{\nu\to\infty}E_{f^\nu}(\Xi^\nu),
$$
then there exists a sequence of gauge transformations 
$u^\nu\in\cG(\R\times Y)$ such that $(u^\nu)^*\Xi^\nu$ 
converges to $\Xi$ in $W^{k,p}_\delta(\R\times Y)$.
\end{thm}

\begin{proof}
Note that, by contradiction, it suffices to prove the convergence statement for a subsequence.
For that purpose we choose $v^\nu\in\cG(\R\times Y)$
such that 
$$
\tilde\Xi^\nu:=(v^\nu)^*\Xi^\nu\in\widetilde\cM(A^-,A^+;X_{f^\nu}).
$$
In particular, $\tilde\Xi^\nu=:\tilde\Phi^\nu\ds + \tilde A^\nu$ 
is in temporal gauge outside of $[-1,1]\times Y$.
Fix a constant $\eps>0$ smaller than the constant in Corollary~\ref{cor:decay}
and note that the exponential $\cC^k$ estimate in Corollary~\ref{cor:decay}
holds with uniform constants $\delta_0:=\delta>0$ and $C_0:=C_k$ 
in a sufficiently small $\cC^{k+1}$ neighborhood of $f$.  
We write $\Xi=\Phi\ds + A$ and choose $T_0>0$ such that 
$$
\int_{-T_0}^{T_0}\int_Y \Abs{\pd_s A-\rd_A\Phi}^2 > E_f(\Xi)-\eps.
$$
Since $\Xi^\nu = \Phi^\nu\ds + A^\nu$ converges in the $\cC^k$ norm 
on compact sets we have 
$$
\int_{-T_0}^{T_0}\int_Y 
\bigl|\pd_s\tilde A^\nu-\rd_{\tilde A^\nu}\tilde\Phi^\nu\bigr|^2
= \int_{-T_0}^{T_0}\int_Y \Abs{\pd_s A^\nu-\rd_{A^\nu}\Phi^\nu}^2 
> E_{f^\nu}(\Xi^\nu)-\eps
$$
and thus $E(\Xi^\nu;(-\infty,T_0])+E(\Xi^\nu;[T_0,\infty))<\eps$
for sufficiently large $\nu\ge\nu_0$. 
Hence it follows from Corollary~\ref{cor:decay} that 
\begin{equation*}
\begin{split}
\bigl\|\tilde A^\nu-A^+\bigr\|_{\cC^k([s,\infty)\times Y)}
&\le C_0e^{-\delta_0(s-T_0)}E(\Xi^\nu;[T_0,\infty)), \\
\bigl\|\tilde A^\nu-A^-\bigr\|_{\cC^k((-\infty,-s]\times Y)}
&\le C_0e^{-\delta_0(s-T_0)}E(\Xi^\nu;(-\infty,-T_0])
\end{split}
\end{equation*}
for $s\ge T_0+1$ and $\nu\ge\nu_0$.  
The same estimate holds with $\tilde A^\nu$ replaced by~$A$.  
Now fix a constant $0<\delta<\delta_0$.  Then there exists a 
constant $C$ (depending on $C_0$, $\delta$, $\delta_0$, $k$, and $p$) 
such that 
$$
\bigl\|\tilde\Xi^\nu-\Xi\bigr\|_{W^{k,p}_\delta((\R\setminus[-T,T])\times Y)}
\le Ce^{-(\delta_0-\delta)(T-T_0)}
$$
for $T\ge T_0+1$ and $\nu\ge\nu_0$.  

Next, fix a sequence $\rho_n \to 0$ and choose 
$T_n\to\infty$ so that $T_n\ge T_0+1$ and
$$
Ce^{-(\delta_0-\delta)(T_n-T_0)} < \frac{\rho_n}{2}.
$$
For fixed $n\in\N$ note that both $\Xi^\nu$ and 
$\tilde\Xi^\nu=(v^\nu)^*\Xi^\nu$ converge to $\Xi$ 
in the $\cC^k$ norm on $[T_n,T_n+1]\times Y$ and 
on $[-T_n-1,-T_n]\times Y$. Using the identity
\begin{equation}\label{vnu}
(v^\nu)^{-1}\rd v^\nu = \tilde\Xi^\nu - (v^\nu)^{-1}\Xi^\nu v^\nu 
\end{equation}
we thus inductively obtain bounds on $v^\nu$ 
in $\cC^{k+1}((\pm[T_n,T_n+1])\times Y)$.
Then, by a compact Sobolev embedding, we find a subsequence 
$\lim_{\ell\to\infty}\nu_n(\ell)=\infty$ such that
$v^{\nu_n(\ell)}|_{(\pm[T_n,T_n+1])\times Y} \to v_n^\pm$
converges in the $\cC^k$ norm. Again using (\ref{vnu}) we see 
that this convergence is in fact in the $\cC^{k+1}$ norm.
On these domains we moreover have
\begin{align*}
\bigl\| (v_n^\pm)^*\Xi - \Xi \bigr\|_{\cC^k}
& =\lim_{\nu=\nu_n(\ell)\to\infty}
\bigl\| (v^\nu)^*\Xi - \Xi \bigr\|_{\cC^k} \\
&\leq \lim_{\nu\to\infty}
\bigl( \bigl\| (v^\nu)^*\Xi^\nu - \Xi \bigr\|_{\cC^k}
+ \bigl\| (v^\nu)^{-1} 
\bigl( \Xi^\nu - \Xi \bigr) v^\nu \bigr\|_{\cC^k} \bigr)
\;=\; 0 .
\end{align*}
First, this implies that $v_n^\pm\in\cG(Y)$ 
is independent of $s\in\pm[T_n,T_n+1]$.
Secondly, by unique continuation (Proposition~\ref{prop:ucon}), 
it implies $(v_n^\pm)^*\Xi=\Xi$ and hence the limits $v_n^\pm\in\cG_{A^\pm}$
must lie in the stabilizer of the limit connections.
Now we can define the gauge transformations $u^\ell_n\in\cG(\R\times Y)$ by
$u^\ell_n=v^{\nu_n(\ell)}(v_n^\pm)^{-1}$ for $\pm s\geq T_n+1$,  
by $u^\ell_n=\one$ for $|s|\le T_n$, 
and, for $s\in\pm[T_n,T_n+1]$, by an interpolation which satisfies 
$d(u^\ell_n,\one)_{\cC^{k+1}((\pm[T_n,T_n+1])\times Y)}\to 0$ 
as $\ell\to\infty$.
With this choice we have
$$
\bigl\|(u^\ell_n)^*\Xi^{\nu_n(\ell)}-\Xi
\bigr\|_{W^{k,p}_\delta((\R\setminus[-T_n-1,T_n+1])\times Y)}
\le \frac{\rho_n}2
$$
from the exponential decay, as before for $(v^\nu)^*\Xi^\nu$, and
$$
\bigl\|(u^\ell_n)^*\Xi^{\nu_n(\ell)}-\Xi
\bigr\|_{W^{k,p}_\delta([-T_n-1,T_n+1])\times Y)}
\leq \frac{\rho_n}2
$$
for all sufficiently large $\ell\geq L_n$,
from the convergence of $\Xi^\nu$ and $u^\ell_n$ on compact subsets.
Now we can pick $\ell_n\geq L_n$ so large that 
$\nu_n:=\nu_n(\ell_n)\to\infty$ and
$\Norm{(u^{\nu_n})^*\Xi^{\nu_n}-\Xi}_{W^{k,p}_\delta(\R\times Y)}
\le \rho_n \to 0$.
This proves the theorem. 
\end{proof}

\begin{cor}\label{cor:compact11}
Let $h_f$ be a regular perturbation 
in the sense of Definition~\ref{def:freg},
and let $A^+,A^-\in\cA(Y,\cL)$ be nondegenerate and
irreducible critical points of $\CS_\cL+h_f$.
Then $\cM^1(A^-,A^+;X_f)/\R$ is compact
and hence is a finite set.
\end{cor}

\begin{proof}
Assume by contradiction that there is a sequence of distinct points 
$[\Xi^\nu]\in\cM^1(A^-,A^+;X_f)/\R$.
These solutions have index $1$ and hence 
fixed energy by Corollary~\ref{cor:index}~(i).
By Theorem~\ref{thm:compact1} we can pick 
a subsequence and representatives $\Xi^{\nu_k}$
that converge to a broken trajectory $(\Xi_1,\dots,\Xi_\ell)$ 
modulo bubbling.
By transversality we do not have solutions of negative index,
so Corollary~\ref{cor:compact} implies that there is no bubbling,
and the index identity in Theorem~\ref{thm:compact1} implies
that $\ell=1$.
Now Theorem~\ref{thm:compact3} implies that $\Xi^{\nu_k}$ 
converges to $\Xi_1$ in the $W^{1,p}_\delta$-norm.
Since $\cM^1(A^-,A^+;X_f)/\R$ is a $0$-manifold this implies 
that $\Xi^{\nu_k}$ is gauge equivalent to a time-shift of $\Xi_1$
in contradiction to the assumption.
\end{proof}

Finally we refine the compactness theorem in the case when bubbling
is excluded but breaking can take place. The precise convergence
statement here will be important for the gluing theory.

\begin{thm}\label{thm:compact2}
Fix a constant $p>1$.
Let $f$, $f^\nu$, $\Xi^\nu$, $s^\nu_i$, $u^\nu_i$, and $\Xi_i$
be as in the conclusion of 
Theorem~\ref{thm:compact1} and suppose that no bubbling 
occurs, i.e.~the curvature of $\Xi^\nu$ is uniformly bounded, 
$((u^\nu_i)^*\Xi^\nu)(s_i^\nu+\cdot)$ converges to $\Xi_i$ 
in the $\cC^k$ topology on compact sets, and 
\begin{equation}\label{eq:energy1}
\sum_{i=1}^\ell E_f(\Xi_i) = \lim_{\nu\to\infty}E_{f^\nu}(\Xi^\nu).
\end{equation}
Then the following holds.

\smallskip\noindent{\bf (i)}
If $\cD_{\Xi_i,\delta}$ is surjective for $i=1,\dots,\ell$ 
then so is $\cD_{\Xi^\nu,\delta}$ for $\nu$ sufficiently large.

\smallskip\noindent{\bf (ii)}
If the set of critical points of $\CS_\cL+h_{f^\nu}$ 
is independent of $\nu$ then, after replacing the broken 
trajectory $(\Xi_1,\ldots,\Xi_\ell)$ by a gauge
equivalent one, and for a subsequence,
there exists a sequence of gauge transformations 
$u^\nu\in\cG(\R\times Y)$ such that
$$
\lim_{\nu\to\infty}
\bigl\|(u^\nu)^*\Xi^\nu-\Xi_i(\cdot-s^\nu_i)\bigr\|_{W^{1,p}(I^\nu_i\times Y)}
= 0,\qquad\text{for}\;\; i=1,\dots,\ell,
$$
$$
I^\nu_i := \left\{\begin{array}{ll}
{(-\infty,\frac34 s^\nu_2 + \frac 14 s^\nu_1]},& i=1, \\
{[\frac 34 s^\nu_{i-1}+ \frac 14 s^\nu_i,
\frac 34 s^\nu_{i+1} + \frac 14 s^\nu_i]},& i=2,\dots,\ell-1,\\
{[\frac 34 s^\nu_{\ell-1}+\frac 14 s^\nu_\ell,\infty)},& i=\ell.
\end{array}\right.
$$
\end{thm}

\begin{proof}
Fix a constant $\eps>0$ smaller than 
the constant of Theorem~\ref{thm:long} and recall that 
the sequences $s^\nu_i$ in Theorem~\ref{thm:compact1} 
are chosen such that
\begin{equation}\label{eq:energy2}
E_{f^\nu}(\Xi^\nu;(-\infty,s^\nu_1]) = \eps/2,
\qquad
E_{f^\nu}(\Xi^\nu;[s^\nu_i,s^\nu_{i+1}])
=E_f(\Xi_i)
\end{equation}
for $\nu$ sufficiently large and $i=1,\dots,\ell-1$.
Since $s^\nu_{i+1}-s^\nu_i\to\infty$ we have for any $T>0$
$$
E_{f^\nu}(\Xi^\nu;[s^\nu_i,s^\nu_i+T])
+E_{f^\nu}(\Xi^\nu;[s^\nu_{i+1}-T,s^\nu_{i+1}])
\le E_{f^\nu}(\Xi^\nu;[s^\nu_i,s^\nu_{i+1}])
$$
for large $\nu$. With $\nu\to\infty$ this gives
$
E_f(\Xi_i;[0,T])+E_f(\Xi_{i+1};[-T,0])\le E_f(\Xi_i)
$
and, by taking the limit $T\to\infty$, 
$
E_f(\Xi_{i+1};(-\infty,0]) 
\le E_f(\Xi_i;(-\infty,0]).
$
Hence
$
E_f(\Xi_i;(-\infty,0])\le\eps/2
$
for all $i$. Choose $\tau_1,\dots,\tau_\ell$ such that
$$
E_f(\Xi_i;[-\tau_i,\tau_i]) = E_f(\Xi_i) - \eps/4.
$$
Then $E_f(\Xi_i;[0,\tau_i])\ge E_f(\Xi_i)-3\eps/4$ and hence 
$E_{f^\nu}(\Xi^\nu;[s^\nu_i,s^\nu_i+\tau_i]) > E_f(\Xi_i)-\eps$
for $\nu$ sufficiently large. 
Moreover, $E_{f^\nu}(\Xi^\nu;[s^\nu_\ell,\infty))$ converges
to $E_f(\Xi_\ell)-\eps/2$, by~(\ref{eq:energy1}) and~(\ref{eq:energy2}).
In summary we have for $i=0,\dots,\ell$ and $\nu$ sufficiently large
\begin{equation}\label{eq:longsmall}
E_{f^\nu}(\Xi^\nu;J^\nu_i) < \eps,\qquad
J^\nu_i := \left\{\begin{array}{ll}
{(-\infty,s^\nu_1]},& i=0, \\
{[s^\nu_i+\tau_i,s^\nu_{i+1}]},& i=1,\dots,\ell-1, \\
{[s^\nu_\ell+\tau_\ell,\infty)},& i=\ell.
\end{array}\right.
\end{equation}
Now choose gauge transformations $v^\nu_i$ on $J_i^\nu\times Y$
such that $(v^\nu_i)^*\Xi^\nu$ is in temporal gauge
on $J_i^\nu\times Y$.  Thus each connection
$(v^\nu_i)^*\Xi^\nu$ is represented by 
a smooth path $\tA^\nu_i:J_i\to\cA(Y,\cL)$.
Then it follows from Theorem~\ref{thm:long}
that there are critical points $B_i^\nu\in\Crit(\CS_\cL+h_{f^\nu})$
and positive constants $C_0$ and $\delta_0$ 
such that, for $i=0,\dots,\ell$, $\tau\ge\tau_i+1$, and $\nu$ 
sufficiently large, we have 
\begin{equation} \label{eq:AB}
\bigl\|\tA^\nu_i - B^\nu_i\bigr\|_{\cC^0(J^\nu_i(\tau)\times Y)} 
+ \bigl\|\tA^\nu_i - B^\nu_i\bigr\|_{W^{1,p}(J^\nu_i(\tau)\times Y),B^\nu_i}
\le C_0 e^{-\delta_0(\tau-\tau_i)} \sqrt{\eps}.
\end{equation}
Here we abbreviate $\tau_0:=0$ and
$$
J^\nu_i(\tau):= \left\{\begin{array}{ll}
{(-\infty,s^\nu_1-\tau]},& i=0, \\
{[s^\nu_i+\tau,s^\nu_{i+1}-\tau]},& i=1,\dots,\ell-1, \\
{[s^\nu_\ell+\tau,\infty)},& i=\ell.
\end{array}\right.
$$
Moreover we use the fact that the constants in  
Theorem~\ref{thm:long} can be chosen uniform for all $f^\nu$. 
Since the estimate is gauge invariant we may modify the gauge 
transformations $v^\nu_i$ so that the sequence $B^\nu_i$ 
converges in the $\cC^k$-norm to the critical point
$B_i$ in the assertion of Theorem~\ref{thm:compact1} 
for every $i$ (see Proposition~\ref{prop:finite}).
Then (\ref{eq:AB}) continues to hold if we drop the subscript 
$B^\nu_i$ in the $W^{1,p}$-norm and replace $C_0$ with a possibly 
larger constant, still denoted by $C_0$. 

Under the assumption of~(ii) we may choose $v^\nu_i$ 
so that $B^\nu_i=B_i$ is independent of $\nu$.  
Now we can argue as in the proof of Theorem~\ref{thm:compact3}.
Combining~(\ref{eq:AB}) with $B^\nu_i=B_i$ and the 
exponential decay of $\Xi_i$ and $\Xi_{i+1}$
we obtain the estimates
\begin{equation}\label{eq:vnu}
\begin{split}
\Norm{(v^\nu_0)^*\Xi^\nu-\Xi_1(\cdot - s^\nu_1) 
}_{W^{1,p}((-\infty,s^\nu_1-\tau]\times Y)} 
&\le C_1 e^{-\delta_0\tau}, \\
\Norm{(v^\nu_i)^*\Xi^\nu-\Xi_i(\cdot - s^\nu_i) 
}_{W^{1,p}([s^\nu_i+\tau,\frac34s^\nu_{i+1}+\frac14s^\nu_i]\times Y)} 
&\le C_1 e^{-\delta_0(\tau-\tau_i)}, \\
\Norm{(v^\nu_i)^*\Xi^\nu-\Xi_{i+1}(\cdot - s^\nu_{i+1}) 
}_{W^{1,p}([\frac34s^\nu_i+\frac14s^\nu_{i+1},s^\nu_{i+1}-\tau]\times Y)} 
&\le C_1 e^{-\delta_0(\tau-\tau_i)}, \\
\Norm{(v^\nu_\ell)^*\Xi^\nu-\Xi_\ell(\cdot - s^\nu_\ell) 
}_{W^{1,p}([s^\nu_\ell+\tau,\infty)\times Y)} 
&\le C_1 e^{-\delta_0(\tau-\tau_\ell)}
\end{split}
\end{equation}
for $\nu$ sufficiently large, some constant $C_1$,
and $i=1,\dots,\ell-1$. Fix a constant $\rho>0$ 
and choose $\tau$ so large that 
$$
C_1 e^{-\delta_0(\tau-\tau_i)}\le\rho/4 \qquad \text{for}\;\; i=1,\ldots,\ell .
$$  
Then, on the interval 
$[s^\nu_{i}-\tau-1,s^\nu_{i}-\tau]\subset J^\nu_{i-1}(\tau)$
the connections $(v^\nu_{i-1})^*\Xi^\nu$ and $(u^\nu_{i})^*\Xi^\nu$ 
are both $W^{1,p}$ close to $\Xi_{i}(\cdot - s^\nu_{i})$
Thus $\bigl((v^\nu_{i-1})^{-1}u^\nu_{i}\bigr)(\cdot+ s^\nu_{i})$ 
is bounded in $W^{2,p}([-\tau-1,-\tau]\times Y)$ and thus, 
for a subsequence, converges to a gauge transformation 
$g_i^-\in\cG^{2,p}([-\tau-1,-\tau]\times Y)$.
For the limit we obtain $(g^-_i)^*\Xi_i=\Xi_i$ on $(-\infty,-1]$
as in Theorem~\ref{thm:compact3}, and we deduce that $g^-_i\in\cG_{B_i}$.
Similarly, we can pick the subsequence such that
$\bigl((v^\nu_i)^{-1} u^\nu_i\bigr) (\cdot+ s^\nu_i)\to g_i^+\in\cG_{B_{i+1}}$
in $W^{2,p}([\tau,\tau+1]\times Y)$ with $(g^+_i)^*\Xi_i=\Xi_i$ on $[1,\infty)$.
With this we can now construct a sequence $u^\nu\in\cG(\R\times Y)$
that satisfies
\begin{itemize}
\setlength{\itemsep}{0ex}
\item
$u^\nu(s) = v^\nu_0(s) g_1^- $ for $s\in (-\infty, s^\nu_1-\tau-1]$,
\item
$u^\nu(s) = u^\nu_1(s)$ for $s\in[s^\nu_1-\tau,s^\nu_1+\tau]$, 
\item
$u^\nu(s) = v^\nu_1(s) g_i^+$ for $s\in[s^\nu_1+\tau+1,s^\nu_{2}-\tau-1]$,  
\item
$u^\nu(s) = u^\nu_i(s) (g_i^-)^{-1} g_{i-1}^+ (g_{i-1}^-)^{-1} \ldots (g_2^-)^{-1}g_1^+$ for 
$s\in[s^\nu_i-\tau ,s^\nu_i+\tau]$ and $i=2,\ldots,\ell$,
\item
$u^\nu(s) = v^\nu_i(s) g_i^+ (g_i^-)^{-1} g_{i-1}^+  \ldots (g_2^-)^{-1}g_1^+$
for $s\in[s^\nu_i+\tau+1,s^\nu_{i+1}-\tau-1]$ and $i=2,\ldots,\ell$.
\item
$\bigl((u^\nu_i)^{-1}u^\nu\bigr)(\cdot+s^\nu_i) \to 
(g_i^-)^{-1} g_{i-1}^+ (g_{i-1}^-)^{-1} \ldots (g_2^-)^{-1}g_1^+$ 
as $\nu\to\infty$ in $W^{2,p}([-\tau-1,\tau+1]\times Y,\rG)$ for $i=1,\ldots,\ell$,
\item
${\rm dist}_{W^{2,p}([s^\nu_i+\tau,s^\nu_{i+1}-\tau])} 
\bigl((v^\nu_i)^{-1}u^\nu , g_i^+ (g_i^-)^{-1} g_{i-1}^+ (g_{i-1}^-)^{-1} 
\hspace{-1pt} \ldots (g_2^-)^{-1}g_1^+\bigr) \to 0$ 
as $\nu\to\infty$ for $i=0,\ldots,\ell$.
\end{itemize}
At the same time we replace the broken trajectory 
$(\Xi_1,\ldots,\Xi_\ell)$ with $\Xi'_1:=\Xi_1$ and 
$\Xi_i':=\bigl((g_i^-)^{-1} g_{i-1}^+ (g_{i-1}^-)^{-1} 
\ldots (g_2^-)^{-1}g_1^+\bigr)^*\Xi_i$ 
for $i=2,\ldots,\ell$.  Note that this again defines a broken 
trajectory $(\Xi'_1,\ldots,\Xi'_\ell)$ between the critical points
\begin{align*}
\lim_{s\to\infty} \Xi_i' 
&= \bigl((g_i^-)^{-1} g_{i-1}^+  \ldots (g_2^-)^{-1}g_1^+\bigr)^* B_{i+1} \\
&= \bigl((g_i^-)^{-1} g_{i-1}^+  \ldots (g_2^-)^{-1}g_1^+\bigr)^* 
\bigl((g_{i+1}^-)^{-1} g_{i}^+\bigr)^* B_{i+1} 
=\lim_{s\to-\infty}\Xi'_{i+1} .
\end{align*}
Here we used the fact that $g_{i+1}^-,g_{i}^+\in\cG_{B_{i+1}}$.
The convergence of $(u^\nu_i)^*\Xi^\nu$ then implies
\begin{equation*}
\bigl\|(u^\nu)^*\Xi^\nu-\Xi_i'(\cdot - s^\nu_i)
\bigr\|_{W^{1,p}([s^\nu_i-\tau,s^\nu_i+\tau]\times Y)} 
\le \rho/3
\end{equation*}
for large $\nu$ and $i=1,\dots,\ell$, 
and from the exponential decay (\ref{eq:vnu})
we obtain
\begin{equation*}
\begin{split}
\bigl\|(u^\nu)^*\Xi^\nu-\Xi_i'(\cdot - s^\nu_i) 
\bigr\|_{W^{1,p}([s^\nu_i+\tau,\frac34s^\nu_{i+1}+\frac14s^\nu_i]\times Y)} 
&\le \rho/3, \\
\bigl\|(u^\nu)^*\Xi^\nu-\Xi_i'(\cdot - s^\nu_i) 
\bigr\|_{W^{1,p}([\frac34s^\nu_{i-1}+\frac14s^\nu_i,s^\nu_i-\tau]\times Y)} 
&\le \rho/3, 
\end{split}
\end{equation*}
for large $\nu$, large $\tau$, and $i=1,\dots,\ell$.
Here we denote $s^\nu_0:=-\infty$ and ${s^\nu_{\ell+1}:=\infty}$, 
and we use the fact that $(g_i^+)^*\Xi_i=\Xi_i$ on $[1,\infty)$ 
and $(g_i^-)^*\Xi_i=\Xi_i$ on $(-\infty,-1]$.
Thus, for every $\rho>0$, we have a subsequence 
$(\nu_n)_{n\in\N}$ and a sequence of gauge 
transformations $u^{\nu_n}_\rho$ such that
$\Norm{(u_\rho^{\nu_n})^*\Xi^{\nu_n}-\Xi_i'(\cdot - s^{\nu_n}_i)
}_{W^{1,p}(I^{\nu_n}_i\times Y)} \le \rho$
holds for all sufficiently large $n\ge N_\rho$.
Assertion~(ii) then follows by taking a diagonal subsequence.

To prove~(i) we can assume by contradiction that, after passing to a subsequence, 
none of the $\cD_{\Xi^\nu,\delta}$ is surjective.
Then we use the $\cC^0$-estimate of~(\ref{eq:AB})
and the same patching construction as for (ii) 
to find a further subsequence and a sequence of 
gauge transformations $u^\nu\in\cG(\R\times Y)$ such that
\begin{equation}\label{eq:argh}
\lim_{\nu\to\infty}
\Norm{(u^\nu)^*\Xi^\nu-\Xi_i'(\cdot-s^\nu_i)
}_{\cC^0(I^\nu_i\times Y)} = 0
\end{equation}
for $i=1,\dots,\ell$. (The $\cC^0$-estimate holds 
on increasingly large domains because $B^\nu_i\to B_i$ 
converges in $\cC^0(\R\times Y)$ - but not in $W^{1,p}(\R\times Y)$.)
By Theorem~\ref{thm:weighted fredholm}~(iv) 
the surjectivity of the linearized operators
$\cD_{\Xi^\nu,\delta}$ is independent of a timeshift 
in the weight function, or equivalently in the connection.
Hence, applying an overall timeshift to each element of the 
sequence $\Xi^\nu$, we may assume w.l.o.g.~that $s^\nu_1=0$ 
and for each $i\ge 2$ we have $s^\nu_i\to\infty$.
By assumption, the linearized operator 
$w\cD_{\Xi_i,\delta}w^{-1}=\Nabla{s}+\cH_{A_i}-I_\lambda$
is surjective on the unweighted Sobolev spaces, 
see Remark~\ref{rmk:halfonto}~(iv),
and so are the operators $\Nabla{s}+\cH_{A_i} - I_\delta$.
(Recall that $\lambda=\pd_s V=w^{-1}\pd_s w$ 
denotes the derivative of the weight function.)
Equivalently, the adjoint operators $-\Nabla{s}+\cH_{A_i}-I_\lambda$ 
resp.\ $-\Nabla{s}+\cH_{A_i}-I_\delta$ are injective.
Hence there is a constant $c$ such that 
\begin{align*}
\Norm{\xi}_{L^p}
\le c\bigl\|-\Nabla{s}\xi +\cH_{A_i}\xi-I_\lambda\xi\bigr\|_{L^p} , \qquad
\Norm{\xi}_{L^p}
\le c\bigl\|-\Nabla{s}\xi +\cH_{A_i}\xi-I_\delta\xi\bigr\|_{L^p} 
\end{align*}
for every $\xi\in W^{1,p}(\R\times Y,\rT^*Y\otimes\cg\oplus\cg)$.
This estimate is stable under $\cC^0$-small perturbations
of $\Xi_i$ and under the action of the gauge group.
Hence, enlarging the constant $c$ if necessary, we obtain 
$$
\supp\xi\subset I^\nu_i\times Y
\qquad\implies\qquad
\Norm{\xi}_{L^p}
\le c\bigl\|-\Nabla{s}\xi+\cH_{A^\nu}\xi-I_\lambda\xi\bigr\|_{L^p}
$$
for all $\xi$ and $i$ and for $\nu$ sufficiently large.  
For $i\ge 1$ this follows directly from (\ref{eq:argh}) 
with $s^\nu_1=0$. For $i\ge 2$ we use the fact that 
$\lambda(s)=\delta$ for all $s\in I^\nu_i$, so we can estimate 
$\bigl\|-\Nabla{s}\xi+\cH_{A^\nu}\xi-I_\lambda\xi\bigr\|_{L^p}$ by
$\bigl\|-\Nabla{s}\xi+\cH_{A_i(\cdot-s^\nu_i)}\xi-I_\delta\xi\bigr\|_{L^p}
- \bigl\|(u^\nu)^*\Xi^\nu-\Xi_i(\cdot-s^\nu_i)\bigr\|_{\cC^0(I^\nu_i\times Y)}
\|\xi\|_{L^p}$
and identify the first term of this with
$\bigl\|(-\Nabla{s}+\cH_{A_i}-I_\delta) 
\xi(\cdot+s^\nu_i)\bigr\|_{L^p}\ge c^{-1}\|\xi\|_{L^p}$.

Now for each $\nu$ we can choose a partition of unity
$h^\nu_i:\R\to[0,1]$ with $\supp h^\nu_i\subset I^\nu_i$ and 
$\sum_i\Norm{\p_sh^\nu_i}_{L^\infty}\to0$.
Then we obtain
\begin{equation*}
\begin{split}
\Norm{\xi}_{L^p}
\le \sum_{i=1}^\ell \Norm{h^\nu_i\xi}_{L^p} 
&\le c\sum_{i=1}^\ell
\bigl\|h^\nu_i(-\Nabla{s}\xi+\cH_{A^\nu}\xi-I_\lambda\xi)
-(\p_sh^\nu_i)\xi\bigr\|_{L^p} \\
&\le \ell c \bigl\|-\Nabla{s}\xi+\cH_{A^\nu}\xi-I_\lambda\xi
\bigr\|_{L^p}
+ \sum_{i=1}^\ell\Norm{\p_sh^\nu_i}_{L^\infty}
\Norm{\xi}_{L^p}.
\end{split}
\end{equation*}
This shows that the operator 
$-\Nabla{s}+\cH_{A^\nu}-I_\lambda$ is injective 
on the unweighted Sobolev spaces for $\nu$ sufficiently large,
and hence its adjoint $\Nabla{s}+\cH_{A^\nu}-I_\lambda$ is surjective.
Since the latter operator is conjugate to $\cD_{\Xi^\nu,\delta}$ this
is a contradiction to the assumption, and the theorem is proved. 
\end{proof}


\section{Transversality}\label{sec:moduli}

Let $(Y,g)$ be a compact oriented Riemannian $3$-manifold with 
metric $g$ and boundary ${\pd Y=\Sigma}$, 
and let $\cL\subset\cA(\Sigma)$ 
be a gauge invariant, monotone, irreducible Lagrangian submanifold 
satisfying (L1-3) on page~\pageref{p:L1}.
Then $\R\times Y$ naturally is a Riemannian $4$-manifold 
with boundary space--time splitting and tubular ends 
in the sense of definition~\ref{def:4g}.
In order to complete the instanton data we must also choose
a perturbation. A detailed construction of holonomy perturbations
is given in Appendix~\ref{app:Xf}. In this section we concentrate on
achieving transversality by the choice of perturbation.

Fix an embedding $\beta:[-1,1]\times\D\to\mathrm{int}(Y)$
and denote by $\Gamma_m$ the set of finite sequences
$\gamma=(\gamma_1,\dots,\gamma_m)$ of embeddings
$\gamma_i:S^1\times\D\to\mathrm{int}(Y)$ that agree 
with $\beta$ in a neighbourhood of $\{0\}\times\D$.
Every $\gamma\in\Gamma_m$ gives rise to a map
$$
\rho=(\rho_1,\dots,\rho_m):\D\times\cA(Y)\to\rG^m
$$ 
where $\rho_i(z,A)$ is the holonomy of the connection $A$ 
around the loop $\gamma_i(\cdot,z)$. Let
$
\cF_m:=\cC^\infty_0(\D\times\rG^m)^\rG
$
denote the space of conjugation invariant real valued 
compactly supported smooth functions on ${\D\times\rG^m}$. 
Each pair ${(\gamma,f)\in\Gamma_m\times\cF_m}$ determines 
a smooth function $h_f:\cA(Y)\to\R$ via
$$
h_f (A) := \int_\D f(z,\rho(z,A))  \,\rd^2 z.
$$
The differential $\rd h_f (A):\rT_A\cA(Y)\to\R$ 
has the form
$$
\rd h_f (A) \alpha = \int_Y \winner{X_f(A)}{\alpha}.
$$
Here $X_f:\cA(Y)\to\Om^2(Y,\cg)$ is a smooth function 
satisfying~(\ref{Xf identities}). We emphasize that 
the tuple $(\one,\dots,\one)$ is a critical point
of every conjugation invariant function $\rG^m\to\R$
and hence the trivial connection $A=0$ is always
a critical point of the perturbed Chern--Simons functional
$\CS_\cL+h_f$; it is nondegenerate by assumption~(L3). 

\begin{dfn}\label{def:Areg}
Fix a perturbation $(\gamma,f)\in\Gamma_m\times\cF_m$ and
two nondegenerate critical points $A^\pm\in{\rm Crit}(\CS_\cL+h_f)$.
A finite energy solution $A:\R\to\cA(Y,\cL)$ 
of the boundary value problem~(\ref{eq:floer})
with limits $\lim_{s\to\pm\infty}A(s)=A^\pm$ 
is called {\bf regular}
if the operator $\cD_{\A,\delta}$ defined in~(\ref{eq:DAd}) 
is surjective for every sufficiently small constant $\delta>0$. 
(This condition is independent of $k$ and $p$.)
\end{dfn}

\begin{dfn}\label{def:freg}
A pair $(\gamma,f)\in\Gamma_m\times\cF_m$ 
is called {\bf regular} (for $(Y,g)$ and $\cL$)  
if it satisfies the following.

\smallskip\noindent{\bf (i)}
Every nontrivial critical point of the perturbed Chern--Simons
functional ${\CS_\cL + h_f}$ 
is irreducible and nondegenerate, i.e. if $A\in\cA(Y,\cL)$ 
is not gauge equivalent to the trivial connection and
satisfies $F_A+X_f(A)=0$ then $H^0_A=0$ and $H^1_{A,f}=0$.

\smallskip\noindent{\bf (ii)}
Let $A:\R\to\cA(Y,\cL)$ be a finite energy solution
of the boundary value problem~(\ref{eq:floer})
with 
$
\delta_f(\A)\le7
$ 
and suppose that at most one of the limits $A^\pm$ is gauge equivalent
to the trivial connection. Then the operator
$\cD_{\A,\delta}$ defined in~(\ref{eq:DAd}) 
is surjective for every integer $k\geq 1$, every $p>1$,
and every sufficiently small constant $\delta>0$. 

\smallskip\noindent
For every $\gamma\in\Gamma_m$ the set of regular elements 
$f\in\cF_m$ will be denoted by $\cF_\mathrm{reg}(\gamma)$. 
\end{dfn}

If $f\in\cF_{\rm reg}(\gamma)$ and $([A^-],[A^+])\ne(0,0)$
then it follows from the discussion in Section~\ref{sec:fredholm} 
that the moduli space $\cM(A^-,A^+;X_f)$, 
introduced in~(\ref{eq:GAA}) and 
the beginning of Section~\ref{sec:compact},
is a smooth manifold of local dimension
$$
\dim_{[\A]} \cM(A^-,A^+;X_f) = \delta_f(\A) .
$$
For every integer $k\ge 1$ we introduce 
the following seminorm on the space of perturbations
$$
\NORM{X_f}_{k} 
:= \sup_{A\in\cA(Y,\cL)}
\left(
\frac{\Norm{X_f(A)}_{\cC^k}}{(1+\Norm{A}_{C^k})^k }
+ \sup_{\alpha\in\rT_A\cA(Y,\cL)}
\frac{\Norm{\rd X_f(A)\alpha}_{\cC^{k-1}}}
{\Norm{\alpha}_{\cC^{k-1}}(1+\Norm{A}_{\cC^{k-1}})^{k-1}}
\right).
$$
We will apply this notation to the difference $X_f - X_{f'}$ 
associated to two pairs $(\gamma,f)\in\Gamma_m\times\cF_m$
and $(\gamma',f')\in\Gamma_{m'}\times\cF_{m'}$. 
This difference can be written as $X_{f-f'}$ associated to
the union $\gamma\cup\gamma'
:=(\gamma_1,\dots,\gamma_m,\gamma'_1,\dots,\gamma_{m'})
\in\Gamma_{m+m'}$,
where $f$ and $f'$ are extended to elements of $\cF_{m+m'}$ 
in the obvious way. Then Proposition~\ref{prop:Xf} implies that
$\NORM{X_{f_\nu}-X_{f_0}}_{k}\to 0$ for
$\|f_\nu-f_0\|_{\cC^{k+1}}\to 0$.

\begin{thm}\label{thm:fcrit}
{\bf (i)}
For every $\gamma\in\Gamma_m$ the set of all 
$f\in\cF_m$ that satisfy condition~(i)
in Definition~\ref{def:freg} is open in $\cF_m$
with respect to the $\cC^2$-topology.

\smallskip\noindent{\bf (ii)}
Let $(\gamma_0,f_0)\in\Gamma_{m_0}\times\cF_{m_0}$ be such that
every nontrivial critical point of ${\CS_\cL+h_{f_0}}$ is irreducible.
Then, for every $\eps>0$ and every $k\in\N$, there exists 
an $n\in\N$ and a pair $(\gamma,f)\in\Gamma_n\times\cF_n$ 
that satisfies condition~(i) in Definition~\ref{def:freg}
and
$
\NORM{X_f-X_{f_0}}_{k} < \eps.
$
\end{thm}

The zero perturbation satisfies the assumptions
of Theorem~\ref{thm:fcrit}~(ii) by~(L3).
Transversality for the critical points near the 
unperturbed equation was established by Taubes~\cite{T}.
The extension to large perturbations requires another proof,
similar to that of the following transversality result for trajectories.

\begin{thm}\label{thm:freg}
{\bf (i)}
The set $\cF_\mathrm{reg}(\gamma)$ is open in $\cF_m$ with 
respect to the $C^2$-topology for every $m\in\N$ and every
$\gamma\in\Gamma_m$. 

\smallskip\noindent{\bf (ii)}
Assume that $(\gamma_0,f_0)\in\Gamma_{m_0}\times\cF_{m_0}$ satisfies
condition~(i) in Definition~\ref{def:freg}.  Then, for every 
$\eps>0$ and $k\in\N$, there exists an $n\in\N$ 
and another pair $(\gamma,f)\in\Gamma_n\times\cF_n$ that is regular,
i.e.\ $f\in\cF_\mathrm{reg}(\gamma)$, and satisfies
$$
\Crit(\CS_\cL+h_{f_0})=\Crit(\CS_\cL+h_f),
$$
$$
A\in\Crit(\CS_\cL+h_{f_0})
\quad\implies\quad
h_f(A)=h_{f_0}(A),
$$
$$
\NORM{X_f-X_{f_0}}_{k} < \eps.
$$
\end{thm}

Note that we do not construct a Banach space of perturbations 
in which regular ones are of Baire second category.
The main reason for this is that the loops in the interior of $Y$
do not form a Banach space.

\begin{rmk}\label{rmk:rho}\rm
Fix a point $y_0\in{\rm int}(Y)$. 
For every based, 
embedded loop $\gamma:[0,1]\to{\rm int}(Y)$ with $\gamma(0)=\gamma(1)=y_0$
denote by $\rho_\gamma:\cA(Y)\to\rG$ the holonomy map.  
For later reference we state two facts that follow from the 
equivalence between connection $1$-forms and parallel transport.
(Note that it suffices to use embedded loops in the interior.)

\smallskip\noindent{\bf (i)}
Two connections $A,B\in\cA(Y)$ are gauge equivalent 
if and only if there is a $g_0\in\rG$ such that 
$$
\rho_\gamma(B) =  g_0^{-1}\rho_\gamma(A)g_0
$$
for every based embedded loop $\gamma$.

\smallskip\noindent{\bf (ii)}
Let $A\in\cA(Y)$ and $\alpha\in\Om^1(Y,\cg)$.
Then $\alpha\in\im\rd_A$ if and only if there 
is a $\xi_0\in\cg$ such that 
$$
\rd\rho_\gamma(A)\alpha 
=  \rho_\gamma(A)\xi_0-\xi_0\rho_\gamma(A)
$$
for every based embedded loop $\gamma$.
\end{rmk}

\begin{proof}[Proof of Theorem~\ref{thm:fcrit}.]
Assertion~(i) follows from the fact that the 
conditions $H^0_A=0$ and $H^1_{A,f}=0$ are open with 
respect to $\cC^2$-variations of $f$ and $A$.
The conditions are moreover gauge invariant, and the
set of nontrivial critical points of $\CS_\cL+h_f$ is 
compact in $\cA(Y,\cL)/\cG(Y)$ for every perturbation~$f$.  
(This follows from Uhlenbeck compactness \cite{U,W}
since $F_A=-X_f(A)$ is $L^\infty$-bounded.)
The proof of~(ii) has three steps.

\medskip\noindent{\bf Step 1.}
{\it Let $(\gamma_0,f_0)\in\Gamma_{m_0}\times\cF_{m_0}$ be given.
Then there is a $\gamma\in\Gamma_m$ 
with $\gamma_i=\gamma_{0i}$ for $i=1,\dots,m_0$
satisfying the following condition.
Define $\sigma:\cA(Y)\to\rG^m$ by 
$$
\sigma(A) := \rho(0,A) = (\rho_1(0,A),\dots,\rho_m(0,A)).
$$
Then, for every critical point $A\in\Crit(\CS_\cL+h_{f_0})$
and every nonzero $1$-form $\eta\in\Om^1(Y,\cg)$ satisfying
\begin{equation}\label{eq:eta0}
\rd_A\eta+\rd X_{f_0}(A)\eta=0,\qquad \rd_A^*\eta=0,\qquad
\eta|_{\pd Y}\in\rT_A\cL,\qquad *\eta|_{\pd Y}=0,
\end{equation}
the vector $[\rd\sigma(A)\eta]\in\rT(\rG^m/\rG)$ is nonzero.}

\medskip\noindent
The trivial connection is nondegenerate by assumption (L3),
so for $\eta\ne 0$ we must have $[A]\ne[0]$, and so by
assumption $A$ is irreducible.
The condition $[\rd\sigma(A)\eta]\neq 0$ is open
with respect to variations of $(A,\eta)$, and it 
is invariant under gauge transformations
$(A,\eta)\mapsto (u^*A,u^{-1}\eta u)$.
Moreover, the set of gauge equivalence classes of pairs 
$(A,\eta)\in\Crit(\CS_\cL+h_{f_0})\times\Om^1(Y,\cg)$
that satisfy $\Norm{\eta}_{L^2}=1$, $[A]\ne[0]$, 
and~(\ref{eq:eta0}) is compact.
(For $\eta$ this follows from elliptic estimates for 
the operator $\rd_A\oplus\rd_A^*$ with boundary condition
$*\eta|_{\pd Y}=0$, see e.g.\ \cite[Theorem~D]{W}.)
Hence it suffices to construct $\gamma$ for a single such pair $(A,\eta)$.
We shall use Remark~\ref{rmk:rho}~(ii) to construct $\gamma$.
In each step it suffices to find the loops $\theta\mapsto\gamma_i(\theta,0)$
(with base point $y_0:=\beta(0,0)$).  Since the condition is open with 
respect to smooth variations of $\gamma$, these loops can be deformed 
and extended to the required
embeddings of $\D\times S^1$ into the interior of $Y$. 

Since $A$ is irreducible we can choose the loops 
$\gamma_{m_0+1}$ and $\gamma_{m_0+2}$ such that 
the matrices $g_1:=\rho_{m_0+1}(0,A)$ and $g_2:=\rho_{m_0+2}(0,A)$ 
do not commute.  Then $\sigma(A)$ lies in the free part of $\rG^m$. 
The tangent space of the $\rG$-orbit through $\sigma(A)$ is
$$
V_0 := \bigl\{v=\bigl( \sigma_i(A)\xi-\xi\sigma_i(A)\bigr)_{i=1,\dots,m}
\,\big|\, \xi\in\cg \bigr\} \subset\rT_{\sigma(A)}\rG^m.
$$
We prove that $\gamma$ can be chosen such that
$\rd\sigma(A)\eta\notin V_0$.

Since $\eta\perp\im\,\rd_A$, it follows from Remark~\ref{rmk:rho}~(ii) 
that for every $\xi\in\cg$ there is a based loop $\gamma$
such that
\begin{equation}\label{eq:xieta}
\rd\rho_\gamma(A)\eta \ne \rho_\gamma(A) \xi - \xi \rho_\gamma(A).
\end{equation}
Since the map
$\xi\mapsto (g_1\xi - \xi g_1 , g_2\xi - \xi g_2 )$
is injective there is a constant $C$ such that
for $|\xi|\geq C$ condition~(\ref{eq:xieta})
holds for one of the loops $\gamma_{m+1}(0,\cdot)$ 
or $\gamma_{m+2}(0,\cdot)$.  The compact set $\{|\xi|\leq C\}$ 
can be covered by finitely many open sets $U_j$, on each of which 
condition (\ref{eq:xieta}) holds with the same loop $\gamma_{m+2+j}$.
Thus we have proved that for every $\xi\in\cg$ there exists an $i$
such that (\ref{eq:xieta}) holds with $\gamma=\gamma_i$.
This implies that $\rd\sigma(A)\eta$ is not contained in 
$V_0$ and hence does not vanish in the tangent space of the
quotient $\rG^m/\rG$.

\medskip\noindent{\bf Step~2.}
{\it Let $\gamma\in\Gamma_m$ be as in Step~1
and fix $p>3$. For $k\in\N$ and $\eps>0$ denote 
\begin{align*}
\cF^{k,\eps}_m 
&:= \bigl\{f\in\cC^{k+1}(\D\times\rG^{m})^\rG \st 
\|f-f_0\|_{\cC^{k+1}}<\eps\bigr\},
\end{align*}
let $\cA^{1,p}(Y,\cL)$ and $\cG^{2,p}(Y)$
denote the $W^{1,p}$- and $W^{2,p}$-closure of
$\cA(Y,\cL)$ and $\cG(Y)$ respectively, 
and 
$$
\widetilde\cM^*(\cF_m^{k,\eps}) 
:= \bigl\{
(A,f)\in\cA^{1,p}(Y,\cL)\times\cF_m^{k,\eps} \,\big|\,
F_A+X_f(A)=0,\,[A]\ne[0]
\bigr\}. 
$$
Then for every $k\in\N$ there is an $\eps>0$ such that
the moduli space 
$$
\cM^*(\cF_m^{k,\eps})
:= \widetilde\cM^*(\cF_m^{k,\eps})/\cG^{2,p}(Y)
$$
is a separable $\cC^k$ Banach manifold.}

\medskip\noindent
We denote
$
W^{1,p}_{\rT_A\cL}(Y,\rT^*Y\otimes\cg)
:=\left\{\alpha\in W^{1,p}(Y,\rT^*Y\otimes\cg)\,\big|\,
\alpha|_{\pd Y}\in\rT_A\cL\right\}
$
and 
${\cF^k_m:=\cC^{k+1}(\D\times\rG^{m})^\rG = \rT_f\cF^{k,\eps}_m}$,
and consider the operator
$$
W^{1,p}_{\rT_A\cL}(Y,\rT^*Y\otimes\cg)\times
W^{1,p}(Y,\cg)\times\cF_m^k
\to L^{p}(Y,\rT^*Y\otimes\cg)\times L^{p}(Y,\cg)
$$
given by 
\begin{equation}\label{eq:Huniv}
(\alpha,\phi,\hat f)\mapsto
\bigl(*\rd_A\alpha+*\rd X_f(A)\alpha 
-\rd_A\phi + *X_{\hat f}(A),-\rd_A^*\alpha\bigr) .
\end{equation}
This operator is $\cH_A \times (*X_{\cdot}(A), 0)$ and hence 
it is the linearized operator of $\widetilde\cM^*(\cF^{k,\eps})$
together with the local slice condition for the $\cG^{2,p}(Y)$-action.
(The nonlinear operator is a $\cC^k$ map since the map
${*X_f:\cA^{1,p}(Y)\to L^p(Y,\rT^*Y\otimes\cg)}$ is $\cC^k$
for $f\in\cC^{k+1}$.)
We must prove that this operator is surjective for every pair
$(A,f)\in\cM^*(\cF_m^{k,\eps})$ when $\eps$ is sufficiently small.
We first prove this for $f=f_0$.  Suppose, by contradiction that there
is a nontrivial critical point $A\in\Crit(\CS_\cL+h_{f_0})$
such that the operator~(\ref{eq:Huniv}) is not onto.  
Then with $q^{-1}=1-p^{-1}$ there is a nonzero element
$$
(\eta,\xi)\in L^q(Y,\rT^*Y\otimes\cg)\times L^q(Y,\cg)
$$
orthogonal to the image of~(\ref{eq:Huniv}). 
Any such element satisfies 
$$
\rd_A\xi=0,\qquad
*\rd_A\eta+*\rd X_{f_0}(A)\eta=0,\qquad
\rd_A^*\eta=0, 
$$
and
\begin{equation}\label{eq:Htrans}
\rd h_{\hat f}(A)\eta = \int_Y\winner{X_{\hat f}(A)}{\eta} = 0
\end{equation}
for every $\hat f\in\cF^k_{m,\eps}$.    
This implies $\xi=0$ because $A$ was assumed to be irreducible.
Since $\eta\ne 0$ it follows from Step~1 that 
$\rd\sigma(A)\eta\ne 0$ and hence the map 
$\R\to\rG^m/\rG:r\mapsto [\rho(0,A+r\eta)]$
is an embedding into the free part of the quotient
near $r=0$. This implies that there exists a map 
$\hat f\in\cF^k_m$ such that
$$
\hat f(z,\rho(z,A+r\eta))=r\beta(r)\beta(|z|),
$$
where $\beta:\R\to[0,1]$ is a smooth cutoff function
that is supported in a sufficiently small neighbourhood
of $0$ and is equal to $1$ near $0$.  Hence 
$$
\rd h_{\hat f}(A)\eta
= \left.\frac{\rd}{\rd r}\right|_{r=0}
\int_\D\hat f(z,\rho(z,A+r\eta))\,\rd^2z
= \int_\D\beta(|z|)\rd^2z >0
$$
in contradiction to~(\ref{eq:Htrans}). 
This proves that the operator~(\ref{eq:Huniv}) is onto
whenever $f=f_0$ and $[A]\ne[0]$.  That this continues to hold for 
$\Norm{f-f_0}_{\cC^{k+1}}$ sufficiently small follows
from compactness and the fact that the trivial connection
is nondegenerate. 

\medskip\noindent{\bf Step~3.}
{\it We prove~(ii).}

\medskip\noindent
By Step~2, the projection $\cM^*(\cF_m^{k,\eps})\to\cF_m^{k,\eps}$ 
is a $C^k$ Fredholm map of Fredholm index zero.
(Its linearization 
$\ker(\cH_A+(*X_\cdot(A),0)\to \rT\cF_m^{k,\eps}$ 
has the same index as the self--adjoint operator $\cH_A$.)
Hence it follows from the Sard--Smale theorem that
the set of regular values of this projection is dense 
in $\cF_m^{k,\eps}$.  
For such a regular value $f\in\cF_m^{k,\eps}$ we have
$\im(*X_\cdot(A),0)\subset\im\cH_A$, so by the surjectivity 
in Step~2, the operator $\cH_A$ itself is surjective 
and hence injective. This shows that $H^1_{A,f}=0$
for all critical points $A\in\Crit(\CS_\cL + h_f)$.
For $\|f-f_0\|_{\cC^2}$ sufficiently small we also
have $H^0_A=0$ by (i), and hence $f$ is `regular'
in the sense that Definition~\ref{def:freg}~(i)
is satisfied. So we have seen that $f_0\in\cF_m$ can be 
approximated by a sequence of `regular'
$\cC^{k+1}$ perturbations $f^\nu\in\cF_m^k$
and due to~(i) also by a sequence of `regular'
smooth perturbations.  This proves the theorem.
\end{proof}

\begin{proof}[Proof of Theorem~\ref{thm:freg}.]
To prove~(i) we suppose by contradiction that there is a $\gamma\in\Gamma_m$
and a sequence $f^\nu\in\cF_m\setminus\cF_{\mathrm{reg}}(\gamma)$ 
converging to some $f\in\cF_{\mathrm{reg}}(\gamma)$ in the 
$C^2$ topology.  By Theorem~\ref{thm:fcrit} we may assume that 
each $f^\nu$ satisfies condition~(i) in Definition~\ref{def:freg}. 
Thus there is a sequence 
${\A^\nu\in\widetilde{\cM}(A^\nu_-,A^\nu_+;X_{f^\nu})}$ 
such that $\delta_{f^\nu}(\A^\nu)\le7$, at most one of the limits $A^\nu_\pm$ 
is gauge equivalent to the trivial connection, and the 
the operator $\cD_{\A^\nu,\delta}$ is not surjective. 
The sequence $\A^\nu$ has bounded energy by Corollary~\ref{cor:index}
and hence a subsequence converges to a 
broken Floer trajectory $(\A_1,\dots,\A_\ell)$
by Theorem~\ref{thm:compact1}.
Since $f\in\cF_{\rm reg}(\gamma)$, all moduli spaces with negative index
and at least one nontrivial limit connection are empty, 
and the assertion of Corollary~\ref{cor:compact} is wrong.
So neither bubbling nor self--connecting trajectories of $[0]$ can occur
in the limit.
Hence $\cD_{\A_j,\delta}$ is surjective for every $j$ and, 
by gluing (see Theorem~\ref{thm:compact2}~(i)), 
the operator $\cD_{\A^\nu,\delta}$ is surjective 
for $\nu$ sufficiently large. This contradiction proves~(i).

We prove~(ii). 
By assumption $\CS_\cL+h_{f_0}$ has only finitely many critical points
in the configuration space $\cA(Y,\cL)/\cG(Y)$.
By Corollary~\ref{cor:index} the energy of a Floer connecting 
trajectory is
$
E = \tfrac 12 \pi^2 \bigl( \eta_{f_0}(A^+) - \eta_{f_0}(A^-) 
- \dim H^1_{A^-} + \dim H^1_{A^+} + j \bigr),
$
where $j$ is the Fredholm index of the linearized operator.
There are finitely many such numbers $E\geq 0$ with $j\leq 7$.
We order them as
$$
0\leq E_0<E_1<\dots<E_\ell .
$$ 

\noindent{\bf Claim.} 
{\it Let $j\in\{0,\ldots,\ell-1\}$ and 
$(\gamma,f)\in\Gamma_m\times\cF_m$ such that
\begin{equation}\label{cond 1}
\left.\begin{array}{c}
\A\in\widetilde\cM(A^-,A^+;X_f), \;
([A^-],[A^+])\neq(0,0), \\
E_f(\A)\leq E_j, \;\delta_f(\A)\leq 7
\end{array}\right\}
\quad\implies\quad \cD_{\A,\delta} \;\text{is onto}
\end{equation}
\begin{equation}\label{cond 2}
\Crit(\CS_\cL+h_f)=\Crit(\CS_\cL+h_{f_0})
\end{equation}
\begin{equation}\label{cond 3}
A\in\Crit(\CS_\cL+h_{f_0})
\quad\implies\quad
h_f(A)=h_{f_0}(A),
\end{equation}
Fix an integer $k\in\N$ and a constant $\eps>0$.
Then there is a perturbation $(\gamma',f')\in\Gamma_{m'}\times\cF_{m'}$
satisfying (\ref{cond 1}) to (\ref{cond 3}) with
$j$ replaced by $j+1$ and }
\begin{equation}\label{cond 4}
\NORM{X_{f'}-X_{f}}_{k} < \eps.
\end{equation}

\noindent
A connection $\A\in\cM(A^-,A^+;X_f)$ with energy 
$E_f(\A)\leq 0$ must be gauge equivalent to the 
constant path $A^-=A^+\not\in[0]$.  By assumption 
these critical points of $\CS_\cL + h_{f_0}$ are 
nondegenerate.  So by Theorem~\ref{thm:iso} 
the hypotheses of the claim are satisfied for
$j=0$ and $(\gamma,f)=(\gamma_0,f_0)$.
Therefore assertion (ii) of the theorem follows 
from the claim by induction on $j$.
We prove the claim in four steps.

\medskip\noindent{\bf Step~1.}
{\it The quotient of the set
$$
\cK := \bigcup_{([A^-],[A^+])\neq(0,0)} 
\left\{ A:\R\to\cA(Y,\cL) \;\left|\;
\begin{array}{l}
\pd_s A + *(F_A+X_f(A))= 0,\\
\lim_{s\to\pm\infty}A(s)\in[A^\pm], \\
E_f(\A)\leq E_{j+1}, \; \delta_f(\A)\leq 7, \\
\cD_{\A,\delta} \;\text{not onto} 
\end{array}
\right.\right\}
$$
by the gauge group $\cG(Y)$ is compact.}

\medskip\noindent
This is proven by the same discussion as in (i).
The argument uses in addition the fact 
that the energy of each limit trajectory $\A_j$
is strictly less than the energy of the $\A^\nu$
if bubbling or breaking of trajectories occurs.
(So the relevant moduli spaces will be transverse
or empty by assumption.)

\medskip\noindent{\bf Step~2.}
{\it There is a $\gamma'\in\Gamma_{m'}$ 
with $\gamma_i'=\gamma_i$ for $i=1,\dots,m$
satisfying the following conditions.
For $z\in\D$ and $A\in\cA(Y)$ let $\rho_i'(z,A)$ be the 
holonomy of $A$ around the loop $\theta\mapsto\gamma_i'(\theta,z)$
and define $\sigma:\cA(Y)\to\rG^{m'}$ by 
$$
\sigma(A) := (\rho_1'(0,A),\dots,\rho'_{m'}(0,A)).
$$
Then, for every $\A\in\cK$, there is an $s_0\in\R$ such that
the following holds.
\begin{description}
\item[(a)]
The tuple $\sigma(A(s_0))$ is not contained in
$\sigma(\Crit(\CS_\cL+h_f))$ and belongs to the free part
of $\rG^{m'}$ for the action of $\rG$ by simultaneous conjugation.
Moreover, $\sigma(A(s))\not\sim \sigma(A(s_0))$ 
for every $s\in\R\setminus\{s_0\}$.
\item[(b)]
For every nonzero section $(\eta,0)\in\ker\cD_{\A,\delta}^*$ the vectors 
$\rd\sigma(A(s_0))\pd_s A(s_0)$ and $\rd\sigma(A(s_0))\eta(s_0)$ 
are linearly independent in $\rT(\rG^{m'}/\rG)$. 
\end{description}
}

\medskip\noindent
For every $s_0\in\R$ and every $\gamma'$ the set of all 
$\A\in\cK$ that satisfy conditions (a) and (b) is open.
Moreover, (a) and (b) are preserved under gauge transformations
and under adding further loops to $\gamma'$.
So it suffices to establish (a) and (b) for a single 
element of $\cK$. (Then $\cK$ is covered by finitely many
gauge orbits of small open sets around such elements,
and the final $\gamma'$ results from taking the union over
all loops that are required by these different elements.)
Hence from now on we fix an element $\A\in\cK$.
Since either $A^+$ or $A^-$ is irreducible, there
is an $s_0\in\R$ such that $A(s_0)$ is irreducible. 
Since the path $s\mapsto(\rd_{A(s)}\xi,0)$ is a solution of 
(\ref{eq:floer-lin}) for every $\xi\in\Om^1(Y,\cg)$, 
it follows from Proposition~\ref{prop:ucon}~(ii) below that
\begin{equation}\label{eq:ds0A}
\pd_s A(s_0) \notin \im \rd_{A(s_0)};
\end{equation}
otherwise we would have $\pd_s A(s)=\rd_{A(s)}\xi$ for all $s\in\R$
and, by partial integration,
$\|\rd_A\xi\|_{L^2(Y)}= -\int_Y \la \rd_A\xi\wedge(F_A+X_f(A))\ra = 0$
which would imply $\pd_s A\equiv 0$ and hence $E_f(\A)=0$.
By Proposition~\ref{prop:ucon}~(i) 
below, we have that
\begin{equation}\label{eq:Ainj}
A(s_0) \notin \bigcup_{s\neq s_0} [A(s)]
\cup \Crit(\CS_\cL+h_f);
\end{equation}
otherwise $A:\R\to\cA(Y,\cL)$ would be constant or periodic
modulo gauge, in contradiction to $0<E_f(\A)<\infty$.
Moreover, for $(\eta,0)\in\ker\cD_{\A,\delta}^*$, we have 
\begin{equation}\label{eq:etaperp}
\eta(s_0) \perp \R \pd_s A(s_0) + \im\rd_{A(s_0)}.
\end{equation}
To see this, fix an element $\xi\in\Om^0(Y,\cg)$.
Then $\alpha(s):=\pd_s A(s) + \rd_{A(s)}\xi$ and $\eta(s)$
satisfy the differential equations
$$
\p_s\alpha + *(\rd_A\alpha+\rd X_f(A)\alpha)=0,\qquad
\p_s\eta + 2\pd_s V \eta - *(\rd_A\eta+\rd X_f(A)\eta)=0.
$$
and the Lagrangian boundary condition 
$\eta(s)|_{\pd Y},\alpha(s)|_{\pd Y}\in T_{A(s)}\cL$.
Hence
\begin{align*}
\frac{\rd}{\rd s}
\exp(2V) \int_Y \inner{\eta}{\alpha} 
=
\exp(2V) \left( \int_Y \inner{\pd_s\eta+ 2\pd_s V \eta}{*\alpha}
+ \int_Y\inner{\eta}{*\pd_s\alpha} \right)
= 0 .
\end{align*}
The last identity uses the fact that the operator
$\alpha\mapsto *(\rd_A\alpha+\rd X_f(A)\alpha)$
with the Lagrangian boundary condition is self-adjoint 
for every $s$.  Since the inner product $e^{2V}\int_Y \inner{\eta}{\alpha}$
converges to zero for $s\to\pm\infty$, this proves~(\ref{eq:etaperp}). 

As in the proof of Theorem~\ref{thm:fcrit} we shall use 
Remark~\ref{rmk:rho} to construct $\gamma'$ and it suffices in each 
step to find the loop $\theta\mapsto\gamma_i'(\theta,0)$. 
Since $A(s_0)$ is irreducible and using (\ref{eq:ds0A})
we can argue exactly as in the proof 
of Step~1 in Theorem~\ref{thm:fcrit}, with $(A,\eta)$ replaced by 
$(A(s_0),\p_sA(s_0))$, to prove that $\gamma'$ can be chosen 
such that $\sigma(A(s_0))$ belongs to the free part of $\rG^{m'}$
and
\begin{equation}\label{eq:BLAH}
\rd\sigma(A(s_0))\pd_s A(s_0)\notin V_0,
\end{equation}
where $V_0\subset\rT_{\sigma(A(s_0))}\rG^{m'}$
is the tangent space of the $\rG$-orbit through $\sigma(A(s_0))$,
namely
$$
V_0 := \bigl\{v=\bigl( 
\sigma_i(A(s_0))\xi_0-\xi_0\sigma_i(A(s_0))
\bigr)_{i=1,\dots,m'}
\,\big|\, \xi_0\in\cg \bigr\}.
$$
This implies that $[\rd\sigma(A(s_0))\pd_s A(s_0)]\ne 0$ 
in the tangent space of the quotient $\rG^{m'}/\rG$.
It follows that the curve 
$[s_0-\delta,s_0+\delta]\to\rG^{m'}/\rG:s\mapsto[\sigma(A(s))]$
is injective for $\delta>0$ sufficiently small.
The set
$$
\cC:=\bigl\{ [A(s)] \st 
|s-s_0|\geq\delta \bigr\} \cup \Crit(\CS_\cL+h_f)/\cG(Y) 
\subset \cA(Y)/\cG(Y)
$$
is compact and, by (\ref{eq:Ainj}), does not contain $[A(s_0)]$.
Now (i) holds if and only if $\sigma(B)\not\sim\sigma(A(s_0))$
for every $[B]\in\cC$. Since this condition is open in $B$, and 
$\cC$ is compact, it suffices to prove this for a fixed element
$[B]\in\cC$.
Given $[B]\in\cC$ it follows from Remark~\ref{rmk:rho}~(i)
that for every $g\in\rG$ there is a based loop $\gamma$
such that
$$
\rho_\gamma(B)
\neq  g^{-1} \rho_\gamma(A(s_0)) g .
$$
For every fixed loop $\gamma$ this condition is open in $g$.
Since $\rG$ is compact there exist finitely many loops $\gamma'_i$
such that the tuple $(\rho_{\gamma'_i}(B))_{i}$ is not simultaneously
conjugate to $(\rho_{\gamma'_i}(A(s_0)))_{i}$.
For this choice of the loops $\gamma_i'$ we have that
${\sigma(B)\not\sim\sigma(A(s_0))}$ as claimed.

To prove~(b) it suffices to consider a fixed nonzero element
${(\eta,0)\in\ker\cD_{\A,\delta}^*}$ because this kernel 
is finite dimensional. Since $\eta(s_0)\neq 0$ (by unique 
continuation as in Proposition~\ref{prop:ucon}~(ii)) 
it follows from~(\ref{eq:etaperp}) that
$$
\eta(s_0) - \lambda\pd_s A (s_0) 
\not\in\im\rd_{A(s_0)}\qquad\forall\lambda\in\R.
$$
By (\ref{eq:BLAH}) we have
$\delta:=\inf_{v\in V_0} |\rd\sigma(A(s_0)\pd_s A(s_0) - v | >0$
and
\begin{equation}\label{eq:slambda}
\rd\sigma(A(s_0)) \bigl( \eta(s_0) - \lambda\pd_s A (s_0) \bigr)
\not\in V_0 
\end{equation}
for $\Abs{\lambda}>\delta^{-1}\Norm{\rd\sigma(A(s_0)\eta(s_0)}=:c$.
We wish prove that~(\ref{eq:slambda}) continues to hold
for all $\lambda\in[-c,c]$ with a suitable choice of $\gamma'$. 
For each fixed $\lambda$ the proof is the same as that 
of Step~1 in the proof of Theorem~\ref{thm:fcrit}.  
Since condition~(\ref{eq:slambda}) is open in
$\lambda$ this proves Step~2. 

\medskip\noindent{\bf Step~3.}
{\it Let 
$
C:=\bigl\{ (z,\rho'(z,A))  \in\D\times\rG^{m'}
 \st  A\in \Crit(\CS_\cL+h_f) \bigr\} .
$ 
For $\eps'>0$ and $k\in\N$ (possibly larger than
the constant in the claim) denote 
$$
\cF^{k,\eps'}_{m'} := \bigl\{ f'\in\cC^{k+1}(\D\times\rG^{m'})^\rG \st 
(f'-f)|_{B_{\eps'}(C)}\equiv 0 ,
\|f'-f\|_{\cC^{k+1}}<\eps'
\bigr\}
$$
and for a fixed $p>4$ let
$$
\widetilde\cM (A^-,A^+,\cF_{m'}^{k,\eps'})
:= \left\{ (\A,f')\in\cA^{1,p}_\delta\times\cF_{m'}^{k,\eps'}
\left|\begin{array}{l}
\A\in\widetilde\cM(A^-,A^+;X_{f'})\\
E_{f'}(\A)\leq E_{j+1} \\
\delta_{f'}(\A)\leq 7 
\end{array}\right.\right\}.
$$
Here we abbreviate 
$\cA^{1,p}_\delta:=\cA^{1,p}_\delta(\R\times Y,\cL;A^-,A^+)$
(see equation~(\ref{eq:Akpdelta})). 
Let $\cG_0^{2,p}(\R\times Y)$ be the $W^{2,p}$-closure of
$\{u:\R\to\cG(Y)\st u(s)=\one \;\;\forall |s|\geq 1 \}$.
Then for every $k\in\N$ there is an $\eps'>0$ such that
the following holds.

Every perturbation $f'\in\cF_{m'}^{k,\eps'}$ 
satisfies conditions (\ref{cond 2}), (\ref{cond 3}), (\ref{cond 4}),
and for every pair of critical points ${([A^-],[A^+])\neq(0,0)}$
the universal moduli space
$$
\cM(A^-,A^+,\cF_{m'}^{k,\eps'})
:=\widetilde\cM (A^-,A^+,\cF_{m'}^{k,\eps'})/\cG_0^{2,p}(\R\times Y)
$$
is a separable $\cC^k$-Banach manifold.}

\medskip\noindent
Conditions (\ref{cond 2}), (\ref{cond 3}), and (\ref{cond 4})
are satisfied for every $f'\in\cF^{k,\eps'}_{m'}$
for  $\eps'>0$ sufficiently small.
The assertion about the universal moduli space
holds whenever the linearized operator
\begin{equation}\label{eq:univop}
(\alpha,\phi,\hat f)\mapsto\cD_{\A,\delta}(\alpha,\phi) 
+ (X_{\hat f}(\A),0)
\end{equation}
is surjective for every pair 
$(\A,f')\in\widetilde\cM (A^-,A^+,\cF_{m'}^{k,\eps'})$.
Here $\cD_{\A,\delta}$ is the operator (\ref{eq:DAd}) with $k=1$.
We first prove that this holds for $f'=f$.
If $\A$ is not gauge equivalent (by $\cG(A^-,A^+)$)
to a connection in $\cK$, 
then the operator $\cD_{\A,\delta}$ is surjective 
by Remark~\ref{rmk:halfonto}~(i),
and hence so is (\ref{eq:univop}).
Let $\A\in\cK$ (after a gauge transformation in $\cG(A^-,A^+)$) 
and $q^{-1}:=1-p^{-1}$, and suppose, 
by contradiction, that there is a nonzero pair
$$
(\eta,\phi)\in L^q_\delta(\R\times Y,\rT^*Y\otimes\cg) 
\times L^q_\delta(\R\times Y,\cg)
$$
orthogonal to the image of (\ref{eq:univop}).
Then we have $\phi=0$ (by the proof of 
Theorem~\ref{thm:weighted fredholm}), 
${\eta\in W^{1,p}_\A(\R\times Y,\rT^*Y\otimes\cg)}$
(by Theorem~\ref{thm:4reg}), $\cD_{\A,\delta}^*(\eta,0)=0$, 
and
\begin{equation}\label{eq:eta}
\int_{-\infty}^\infty \exp(2V(s)) \rd h_{\hat f}(A(s))\eta(s)  \ds 
= 0 
\end{equation}
for every $\hat f\in\rT_f\cF^{k,\eps'}_{m'}$.
By Step~2 there is $s_0\in\R$ such that 
${\sigma(A(s))\ne \sigma(A(s_0))}$ for $s\ne s_0$ and
the tangent vectors $\rd\sigma(A(s_0))\pd_sA(s_0)$,
$\rd\sigma(A(s_0))\eta(s_0)$ are linearly independent.
Hence the map
$$
(r,s)\mapsto \rho(z,A(s)+r\eta(s))
$$
is an embedding in a neighbourhood of $(0,s_0)\in\R^2$ 
for every sufficiently small $z\in\D$.  It follows that
there exists a smooth $\rG$-invariant map 
$\hat f:\D\times\rG^{m'}\to\R$ vanishing in a neighbourhood
of $C$ and satisfying 
$$
\hat f(z,\rho(z,A(s)+r\eta(s))) = r\beta(r)\beta(s-s_0)\beta(|z|)
$$
for a suitable cutoff function $\beta:\R\to[0,1]$
that is supported in a neighbourhood of $0$ 
and is equal to $1$ near $0$. This implies
\begin{align*}
dh_{\hat f}(A(s))\eta(s) 
&= \int_\D\left.\frac{\pd}{\pd r}\right|_{r=0}
   \hat f(z,\rho(z,A(s)+r\eta(s)))\rd^2z  \\
&= \beta(s-s_0)\int_\D\beta(|z|)\rd^2z \;\ge\; 0
\end{align*}
for every $s\in\R$.  Hence the integral on the right hand side
of~(\ref{eq:eta}) does not vanish, contradiction. 
Thus we have proved that the operator~(\ref{eq:univop}) 
is onto whenever $f'=f$.  

We must prove that~(\ref{eq:univop}) is onto
when $\|f'-f\|_{\cC^{k+1}}$ is sufficiently small.
Otherwise there are sequences
$\cF^{k,\eps'}_{m'}\ni f^\nu\to f$
and ${\A^\nu\in\widetilde\cM(A^-,A^+;X_{f^\nu})}$
such that the operator~(\ref{eq:univop}),
with $(\A,f')$ replaced by $(\A^\nu,f^\nu)$, is not onto.
If $\A^\nu$ converges (modulo gauge) to $\A\in\cK$ then
(\ref{eq:univop}) is surjective for the pair $(\A,f)$
and hence for $(A^\nu,f^\nu)$ when $\nu$ is
sufficiently large.
Otherwise it follows from the compactness and gluing 
theorems as in the proof of (i) 
that $\cD_{\A^\nu,\delta}$ is surjective
for $\nu$ sufficiently large.
This contradiction finishes the proof of Step~3.

\medskip\noindent{\bf Step 4.}
{\it We prove the claim.}

\medskip\noindent
By Step~3 the projection
$\cM (A^-,A^+,\cF_{m'}^{k,\eps'})\to \cF_{m'}^{k,\eps'}$
is a Fredholm map of index at most $7$
for every pair $A^\pm\in\Crit(\CS_\cL+h_f)$ 
with $([A^-],[A^+])\neq(0,0)$.
(The index at $(\A,f)$ is the same as that of the linearized
operator $\cD_{\A,\delta}$.)
Hence it follows from the Sard--Smale theorem that, 
for $k\geq 8$, the set of regular values is of the
second category in the sense of Baire.
Any such regular value $f\in\cF_{m'}^{k,\eps'}$ 
satisfies~(\ref{cond 1}). To prove the claim, pick a 
regular value of the projection and approximate it 
by a smooth perturbation $f'$. 
In the last step we use the fact that the
set of all perturbations that satisfy the 
requirements of the claim is open 
in the $\cC^{k+1}$-topology. 
(The proof is analogous to the proof of (i).)
This proves the theorem.
\end{proof}

\label{page:bubble}
The main difference between our proof of 
Theorem~\ref{thm:freg} and the argument in 
Donaldson's book~\cite[p 144]{Donaldson book}
for the closed case is that we do not have a gluing 
theorem converse to bubbling on the boundary 
and hence cannot work on a compact part of
the moduli space in the presence of bubbling
on the boundary.  To circumvent this difficulty we have 
restricted the discussion to the monotone case and
to Floer connecting trajectories of index less than or equal 
to seven. We also made use of a unique continuation result
for perturbed anti-self-dual connections with 
Lagrangian boundary conditions, which is established next.

\subsection*{Unique Continuation}

\begin{prp}\label{prop:ucon}
Let $(\gamma,f)\in\Gamma_m\times\cF_m$ and fix an
open interval $I\subset\R$. 

\smallskip\noindent{\bf (i)}
Let $A,B:I\to\cA(Y)$ be two solutions
of the Floer equation
\begin{equation}\label{eq:floer1}
\pd_s A + *F_A + *X_f(A) = 0, \qquad
A(s)|_\Sigma \in\cL.
\end{equation}
If $A(s_0)=B(s_0)$ for some $s_0\in I$
then $A(s)=B(s)$ for all $s\in I$. 

\smallskip\noindent{\bf (ii)}
Let $A:I\to\cA(Y,\cL)$ and $\xi=(\alpha,\phi):I\to\Om^1(Y,\cg)\times\Om^0(Y,\cg)$
be smooth maps satisfying the (augmented) linearized Floer equation
\begin{equation}\label{eq:floer-lin}
\pd_s\xi + \cH_A\xi = 0 , \qquad
\alpha(s)|_\Sigma \in \rT_{A(s)}\cL, \quad *\alpha(s)|_\Sigma=0 .
\end{equation}
If $\xi(s_0)=0$ for some $s_0\in I$
then $\xi(s)=0$ for all $s\in I$. 
\end{prp}

The proof will use the following local continuation result
in the interior. This was proven by Taubes \cite{T:ucon}
in a slightly different formulation; we include the proof
for the sake of completeness.

\begin{lemma}\label{le:ucon-loc}
Let $U$ be a (not necessarily compact) $3$-manifold without boundary
and $I\subset\R$ be an open interval. 

\smallskip\noindent{\bf (i)}
Let $A,B:I\to\cA(U)$ be two solutions of the unperturbed 
Floer equation~(\ref{eq:floer1}) with $f=0$.
If $A(s_0)=B(s_0)$ for some $s_0\in I$
then $A(s)=B(s)$ for all $s\in I$. 

\smallskip\noindent{\bf (ii)}
Let $A:I\to\cA(U)$ and $\xi=(\alpha,\phi):I\to\Om^1(U,\cg)\times\Om^0(U,\cg)$
be smooth maps satisfying the unperturbed linearized 
Floer equation~(\ref{eq:floer-lin}) with $f=0$.
If $\xi(s_0)=0$ for some $s_0\in I$ then $\xi(s)=0$ for all $s\in I$. 
\end{lemma}

\begin{proof}
To prove~(i) assume by contradiction that 
$A(s',y_1)\neq B(s',y_1)$ 
for some $(s',y_1)\in I\times U$.
Let $D_r(y_1)\subset U$ be a geodesic ball
of radius $r>0$ around $y_1$ and denote
$$
J:=\bigl\{s\in I \st A(s)|_{D_{r/2}(y_1)}
=B(s)|_{D_{r/2}(y_1)} \bigr\} \subset I .
$$
This set contains $s_0$ by assumption and it 
is a closed subset of $I$ because $A-B$ is continuous.
We claim that $J\subset I$ is open and hence 
$J=I$ in contradiction to the assumption.

To prove that $J$ is open we fix an element $s_1\in J$.
Then $A-B$ vanishes to infinite order 
(i.e.~with all derivatives) at $x_1:=(s_1,y_1)$.
For the derivatives in the direction of $I$ 
this follows from the Floer equation.
Let $D_r(x_1)\subset I\times U$ denote the geodesic ball 
centred at $x_1$.
We fix gauge transformations $u_A,u_B\in\cG(D_r(x_1))$ 
with $u_A(x_1)=u_B(x_1)=\one$ such that
$u_A^*A$ and $u_B^*B$ are in radial gauge on $D_r(x_1)$.
Then these can be pulled back to connections in temporal
gauge $A',B': (-\infty,\log r)\to \cA(S^3)$
by geodesic polar coordinates
$(-\infty,\log r)\times S^3 \overset{\sim}{\rightarrow} 
D_r(x_1)\setminus\{x_1\}$.
The fact that $u_A^*A-u_B^*B$ vanishes to infinite order at $x_1$
translates into superexponential convergence
${A'(s)-B'(s)\to 0}$ as $s\to-\infty$.
In particular, for every $K>0$, we have
\begin{equation}\label{eq:superexp}
\lim_{s\to-\infty}e^{-Ks} \|A'(s)-B'(s)\|_{L^2(S^3)} = 0 .
\end{equation}
The pullback metric on $(-\infty,\log r)\times S^3$ has the form
$e^{2s}(\ds^2 + g_s )$, where $g$ is a smooth family of metrics on $S^3$
that converges exponentially to the standard metric on $S^3$ as $s\to-\infty$.
Since the anti-self-duality equation is conformally invariant, 
the connections $A'$ and $B'$ also satisfy~(\ref{eq:floer1}) with respect to
the metric $\ds^2 + g_s$ on $(-\infty,\log r)\times S^3$.
We now denote 
$
\alpha:=B'-A':(-\infty,\log r)\to\Om^1(S^3,\cg)
$
and use the technique of Agmon--Nirenberg in Appendix~\ref{app:uc}
to prove that $\alpha\equiv 0$.
The Floer equations (i.e.~the anti-self-duality 
of $A'$ and $B'$ w.r.t.\ the conformally rescaled metric) 
imply that $\alpha$ satisfies
$$
\pd_s\alpha + *\rd_{A'+\frac 12 \alpha} \alpha = 0.
$$
We shall use the operator 
$\bF:=-*\rd_{A'+\frac 12\alpha}$ 
(corresponding to $A(s)$, appropriately shifted, 
in the notation of Appendix~\ref{app:uc}) 
which is self--adjoint with respect
to the time dependent inner product 
$$
\Inner{\alpha}{\beta}_s 
:= \int_{S^3} \winner{\alpha}{*_s \beta}
= \Inner{Q(s)\alpha}{Q(s)\beta}_E .
$$ 
Here $*_s$ is the Hodge operator for the metric $g_s$ on $S^3$, 
and the subscript ${\scriptstyle E}$ indicates the 
use of the standard metric on $S^3$. 
The operator $Q(s):\Om^1(S^3,\cg)\to\Om^1(S^3,\cg)$ is defined
as in~\cite[p.151]{DK}, as a self--adjoint operator 
such that $Q(s)^2=*_E*_s$.
This square root exists since $*_E*_s$ is 
positive definite.
These operators satisfy~$(Q1)$ in Appendix~\ref{app:uc} 
by the exponential convergence of $g_s$ as $s\to-\infty$. 
Moreover,
\begin{align*}
- \frac\rd{\rd s} \langle{\alpha} , {*\rd_{A'+\frac 12\alpha}\alpha}\rangle_s
+ 2 \langle{\pd_s\alpha} , {*\rd_{A'+\frac 12\alpha}\alpha}\rangle_s 
&= - \int_{S^3} 
\winner{\alpha}{[\pd_s (A'+\tfrac 12 \alpha)\wedge\alpha]} \\
&\leq  \Norm{\pd_s A'+\tfrac 12\pd_s\alpha}_{L^\infty(S^3)} 
\Norm{\alpha}_s^2 .
\end{align*}
Hence the function $x(s):=\alpha(s_2-s)$, 
with $s_2\in(-\infty,\log r)$, satisfies the assumptions 
of Theorem~\ref{thm:agni2} with ${c_1=c_2=0}$ and 
$c_3(s)=\|\pd_sA' + \frac 12\pd_s\alpha\|_{L^\infty(S^3)}$.
The constant $c$ in Theorem~\ref{thm:agni2} 
is finite because
$
\int_{-\infty}^{s_2}
\|\pd_sA'+\tfrac 12\pd_s\alpha\|_{L^\infty(S^3)}<\infty,
$
by the exponential decay of $A'$ and $B'$ 
(see Theorem~\ref{thm:decay}). We thus obtain
$$
{\Norm{\alpha(s)}_s \ge e^{-c(s_2-s)}\Norm{\alpha(s_2)}}_{s_2}
$$
for all $s\in (-\infty,s_2]$. This estimate contradicts 
the superexponential convergence in~(\ref{eq:superexp})
unless $\alpha(s_2)=0$.  Since $s_2$ is any element
of the interval $(-\infty,\log r)$ we have shown that
$\alpha\equiv 0$ and hence ${u_A^*A=u_B^*B}$ on the
geodesic ball $D_r(x_1)$ around $x_1=(s_1,y_1)$.
This ball contains the set 
$[s_1-\frac r2,s_1+\frac r2]\times D_{r/2}(y_1)$.
From the construction of the gauge transformations 
with $A=B$ on $\{s_1\}\times D_{r/2}(y_1)$ we know that
$
u_A|_{s=s_1}=u_B|_{s=s_1}.
$
Now there is a unique gauge transformation $v$
on $[s_1-\frac r2,s_1+\frac r2]\times D_{r/2}(y_1)$
with 
$
v|_{s=s_1}=u_A^{-1}|_{s=s_1}=u_B^{-1}|_{s=s_1}
$ 
that puts $u_A^*A=u_B^*B$ back into temporal gauge.
By the uniqueness of the temporal gauge with
$u_A v|_{s=s_1}=u_B v|_{s=s_1}=\one$ this implies
$$
A=(u_A v)^*A=(u_B v)^*B=B
\qquad\text{on}\quad [s_1-\tfrac r2,s_1+\tfrac r2]\times D_{r/2}(y_1)
$$
and hence $[s_1-\frac r2,s_1+\frac r2]\subset J$.
This proves that $J$ is open as claimed.

The proof of (ii) is analogous to (i).
In conformal polar coordinates near $x_1$ we choose 
the radial gauge $u_A^*A$ as before.
The pullback $\xi':(-\infty,\log r)\to\Om^1(S^3,\cg)\times\Om^0(S^3,\cg)$
then satisfies the linearized Floer equation with respect to $A'$.
Now the Agmon-Nirenberg technique for $x=\xi'$ 
(with the Hessian $\cH_{A'(s)}$ as self-adjoint operator)
shows that $\xi'\equiv 0$ and hence $\xi=0$ on $D_r(x_1)$.
The relevant estimate is
\begin{align*}
& - \frac\rd{\rd s} \Inner{\xi'}{\cH_{A'}\xi'}_s
+ 2 \Inner{\pd_s\xi'}{\cH_{A'}\xi'}_s \\
&= - \frac\rd{\rd s}\biggl( \int_{S^3} \winner{\alpha'}{\rd_{A'}\alpha'} 
 -2\int_{S^3} \winner{\alpha'}{*\rd_{A'}\phi'} \biggr) 
 + 2 \Inner{\pd_s(\alpha',\phi')}{\cH_{A'}(\alpha',\phi')}_s \\
&=  - \int_{S^3} \winner{\alpha'}{[\pd_s A',\alpha']} 
 + 2\int_{S^3} \winner{\alpha'}{*[\pd_s A',\phi']} 
+ 2\int_{S^3} \winner{\alpha'}{(\pd_s*)\rd_{A'}\phi'} \\
& \leq 2c_2(s)\Norm{\cH_{A'}\xi'}_s \Norm{\xi'}_s + c_3(s)\Norm{\xi'}_s^2 ,
\end{align*}
where $c_2(s)=2\delta^{-1}c_Q(s)$ and
$c_3(s)=2\Norm{\pd_s A'}_{L^\infty(S^3)} 
+ 8\delta^{-1}c_Q(s)\Norm{F_{A'}}_{L^\infty(S^3)}^{1/2}$
with $\delta$ and $c_Q$ as in (Q1) in Appendix~\ref{app:uc}.
We have used the identity $\pd_s*=*_{E}\pd_s Q^2$,
which implies $\|\pd_s*\|_s\leq 2\delta^{-3} c_Q$,
and 
$$
\Norm{\rd_{A'}\phi'}_s^2 + \Norm{\rd_{A'}\alpha'}_s^2
= \Norm{*\rd_{A'}\alpha' - \rd_{A'}\phi'}_s^2
+ 2\tint_{S^3} \winner {\alpha'}{[F_{A'},\phi']} 
$$
which implies
$\Norm{\rd_{A'}\phi'}_s \leq  \Norm{\cH_{A'}\xi'}_s
+ 2 \Norm{F_{A'}}^{1/2}_{L^\infty(S^3)} \Norm{\xi'}_s$. 
\end{proof}

\begin{proof}[Proof of Proposition~\ref{prop:ucon}.]
The proof of~(i) is similar to that of Lemma~\ref{le:ucon-loc}
except for the presence of boundary terms.
To control these we first use Lemma~\ref{le:ucon-loc}~(i)
on $U:=N\setminus\pd Y$ for a neighbourhood $N\subset Y$ 
of $\pd Y$ on which $X_f\equiv 0$.
It implies that $A$ and $B$ agree on $I\times U$
and hence by continuity on $I\times N$. 
In particular, the $1$-form
$\alpha(s) := B(s) - A(s) \in\Om^1(Y,\cg)$
vanishes near $\pd Y$ and hence 
belongs to the space $\Om^1_{A(s)}(Y,\cg)$ for every~$s$.
To establish unique continuation in the interior we 
assume, by contradiction, that $\alpha(s_1)\ne 0$ for some ${s_1<s_0}$.
We will apply Theorem~\ref{thm:agni1} to $x(s)=\alpha(s_1-s)$ and
the symmetric operator
$$
\bF(s) := *\rd_{A(s)} + *\rd X_f(A(s)) \,:\;
\Om^1_{A(s)}(Y,\cg)\to\Om^1_{A(s)}(Y,\cg)
$$
for $s\in I$.  We have $\alpha(s_0)=0$ and
\begin{align*}
& \p_s\alpha+\bF\alpha 
= -\tfrac12 *[\alpha\wedge\alpha] 
- *\bigl(X_f(A+\alpha)-X_f(A)-\rd X_f(A)\alpha\bigr), \\
&
\pd_s \inner{\alpha}{\bF\alpha}
- 2 \inner{\pd_s\alpha}{\bF\alpha} 
= \tint_Y \winner{\alpha}{[\p_s A\wedge\alpha]}
+\tint_Y\inner{\alpha}{\rd^2X_f(A)(\p_sA,\alpha)}.
\end{align*}
Hence it follows from Proposition~\ref{prop:Xf}~(v)
that
\begin{align*}
\Norm{\p_s\alpha(s)+\bF(s)\alpha(s)}_{L^2(Y)}
&\le c_1\Norm{\alpha(s)}_{L^2(Y)}, \\
\pd_s\inner{\alpha(s)}{\bF(s)\alpha(s)}
- 2 \inner{\pd_s\alpha(s)}{\bF(s)\alpha(s)} 
&\le c_3 \Norm{\alpha(s)}_{L^2(Y)}^2
\end{align*}
for $s_1\le s\le s_0$ and suitable constants 
$c_1$ and $c_3$. This shows that the path 
$s\mapsto\alpha(s)$ and the operator family $\bF(s)$ 
satisfy the hypotheses of Theorem~\ref{thm:agni1} with $c_2=0$. 
Hence $\alpha(s)=0$ for $s_1 < s\le s_0$ and $\alpha(s_1)=0$ follows by continuity,
in contradiction to the assumtion.
The argument for $s_1>s_0$ is simlar and this proves~(i).
Assertion~(ii) follows from Lemma~\ref{le:ucon-loc}~(ii)
and the analogous estimates for the solutions of~(\ref{eq:floer1}).  
This proves the proposition.
\end{proof}


\section{Gluing}\label{sec:gluing}

Let $Y$ be a compact oriented Riemannian $3$-manifold 
with boundary $\pd Y=\Sigma$ and $\cL\subset\cA(\Sigma)$ 
be a gauge invariant, monotone, irreducible Lagrangian submanifold 
satisfying (L1-3) on page~\pageref{p:L1}.
Fix a regular perturbation $h_f:\cA(Y)\to\R$ 
in the sense of Definition~\ref{def:freg}.

Let $B_0,B_1,B_2\in\cA(Y,\cL)$ be nondegenerate and
irreducible critical points of $\CS_\cL+h_f$.
We denote by $\cA(\R\times Y,\cL;B_0,B_2)$ the space 
of smooth connections on $\R\times Y$ with boundary 
values in $\cL$ and $\cC^\infty$-limits 
$B_0$ and $B_2$ as in~(\ref{eq:limX});
this is a special case of the notation~(\ref{eq:Akpdelta}).
Also recall the notation $\widetilde\cM(B_0,B_1;X_f)$
from chapter~\ref{sec:compact} 
for the space of solutions that are in temporal gauge over the ends,
and $\cM(B_0,B_1;X_f)$ for this space modulo gauge equivalence.
For $T>1$ we define a pregluing map
\begin{equation}\label{eq:preg}
\begin{split}
\widetilde\cM(B_0,B_1;X_f) \times 
\widetilde\cM(B_1,B_2;X_f)
&\to \cA(\R\times Y,\cL;B_0,B_2) \\
(\Xi_1,\Xi_2) &\mapsto \Xi_1\#_T\Xi_2
\end{split}
\end{equation}
as follows.  The connections $\Xi_i=A_i+\Phi_i\ds$ 
are in temporal gauge outside the compact set $[-1,1]\times Y$
and have limits 
$$
\lim_{s\to-\infty}A_1(s)=B_0,\qquad
\lim_{s\to\infty}A_1(s)=B_1=\lim_{s\to-\infty}A_2(s),\qquad
\lim_{s\to\infty}A_2(s)=B_2 .
$$
Define $\Xi_1\#_T\Xi_2:=A+\Phi\ds$ by
$$
\Phi(s) := \left\{\begin{array}{ll}
\Phi_1(s+T),&s\le 0,\\
\Phi_2(s-T),&s\ge 0,
\end{array}\right.
$$
$$
A(s) := \left\{\begin{array}{ll}
A_1(\tfrac T2-\phi(-\tfrac T2-s)),& s< -\tfrac T2,\\
B_1,& s\in[-\tfrac T2, \tfrac T2] , \\
A_2(-\tfrac T2+\phi(-\tfrac T2+s)),& s > \tfrac T2 ,
\end{array}\right.
$$
where $\phi:(0,\infty)\to\R$ is a smooth function
satisfying
$$
\phi(s)=\left\{\begin{array}{ll}
s,&s\ge 2,\\
- \tfrac 1s,&s\le \tfrac 12,
\end{array}\right.\qquad
\pd_s\phi>0.
$$
This connection is smooth because $A_1$ and $A_2$ converge 
exponentially as $s$ tends to $\pm\infty$. 
It satisfies the limit conditions and the 
Lagrangian boundary conditions by construction. 
In fact, this is why we use rescaling in time rather 
than convex interpolation in space.  The map
$(\Xi_1,\Xi_2)\mapsto\Xi_1\#_T\Xi_2$ 
is gauge equivariant in the sense that
$$
(u_1^*\Xi_1)\#_T(u_2^*\Xi_2)=u^*(\Xi_1\#_T\Xi_2),
\qquad
u(s) := \left\{\begin{array}{ll}
u_1(s+T),&s\le 0,\\
u_2(s-T),&s\ge 0 
\end{array}\right.
$$
for each pair 
$(u_1,u_2)\in\cG(B_0,B_1)\times\cG(B_1,B_2)$. 
Recall from the beginning of Section~\ref{sec:compact}
that each $u_1\in\cG(B_0,B_1)$ satisfies
$\p_su_1(s)=0$ for $|s|\ge1$, $u_1(s)\in\cG_{B_1}$
for $s\ge 1$, and $u_1(s)\in\cG_{B_0}$ for $s\le -1$;
similarly for $u_2$. Since $B_1$ is irreducible 
we have $u_1(s)=u_2(-s)=\one$ for $s\ge 1$.

\begin{thm} \label{thm:glue}
Let $B_0,B_1,B_2\in\cA(Y,\cL)$ be nondegenerate and
irreducible critical points of $\CS_\cL+h_f$,
and fix $\Xi_1\in\widetilde\cM(B_0,B_1;X_f)$ and 
$\Xi_2\in\widetilde\cM(B_1,B_2;X_f)$ 
with $\delta_f(\Xi_1)=\delta_f(\Xi_2)=1$.
Then, for every $p>2$,
there exist positive constants $\kappa$, $T_0$ and a map
$$
\tau :
(T_0,\infty) \to \cM^2(B_0,B_2;X_f)/\R , \qquad
T \mapsto \tau_T(\Xi_1,\Xi_2) 
$$
with the following properties:
\begin{description}
\item[(i)]
$\tau$ is a diffeomorphism onto its image.
\item[(ii)]
The connections $\tau_T(\Xi_1,\Xi_2)$ converge without bubbling 
(as in Theorem~\ref{thm:compact1}) to the broken trajectory 
$(\Xi_1,\Xi_2)$  as $T\to\infty$.
\item[(iii)]
If $\Xi$ is a solution of the Floer equation (\ref{eq:floerPhi}) and
$$
\bigl\| \Xi - (\Xi_1 \#_T \Xi_2) \big\|_{W^{1,p}(\R\times Y)} \leq \kappa
$$
for some $T\ge T_0+1$, then its gauge and time-shift equivalence class 
$[\Xi]$ lies in the image of $\tau$.
\end{description}
\end{thm}

\begin{proof}
The preglued connection 
$$
\Xi_1\#_T\Xi_2=:\Xi_T=A_T+\Phi_T\ds
$$
is an approximate solution of the Floer equation 
and $\tau_T(\Xi_1,\Xi_2)$ will
be constructed as a nearby true solution.
More precisely, we have
\begin{equation}\label{eq:approx solution}
\bigl\| \pd_s A_T - \rd_{A_T}\Phi_T 
+ *\bigl(F_{A_T} + X_f(A_T) \bigr) \bigr\|_{L^p(\R\times Y)}
\leq C e^{-\delta T}
\end{equation}
for some constants $C$ and $\delta>0$ by exponential decay, Theorem~\ref{thm:decay}.
We will use the inverse function theorem to find near
the approximate solution $\Xi_T$ a true solution 
$\Tilde{\Xi}_T\in\widetilde\cM(B_0,B_2;X_f)$.
For that purpose we use the Banach manifold structure
of the space $\cA^{1,p}(\R\times Y,\cL;B_0,B_2)$, see (\ref{eq:Akpdelta}).
Its tangent space $\rT_{\Xi_T}\cA^{1,p}(\R\times Y,\cL;B_0,B_2)$
is the space of all $1$-forms $\xi=\alpha+\phi\ds$ with
$\alpha\in W^{1,p}(\R\times Y,\rT^*Y\otimes\cg)$
and $\phi\in W^{1,p}(\R\times Y,\cg)$ satisfying
the boundary condition $\alpha(s)\in \rT_{A_T(s)}\cL$.
Using the exponential map of Theorems~\ref{thm:expmap} and Corollary~\ref{cor:expmap4} 
we obtain a continuously differentiable map
$$
\rT_{\Xi_T}\cA^{1,p}(\R\times Y,\cL;B_0,B_2) \supset\tilde{\cU}
\to \cA^{1,p}(\R\times Y,\cL;B_0,B_2):
\xi\mapsto \tilde E(\Xi_T;\xi)
$$
defined on a neighbourhood $\tilde{\cU}$ of zero by
$$
\tilde E(\Xi_T;\xi):=E_{A_T(s)}(\alpha(s)) + (\Phi_T(s)+\phi(s))\ds .
$$
We now look for a solution of the form
$\Tilde{\Xi}_T=\Tilde{A}_T+\Tilde{\Phi}_T\ds = \tilde E(\Xi_T;\xi)$,
where $\xi \in \tilde{\cU}$
satisfies
\footnote{
Here $*$ denotes the Hodge $*$ operator on the four-manifold 
$\R\times Y$ unlike in (\ref{eq:floer-true}) below.
The first two conditions fix the gauge whereas the third
condition fixes a complement of the kernel of the linearized
operator for combined anti-self-duality and gauge fixing.
} 
\begin{equation}\label{eq:Xitilde}
\rd_{\Xi_T}^*\xi=0, \qquad *\xi|_{\R\times\pd Y}=0  , 
\qquad \xi\in \im\cD_T^* .
\end{equation}
Note that $\tilde{\Xi}_T$ automatically satisfies the boundary conditions
$\Tilde{A}_T(s)|_{\pd Y} \in \cL$ and has the limits
$\lim_{s\to-\infty}\Tilde{A}_T(s)=B_0$,
$\lim_{s\to\infty}\Tilde{A}_T(s)=B_2$.
So it remains to solve the Floer equation
\begin{equation}\label{eq:floer-true}
\pd_s\Tilde{A}_T-\rd_{\Tilde{A}_T}\Tilde{\Phi}_T
+ *\bigl(F_{\Tilde{A}_T}+X_f(\Tilde{A}_T) \bigr) = 0
\end{equation}
for $\xi$ subject to (\ref{eq:Xitilde}).
The precise setup for the inverse function theorem is as follows:
In order to keep track of the $T$-dependence we use the version
\cite[Proposition~A.3.4.]{MS} which provides explicit constants.
We apply this version of the inverse function theorem to the
$\cC^1$-map 
$$
f_T:X_T\to Z , \qquad
f_T(\xi):= (F_{\tilde E(\Xi_T;\xi)}^+,\rd^*_{\Xi_T}\xi) .
$$
Its domain is a neighbourhood of zero in the Banach space $X_T$ consisting of 
$\xi \in \rT_{\Xi_T}\cA^{1,p}(\R\times Y,\cL;B_0,B_2)$
that satisfy the boundary condition $*\xi|_{\R\times\pd Y}=0$.
(Note that the domain depends on $T$.
One could also work with a $T$-independent domain 
by using simple reparametrizations in $s\in\R$ to identify $X_T\cong X_{T_0}$ 
for a fixed $T_0$.
This gives rise to a continuous family of
inverse function problems $\tilde f_T:X_{T_0}\to Z$ for $T\in[T_0,\infty)$.)
The first component, $F_{\tilde E(\Xi_T;\xi)}^+$, is identified with the left hand
side of (\ref{eq:floer-true}), so the target space of $f_T$ is the Banach space
$$
Z=L^p(\R\times Y,\rT^*Y\otimes\cg)\times L^p(\R\times Y,\cg) .
$$
The differential $\rd f_T(0)$ at $x_0=0$ then is the 
linearized operator $\cD_T:=\cD_{\Xi_T}$.
To check that the differential $\rd f_T$ is uniformly continuous at $0\in X_T$ 
we calculate for all $\xi,\zeta\in X_T$ 
\begin{align}\label{eq:mightneedthis}
& \bigl\| \bigl(\rd f_T (\xi) - \cD_T \bigr)\zeta \bigr\|_{L^p(\R\times Y)}
= \bigl\| * \bigl[ \bigl(\tilde E(\Xi_T;\xi) - \Xi_T \bigr) \wedge * \zeta \bigr] \bigr\|_{L^p(\R\times Y)}
\nonumber\\
& \leq C \sup_{s\in\R} \bigl\| E_{A_T(s)}(\alpha(s)) - A_T(s) + \phi(s)\rd s \bigr\|_{L^{2p}(Y)} 
\| \zeta \|_{W^{1,p}(\R\times Y)} .
\end{align}
Here $C$ is the constant from the Sobolev embedding $W^{1,p}(Y)\hookrightarrow L^{2p}(Y)$ and
the second factor converges to zero uniformly in $T$ as 
$\|\xi\|_{W^{1,p}}=\|\alpha+\phi\ds\|_{W^{1,p}} \to 0$.
Indeed, given $\eps>0$ there is $\delta_{T,s}>0$ such that 
$\| E_{A_T(s)}(\alpha) - A_T(s) \|_{L^{2p}(Y)}\leq \eps$ for all 
$\alpha\in\rT_{A_T(s)}\cA(Y,\cL)$ with $\|\alpha\|_{L^{2p}(Y)}\leq\delta_{T,s}$.
We can choose $\delta_{T,s}=\delta>0$ uniform for all $T>1$, $s\in\R$ 
because the image of $A_T$ in $\cA(Y,\cL)$ is compact and independent of $T$.

That the linearized operator is surjective for sufficiently large $T$ 
with a uniform bound for its right inverse
$Q_T:=\cD_T^*(\cD_T\cD_T^*)^{-1}$ follows from the estimates
\begin{align}
\|\eta\|_{W^{1,p}(\R\times Y)} 
&\leq C \|\cD_T^*\eta\|_{L^p(\R\times Y)}, \label{eq:DT1} \\
\|\cD_T^*\eta\|_{W^{1,p}(\R\times Y)} 
&\leq C \|\cD_T\cD_T^*\eta\|_{L^p(\R\times Y)}. \label{eq:DT2} 
\end{align}
These estimates hold for $T$ sufficiently large, 
and the constant $C$ is independent of $T$.
The inequality~(\ref{eq:DT1}) implies that $\cD_T$ 
is surjective and $Q_T:Y\to X_T$ is defined,
and~(\ref{eq:DT2}) gives a uniform bound for $Q_T$.
The proof of the estimates is as 
in~\cite[Proposition~3.9]{Donaldson book},
\cite[Proposition~3.9]{Sal}, or Theorem~\ref{thm:compact2}. 
It rests on the fact that the connections 
$\Xi_{1,T}:=\Xi_1\#_T B_1$ and $\Xi_{2,T}:=B_1\#_T \Xi_2$
(which coincide with $\Xi_T$ for $s\le\tfrac T2$ 
and $s\ge -\tfrac T2$ respectively)
satisfy exponential estimates of the form 
$\|\Xi_{i,T}-\Xi_i(\cdot\pm T)\|_{\cC^k}\le C_ke^{-\delta T}$,
and hence their linearized operators are surjective 
with uniform estimates. Here we use the fact that $\Xi_1$ and $\Xi_2$ 
are regular in the sense of Definition~\ref{def:Areg}.

We have thus checked that the assumptions of \cite[Proposition~A.3.4.]{MS}
are satisfied with uniform constants for all $T\geq T_0$, where
$T_0>1$ is determined by comparing (\ref{eq:approx solution}) with
\cite[(A.3.5)]{MS}.
Hence the inverse function theorem provides unique solutions
$\xi_T\in\im Q_T\subset X_T$ of $f_T(\xi_T)=0$.
In other words, we can define 
$
\tau_T(\Xi_1,\Xi_2) := \Tilde{\Xi}_T = \tilde E(\Xi_T,\xi_T),
$
where $\Tilde{\Xi}_T\in\widetilde\cM(B_0,B_2;X_f)$ is the unique solution
of the form~(\ref{eq:Xitilde}) with $\Xi_T=\Xi_1\#_T\Xi_2$. 
This map is gauge equivariant and induces a map to the moduli space.
Note moreover that $\xi_T$ will be continuous with respect to $T$
in the $W^{1,p}$-norm and hence $\tilde\Xi_T$ as well as $\tau$ will
depend continuously on $T\in[T_0,\infty)$.
In the following we sketch the proof of properties (i)--(iii).

The convergence in~(ii) follows from the fact 
that the infinitesimal connection $\xi_T$ obtained 
in the inverse function theorem satisfies an
estimate of the form 
$\|\xi_T\|_{W^{1,p}}\leq C \| f_T(0) \|_{L^p} \leq C' e^{-\delta T}$ 
for uniform constants $C,C'$.

The index of $\tau_T(\Xi_1,\Xi_2)$ is given 
by~(\ref{eq:deltamu}), i.e.
\begin{align*}
\delta_f(\tau_T(\Xi_1,\Xi_2))
&= \mu_f(B_0,\tB_0) - \mu_f(B_1,\tB_1) + \mu_f(B_1,\tB_1) 
- \mu_f(B_2,\tB_2) \\
&= \delta_f(\Xi_1) + \delta_f(\Xi_2) = 2.
\end{align*}
Here $\tB_i:[0,1]\to\cL$ are paths from $\tB_i(0)=B_i$
to $\tB_i(1)=0$, where we pick any $\tB_1$ and pick 
the other paths such that $\tB_0$ is homotopic 
to the catenation of $\Xi_1|_{\R\times\Sigma}$ 
with $\tB_1$ and $\tB_1$ is homotopic to the catenation 
of $\Xi_2|_{\R\times\Sigma}$ with $\tB_2$.
Then, by construction, $\tB_0$ is homotopic to the 
catenation of $(\Xi_1\#_T \Xi_2)|_{\R\times\Sigma}$ 
with $\tB_2$. Moreover, $\tau_T(\Xi_1,\Xi_2)|_{\R\times\Sigma}$ 
is homotopic to $(\Xi_1\#_T \Xi_2)|_{\R\times\Sigma}$.

To see that $\tau$ is a diffeomorphism note first that both
domain and target are $1$-dimensional manifolds
(by the regularity and additivity of the indices).
Hence it suffices to show that $\tau$ is an injective immersion
by following the argument in~\cite[p.96]{Donaldson book}.
In fact, since the domain of $\tau$ is connected, it suffices to show that $\rd\tau$ 
is nonzero for all sufficiently large $T$.
We will show below that $\tau$ is $\cC^1$-close to the pregluing
$T\mapsto\Xi_T=\Xi_1\#_T\Xi_2$ as a map $[T_0,\infty)\to\cA^{1,p}(\R\times Y,\cL;B_0,B_2)$,
i.e.\ 
\begin{equation} \label{eq:dtau claim}
\bigl\| \tfrac\rd{\rd T}{\tilde{\Xi}_T} - \tfrac\rd{\rd T}\Xi_T \bigr\|_{W^{1,p}(\R\times Y)}
\underset{T\to\infty}{\longrightarrow} 0 .
\end{equation}
With this, the immersion condition
$\tfrac\rd{\rd T}\tau \neq 0 \in \rT_{\tau(T)}\cM(B_0,B_2;X_f)/\R$
follows if we can prove that the pregluing map is an immersion modulo 
gauge and time-shift with a uniform estimate.
Indeed, taking the infimum over all
$\psi\in\cC^\infty(\R\times Y,\cg)$, $\lambda\in\R$ we have
\begin{align*}
& \inf_{\psi,\lambda} \;
\bigl\| \tfrac\rd{\rd T}{\Xi_T} - \rd_{\Xi_T}\psi - \lambda\cdot\pd_s\Xi_T \bigr\|_{W^{1,p}(\R\times Y)} \\
&\geq 
\inf_{\lambda} \Bigl(
\inf_{\psi}\bigl\| \pd_s A_1 - \rd_{A_1}\psi - \lambda\cdot\pd_s A_1 \bigr\|_{W^{1,p}((-\infty,-1]\times Y)} \\
&\qquad\quad 
+ \inf_{\psi}\bigl\| - \pd_s A_2 - \rd_{A_2}\psi - \lambda\cdot\pd_s A_2 \bigr\|_{W^{1,p}([1,\infty)\times Y)}
\Bigr)
\; \geq \Delta >0 .
\end{align*}
Here we restricted the $W^{1,p}$-norm to the half cylinders $s\leq -T-1$ resp.\ $s\geq T+1$, where
$\Xi_T(s)=A_1(s+T)$ resp.\ $\Xi_T(s)=A_2(s-T)$. We also dropped the $\ds$-terms and applied various shifts.
The constant $\Delta>0$ is obviously independent of $T$.
It is positive since otherwise one could pick a minimizing sequence converging to limits 
$\lambda,\psi_1,\psi_2$ such that $(1-\lambda) \pd_s A_1=\rd_{A_1}\psi_1$
and $(1+\lambda)\pd_s A_2=-\rd_{A_2}\psi_2$.
However, from unique continuation (Proposition~\ref{prop:ucon}~(ii)) we know that
$\pd_s A_i(s) \not\in\im\rd_{A_i(s)}$, so $\rd_{A_i}\psi_i$ vanishes on both half cylinders,
which leaves the contradiction $1=\lambda=-1$.

It remains to establish (\ref{eq:dtau claim}).
We write $\dot{(\ldots)}$ for $\tfrac\rd{\rd T}(\ldots)$ and claim that
$$
\bigl\| \tfrac\rd{\rd T}{\tilde{\Xi}_T} - \dot\Xi_T \bigr\|_{W^{1,p}}
\leq \bigl\|\partial_1 \tilde E(\Xi_T,\xi_T) - {\rm Id} \bigr\|\, \bigl\| \dot\Xi_T \bigr\| 
+ \bigl\| \partial_2 \tilde E(\Xi_T,\xi_T) \dot\xi_T \bigr\|
\underset{T\to\infty}{\longrightarrow} 0
$$ 
due to the identities $\tilde E(\cdot,0)={\rm Id}$ and
$\pd_2 \tilde E(\Xi_T,0)={\rm Id}$, the boundedness of $\|\dot\Xi_T\|_{W^{1,p}}$ 
(due to exponential decay), 
and the convergence $\xi_T\to 0$ and $\|\dot\xi_T\|_{W^{1,p}}\to 0$.
To check the latter recall the abstract setup for the inverse function theorem.
Taking the $T$-derivative of $f_T(\xi_T)=0$ we obtain
\begin{align*}
&\bigl\|\rd f_T(\xi_T) \dot{\xi}_T\bigr\|_{L^p} = \bigl\|\dot{f_T}(\xi_T)\bigr\|_{L^p}  \\
&=\bigl\| \rd^+_{\tilde E(\Xi_T,\xi_T)} \partial_1 \tilde E(\Xi_T,\xi_T) \dot\Xi_T \bigr\|_{L^p}
+ \bigl\| [\dot\Xi_T \wedge * \xi_T ] \bigr\|_{L^p} 
\underset{T\to\infty}{\longrightarrow} 0 .
\end{align*}
This convergence uses the same estimates as before and the fact that $\rd^+_{\Xi_T}\dot\Xi_T$
vanishes except for near $s=\pm\frac T2$, where it is exponentially small.
Now write $\xi_T=Q_T\eta_T$ with $\eta_T=\cD_T\xi_T$, then
$$
\dot{\xi}_T = \zeta_T + \dot{Q}_T \cD_T \xi_T
\qquad\text{with}\quad \zeta_T=Q_T\dot{\eta}_T\in\im Q_T .
$$
We have $\|\dot{Q}_T \cD_T \xi_T\|_{W^{1,p}}\to 0$ since $\xi_T\to 0$
and the operators $\cD_T:W^{1,p}\to L^p$ and $\dot{Q}_T:L^p\to\dom\cD_T\subset W^{1,p}$
are uniformly bounded. The first bound is due to
$\|\cD_T-\cD_{T_0}\| \leq \| \Xi_T - \Xi_{T_0} \|_{\cC^0}$;
similarly $\dot\cD_T:W^{1,p}\to L^p$ and $\dot\cD_T^*:W^{2,p}\to W^{1,p}$
are bounded in terms of $\|\dot\Xi_T\|_{\cC^0}$ resp.\ $\|\dot\Xi_T\|_{\cC^1}$,
and we have the identity
$\dot{Q}_T=\dot\cD_T^*(\cD_T\cD_T^*)^{-1} 
- Q_T (\dot\cD_T\cD_T^* + \cD_T\dot\cD_T^*)(\cD_T\cD_T^*)^{-1}$.
Here the uniform bound on $(\cD_T\cD_T^*)^{-1}$, that is 
$\|\eta\|_{W^{2,p}}\leq C \|\cD_T\cD_T^*\eta\|_{L^p}$, 
follows from combining (\ref{eq:DT2}) with the $W^{2,p}$-version of (\ref{eq:DT1}).

Finally, we can prove that $\|\zeta_T\|_{W^{1,p}}\to 0$ because, starting from (\ref{eq:DT2}),
\begin{align*}
&\|\zeta_T\|_{W^{1,p}}\leq C \|\cD_T\zeta_T\|_{L^p} \\
&\leq C \bigl(\| (\rd f_T(\xi_T)-\cD_T ) \zeta_T\|_{L^p} +  \|\rd f_T(\xi_T)\dot\xi_T\|_{L^p} 
+  \|\rd f_T(\xi_T) \dot{Q}_T \cD_T \xi_T\|_{L^p} \bigr) .
\end{align*}
Here the first term can be absorbed into the left hand side
by (\ref{eq:mightneedthis}) for sufficiently large $T>T_0$ and the other terms
converge to zero as $T\to\infty$, using a uniform bound on $\rd f_T(\xi_T)$
from $\|\rd f_T(\xi_T)-\cD_{T}\| \leq \| \tilde E(\Xi_T,\xi_T) - \Xi_{T} \|_{\cC^0}$.
This finishes the proof that $\dot\xi_T\to 0$, hence (\ref{eq:dtau claim}) holds
and (i) is proven.

Assertion~(iii) follows from the uniqueness statement in the inverse function theorem
if we can find $u\in\cG(\R\times Y)$, $\sigma\in\R$, and $T'>T_0$
such that $u^*\Xi(\cdot+\sigma)=\tilde E(\Xi_{T'},\xi)$ with
$\xi$ satisfying (\ref{eq:Xitilde}) and $W^{1,p}$-small.
For each $(\sigma,T')$ close to $(0,T)$ we can 
use the local slice theorem to find 
$u_{\sigma,T'}$ and $\xi_{\sigma,T'}$ satisfying
$$
u_{\sigma,T'}^*\Xi(\cdot+\sigma)=\tilde E(\Xi_{T'},\xi_{\sigma,T'}),
\qquad
\rd_{\Xi_{T'}}^*\xi_{\sigma,T'}=0, \qquad 
*\xi_{\sigma,T'}|_{\R\times\pd Y}=0   .
$$
One then finds $(\sigma,T')$ satisfying
$\xi_{\sigma,T'}\in\im\cD_{T'}^*=(\ker\cD_{T'})^\perp$
by a further implicit function theorem.
Namely, there is a basis $(\eta_{1,T'},\eta_{2,T'})$ 
of $\ker\cD_{T'}$ close to 
$(\pd_s\Xi_1 \#_{T'} 0 \,,\, 0 \#_{T'} \pd_s\Xi_2)$.
Then the map
$(\sigma,T')\mapsto 
(\langle \xi_{\sigma,T'} , \eta_{1,T'} \rangle ,
\langle \xi_{\sigma,T'} , \eta_{2,T'} \rangle )$
is invertible and has a zero close to $(0,T)$.
\end{proof}

\begin{rmk}\rm \label{rmk:redglue}
In Theorem~\ref{thm:glue} we can allow $B_1$ to be reducible 
(but still nondegenerate). Then we obtain a gluing map
$$
\tau :\; (T_0,\infty) \times \bigl(\cG_{B_1}/\{\pm\one\}\bigr) 
\;\to\; \cM^{2+\dim H^0_{B_1}}(B_0,B_2;X_f)/\R
$$
with the same properties as in Theorem~\ref{thm:glue}.
This map is constructed by starting from a preglued connection
$\Xi_1\#_{g,T} \Xi_2$ that takes $g\in\cG_{B_1}/\{\pm\one\}$ into account by
$$
A(s) := \left\{\begin{array}{ll}
A_1(\tfrac T2-\phi(-\tfrac T2-s)),& s\le -\tfrac T2,\\
B_1=g^*B_1,& s\in[-\tfrac T2 ,\tfrac T2] , \\
g^*A_2(-\tfrac T2+\phi(-\tfrac T2+s)),& s\ge \tfrac T2 .
\end{array}\right.
$$
The index identity again follows from (\ref{eq:deltamu})
and the uniformly bounded right inverse can be constructed 
using weighted spaces, as described in \cite[4.4.1]{Donaldson book}.

This shows that the breaking of trajectories at the zero connection 
can be excluded in low dimensional moduli spaces since the stabilizer 
$\cG_0\subset\cG(Y)$ adds $3$ to the index of the glued 
connection. However, this argument is not needed for the construction 
of Floer homology. In the proof of Corollary~\ref{cor:glue} below, 
we use simpler index bounds to exclude breaking at the zero connection.
\end{rmk}

Theorem~\ref{thm:glue} gives rise to maps
$$
\tau_T : \cM^1(B_0,B_1)/\R \times \cM^1(B_1,B_2)/\R \to \cM^2(B_1,B_2)/\R 
$$
defined by choosing one representative for each gauge and 
shift equivalence class in each moduli space $\cM^1(A^-,A^+)/\R$ 
with $[A^+],[A^-]\in\cR_f\setminus[0]$.

\begin{cor} \label{cor:glue}
Let $A^+,A^-\in\Crit(\CS_\cL+h_f)\setminus[0]$.
Then, for $T_0$ sufficiently large, the sets
$\tau_{(T_0,\infty)}([\Xi_1],[\Xi_2])\subset\cM^2(A^-,A^+)/\R$,
indexed by ${[B]\in\cR_f\setminus[0]}$ and
$([\Xi_1],[\Xi_2])\in \cM^1(A^-,B)/\R \times \cM^1(B,A^+)/\R $,
are pairwise disjoint. Moreover, their complement
$$
\cM^2(A^-,A^+)/\R \quad\setminus 
\bigcup_{[0]\neq[B]\in\cR_f} \;\bigcup_{T> T_0} 
\tau_T\bigl(\cM^1(A^-,B)/\R \times \cM^1(B,A^+)/\R \bigr)
$$
is compact.
\end{cor}
\begin{proof}
The sets $\tau_{(T_0,\infty)}([\Xi_1],[\Xi_2])$ are disjoint
for $T_0$ sufficiently large since they converge to different
broken trajectories for $T_0\to\infty$, see 
Theorem~\ref{thm:glue}~(ii).

To prove compactness we assume by contradiction that there exists
a sequence $[\Xi^\nu]\in\cM^2(A^-,A^+;X_f)/\R$ in the complement
of the image of $\tau$ as above, and that has no convergent subsequence.
These solutions have index $2$ and hence 
fixed energy by Corollary~\ref{cor:index}~(i).
By Theorem~\ref{thm:compact1} we can pick 
a subsequence and representatives, still denoted by $\Xi^\nu$,
that converge to a broken trajectory $(\Xi_1,\dots,\Xi_\ell)$ 
modulo bubbling.
By transversality we do not have solutions of negative index,
so Corollary~\ref{cor:compact} implies that there is no bubbling,
and the index identity in Theorem~\ref{thm:compact1} implies
$\ell\leq 2$.
In the case $\ell=1$ we would obtain a convergent subsequence 
from Theorem~\ref{thm:compact3}, hence the limit must be
a broken trajectory with two index $1$ solutions and an
irreducible intermediate critical point $B$.
The time-shifts and gauge transformations in 
Theorem~\ref{thm:compact1} can be chosen such that the
limit $(\Xi_1,\Xi_2)$ consists of the fixed representatives
used in the definition of $\tau_T$.
Now the assertion of Theorem~\ref{thm:compact2}~(ii) 
can be reformulated as
$$
\bigl\| {v^\nu}^*\Xi^\nu (\cdot + \tfrac 12(s^\nu_1+s^\nu_2))
- \Xi_1\#_{T^\nu} \Xi_2 \bigr\|_{W^{1,p}(\R\times Y)} \to 0
$$
for $T^\nu:=\frac 12(s^\nu_2-s^\nu_1)\to\infty$.
Then, by Theorem~\ref{thm:glue}~(iii), $[\Xi^\nu]$ lies in 
the image of $\tau$ for sufficiently large $\nu$,
in contradiction to the assumption.
\end{proof}


\section{Coherent orientations}\label{sec:orient}

Let $Y$ be a compact oriented Riemannian $3$-manifold 
with boundary $\pd Y=\Sigma$ and $\cL\subset\cA(\Sigma)$ 
be a gauge invariant, monotone, irreducible Lagrangian submanifold 
satisfying (L1-3) on page~\pageref{p:L1}.
In this section it is essential that we restrict to the 
case of $Y$ being connected with nonempty boundary, 
so that the gauge group $\cG(Y)$ is connected.
The construction of orientations for closed $Y$ can
be found in~\cite[5.4]{Donaldson book}. 
Fix a perturbation $h_f$ such that
every critical point of $\CS_\cL + h_f$ is nondegenerate
and every nontrivial critical point is irreducible
(see Definition~\ref{def:freg}).
For every pair of irreducible critical points  
${A^-,A^+\in\Crit(\CS_\cL+h_f)}$ we consider the space
$$
\cA(A^-,A^+):= \bigl\{ \A\in\cA(\R\times Y,\cL) \st 
\A|_{[s,s+1]\times Y} 
\underset{s\to\pm\infty}{\longrightarrow}
0\ds + A^\pm \;\text{exponent.} \bigr\},
$$ 
which consists of smooth connections $\A=\Phi\ds+A$ on $\R\times Y$
that are given by paths $\Phi:\R\to\Om^0(Y,\cg)$ and $A:\R\to\cA(Y,\cL)$ 
that converge exponentially with all derivatives to 
$0$ and $A^\pm$, respectively, as $s\to\pm\infty$.
If we allow the limits $A^\pm$ to vary within 
gauge orbits of critical points, we obtain the spaces
$$
\cA([A^-],[A^+]):=
\bigcup_{u^\pm\in\cG(Y)}
\cA((u^-)^*A^-,(u^+)^*A^+) .
$$
We denote by 
$$
\Or([A^-],[A^+])
:=\bigsqcup_{\A\in\cA([A^-],[A^+])}\Or(\cD_\A) 
\;\to\;\cA([A^-],[A^+])
$$ 
the principal $\Z_2$-bundle whose fibre over $\A\in\cA([A^-],[A^+])$
is the set $\Or(\cD_\A)$ of orientations of the determinant line
$$
\det(\cD_\A) := \Lambda^{\text{max}} \bigl( \ker\cD_\A \bigr) 
\otimes \Lambda^{\text{max}} \bigl( \text{coker}\cD_\A \bigr)^* .
$$
Here $\cD_\A$ is the linearized operator~(\ref{eq:DA}).
Any homotopy $[0,1]\to\cA([A^-],[A^+])$, $\lambda\mapsto A_\lambda$
induces an isomorphism
$$
\Or(\cD_{\A_0}) \to 
\Or(\cD_{\A_1})
$$ 
by path lifting.
A gauge transformation ${u\in\cG(\R\times Y)}$ 
which converges exponentially to 
$u^\pm\in\cG(Y)$
as $s\to\pm\infty$ gives rise to a bundle isomorphism
$$
u^*:\Or(A^-,A^+) \to 
\Or((u^-)^*A^-,(u^+)^*A^+)
$$ 
induced by the conjugate action of $u$ on kernel and cokernel.
The pregluing construction in~(\ref{eq:preg})
for $\A_1\in\cA(B_0,B_1)$ and $\A_2\in\cA(B_1,B_2)$
induces a natural isomorphism 
$$
\sigma_T :
\Or(\cD_{\A_1})\otimes \Or(\cD_{\A_2}) 
\to \Or(\cD_{\A_1\#_T\A_2})
$$
for sufficiently large $T$.
If both $\cD_{\A_1}$ and $\cD_{\A_2}$ are surjective, then 
$\cD_{\A_1\#_T\A_2}$ is surjective for $T$ sufficiently large, 
by estimates as in the proof of Theorem~\ref{thm:compact2}, and
$\sigma_T$ is induced by the isomorphism
$\ker(\cD_{\A_1})\times \ker(\cD_{\A_2}) 
\to \ker(\cD_{\A_1\#_T\A_2})$.
The general case is reduced to the surjective case by the
method of stabilizations as in \cite[Section~3(a)]{Donaldson orient}.

We will also have to glue connections
over $S^4$ to connections over $\R\times Y$.
For that purpose we denote by 
$\cA(P_u)$ the space of connections
on the bundle $P_u$ that is obtained by gluing two copies
of $\C^2\times B^4$ with the transition function $u\in\cG(S^3)$.
Then for every $\A\in\cA(A^-,A^+)$ and 
$\Xi_u\in\cA(P_u)$ we can construct a preglued 
connection $\A\#_T\Xi_u\in\cA(A^-,\bar u^*A^+)$
by taking the connected sum $(\R\times Y)\#_{\pd D_T} S^4$
and trivializing the induced bundle over $\R\times Y$.
Here we denote by $D_T\subset \R\times Y$ the ball of radius $T^{-1}$ 
centred at $(0,y)$ for some $y\in{\rm int}(Y)$, and 
after the trivialization we have
$$
\bigl(\A\#_T\Xi_u\bigr) |_{(\R\times Y) \setminus D_T} = \tilde u^*\A
$$
for a gauge transformation $\tilde u$ on $(\R\times Y)\setminus D_T$ 
with $\tilde u|_{\pd D_T}\cong u$.
We fix these extensions such that
$\tilde u |_{(-\infty,-1]\times Y} \equiv \one$
and $\tilde u |_{[1,\infty)\times Y} \equiv \one$,
and hence $\tilde u|_{\R\times\pd Y}$ defines
a path $\bar u:\R\to\cG(\Sigma)$ with $\bar u(s)=\one$ for $|s|\ge 1$.
A partial integration on $[-1,1]\times Y$ then shows that the degree of
this loop is $\deg(\bar u)=\deg(u)$.
So we have both $\A,\,\A\#_T\Xi_u\in\cA(A^-,A^+)$, but
the homotopy classes (of paths in $\cL$ with fixed endpoints)
of $\A|_{\pd Y}$ and $(\A\#_T\Xi_u)|_{\pd Y}$ differ by $\deg(u)$.
The determinant line bundle over the contractible space $\cA(P_u)$
is canonically oriented (compatible with gauge transformations, 
homotopies, and gluing, see e.g.\ \cite[Proposition~5.4.1]{DK}),
and as before pregluing induces an isomorphism
$$
\sigma_T :
\Or(\cD_{\A})\otimes \Or(\cD_{\Xi_u}) 
\to \Or(\cD_{\A\#_T\Xi_u})
$$
for $T$ sufficiently large. The various isomorphisms, 
induced by homotopies, gauge transformations, 
and pregluing, all commute in the appropriate sense.

\begin{dfn} \label{dfn:orient}
A {\bf system of coherent orientations}
is a collection of sections 
$$
\cA([A^-],[A^+]) \to \Or([A^-],[A^+]) : \A \mapsto o_\A ,
$$
one for each pair $[A^-],[A^+]\in\Crit(\CS_\cL+h_f)/\cG(Y)\setminus[0]$
of nontrivial gauge equivalence classes of critical points,
satisfying the following conditions.
\begin{description}
\item[(Homotopy)]
The sections $o:\cA([A^-],[A^+]) \to \Or([A^-],[A^+])$ are continuous.
In other words, if $[0,1]\to\cA([A^-],[A^+]):\lambda\mapsto\A_\lambda$ 
is a continuous path, then the induced isomorphism 
$\Or(\cD_{\A_0})\to\Or(\cD_{\A_1})$ sends
$o_{\A_0}$ to $o_{\A_1}$.
\item[(Equivariance)]
For every $\A\in\cA(A^-,A^+)$ and every $u\in\cG(\R\times Y)$ 
that converges exponentially to $u^\pm\in\cG(Y)$
as $s\to\pm\infty$ we have 
$$
o_{u^*\A}=u^* o_\A .
$$
\item[(Catenation)]
Let $\A\in\cA(B_0,B_1)$ and $\A'\in\cA(B_1,B_2)$,
then for $T$ sufficiently large we have
$$
o_{\A\#_T\A'} = \sigma_T(o_\A\otimes o_{\A'}).
$$
\item[(Sum)]
Let $\A\in\cA(A^-,A^+)$, $u\in\cG(S^3)$, and $\Xi_u\in\cA(P_u)$,
then for $T$ sufficiently large we have
$$
o_{\A\#_T\Xi_u} = \sigma_T(o_\A\otimes o_{\Xi_u}) .
$$
\item[(Constant)]
If $\A\equiv A^- = A^+$, then $o_\A$ 
is the orientation induced by the canonical isomorphism 
$\det(\cD_\A)\to\R$. (Under this assumption $\cD_\A$ 
is bijective.)
\end{description}
\end{dfn}

\begin{rmk}\label{rmk:or}\rm
{\bf (i)}
The {\it (Equivariance)} axiom follows from the 
{\it (Homotopy)} axiom.  To see this note that,
since $Y$ is connected with nonempty boundary,
the gauge groups $\cG(Y)$ and hence $\cG(\R\times Y)$
are connected. (Here we do not fix the boundary values
or limits of the gauge transformations.) 
The claim then follows from the following observation. 

\smallskip\noindent{\bf (ii)}
For every continuous path $[0,1]\to\cG(\R\times Y):\lambda\mapsto u_\lambda$ 
with $u_0=\one$ the isomorphism $u_1^*:\det(\cD_\A)\to\det(\cD_{u_1^*\A})$ 
coincides with the isomorphism induced by 
the homotopy $\lambda\mapsto u_\lambda^*\A$.
To see this consider the continuous family of paths
$[0,1]\to\cG(\R\times Y):\lambda\mapsto u_{\tau\lambda}$ for $\tau\in[0,1]$.
Then the assertion holds obviously for $\tau=0$
(both maps are the identity) and hence,
by continuity, for all $\tau$.
\end{rmk}

\begin{thm} \label{thm:orient}
Fix representatives $B_1,\ldots,B_N$, one for each nontrivial 
gauge equivalence class in $\Crit(\CS_\cL+h_f)/\cG(Y)\setminus [0]$, 
connections $\A_i\in\cA(B_i,B_{i+1})$, and orientations 
$o_i\in\Or(\cD_{\A_i})$ for $i=1,\ldots,N-1$.
Then there is a unique system of coherent orientations 
$o_\A\in\Or(\cD_{\A})$ such that $o_{\A_i}=o_i$ for all $i$.
\end{thm}

The proof of this theorem will make use of the following lemma.

\begin{lemma}\label{le:orient}
Fix a pair $A^\pm\in\cA(Y,\cL)$ of irreducible
and nondegenerate critical points of $\CS_\cL+h_f$.
Let $[0,1]\to\cA([A^-],[A^+]):\lambda\mapsto\A_\lambda$
be a smooth path and $u\in\cG(\R\times Y)$
such that $\A_1=u^*\A_0$. Then the isomorphism
$$
u^*:\Or(\cD_{\A_0})\to\Or(\cD_{\A_1})
$$
agrees with the one induced by the path $\lambda\mapsto\A_\lambda$.
In particular, the orientation bundle
$
\Or([A^-],[A^+])\to \cA([A^-],[A^+])
$
admits a trivialization.
\end{lemma}

\begin{proof}
By continuity, it suffices to prove the identity under the asumption
$\p_sA_\lambda(s)=0$, $\Phi_\lambda(s)=0$,
and $\p_su(s)=0$ for $\Abs{s}\ge1$.
Then there are paths $[0,1]\to\cG(Y), \lambda\mapsto v_\lambda^\pm$
such that $(v_\lambda^-)^*A_\lambda(s)=A^-$ for $s\le-1$ and
$(v_\lambda^+)^*A_\lambda(s)=A^+$ for $s\ge1$.
We can replace $A^\pm$ by $((v_0^\pm)^{-1})^*A^\pm$ 
and thus assume in addition that $v_0^\pm=\one$.
Now there is a smooth map
$[0,1]\times\R\to\cG(Y):(\lambda,s)\mapsto u_\lambda(s)$
such that $u_0\equiv\one$, $u_\lambda(s)=v_\lambda^-$ for $s\le-1$
and $u_\lambda(s)=v_\lambda^+$ for $s\ge1$. Define  
$$
\A_\lambda^\tau:=u_{\lambda\tau}^*\A_\lambda,\qquad
u^\tau := u_0^{-1}uu_\tau 
$$
for every $\tau\in[0,1]$. Then we have $\A^\tau_1=(u^\tau)^*\A^\tau_0$.
By continuity, the assertion now holds for $\tau=1$ if and only if 
it holds for $\tau=0$, that is for the original pair $(\{\A_\lambda\},u)$.  
For $\tau=1$ we have $A_\lambda^1(s)=A^\pm$ and $u^1(s)=\one$ for $\pm s\ge1$. 

Finally, we prove the lemma in the case 
$A_\lambda(s)=A^\pm$ for $\pm s\ge 1$.
For $T\ge 2$ we define the catenation 
$\B^T_\lambda:=\Phi^T_\lambda\ds + B^T_\lambda \in\cA(\R/2T\Z\times Y,\cL)$ 
and $u^T\in\cG(\R/2T\Z\times Y)$ by 
\begin{equation*}
\begin{split}
B^T_\lambda(s) &:= \left\{\begin{array}{ll}
A_\lambda(s),&-T/2\le s\le T/2,\\
A_0(T-s),&T/2\le s\le 3T/2,
\end{array}\right. \\
\Phi^T_\lambda(s) &:= \left\{\begin{array}{ll}
\Phi_\lambda(s),&-T/2\le s\le T/2,\\
\Phi_0(T-s),&T/2\le s\le 3T/2,
\end{array}\right. \\
u^T(s) &:= \left\{\begin{array}{ll}
u(s),&-T/2\le s\le T/2,\\
\one,&T/2\le s\le 3T/2.
\end{array}\right.
\end{split}
\end{equation*}
Then $B^T_\lambda(s)=A^\pm$ and $\Phi^T_\lambda=0$ for $\pm s\in[1,T-1]$.
Moreover we have $\B^T_1=(u^T)^*\B^T_0$.  
For $T$ sufficiently large the linear gluing theory gives rise 
to a continuous family of isomorphisms
$$
\phi^T_\lambda:\Or(\cD_{\A_\lambda})\to \Or(\cD_{\B^T_\lambda}) ,
$$
where $\cD_{\B^T_\lambda}$ denotes the anti-self-duality
operator on $\R/2T\Z\times Y$ introduced in Section~\ref{sec:S1Y}. 
The gluing operators commute with the gauge 
transformations, i.e.
$$
\phi^T_1\circ u^* = (u^T)^*\circ\phi^T_0:
\Or(\cD_{\A_0})\to \Or(\cD_{\B^T_1}).
$$
The isomorphisms induced by the homotopies $\lambda\mapsto \A_\lambda$
and $\lambda\mapsto \B_\lambda^T$ satisfy the same relation. 
By Theorem~\ref{thm:S1Y}~(iv) (with $v=\one$), the isomorphism
$(u^T)^*:\Or(\cD_{\B^T_0})\to\Or(\cD_{\B^T_1})$
agrees with the one induced by the 
path $\lambda\mapsto B^T_\lambda$. 
Hence the same holds for $u^*$ and this proves the 
desired identity.

To see that $\Or([A^-],[A^+])\to \cA([A^-],[A^+])$ admits a trivialization
we only need to check that parallel transport around loops induces the
identity isomorphism on the fibre. 
This follows immediately from the identification
of the homotopy induced isomorphism 
with $u^*:\Or(\cD_{\A_0})\to\Or(\cD_{\A_0})$
for $u=\one$.
\end{proof}

\begin{proof}[Proof of Theorem~\ref{thm:orient}.]
The orientation bundle over the constant component of 
$\cA([B_i],[B_i])$ is canonically oriented by the {\it (Homotopy)} 
and {\it (Constant)} axioms.  The orientation on the 
other components of $\cA([B_i],[B_i])$ is determined 
by the {\it (Sum)} axiom because any connection $\A\in\cA([B_i],[B_i])$  
is homotopic to $\B_i\#_T\Xi_u$ for the constant solution 
$\B_i\equiv B_i$, a connection $\Xi_u$ over $S^4$ 
associated to a nontrivial $u\in\cG(S^3)$, and any $T>0$.
Indeed, since $\cG(Y)$ is connected, $\A$ can be homotoped 
to a connection with fixed limits in $\cA(B_i,B_i)$.
Moreover, there is a homotopy equivalence
$\cA(B_i,B_i)\to\cC^\infty(S^1,\cL)$
which assigns to each connection $\A\in\cA(B_i,B_i)$
a based loop in $\cL$ obtained from the path 
$\A|_{\pd Y}:\R\to\cL$ with endpoints $B_i|_{\pd Y}$.
Now, by (L2), the loop $\A|_{\pd Y}$ in $\cL$ is homotopic to 
$\hat u^* B_1|_{\pd Y}$ for some loop $\hat u:S^1\to\cG(\Sigma)$.
Hence $\A$ is homotopic to $\B_i\#_T\Xi_u$ 
for the associated $u\in\cG(S^3)$.
Similarly, the orientation bundle over 
$\cA([B_i],[B_{i+1}])$ is oriented by $o_{\A_i}$ 
and the {\it (Homotopy)} and {\it (Sum)} axioms, 
because any connection in $\cA([B_i],[B_{i+1}])$ 
is homotopic to $\A_i\#_T\Xi_u$ for some $u\in\cG(S^3)$.  
Finally, the orientation bundles over general spaces
$\cA([B_i],[B_j])$ are oriented by the {\it (Catenation)} 
axiom and the previously fixed orientations. 
This proves uniqueness.  

To establish existence note that, by Lemma~\ref{le:orient}, 
we have a choice of two possible orientations over every 
component of each $\cA([A^-],[A^+])$. 
Each of the possible combinations of choices 
satisfies the {\it (Homotopy)} axiom by construction.
To see that the choices can be made such that 
the {\it (Constant)}, {\it (Catenation)}, and {\it (Sum)} 
axioms are satisfied (and so the {\it (Equivariance)} 
axiom follows from Remark~\ref{rmk:or}),
one needs to check that the isomorphisms in the
{\it (Catenation)}, {\it (Sum)}, and {\it (Homotopy)} 
axioms all commute. For example, let $\A_\lambda\in\cA(B_i,B_j)$ 
and ${\A'_\lambda\in\cA(B_j,B_k)}$ be smooth families
parametrized by $\lambda\in[0,1]$ and denote~by
\begin{align*}
\rho:\Or(\cD_{\A_0})\to\Or(\cD_{\A_1}),\qquad
&\rho':\Or(\cD_{\A'_0})\to\Or(\cD_{\A'_1}), \\
\rho^T:\Or(\cD_{\A_0\#_T\A'_0})&\to\Or(\cD_{\A_1\#_T\A'_1})
\end{align*}
the isomorphisms induced by the homotopies $\lambda\mapsto\A_\lambda$,
$\A'_\lambda$, and $\A_\lambda\#\A'_\lambda$.
Let
$$
\sigma^T_\lambda:\Or(\cD_{\A_\lambda})\otimes\Or(\cD_{\A'_\lambda})
\to\Or(\cD_{\A_\lambda\#_T\A'_\lambda})
$$
denote the catenation isomorphisms for $T$ sufficiently large. 
A parametrized version of the linear gluing construction then proves that
$$
{\sigma^T_1\circ(\rho\otimes\rho')=\rho^T\circ\sigma^T_0}.
$$
A similar statement holds for the {\it (Homotopy}) and {\it (Sum)} isomorphisms. 
That two {\it (Catenation)} isomorphisms commute
is a kind of associativity rule modulo homotopy
and the proof involves a simultaneous gluing construction
for three connecting trajectories; similarly for the 
commutation rules of the {\it (Sum)} and {\it (Catenation}) 
isomorphisms.  All these arguments are exactly 
as in the standard theory and the details will be omitted. 
\end{proof}


\section{Floer homology}\label{sec:floer}

Let $Y$ be a compact connected oriented $3$-manifold 
with boundary $\p Y=\Sigma$ and $\cL\subset\cA(\Sigma)$ 
be a gauge invariant, monotone, irreducible Lagrangian submanifold 
satisfying (L1-3) on page~\pageref{p:L1}.
Fix a Riemannian metric $g$ on $Y$, a regular perturbation
$(\gamma,f)\in\Gamma_m\times\cF_m$ as in 
Theorem~\ref{thm:freg}, and a system $\co=\{o_\A\}_\A$
of coherent orientations as in Theorem~\ref{thm:orient}.
Associated to these data we define a Floer homology group
$\HF(Y,\cL;g,f,\co)$ as follows. 

Since the trivial connection is nondegenerate by (L3), the set
$$ 
\cR_f
:= \bigl\{ A\in\cA(Y) \st F_A + X_f(A)= 0 ,\; A|_{\pd Y}\in\cL \bigr\}/\cG(Y) 
$$ 
of gauge equivalence classes of critical points of $\CS_\cL+h_f$
is finite, by Proposition~\ref{prop:finite}. 
The nontrivial critical points determine a chain complex 
$$
\CF(Y,\cL;f) 
:= \bigoplus_{[A]\in\cR_f\setminus[0]} \Z \, \langle A \rangle .
$$
with a $\Z/8\Z$-grading $\mu_f:\cR_f\to\Z/8\Z$ defined by
the spectral flow (see Corollary~\ref{cor:eta mu}).
We emphasize that the spectral flow is 
invariant under homotopies of the metric and 
of the perturbation with fixed critical points.
To define the boundary operator we consider the space 
$$
\widetilde\cM(A^-,A^+;g,X_f)
:= \left\{ \begin{array}{l}
A + \Phi\ds \\
\in \cA(\R\times Y)
\end{array}
\left|
\begin{array}{l}
\pd_s A - \rd_A\Phi + * (F_A + X_f(A)) = 0 \\
A(s)|_\Sigma\in\cL \quad\forall s\in\R \\
\lim_{s\to\pm\infty}A(s)= A^\pm \\
\Phi|_{\{|s|\geq 1\}}\equiv 0 \\
E_f(\A)<\infty 
\end{array}\right.\right\} .
$$
This space is invariant under the group $\cG(A^-,A^+)$ of 
gauge transformations $u\in\cG(\R\times Y)$ that satisfy
$u(s)=u^\pm\in\cG_{A^\pm}$ for $\pm s\geq 1$. The quotient spaces 
$
\widetilde\cM(A^-,A^+;g,X_f)/\cG(A^+,A^-)
$
are canonically isomorphic for different choices of 
representatives $A^\pm$ of critical points. 
The index of the linearized operator at $[\A]$ 
is $\delta_f(\A) \equiv \mu_f(A^-)-\mu_f(A^+)$ (modulo $8$).
For $k\in\Z$ we denote the index $k$ part of the Floer moduli 
space by
$$
\cM^k(A^-,A^+;g,X_f)
:=\bigl\{[\A]\in\widetilde\cM(A^-,A^+;g,X_f)/\cG(A^-,A^+) \st
\delta_f(\A)=k \bigr\}.
$$
For $k\leq 7$ this is a smooth $k$-dimensional manifold 
(see Section~\ref{sec:fredholm} and Definition~\ref{def:freg}).
The energy of a solution in this space is
$
E_f(\A)=\frac 12\pi^2 (k+\eta_f(A^+)-\eta_f(A^-))
$
by Corollary~\ref{cor:index}~(i), and hence is independent of $\A$. 
Moreover, $\R$ acts on $\cM^k(A^-,A^+;g,X_f)$ by time--shift,
and the action is proper and free unless $A^-= A^+$ and $k=0$.
For $k=1$ the quotient space $\cM^1(A^-,A^+;X_f)/\R$ is a 
finite set, by Corollary~\ref{cor:compact11}.
Counting the elements with signs gives rise to a boundary 
operator on $\CF(Y,\cL;f)$ via
\begin{equation}\label{eq:bd}
\pd \langle A^-\rangle:= 
\sum_{[A^+]\in \cR_f\setminus[0]}\left(
\sum_{[\A]\in\cM^1(A^-,A^+;g,X_f)/\R}\nu(\A)\right)
\langle A^+ \rangle.
\end{equation}
Here $\nu(\A):=1$ whenever the element
$\p_s\A\in\ker\cD_\A=\det(\cD_\A)$ is positively oriented with 
respect to $o_\A$ and $\nu(\A):=-1$ otherwise. 
The next two theorems are the main results of this paper;
their proofs take up the rest of this section.

\begin{thm}\label{thm:floer1}
The operator $\p:\CF(Y,\cL;f)\to\CF(Y,\cL,f)$
defined by~(\ref{eq:bd}) satisfies $\p\circ\p=0$.
\end{thm}

The {\bf Floer homology group} of the pair 
$(Y,\cL)$ equipped with the regular data 
$(g,f,\co)$ is defined by 
$$
\HF(Y,\cL;g,f,\co) := \frac{\ker\p:\CF(Y,\cL;f)\to\CF(Y,\cL;f)}
{\im\p:\CF(Y,\cL;f)\to\CF(Y,\cL;f)}.
$$
The next theorem shows that it is independent of the 
choices of metric, perturbation, and coherent orientations.

\begin{thm}\label{thm:floer2}
There is a collection of isomorphisms
$$
\Phi^{\beta\alpha}:\HF(Y,\cL;g^\alpha,f^\alpha,\co^\alpha)
\to \HF(Y,\cL;g^\beta,f^\beta,\co^\beta),
$$
one for any two regular triples $(g^\alpha,f^\alpha,\co^\alpha)$
and $(g^\beta,f^\beta,\co^\beta)$, such that
\begin{equation}\label{eq:abc}
\Phi^{\gamma\beta}\circ\Phi^{\beta\alpha}=\Phi^{\gamma\alpha},\qquad
\Phi^{\alpha\alpha} = \Id
\end{equation}
for any three regular triples $(g^\alpha,f^\alpha,\co^\alpha)$,
$(g^\beta,f^\beta,\co^\beta)$, and $(g^\gamma,f^\gamma,\co^\gamma)$.
\end{thm}

\begin{proof}[Proof of Theorem~\ref{thm:floer1}.]
For $A^\pm\in\cR_f\setminus[0]$ denote
$$
n(A^-,A^+):=\sum_{[\A]\in\cM^1(A^-,A^+;g,X_f)/\R}\nu(\A) .
$$
Then the equation $\pd\circ\pd=0$ is equivalent
to the formula 
\begin{equation}\label{eq:d20}
\sum_{[B]\in\cR_f\setminus[0]}
n(A^-,B) \, n(B,A^+) = 0
\end{equation}
for all $A^\pm\in\cR_f\setminus [0]$.
The proof of~(\ref{eq:d20}) is exactly as in the standard
case. One studies the moduli space $\cM^2(A^-,A^+)/\R$.  
This is a $1$-manifold, oriented by the coherent 
orientations of Theorem~\ref{thm:orient}.
By Corollary~\ref{cor:glue} its ends are in one-to-one
correspondence with pairs of trajectories in 
$\cM^1(A^-,B;X_f)/\R\times\cM^1(B,A^+;X_f)/\R$
for any critical point $[B]\in\cR_f\setminus[0]$,
which are exactly what is counted on the left hand side
of~(\ref{eq:d20}). By the {\it (Catenation}) 
axiom in Section~\ref{sec:orient} 
the signs agree with the orientation of the 
boundary of $\cM^2(A^-,A^+;X_f)/\R$.
Hence the sum must be zero and this 
proves $\p\circ\p=0$.
\end{proof}

\begin{proof}[Proof of Theorem~\ref{thm:floer2}.]
That the Floer homology groups are independent of the 
choice of the system of coherent orientations is obvious;
two such systems give rise to isomorphic boundary 
operators via a sign change isomorphism (with $\pm1$
on the diagonal). To prove the independence of 
metric and perturbation, we fix two Riemannian metrics 
$g^\pm$ on $Y$ and two sets of regular perturbation 
data $(\gamma^\pm,f^\pm)$.
We will construct a chain map
from $\CF(Y,\cL;g^-,f^-)$ to $\CF(Y,\cL;g^+,f^+)$
following the familiar pattern. 
As in the closed case we choose a metric $\tilde g$ 
on $\R\times Y$ such that $\tilde g = g^\pm$ for $\pm s$ 
sufficiently large. However, unlike the closed case
this metric cannot necessarily be chosen
in split form since it is required to be
compatible with the boundary space-time splitting 
in the sense of Definition~\ref{def:4g}
(see Example~\ref{ex:RY} or~\cite[Example~1.4]{W elliptic}).
Next we choose a holonomy perturbation 
$\tilde X:\cA(\R\times Y)\to\Om^2(\R\times Y,\cg)$
of the form 
$\tilde X=\beta X_{f^-} + (1-\beta)X_{f^+} + X_{f'}$
for some cutoff function $\beta\in\cC^\infty(\R,[0,1])$
and a further holonomy perturbation $X_{f'}$ as in 
Definition~\ref{def:data}. This uses thickened loops 
$\gamma'_i:S^1\times B^3\hookrightarrow \R\times{\rm int}(Y)$
in a compact part of $\R\times Y$, so that we have
$\tilde X=X_{f^\pm}$ for $\pm s$ sufficiently large.
This perturbation is still gauge equivariant
but no longer translation invariant.
We use these interpolation data to set up 
the $4$-dimensional version of the perturbed 
anti-self-duality equation on $\R\times Y$
as described in Section~\ref{sec:fredholm}.
For critical points $A^\pm\in\cR_{f^\pm}\setminus[0]$
from the two Floer chain complexes we consider the
space of generalized Floer trajectories
$$
\widetilde\cM(A^-,A^+;\tilde g,\tilde X)
:= \left\{ \begin{array}{l}
\Xi=A + \Phi\ds \\
\in \cA(\R\times Y)
\end{array}
\left|
\begin{array}{l}
F_\Xi + \tilde X(\Xi) + *_{\scriptscriptstyle\tilde g} 
(F_\Xi + \tilde X(\Xi)) = 0 \\
A(s)|_\Sigma\in\cL \quad\forall s\in\R \\
\lim_{s\to\pm\infty}A(s)= A^\pm \\
\Phi|_{\{|s|\geq 1\}}\equiv 0 \\
E_f(\Xi)<\infty 
\end{array}\right.\right\} .
$$
Here $*_{\scriptscriptstyle\tilde g}$ denotes the Hodge operator 
on $\R\times Y$ with respect to the metric $\tilde g$.
This space is invariant under the gauge group $\cG(A^-,A^+)$
as before, and if the perturbation $\tilde X$ is regular, 
then the quotient $\widetilde\cM(A^-,A^+;\tilde g, \tilde X)/\cG(A^+,A^-)$
will be a smooth manifold whose local dimension near $[\A]$ is given
by the Fredholm index 
${\delta(\A) \equiv \mu_{f^-}(A^-)-\mu_{f^+}(A^+)}$ (modulo $8$).
By transversality arguments similar to Section~\ref{sec:moduli}
we can find a perturbation $X_{f'}$ (and thus $\tilde X$) such 
that the linearized operators of index less than or equal to $7$ 
are indeed surjective. Thus we obtain smooth $k$-dimensional 
moduli spaces
$$
\cM^k(A^-,A^+;\tilde g,\tilde X)
:=\bigl\{[\A]\in\widetilde\cM(A^-,A^+;\tilde g,\tilde X)/\cG(A^-,A^+) \st
\delta(\A)=k \bigr\}
$$
for $k\le7$.  The $0$-dimensional moduli spaces are compact 
by the same analysis as in Section~\ref{sec:compact}.
Namely, the main component will converge to a new solution 
without time-shift; energy cannot be lost by bubbling or 
by shift to $\pm\infty$ since the remaining solution would 
have negative index.  So -- again using the orientations from 
Section~\ref{sec:orient} -- we can define a homomorphism
$$
\Phi : \CF(Y,\cL;g^-,f^-)\to \CF(Y,\cL;g^+,f^+) ,
$$
which preserves the grading and is given by
$$
\Phi \langle A^-\rangle := 
\sum_{[A^+]\in \cR_{f^+}\setminus[0]} 
\left(\sum_{\A\in\cM^0(A^-,A^+;\tilde g,\tilde X)}\nu(\A)\right) 
\langle A^+ \rangle.
$$
This time the linearized operator is bijective,
so $\det(\cD_\A)$ is canonically isomorphic to $\R$,
and the sign $\nu(\A)=\pm1$ is obtained by comparing the
coherent orientation $o_\A$ with the standard
orientation of $\R$. 

As in the standard theory 
there are three identities to verify
(e.g.~\cite[Section~3.2]{Sal}).
First, we must prove that $\Phi$ is a chain map, i.e.
\begin{equation}\label{eq:phid}
\p^+\circ\Phi=\Phi\circ\p^-.
\end{equation}
This is proved just like the formula $\p\circ\p=0$
in Theorem~\ref{thm:floer1}. In this case the relevant
$1$-manifold is the moduli space 
$\cM^1(A^-,A^+;\tilde g,\tilde X)$.
A compactness and gluing theory similar to 
Corollary~\ref{cor:glue} identifies the ends of this 
moduli space with the pairs of trajectories in 
$\cM^1(A^-,B^-;g^-,X_{f^-})\times\cM^0(B^-,A^+;\tilde g,\tilde X)$
for $[B^-]\in\cR_{f^-}\setminus[0]$ and in
$\cM^0(A^-,B^+;\tilde g,\tilde X)\times\cM^1(B^+,A^+;g^+,X_{f^+})$
for $[B^+]\in\cR_{f^+}\setminus[0]$.
Summing over these oriented ends of a $1$-manifold then proves
that $\Phi$ satisfies~(\ref{eq:phid})
and hence descends to a morphism on Floer homology.

Second, we must prove that the induced map on homology
is independent of the choices. Given two such maps 
$
\Phi_0,\Phi_1:\CF(Y,\cL;g^-,f^-)\to \CF(Y,\cL;g^+,f^+)
$
associated to $(\tilde g_0,\tilde X_0)$ and $(\tilde g_1,\tilde X_1)$ 
we must find a chain homotopy equivalence
$
H:\CF(Y,\cL;g^-,f^-)\to \CF(Y,\cL;g^+,f^+) 
$
satisfying 
\begin{equation}\label{eq:t}
\Phi_1-\Phi_0 = \p^+\circ H + H\circ\p^-.
\end{equation}
To construct $H$ we choose a $1$-parameter family 
$\{\tilde g_\lambda,\tilde X_\lambda\}_{0\le\lambda\le1}$
of interpolating pairs of metric and perturbation.
By Lemma~\ref{le:metric} the metrics can be interpolated within 
the space of metrics that are equal to $g_\pm$ over the ends and are 
compatible with the space-time splitting of the boundary.
The perturbations~$\tilde X_\lambda$ can be chosen as 
convex combinations.  We then add further compactly supported 
holonomy perturbations for $0<\lambda<1$ to achieve transversality 
of the parametrized moduli spaces
$$
\cM^k\bigl(A^-,A^+;\{\tilde g_\lambda,\tilde X_\lambda\}\bigr)
:= \bigl\{ (\lambda,[\A]) \st 
[\A]\in\cM^k(A^-,A^+;\tilde g_\lambda,\tilde X_\lambda) \bigr\} .
$$
For $k=-1$ these are compact oriented $0$-manifolds
which we use to define $H$:
$$
H \langle A^-\rangle 
:= \sum_{[A^+]\in \cR_{f_+}\setminus[0]}\left(
\sum_{(\lambda,\A)\in
\cM^{-1}(A^-,A^+;\{\tilde g_\lambda,\tilde X_\lambda\})}
\nu(\lambda,\A)\right) 
\langle A^+ \rangle.
$$
The linearized operator has a
$1$-dimensional cokernel which projects isomorphically 
to $\R$ and $\nu(\lambda,\A)$ is the sign 
of this projection. To prove~(\ref{eq:t})
one studies the $1$-dimensional moduli space
$\cM^0(A^-,A^+;\{\tilde g_\lambda,\tilde X_\lambda\})$
in the usual fashion with the contributions of 
$\Phi_0$ corresponding to the boundary at $\lambda=0$,
the contributions of $\Phi_1$ to the boundary at $\lambda=1$,
and the contributions on the right in~(\ref{eq:t})
to the noncompact ends with $0<\lambda<1$. 
These ends have either the form of a pair
in 
$
\cM^{-1}(A^-,B^+;\{\tilde g_\lambda,\tilde X_\lambda\})
\times \cM^1(B^+,A^+;g^+,X_{f^+})/\R
$
with $[B^+]\in\cR_{f^+}\setminus\{[0]\}$ or in
$
\cM^1(A^-,B^-;g^-,X_{f^-})/\R 
\times\cM^{-1}(B^-,A^+;\{\tilde g_\lambda,\tilde X_\lambda\})
$
with $[B^-]\in\cR_{f^-}\setminus\{[0]\}$.
Counting all the ends and boundary points
with appropriate signs proves that $H$ 
satisfies~(\ref{eq:t}).  

Third, we must establish the composition rule 
in~(\ref{eq:abc}) for three sets of regular data
$(g^\alpha,f^\alpha)$, $(g^\beta,f^\beta)$, $(g^\gamma,f^\gamma)$.
We choose regular interpolating metrics and perturbations
to define $\Phi^{\beta\alpha}$ and $\Phi^{\gamma\beta}$
on the chain level. The catenation 
(with gluing parameter $T$) of these data gives 
rise to a regular interpolation from $(g^\alpha,f^\alpha)$
to $(g^\gamma,f^\gamma)$ for $T$ sufficiently large.
The resulting morphism $\Phi_T^{\gamma\beta}$
will then, for large $T$, agree with 
$\Phi^{\gamma\beta}\circ\Phi^{\beta\alpha}$ 
on the chain level.  This follows from a gluing theorem
as in Section~\ref{sec:gluing} and compactness arguments 
as in Theorem~\ref{thm:compact2} and Corollary~\ref{cor:glue}. 
In particular, the breaking of connecting trajectories in the limit $T\to\infty$ 
at the zero connection is excluded since the stabilizer 
$\cG_0\subset\cG(Y)$ adds $3$ to the index of the glued connection 
(compare with Remark~\ref{rmk:redglue} or use index
inequalities as in Corollary~\ref{cor:compact}.).
Again, the orientations are compatible with the gluing by the
{\it (Catenation)} axiom. The upshot is that, for suitable
choices of interpolating data, equation~(\ref{eq:abc})
already holds on the chain level. 

Once these three relations have been established
one just needs to observe that $\Phi^{\alpha\alpha}$
is the identity on the chain level for the obvious 
product metric and perturbation on $\R\times Y$.
It follows that each $\Phi$ induces an isomorphism
on Floer homology.  This proves Theorem~\ref{thm:floer2}.
\end{proof}


\appendix

\section{The spectral flow}\label{app:spec}

In this appendix we adapt the results of~\cite{RS}
to families of self-adjoint operators with varying domains.
Similar results have appeared in various forms
(see~\cite{BossZhu,KirkLesch,BBC}).

Let $H$ be a separable real Hilbert space.
Throughout we identify $H$ with its dual space. 
We consider a family of bounded linear operators
$$
A(s):W(s)\to H
$$
indexed by $s\in\R$.  Here $W(s)$ is a Hilbert space 
equipped with a compact inclusion
$
W(s)\subset H
$
with a dense image. 
We formulate conditions under which the
unbounded operator 
$$
\cD := \pd_s + A
$$ 
on $L^2(\R,H)$ is Fredholm and its index is
the spectral flow of the operator family $s\mapsto A(s)$.
In contrast to \cite{RS} the domain of $A(s)$ varies with
$s\in\R$. Our axioms give rise to an isomorphic
family of operators with constant domain but which are self--adjoint
with respect to inner products which vary with $s\in\R$.
More precisely, we assume that the disjoint union $\bigsqcup_{s\in\R}W(s)$ 
is a Hilbert space subbundle of $\R\times H$ in the following 
sense.
\begin{description}
\item[(W1)]
There is a dense subspace $W_0\subset H$ with a compact 
inclusion and a family of isomorphisms $Q(s):H\to H$ such that
$
Q(s)W_0=W(s)
$
for every $s\in\R$.  
\item[(W2)]
The map $Q:\R\to\cL(H)$ is continuously differentiable 
in the weak operator topology and there is a $c_0>0$ 
such that, for all $s\in\R$ and $\xi\in W_0$,
$$
c_0^{-1}\|\xi\|_{W_0}\le \|Q(s)\xi\|_{W(s)} \le c_0 \|\xi\|_{W_0},
$$
$$
\| Q(s)\xi\|_H + \|\pd_s Q(s)\xi\|_H \le c_0 \|\xi\|_H.
$$
\item[(W3)]
There exist Hilbert space isomorphisms $Q^\pm\in\cL(H)$ 
such that
$$
\lim_{s\to\pm\infty}\|Q(s)-Q^\pm\|_{\cL(H)} = 0. 
$$
\end{description}
Two trivializations $Q_1,Q_2:\R\to\cL(H)$ satisfying~(W1-3)
with $W_{01},W_{02}$, respectively, are called {\bf equivalent}
if there is a family of Hilbert space isomorphisms
$\Phi(s)\in\cL(H)$ such that
$$
\Phi(s)W_{01}=W_{02},\qquad Q_2(s)\Phi(s)=Q_1(s)
$$
for every $s$, the map $\Phi:\R\to\cL(H)$ is continuously 
differentiable in the weak operator topology, the map 
$\Phi:\R\to\cL(W_{01},W_{02})$ is continuous in the norm 
topology, $\sup_{s\in\R}\|\pd_s\Phi(s)\|_{\cL(H)}<\infty$,
and there exist Hilbert space isomorphisms 
$\Phi^\pm\in\cL(H)\cap\cL(W_{01},W_{02})$ such that 
$$
\lim_{s\to\pm\infty}\|\Phi(s)-\Phi^\pm\|_{\cL(H)} = 0.
$$

\begin{remark}\label{rmk:W}\rm
To verify~(W1-3) it suffices to construct local 
trivializations on a finite cover $\R=\bigcup U_\alpha$
that satisfy these conditions (where condition~(W3) is only 
required near the ends) and that are equivalent 
over the intersections $U_\alpha\cap U_\beta$. 
\end{remark}

We now impose the following conditions on the operator family~$A$.
Again, it suffices to verify these in the local trivializations 
of Remark~\ref{rmk:W}.

\begin{description}
\item[(A1)]
The operators $A(s)$ are uniformly self-adjoint.  
This means that for each $s\in\R$ the operator $A(s)$ 
when considered as an unbounded operator on $H$ 
with $\mathrm{dom}\,A(s) = W(s)$ is self-adjoint 
and that there is a constant $c_1$
such that
$$
\|\xi\|_{W(s)}^2 \le 
c_1\left(\|A(s)\xi\|_H^2 + \|\xi\|_H^2\right).
$$
for every $s\in\R$ and every $\xi\in W(s)$.
\item[(A2)]
The map $B:=Q^{-1}AQ:\R\to\cL(W_0,H)$ is continuously differentiable 
in the weak operator topology and there exists a constant
$c_2>0$ such that
$$
\|B(s)\xi\|_H + \|\pd_s B(s)\xi\|_H \le c_2 \|\xi\|_{W_0}.
$$
for every $s\in\R$ and every $\xi\in W_0$.
\item[(A3)]
There are invertible operators $B^\pm\in\cL(W_0,H)$ 
such that
$$
\lim_{s\to\pm\infty} \|B(s) - B^\pm\|_{\cL(W_0,H)} = 0.
$$
\end{description}

Given a differentiable curve $\xi:\R\to H$ with
$\xi(s)\in W(s)$ for all $s\in\R$ we define
$\cD\xi:\R\to H$ by
$$
(\cD\xi)(s) = \pd_s{\xi}(s)+A(s)\xi(s) .
$$
This map extends to a bounded linear operator 
$$
\cD:W^{1,2}(\R,H)\cap L^2(\R,W)\to L^2(\R,H).
$$ 
Here 
$
L^2(\R,W):=\bigl\{ Q\eta_0 \st \eta_0\in L^2(\R,W_0) \bigr\}
$
is a Hilbert space with the norm
$$
\Norm{\eta}_{L^2(\R,W)}^2 
= \int_{-\infty}^\infty \Norm{\eta(s)}_{W(s)}^2 \rd s .
$$
By (W2) this norm is equivalent to the norm on $L^2(\R,W_0)$ 
under the isomorphism $\eta\mapsto Q^{-1}\eta$.
We will prove the following 
estimate, regularity, and index identity.

\begin{lemma}\label{lem:spec est}
There exist constants $c$ and $T$ such that
$$
\int_{-\infty}^\infty 
\Bigl(\Norm{\pd_s\xi(s)}_H^2 + \Norm{\xi(s)}_{W(s)}^2  \Bigr)\rd s
\leq c^2 \biggl( \int_{-\infty}^\infty \Norm{\cD\xi (s)}_H^2 \rd s
+ \int_{-T}^T \Norm{\xi(s)}_H^2 \rd s \biggr)
$$
for every $\xi\in W^{1,2}(\R,H)\cap L^2(\R,W)$.
\end{lemma}

\begin{thm}\label{thm:spec reg}
Suppose that $\xi,\eta\in L^2(\R,H)$ satisfy
$$
\int_{-\infty}^\infty \Bigl( \inner{\pd_s\phi(s)-A(s)\phi(s)}{\xi(s)} 
+ \inner{\phi(s)}{\eta(s)} \Bigr) \rd s = 0
$$
for every test function $\phi:\R\to H$ such that 
$Q^{-1}\phi\in\cC^1_0(\R,W_0)$.   Then 
$$
\xi\in W^{1,2}(\R,H)\cap L^2(\R,W),\qquad
\cD\xi=\eta.
$$
\end{thm}

\begin{thm}\label{thm:Fredholm}
The operator $\cD$ is Fredholm and its index is equal to the
upward spectral flow of the operator family $s\mapsto A(s)$.
\end{thm}

As in the case of constant domain the spectral flow can be
defined as the sum of the crossing indices
\begin{equation}\label{eq:spec flow}
\mu_\spec(A) := \sum_s {\rm sign}\; \Gamma(A,s).
\end{equation}
In the present case the crossing form 
$
\Gamma(A,s):\ker A(s)\to\R
$ 
is defined by
$$
\Gamma(A,s)(\xi):= \left.\frac{d}{dt}\right|_{t=0}
\inner{\xi(t)}{A(s+t)\xi(t)},
$$
where $\xi(t)\in W(s+t)$ is chosen such $\xi(0)=\xi$ and the path
$t\mapsto A(s+t)\xi(t)\in H$ is differentiable
(for example $\xi(t):=Q(s+t)Q(s)^{-1}\xi$);
the value of the crossing form at $\xi$ is independent 
of the choice of the path $t\mapsto\xi(t)$. 
We assume that the crossings are all regular, i.e.\ 
$\Gamma(A,s)$ is nondegenerate for every $s\in\R$ with
$\ker A(s)\neq\{0\}$. Under this assumption the sum in
(\ref{eq:spec flow}) is finite.

Two operator families $A_1(s):W_1(s)\to H$ and $A_2(s):W_2(s)\to H$
with the same endpoints $A^\pm$ are called {\bf homotopic} if they can be 
connected by an operator family
$A_\lambda(s):W_\lambda(s)\to H$, $1\le\lambda\le 2$,
with the following properties.  There is a family of Hilbert space 
isomorphisms $Q_\lambda(s):H\to H$ that is continuously 
differentiable in $\lambda$ and $s$ with respect to the 
weak operator topology and satisfies $Q_\lambda(s)W_0=W_\lambda(s)$
as well as conditions~(W2-3) uniformly in $\lambda$. 
Moreover $A_\lambda(s)$ satisfies~(A1-3) with constants independent
of $\lambda$ and the map $[1,2]\times\R\to\cL(W_0,H):(\lambda,s)\mapsto
Q_\lambda(s)^{-1}A_\lambda(s)Q_\lambda(s)$ is continuously differentiable
in the weak operator topology. 

\smallbreak

The spectral flow has the following properties: \label{axioms}
\begin{description}
\item[(Homotopy)]
The spectral flow is invariant under homotopy. 
\item[(Constant)] 
If $W(s)$ and $A(s)$ are independent of $s\in\R$ then $\mu_\spec(A)=0$. 
\item[(Direct sum)]
The spectral flow of a direct sum of two operator families $A$ and $B$
is the sum of their spectral flows, i.e.
$$
\mu_\spec(A\oplus B)=\mu_\spec(A)+\mu_\spec(B).
$$
\item[(Catenation)]
The spectral flow of the catenation of two operator families
$A_{01}$ from $A_0$ to $A_1$ and $A_{12}$ from $A_1$
to $A_2$ is the sum of their spectral flows, i.e.
$$
\mu_\spec(A_{01}\#A_{12})=\mu_\spec(A_{01})_+\mu_\spec(A_{12}).
$$
\item[(Normalization)]
For $W=H=\R$, $A(s)=\arctan(s)$ we have ${\mu_\spec(A)=1}$.
\end{description}
The spectral flow is uniquely determined by the homotopy, constant,
direct sum, and normalization axioms.  The proof is the same
as that of~\cite[Theorem~4.23]{RS} and will be omitted.

\begin{proof}[Proof of Lemma~\ref{lem:spec est}.]
The proof is analogous to that of \cite[Lemma~3.9]{RS}.
The only difference is in the first step where we 
prove the estimate with $T=\infty$.
For every $\xi:\R\to H$ such that $\eta:=Q^{-1}\xi\in\cC^1_0(\R,W_0)$
we have
$$
\int_{-\infty}^\infty \|\cD\xi\|_H^2 \rd s
= \int_{-\infty}^\infty \Bigl( \|\pd_s\xi\|_H^2 + \|A\xi\|_H^2 
 + 2\inner{\pd_s\xi}{A\xi} \Bigr) \rd s .
$$
The last summand can be estimated by
\begin{align*}
&2\int_{-\infty}^\infty \inner{\pd_s\xi}{A\xi} \rd s \\
&=
\int_{-\infty}^\infty \Bigl( 2\inner{(\pd_s Q)\eta}{AQ\eta}  
+ \inner{Q\pd_s\eta}{AQ\eta} 
+ \inner{Q\eta}{AQ\pd_s\eta} \Bigr)\rd s \\
&=
\int_{-\infty}^\infty \Bigl( \inner{(\pd_s Q)\eta}{QB\eta}  
- \inner{Q\eta}{(\pd_s Q)B\eta} 
- \inner{Q\eta}{Q(\pd_s B)\eta} \Bigr)\rd s \\
&\leq
3 c_0^2 c_2\int_{-\infty}^\infty  \|\eta\|_H \|\eta\|_{W_0} \rd s \\
&\leq
c  \|\xi\|_{L^2(\R,H)} \|\xi\|_{L^2(\R,W)} 
\end{align*}
with $c:=3 c_0^4 c_2$. Here we used partial integration 
and the identity $AQ=QB$. Now use (A1) to obtain
\begin{align*}
& \|\cD\xi\|_{L^2(\R,H)}^2 \\
&\ge 
\|\pd_s\xi\|_{L^2(\R,H)}^2 
+ c_1^{-1}\|\xi\|_{L^2(\R,W)}^2 - \|\xi\|_{L^2(\R,H)}^2
   - c \|\xi\|_{L^2(\R,H)} \,\|\xi\|_{L^2(\R,W)} \\
&\ge 
   \|\pd_s\xi\|_{L^2(\R,H)}^2 + (2c_1)^{-1}\|\xi\|_{L^2(\R,W)}^2
   - \left(1+\tfrac12 c^2 c_1\right) \|\xi\|_{L^2(\R,H)}^2 .
\end{align*}
This proves the estimate for $T=\infty$.
\end{proof}

\begin{proof}[Proof of Theorem~\ref{thm:spec reg}.]
We follow the line of argument in \cite[Thm.~3.10]{RS}.

\medskip

\noindent{\bf Step 1:}
{\it 
Define $\xi_0,\eta_0\in L^2(\R,H)$ by
$$
\xi_0(s):= Q(s)^*\xi(s) , \qquad 
\eta_0(s):= Q(s)^*\eta(s) + (\pd_s Q(s)^*) \xi(s) .
$$
Then $\xi_0\in W^{1,2}(\R,W_0^*)$ and
\begin{equation}\label{eq:weak}
\pd_s\xi_0(s) = - B(s)^*\xi_0(s)  + \eta_0(s) .
\end{equation}
}
To see this we calculate for $\phi_0\in\cC^\infty_0(\R,W_0)$
\begin{align*}
\int_{-\infty}^\infty \inner{\pd_s\phi_0}{\xi_0}_H\,\ds
&=
\int_{-\infty}^\infty \inner{\pd_s(Q \phi_0) 
- (\pd_s Q)\phi_0}{\xi}_H\,\ds \\
&=
\int_{-\infty}^\infty \Bigl( 
\inner{AQ\phi_0}{\xi}_H 
- \inner{Q\phi_0}{\eta}_H 
- \inner{(\pd_s Q)\phi_0}{\xi}_H \Bigr)\,\ds \\
&=
\int_{-\infty}^\infty \Bigl( 
\inner{\phi_0}{B^*\xi_0 - \eta_0 }_{W_0,W_0^*} \Bigr)\,\ds .
\end{align*}
Here the self-adjoint operator $A(s)$ 
extends to an operator in $\cL(H,W(s)^*)$
which we also denote by $A(s)$.
We denote the dual of the trivialization $Q(s)$
by ${Q(s)^*\in\cL(H)}$, which extends 
to an isomorphism $W(s)^*\to W_0^*$.
With this we can write $B^*=Q^*A(Q^*)^{-1}$ for 
the dual operator family of $B=Q^{-1}AQ$, 
which is continuously differentiable in $\cL(H,W_0^*)$
with a uniform estimate dual to that in (A2).
So we have $B^*\xi_0 - \eta_0\in L^2(\R,W_0^*)$, 
and since the derivatives of test functions $\phi_0$ 
are dense in $L^2(\R,W_0)$ this implies Step~1.

\medskip\noindent{\bf Step 2.}
{\it
Suppose that $\xi$ and $\eta$ are supported 
in an interval $I$ such that for all $s\in I$
the operator $B(s):W_0\to H$ is bijective and
satisfies a uniform estimate 
$$
        \|B(s)^{-1}\|_{\cL(H,W_0)}\le c .
$$
Fix a smooth function ${\rho:\R\to[0,\infty)}$ with
support in $(-1,1)$ and $\int\rho = 1$ and denote by
$\rho_\delta(s) = \delta^{-1}\rho(\delta^{-1}s)$
for $\delta>0$ the standard mollifier.
Then we find a constant $C$ such that 
$\rho_\delta * (Q^{-1}\xi) \in W^{1,2}(\R,H)\cap L^2(\R,W_0)$ 
for all $\delta > 0$ and}
$$
\|\cD Q (\rho_\delta * (Q^{-1}\xi))\|_{L^2(\R,H)} \le C .
$$

\medskip\noindent
Multiply equation~(\ref{eq:weak}) by $(B^*)^{-1}$ to obtain
$\xi_0 = (B^*)^{-1}\bigl(\eta_0 - \pd_s\xi_0 \bigr)$ 
and note that $Q^{-1}\xi = (Q^*Q)^{-1}\xi_0$.
Then convolution gives
\begin{align*}
&\rho_\delta * (Q^{-1}\xi) \\
&=
\rho_\delta * \bigl( \pd_s\bigl((B^*Q^*Q)^{-1}\bigr) \xi_0
                     + (B^*Q^*Q)^{-1} \eta_0 \bigr)  
- \dot\rho_\delta * \bigl( (B^*Q^*Q)^{-1} \xi_0 \bigr) \\
&=
\rho_\delta * \zeta_0
- \dot\rho_\delta * \bigl( (Q^*QB)^{-1} \xi_0 \bigr)
\end{align*}
with 
$\zeta_0 =  \pd_s\bigl((Q^*Q B)^{-1}\bigr) \xi_0 + (Q^*Q B)^{-1} \eta_0
\in L^2(\R,W_0)$. This takes values in $W_0$ since 
\begin{equation}\label{eq:BQ}
(B^*Q^*Q)^{-1} = Q^{-1}A^{-1}(Q^*)^{-1} = B^{-1} (Q^*Q)^{-1}
\end{equation}
and its derivative are uniformly bounded in $\cL(H,W_0)$.

So, after convolution, $Q\bigl(\rho_\delta * (Q^{-1}\xi)\bigr)$ 
lies in the domain of $\cD$ and
\begin{align*}
& Q^{-1}\cD Q \bigl( \rho_\delta * (Q^{-1}\xi) \bigr) \\
&= \dot\rho_\delta * (Q^{-1}\xi) 
+ Q^{-1}(\pd_s Q) \bigl( \rho_\delta * (Q^{-1}\xi) \bigr)
+ B \bigl(\rho_\delta * (Q^{-1}\xi)\bigr) \\
&= B \bigl( B^{-1} \bigl(\dot\rho_\delta * (Q^{-1}\xi) \bigr) -
\dot\rho_\delta * \bigl( (Q^*QB)^{-1} \xi_0 \bigr) \bigr) \\
&\quad 
+ Q^{-1}(\pd_s Q) \bigl( \rho_\delta * (Q^{-1}\xi) \bigr)
+ B \bigl( \rho_\delta * \zeta_0 \bigr)
\end{align*}
The second line is unifomly bounded in $L^2(\R,H)$.
For the first term we have
\begin{align*}
&\int_{-\infty}^\infty 
\bigl\|  B^{-1} \bigl(\dot\rho_\delta * (Q^{-1}\xi) \bigr)(s) -
\dot\rho_\delta * \bigl( (Q^*QB)^{-1} \xi_0 \bigr)(s) \bigr\|_{W_0} \,\ds \\
&= 
\int_{-\infty}^\infty 
\biggl\| \int_{s-\delta}^{s+\delta} \tfrac 1\delta \dot\rho(\tfrac{t-s}\delta) 
\frac{B(s)^{-1} - B(t)^{-1}}\delta 
\bigl(Q(t)^*Q(t)\bigr)^{-1} \xi_0(t) \,\rd t \Bigr\|_{W_0} \,\ds \\
&\leq
C \int_{-\infty}^\infty\int_{-\infty}^\infty
\left|\tfrac{1}{\delta}\dot\rho(\tfrac{t-s}{\delta})\right|\,
     \left\|\xi_0(t)\right\|_H\,\dt\,\ds \\
&\leq
C \|\dot\rho\|_{L^1(\R)} \int_{-\infty}^\infty \left\|\xi_0(t)\right\|_H\,\dt .
\end{align*}
Here the constant $C$ contains a uniform bound for 
$\pd_s B^{-1}=-B^{-1}(\pd_s B) B^{-1}$ on $I$.
This proves Step 2.

\medskip\noindent{\bf Step 3.}  
{\it $\xi\in  W^{1,2}(\R,H)\cap L^2(\R,W)$ and $\cD\xi=\eta$.}

\medskip\noindent
Under the assumptions of Step~2 it follows from Lemma~\ref{lem:spec est} 
that $\rho_\delta * (Q^{-1}\xi)$ is uniformly bounded in 
$\cW:=L^2(\R,W_0)\cap W^{1,2}(\R,H)$ for all $\delta>0$.
So there is a sequence $\delta_\nu\to 0$ such that 
$\rho_{\delta_\nu} * (Q^{-1}\xi)$ converges weakly in $\cW$.
The limit has to coincide with the strong $L^2(\R,H)$-limit $Q^{-1}\xi$.
Thus we have $\xi\in L^2(\R,W)\cap W^{1,2}(\R,H)$.  
Now it follows from~(\ref{eq:weak}) and (\ref{eq:BQ}) that
\begin{align*}
\cD\xi
&= (Q^*)^{-1} \pd_s\xi_0 - (Q^*)^{-1}(\pd_s Q^*)(Q^*)^{-1} \xi_0 
+ B (Q^*)^{-1} \xi_0 \\
&= \eta  - (Q^*)^{-1} B^*\xi_0 + B (Q^*)^{-1} \xi_0 
\;=\;\eta .
\end{align*}
This proves the theorem under the assumption that $\xi$
and $\eta$ are supported in an interval on which $B$ is bijective.
In general, one can cover the real axis by finitely many open intervals 
on which $\lambda\one+B(s):W_0\to H$ has uniformly bounded inverses 
for some $\lambda\in\R$.
Then one can use a partition of unity argument to deduce the regularity 
and equation for $\xi$ on each interval.
\end{proof}

\begin{proof}[Sketch of proof of Theorem~\ref{thm:Fredholm}.]
By Lemma~\ref{lem:spec est} the operator $\cD$ has a finite 
dimensional kernel and a closed image.  
By Theorem~\ref{thm:spec reg} the cokernel of $\cD$ is the kernel 
of the operator $\cD'$ with $A$ replaced by $-A$.
Hence the cokernel of $\cD$ is also finite dimensional and thus 
$\cD$ is Fredholm.

To prove the index identity one verifies as in~\cite[Theorem~4.1]{RS} 
that the Fredholm index satisfies the axioms on page~\pageref{axioms},
which characterize the spectral flow.
For the homotopy and the direct sum property one can extend the 
proofs in~\cite{RS} without difficulty to nonconstant domains; 
the constant and normalization properties are immediate since they
only refer to constant domains.
\end{proof}

We conclude this appendix with a version of the index identity
for twisted loops of self-adjoint operators.

\begin{thm}\label{thm:indexS1}
Let $A(s):W(s)\to H$ be an operator family
that satisfies the conditions~$(W1-2)$, $(A1-2)$, and
$$
W(s+1)=Q^{-1}W(s),\qquad A(s+1) = Q^{-1}A(s)Q
$$
for every $s\in\R$ and a suitable Hilbert space 
isomorphism $Q:H\to H$.  Then $A$ induces a Fredholm
operator
$
\cD = \p_s+A:\cW\to\cH,
$
where
$$
\cH:=\left\{\xi\in L^2_\loc(\R,H)\,\big|\,\xi(s+1)=Q^{-1}\xi(s)\right\}, 
$$
$$
\cW:=\bigl\{\xi\in L^2_\loc(\R,W)\cap W^{1,2}_\loc(\R,H)\,\big|\,
\xi(s+1)=Q^{-1}\xi(s)\bigr\}.
$$
Its Fredholm index is equal to the upward spectral flow 
of the operator family $A$ on a fundamental domain $[s_0,s_0+1]$.
\end{thm}

\begin{proof}
The Fredholm property follows from Lemma~\ref{lem:spec est} 
and Theorem~\ref{thm:spec reg}. 
The proof of the index formula can be reduced to 
Theorem~\ref{thm:Fredholm} by using the homotopy invariance 
of spectral flow and Fredholm index, stretching the fundamental domain,
and comparing kernel and cokernel with a corresponding 
operator over $\R$ via a gluing argument.  We omit the details.
For a version of the relevant linear gluing theorem 
see~\cite[Propositions~3.8, (3.2)]{Donaldson book}.
\end{proof}


\section{The Gelfand--Robbin quotient} \label{app:GR}

In this appendix we collect various results on the Gelfand-Robbin 
quotient, associated to an unbounded symmetric operator, 
whose Lagrangian subspaces correspond to self-adjoint extensions.  
Related results concerning the spectral flow and the Maslov index 
for Fredholm Lagrangian pairs can be found in various places
(for example~\cite{BBF,BossZhu,RS}). However, the existing literature 
on this subject does not seem to fully cover what is needed in this paper.

Let $H$ be a Hilbert space and $D:\mathrm{dom}\,D\to H$
be an injective, symmetric, but not necessarily self-adjoint, 
operator with a dense domain and a closed image.  
Then the domain of the adjoint operator $D^*:\mathrm{dom}\,D^*\to H$ 
contains the domain of $D$ and the restriction of $D^*$ to the domain
of $D$ agrees with $D$.  The Gelfand--Robbin quotient 
$$
V := \mathrm{dom}\,D^*/\mathrm{dom}\,D
$$
carries a natural symplectic form 
$$
\om([x],[y]) := \inner{D^*x}{y} - \inner{x}{D^*y}.
$$
The Lagrangian subspaces $\Lambda\subset V$ are in one-to-one
correspondence to self-adjoint extensions $D_\Lambda$ of $D$ with 
$$
\mathrm{dom}\,D_\Lambda:=
\left\{x\in\mathrm{dom}\,D^*\,|\,[x]\in\Lambda\right\}.
$$
Moreover, the kernel of $D^*$ determines a Lagrangian subspace
\begin{equation}\label{eq:ker}
\Lambda_0:=\left\{[x]\in V\,|\,x\in\mathrm{dom}\,D^*,\,D^*x=0\right\} .
\end{equation}
The operator $D_\Lambda$ is bijective if and only if 
$V=\Lambda_0\oplus\Lambda$. (See Lemma~\ref{lem:Fredtrip} below.)

The domain of $D^*$ is a Hilbert space with the graph 
inner product 
$$
\Inner{x}{y}_{D^*} := \Inner{x}{y}_H + \Inner{D^*x}{D^*y}_{H}.
$$
The domain of $D$ is a closed subspace because $D$
has a closed graph. Hence both $\dom\,D$ and the quotient space 
$V=\mathrm{dom}\,D^*/\mathrm{dom}\,D$ inherit a Hilbert space 
structure from $\mathrm{dom}\,D^*$. 
One can now check (using the next remark) that 
$(V,\omega)$ is a {\em symplectic Hilbert space}
in the sense that the symplectic form is bounded and 
the linear map $V\to V^*:v\mapsto I_\omega(v):=\om(v,\cdot)$ 
is an isomorphism. If $\Lambda\subset V$ is a Lagrangian subspace,
i.e.~the annihilator $\Lambda^\perp\subset V^*$ is given by
$\Lambda^\perp = I_\omega(\Lambda)$,  then $\Lambda$ 
is closed and hence inherits a Hilbert space structure from $V$.

\begin{rmk}\label{rmk:V}\rm
{\bf (i)}
The graph norm on $\dom D$ is equivalent to the norm
$$
\Inner{x}{y}_{D} := \Inner{Dx}{Dy}_{H}.
$$
because $D$ is injective and has a closed image.

\smallskip\noindent{\bf (ii)}
It is sometimes convenient to identify 
the Gelfand--Robbin quotient $V=\dom D^*/\dom D$ 
with the orthogonal complement
$$
V = (\dom D)^\perp = \left\{x\in\dom D^*\,|\,
D^*x\in\dom D^*,\,D^*D^*x+x=0\right\}.
$$
The orthogonal projection of $\dom D^*$ onto $V$
along $\dom D$ is given by
$$
\dom D^*\to V:x\mapsto x - (\one+D^*D)^{-1}(x+D^*D^*x),
$$
where $\one+D^*D$ is understood as an operator from
$\dom D$ to $(\dom D)^*$. 
The graph inner product on $V$ is compatible with the symplectic
form and the associated complex structure is $x\mapsto Jx:=D^*x$, that is
$\omega(x,Jy)=\la x, y \ra_{D^*}$. This shows that $(V,\om)$
is indeed a symplectic Hilbert space. 

\smallskip\noindent{\bf (iii)}
In the formulation of (ii) the subspace $\Lambda_0$ and its orthogonal 
complement are given by 
$$
\Lambda_0 = \left\{x\in V\,|\,\exists \xi\in\dom D\,
s.t.~D^*(x+\xi)=0\right\}
= \left\{x\in V\,|\,D^*x\in\im D\right\}
$$
and
$$
\Lambda_0^\perp = D^*\Lambda_0 = V\cap\im D. 
$$
\end{rmk}

\begin{dfn}\label{dfn:Lfred}
A triple $(V,\Lambda_1,\Lambda_2)$ consisting of a Hilbert space
$V$ and two closed subspaces $\Lambda_1,\Lambda_2\subset V$ 
is called {\bf Fredholm} if $\Lambda_1\cap\Lambda_2$ is finite 
dimensional, $\Lambda_1+\Lambda_2$ is a closed subspace of $V$, 
and the cosum $V/(\Lambda_1+\Lambda_2)$ is finite 
dimensional (see~\cite{RRS}); equivalently 
the linear operator $S:\Lambda_1\times\Lambda_2\to V$ 
given by $S(x_1,x_2):=x_1+x_2$ is Fredholm.
The {\bf Fredholm index} of a Fredholm triple 
$(V,\Lambda_1,\Lambda_2)$ is defined by 
$$
\mathrm{index}(V,\Lambda_1,\Lambda_2)
:= \dim(\Lambda_1\cap\Lambda_2) - \dim(V/(\Lambda_1+\Lambda_2))
= \mathrm{index}(S).
$$
\end{dfn}

\begin{lemma} \label{lem:Fredtrip}
Let $\Lambda\subset V$ be a Lagrangian subspace.
Then $D_\Lambda:\dom D_\Lambda \to H$ is a Fredholm operator if and only if
$(V,\Lambda_0,\Lambda)$ is a Fredholm triple.
\end{lemma}

\begin{proof}
This follows from the definition and the fact that the homomorphisms
$$
\ker D_\Lambda\to\Lambda_0\cap\Lambda:x\mapsto[x] ,
$$ 
$$
\frac{V}{\Lambda_0+\Lambda}
\to\frac{H}{\im\,D_\Lambda}:[x]\mapsto[D^*x]
$$
are bijective. For the second map this uses Lemma~\ref{le:D*D} below.
\end{proof}

\begin{lemma}\label{le:D*D}
Let $D:\dom D\to H$ be an injective symmetric operator 
with a closed image and a dense domain.  Then 
$$
Y:=\left\{\xi\in\dom D\,|\,D\xi\in\dom D^*\right\}
$$
is a Hilbert space with the inner product 
$$
\inner{\xi}{\eta}_Y := 
\inner{\xi}{\eta}_H+\inner{D\xi}{D\eta}_H+\inner{D^*D\xi}{D^*D\eta}_H
$$
and the operator $D^*D:Y\to H$ is an isomorphism.  

Moreover, if the inclusion $\dom D\to H$ is a compact operator
then the operator 
$
D(D^*D)^{-1}:H\to H
$
is compact.
\end{lemma}

\begin{proof}
We prove that $Y$ is complete.  Let $\xi_i\in Y$ be a Cauchy sequence.
Then $\xi_i$, $D\xi_i$, $D^*D\xi_i$ are Cauchy sequences in $H$.  
Define $\xi:=\lim\xi_i$, $x:=\lim D\xi_i$, $y:=\lim D^*D\xi_i$.
Since $D$ and $D^*$ have closed graphs we have 
$\xi\in\dom D$, $x\in\dom D^*$, $D\xi=x$, and $D^*x=y$.
Hence $\xi\in Y$ and $\xi_i$ converges to $\xi$ in $Y$. 

That $D^*D:Y\to H$ is injective follows from the fact that 
$D$ is injective and $\inner{D^*D\xi}{\xi}_H=\Norm{D\xi}_H^2$ 
for $\xi\in Y$.  Now consider the  Gelfand triple
$$
Z\subset H\subset Z^*,
$$
where $Z:=\dom D$ and $\inner{\xi}{\eta}_Z=\inner{D\xi}{D\eta}_H$.
We identify $H$ with its dual space and define the inclusion
$H\to Z^*$ as the dual operator of the inclusion $Z\to H$.
We can think of $D:Z\to H$ as a bounded linear operator 
and of its adjoint as bounded linear operator $D^*:H\to Z^*$. 
Then
$
\dom D^* = \left\{x\in H\,|\,D^*x\in H\right\}.
$
Since $D:Z\to H$ is injective and has a closed 
image the dual operator $D^*:H\to Z^*$ is surjective. 
Now let $y\in H$.   Then $y\in Z^*$ and hence there 
exists an element $x\in H$ with $D^*x=y$.
Since $D^*x\in H$ we have $x\in\dom D^*$. 
Now it follows from the definitions that the kernel
of $D^*$ is the orthogonal complement of the image 
of $D$.  Since the image of $D$ is closed this implies
$
H = \ker D^*\oplus\im\,D.
$
Hence there is a vector $x_0\in\ker D^*$ such that 
$x-x_0\in\im\,D$.  Choose $\xi\in\dom D$ such that
$D\xi=x-x_0$.  Then $D\xi\in\dom D^*$ and $D^*D\xi=D^*x=y$.
This proves that $D^*D$ is surjective.

Now assume that the inclusion $Z\to H$ is compact.
To prove that the operator $D(D^*D)^{-1}:H\to H$ is compact 
we observe that 
\begin{equation}\label{eq:Z*}
\Norm{x}_{Z^*} 
:= \sup_{0\ne\xi\in\dom D}\frac{\inner{x}{\xi}_H}{\Norm{D\xi}_H}
= \Norm{D(D^*D)^{-1}x}_H
\end{equation}
for every $x\in H\subset Z^*$. 
Here the last equation follows from the fact that the supremum
in the second term is attained at the vector $\xi_0=(D^*D)^{-1}x$
with $x=D^*D\xi_0$. Now let $x_i$ be a bounded sequence in $H$.
Since the inclusion $H\to Z^*$ is compact, there exists a subsequence
$x_{i_k}$ which converges in $Z^*$ and it follows from~(\ref{eq:Z*})
that the sequence $D(D^*D)^{-1}x_{i_k}$ converges in $H$.
This proves the lemma. 
\end{proof}

\begin{rmk} \label{rm:llv} \rm
{\bf (i)}
By Lemma~\ref{le:D*D} the subspaces $\Lambda_0$ 
and $\Lambda_0^\perp$ of $V$ in Remark~\ref{rmk:V}
can also be written in the form
$$
\Lambda_0 = \left\{x\in\dom D^*\,|\,D^*x + D(D^*D)^{-1}x = 0\right\}
= \ker(D^*+T),
$$
$$
\Lambda_0^\perp =
\left\{x\in V \,|\,x = D(D^*D)^{-1}D^*x\right\}
= T\ker(D^*+T) ,
$$
where $T:=D(D^*D)^{-1}:H\to H$ maps to $\im T=\dom D^*\cap\im D$.

\smallskip\noindent{\bf (ii)}
The orthogonal projection of $V$ onto $\Lambda_0^\perp$ extends
to a bounded linear operator $\Pi_0: H\to H$ given by 
$$
\Pi_0 x = D(D^*D+TD)^{-1}(D^*x+Tx) .
$$
Here $D^*D+TD : \dom D \to (\dom D)^*$ is an isomorphism 
because 
$$
\la x, D^*Dx + TD x \ra = \|Dx\|^2 + \|TDx\|^2\ge \delta \|x\|^2 .
$$
In fact, $\Pi_0$ is a projection on all of $H$, its kernel is $\Lambda_0$, 
and its image is equal to the image of $D$. In particular, $\Pi_0|_{\Lambda_0^\perp}=\one$.

\smallskip\noindent{\bf (iii)}
In all our applications the inclusion $\dom D\to H$ is a compact
operator. Then, by Lemma~\ref{le:D*D}, $T:H\to H$ is compact, and thus
the inclusion $\Lambda_0^\perp\to H$ is compact.  
Indeed, the inclusion is given by the composition $x\mapsto TD^*x$ 
of a compact and a bounded operator.
\end{rmk}

The inclusions $\dom D^*\to H$ and $\Lambda_0\to H$, however, 
are not compact unless $V$ is finite dimensional. 
Namely, if $V$ is infinite dimensional then so is the kernel of $D^*$
(since $\Lambda_0\subset V$ is Lagrangian)
and the inclusion $\dom D^*\supset\ker D^*\to H$ is an isometric embedding. 
Lemma~\ref{le:DLcpct} below gives a condition under which the domain
of a self-adjoint extension of $D$ has a compact embedding into $H$.
This requires the notion of a compact perturbation of a closed
subspace of $V$. 

\begin{dfn}\label{dfn:Lcpct}
Let $V$ be a Hilbert space and $\Lambda\subset V$ 
be a closed subspace.  
A closed subspace $\Lambda'\subset V$ is called
a {\bf compact perturbation of $\Lambda$} if 
the projection of $\Lambda'$ onto some (and hence 
every) complement of $\Lambda$ in $V$ is a compact 
operator and vice versa. 
\end{dfn}

\begin{rmk}\label{rmk:Lcpct ?} \rm
The notion of {\em compact perturbation} defines an equivalence 
relation on the set of closed subspaces of $V$.
To see this denote by $\Pi:V\to\Lambda$
and $\Pi':V\to\Lambda'$ the orthogonal projections.  
If $\Lambda'$ is a compact perturbation of $\Lambda$ 
and $\Lambda''$ is a compact perturbation of $\Lambda'$ 
then the operators $\one-\Pi:\Lambda'\to\Lambda^\perp$
and $\one-\Pi':\Lambda''\to(\Lambda')^\perp$
are compact.  Hence the operator 
$
(\one-\Pi)|_{\Lambda''}
= (\one-\Pi)(\one-\Pi')|_{\Lambda''}
  + (\one-\Pi)\Pi'|_{\Lambda''}
$
is compact. Repeating this argument with $\Lambda$ and
$\Lambda''$ interchanged we see that $\Lambda''$ 
is a compact perturbation of $\Lambda$.
\end{rmk}

\begin{lemma}\label{le:fredstable}
Let $V$ be a Hilbert space and $\Lambda_1,\Lambda,\Lambda'\subset V$
be closed subspaces such that $\Lambda'$ is a compact 
perturbation of $\Lambda$. If $(V,\Lambda_1,\Lambda)$ 
is a Fredholm triple then so is $(V,\Lambda_1,\Lambda')$.
\end{lemma}

\begin{proof}
Let $\Pi:V\to\Lambda$ and $\Pi':V\to\Lambda'$ 
be the orthogonal projections.  Then 
$(\one-\Pi\Pi')|_{\Lambda}=\Pi(\one-\Pi')|_{\Lambda}:\Lambda\to\Lambda$ 
and
$(\one-\Pi'\Pi)|_{\Lambda'}:\Lambda'\to\Lambda'$ 
are compact operators.  This implies that 
$\Pi|_{\Lambda'}:\Lambda'\to\Lambda$ and 
$\Pi'|_{\Lambda}:\Lambda\to\Lambda'$ are Fredholm 
operators with opposite indices; see e.g.\ \cite[Chapter III.3]{KG}.

Now suppose that $(V,\Lambda_1,\Lambda)$ is a Fredholm triple, i.e.\ 
the map $S:\Lambda_1\times\Lambda\to V$ given by $S(v_1,v)=v_1+v$ 
is Fredholm.  Then the operator
$$
S'':=S\circ(\one\times\Pi):\Lambda_1\times\Lambda'\to V
$$
is Fredholm. Define the map $S':\Lambda_1\times\Lambda'\to V$ by 
$S'(v_1,v')=v_1+v'$.  Since $ S'(v_1,v')-S''(v_1,v') = (\one-\Pi\Pi') v' $
the operator $S'-S''$ is compact. Hence $S'$ is a Fredholm operator 
and so $(V,\Lambda_1,\Lambda')$ is a Fredholm triple.
\end{proof}

\begin{lemma}\label{le:LL'}
Let $(V,\om)$ be a symplectic Hilbert space.  
Let $\Lambda,\Lambda'\subset V$ be Lagrangian subspaces. 
Then the following are equivalent.
\begin{description}
\item[(i)]
$\Lambda'$ is a compact perturbation of $\Lambda$. 
\item[(ii)]
The projection of $\Lambda'$ onto $\Lambda^\perp$ is a compact operator.
\item[(iii)]
The operator 
$\Lambda'\to\Lambda^*:v'\mapsto\om(v',\cdot)$ 
is compact.
\end{description}
\end{lemma}

\begin{proof}
By definition, (i) implies (ii).
The Lagrangian condition asserts that the orthogonal complement 
$\Lambda^\perp$ is isomorphic to $\Lambda^*$ via the isomorphism 
$\Lambda^\perp\to\Lambda^*:v\mapsto\om(v,\cdot)$. 
Under this isomorphism the orthogonal projection 
$\Lambda'\to\Lambda^\perp$ corresponds to the operator
$\Lambda'\to\Lambda^*:v'\mapsto\om(v',\cdot)$, hence
(ii) and (iii) are equivalent.
To see that (iii) implies (i) note that the operators 
$\Lambda'\to\Lambda^*:v'\mapsto\om(v',\cdot)$ 
and $\Lambda\to(\Lambda')^*:v\mapsto-\om(v,\cdot)$
are dual to each other. Using ``$\mathrm{(iii)}\Leftrightarrow\mathrm{(ii)}$''
we see that (iii) implies compactness of both projections $\Lambda'\to\Lambda^\perp$ 
and $\Lambda\to(\Lambda')^\perp$. This proves the lemma.
\end{proof}

\begin{lemma}\label{le:DLcpct}
Let $D:\dom D\to H$ be an injective symmetric operator 
with a closed image and a dense domain and suppose that 
the inclusion $\dom D\to H$ is a compact operator. 
Let $V=(\dom D)^\perp$ be the Gelfand--Robbin quotient,
$\Lambda\subset V$ be a Lagrangian subspace,
and $\Lambda_0,\Lambda_0^\perp$ be as in Remark~\ref{rmk:V}. 
Then the following are equivalent. 
\begin{description}
\item[(i)]
The inclusion $\dom D_\Lambda\to H$ is compact.
\item[(ii)]
The inclusion $\Lambda\to H$ is compact.
\item[(iii)]
$\Lambda$ is a compact perturbation of $\Lambda_0^\perp$. 
\end{description}
\end{lemma}

\begin{proof}
Let $\Pi_0:V\to V$ denote the orthogonal projection onto $\Lambda_0^\perp$. 
Then $\Pi_0:V\to H$ is compact since the inclusion of the image 
$\Pi_0(V)=\Lambda_0^\perp$ into $H$ is compact by Remark~\ref{rm:llv}~(ii).
By Lemma~\ref{le:LL'}, (iii) holds if and only if the operator 
$(\one-\Pi_0)|_{\Lambda}:\Lambda\to\Lambda_0$ is compact.
Moreover, the graph norm of $D^*$ on $\Lambda_0=\ker(D^*+T)$ 
is equivalent to the norm of $H$ so, in fact, (iii) holds if and only if
the operator $(\one-\Pi_0)|_{\Lambda}:\Lambda\to H$ is compact.
We deduce that (iii) is equivalent to (ii) because 
the inclusion $\Lambda\to H$ is given by the sum
$\one|_{\Lambda}=(\one-\Pi_0)|_{\Lambda}+\Pi_0|_{\Lambda}$, 
where $\Pi_0|_{\Lambda}:\Lambda\to H$ is compact.  

That~(i) is equivalent to~(ii) follows from the 
fact that the inclusion of $\dom D$ into $H$ is compact, 
by assumption, and $\dom D_\Lambda=\dom D\oplus\Lambda$.
\end{proof}

\begin{lemma}\label{le:DPL}
Let $D:\dom D\to H$ be an injective symmetric operator 
with a closed image and a dense domain and suppose that 
the inclusion $\dom D\to H$ is a compact operator. 
Let $V=\dom D^*/\dom D$ be the Gelfand--Robbin quotient
and $\Lambda_0=\{[x]\in V\,|\, D^*x=0\}$ as in~(\ref{eq:ker}). 
Let $P:H\to H$ be a self-adjoint bounded linear operator
such that $D+P:\dom D\to H$ is injective.  
Then the following are equivalent. 
\begin{description}
\item[(i)]
The composition of $P$ with the inclusion $\dom D^*\to H$ is a compact 
operator. 
\item[(ii)]
The operator $P|_{\ker(D^*+P)}:\ker(D^*+P)\to H$ is compact. 
\item[(iii)]
$\Lambda_P:= \left\{[x]\in V\,|\,D^*x+Px=0\right\}$ 
is a compact perturbation of $\Lambda_0$. 
\end{description}
\end{lemma}

\begin{proof}
Abbreviate ${D_P:=D+P}$. Then $\dom D_P^*=\dom D^*$ and the 
graph norm of $D^*$ is equivalent to the graph norm of $D_P^*$.
Moreover, on $\ker D_P^*$ and $\ker D^*$ both graph norms are 
equivalent to the norm of $H$.
For $[x]\in\Lambda_P$, represented by $x\in\ker D_P^*$, and 
$[x_0]\in\Lambda_0$, represented by $x_0\in\ker D^*$, we have 
$$
\om([x],[x_0]) = \inner{D^*x}{x_0} - \inner{x}{D^*x_0}
= - \inner{Px}{x_0} = \inner{TPx-Px}{x_0} ,
$$
where $T:=D(D^*D)^{-1}:H\to\dom D^*$.  
Using Lemma~\ref{le:LL'} and the compactness of 
the inclusion $\dom D\to H$, we see that 
$\Lambda_P$ is a compact perturbation of $\Lambda_0$ if and only if 
${(P-TP)|_{\ker D_P^*}:\ker D_P^*\to\ker D^*}$
is a compact operator.
Since $T$ is compact, by Lemma~\ref{le:D*D},
this shows that~(ii) is equivalent to~(iii).
That~(i) implies~(ii) is obvious.  
To prove that~(ii) implies~(i) note that, 
by Remark~\ref{rm:llv} with $D$ replaced
by $D_P$, the inclusion of ${\dom D^*\cap\im D_P}$
into $H$ is compact.  Since the decomposition 
$
\dom D^* = (\dom D^*\cap\im D_P)\oplus\ker D_P^*
$
is orthogonal with respect to the graph norm of $D_P^*$,
this shows that~(ii) implies~(i).
\end{proof}

\begin{rmk}\label{rmk:P}\rm
Let $D$, $V$, $\Lambda_0$ be as in Lemma~\ref{le:DPL},
$P:H\to H$ be a bounded self-adjoint operator, and 
denote $\Lambda_P:= \left\{[x]\in V\,|\,D^*x+Px=0\right\}$.

\smallskip\noindent{\bf (i)}
Let $\Lambda_P^{\perp,{\scriptscriptstyle P}}$ denote the 
orthogonal complement of $\Lambda_P$ with respect to 
the graph inner product of $D^*+P$.  
Then it always is a compact perturbation of~$\Lambda_0^\perp$.  
Namely, by Remark~\ref{rm:llv} with $D$ replaced by $D+P$,
the inclusion $\Lambda_P^{\perp,{\scriptscriptstyle P}}\to H$ is compact.
Hence, by Lemma~\ref{le:DLcpct} with $D$ replaced by $D+P$ and
$\Lambda:=\Lambda_P^{\perp,{\scriptscriptstyle P}}$,
the inclusion $\left\{v\in\dom D^*\,|\,[v]\in \Lambda\right\}\to H$
is compact.  Using Lemma~\ref{le:DLcpct} again we deduce that
$\Lambda$ is a compact perturbation of $\Lambda_0^\perp$.

\smallskip\noindent{\bf (ii)}
The orthogonal complement $\Lambda_P^\perp$ 
with respect to the graph inner product of $D^*$
is a compact perturbation of $\Lambda_0^\perp$ 
if and only if the restriction of $P$ to $\dom D^*$ 
is a compact operator.
This follows from Lemma~\ref{le:DPL} and the fact that 
$\Lambda_P^\perp=D^*\Lambda_P$ and $\Lambda_0^\perp=D^*\Lambda_0$
in the notation of Remark~\ref{rmk:V}, where $D^*$ is a compatible
complex structure on $V$.

\smallskip\noindent{\bf (iii)}
It follows from~(i) and~(ii) that $\Lambda_P^{\perp,{\scriptscriptstyle P}}$
is a compact perturbation of $\Lambda_P^\perp$ if and only if 
the restriction of $P$ to the domain of $D^*$ is a compact operator.

\smallskip\noindent{\bf (iv)}
If $\Lambda$ is a compact perturbation of~$\Lambda_0^\perp$
then $(V,\Lambda_P,\Lambda)$ is a Fredholm triple.
Since $(V,\Lambda_P,\Lambda_P^{\perp,{\scriptscriptstyle P}})$ 
is a Fredholm triple, this follows from~(i) and Lemma~\ref{le:fredstable}.
\end{rmk}

\begin{lemma}\label{le:maslov}
Let $D,V,\Lambda_0$ be as in Lemma~\ref{le:DPL}
and let $P(s):H\to H$ for $s\in\R$ be a continuously differentiable
family of self-adjoint bounded linear operators.
Assume that $P(s)$ converges to $P(\pm\infty)=:P^\pm$ in the 
operator norm as $s$ tends to~$\pm\infty$, that
$(D+P(s))|_{\dom D}$ is injective for every~$s\in\R\cup\{\pm\infty\}$,
and that
$$
\frac{\ker (D^*+P^-)\oplus\mathrm{dom}\, D} 
{\mathrm{dom}\, D} 
= \frac{\ker (D^*+P^+)\oplus\mathrm{dom}\, D} 
{\mathrm{dom}\, D}
=: \Lambda'_0\subset V.
$$
Then the spectral flow of the operator family 
$s\mapsto\left(D+P(s)\right)_\Lambda$ is independent of
the Lagrangian subspace $\Lambda\subset V$ such that
$V = \Lambda'_0\oplus\Lambda$ and 
$\Lambda$ is a compact perturbation of $\Lambda_0^\perp$.
\end{lemma}

\begin{rmk}\label{rmk:loopmaslov}\rm
Let $D,V,\Lambda_0$ be as in Lemma~\ref{le:DPL}, 
and let $Q:H\to H$ be a Hilbert space isomorphism such that 
$$
x\in\dom D^*\qquad\implies\qquad x-Qx\in\dom D.
$$
Then $Q$ induces the identity on $V$. 
Let $P(s):H\to H$ for $s\in\R$ be a continuously differentiable
family of self-adjoint bounded linear operators such that 
$$
D+P(s+1) = Q^{-1}(D+P(s))Q.
$$
Assume $(D+P(s))|_{\dom D}$ is injective for every~$s\in\R$ and
denote 
$$
\Lambda_0':=\frac{\ker(D^*+P(0))\oplus\dom D}{\dom D} 
=\frac{\ker(D^*+P(1))\oplus\dom D}{\dom D} .
$$
Then the spectral flow of the operator family 
$s\mapsto\left(D+P(s)\right)_\Lambda$ 
on the fundamental domain $[0,1]$ is independent of
the Lagrangian subspace $\Lambda\subset V$ such that
$V = \Lambda_0'\oplus\Lambda$ and 
$\Lambda$ is a compact perturbation of $\Lambda_0^\perp$.
The proof is the same as that of Lemma~\ref{le:maslov}.
\end{rmk}

\begin{proof}[Proof of Lemma~\ref{le:maslov}.]
The operators $D+P(s):\dom D \to H$ satisfy the assumptions of this section
and give rise to the constant Gelfand-Robbin quotient $V=\dom D^*/\dom D$
since $\dom (D+P(s))^*=\dom D^*$.
Hence any Lagrangian subspace $\Lambda\subset V$ gives rise to a family
of self-adjoint operators $A(s):=(D+P(s))_\Lambda : \dom D_\Lambda \to H$, 
which satisfies the conditions (A1--3) of Section~\ref{app:spec}
whenever $V=\Lambda_0'\oplus\Lambda$.
In particular, the estimate in (A1) holds for $s=\pm\infty$, 
i.e.~$\Norm{x}_{D^*}\leq C\Norm{(D+P^\pm)_\Lambda x}_H$
for $x\in\dom D_\Lambda$, because
$(V,\Lambda_0',\Lambda)$ is a Fredholm triple and
$(D+P^\pm)_\Lambda$ is injective.
The estimate for $s\in\R$ follows from a uniform bound of the 
form $\|P(s)-P^\pm\|\leq C$ for the operator norm on~$H$.
The assumptions (W1--3) are satisfied with the trivial map $Q\equiv\one$
and the constant domain $W_0=\dom D_\Lambda$. 
In particular, the domain embeds compactly
to $H$, by Lemma~\ref{le:DLcpct}, whenever 
$\Lambda$ is a compact perturbation of $\Lambda_0^\perp$.
Hence the spectral flow is well defined under our assumptions 
(see Appendix~\ref{app:spec}). 

We prove that the set $\cS$ of Lagrangian subspaces of $V$ 
that are transverse to $\Lambda_0'$ and
are compact perturbations of $\Lambda_0^\perp$ is connected.
For that purpose let $\Lambda_1\subset V$ denote the orthogonal complement of $\Lambda_0'$
with respect to the graph inner product of $D^*+P^+$ and let 
$I_\omega:\Lambda_0'\to\Lambda_1^*$ be the isomorphism 
given by $v\mapsto\om(v,\cdot)$.  Then a subspace 
$\Lambda\subset V=\Lambda_1\oplus\Lambda_0'$ is a complement of 
$\Lambda_0'$ if and only if it is the graph of a linear operator
from $\Lambda_1$ to $\Lambda_0'$ or, equivalently, 
$\Lambda=\Lambda_A:=\mathrm{graph}(I_\omega^{-1}\circ A)$
for some linear operator $A:\Lambda_1\to\Lambda_1^*$.
One can check that the subspace $\Lambda_A$ is Lagrangian
if and only if $A$ is self-adjoint and that it is 
a compact perturbation of $\Lambda_0^\perp$
if and only if $A$ is compact.  
The last assertion uses the explicit formula 
$x+I_\omega^{-1} A x \mapsto I_\omega^{-1} A x$
for the projection $\Lambda_A\to \Lambda_0'$ along $\Lambda_1$
and the fact that $\Lambda_A$ is a compact perturbation 
of $\Lambda_0^\perp$ if and only if it is a compact perturbation of $\Lambda_1$,
by Remark~\ref{rmk:P}~(i) and Remark~\ref{rmk:Lcpct ?}.
Thus we have identified $\cS$
with the vector space of compact self-adjoint operators 
$A:\Lambda_1\to \Lambda_1^*$ and so $\cS$ is contractible, as claimed.
 
Now the result follows from the homotopy invariance of the spectral flow. 
The homotopies of Lagrangian subspaces do not directly translate into
homotopies in the sense of Section~\ref{app:spec}, see the proof of 
Lemma~\ref{lem:kernel} below.
However, the homotopy invariance of the spectral flow of the family
$s\mapsto\left(D+P(s)\right)_\Lambda$ follows from Remark~\ref{rmk:smaslov},
where the spectral flow is identified with a Maslov index, which in turn
is invariant under homotopies of $\Lambda$.
\end{proof}

\begin{rmk}\label{rmk:smaslov}\rm  
{\bf (i)}
Let $[0,1]\ni s\mapsto(\Lambda_0(s),\Lambda_1(s))$ be 
a smooth path of pairs of Lagrangian subspaces
of $V$ such that $(V,\Lambda_0(s),\Lambda_1(s))$ 
is a Fredholm triple for every $s$.  
For each $s$ define the crossing form 
$
\Gamma(\Lambda_0,\Lambda_1,s):
\Lambda_0(s)\cap\Lambda_1(s)\to\R
$
by 
$$
\Gamma(\Lambda_0,\Lambda_1,s)(v) 
:= \left.\frac{\rd}{\rd t}\right|_{t=0}
\bigl(\om(v,v_0'(t)) - \om(v,v_1'(t))\bigr)
$$
for $v\in\Lambda_0(s)\cap\Lambda_1(s)$, 
where $\Lambda_0',\Lambda_1'\subset V$ are Lagrangian
subspaces such that 
$
V=\Lambda_0(s)\oplus\Lambda_0'=\Lambda_1(s)\oplus\Lambda_1'
$
and $v_0'(t)\in\Lambda_0'$, $v_1'(t)\in\Lambda_1'$
are chosen such that $v+v_0'(t)\in\Lambda_0(s+t)$
and $v+v_1'(t)\in\Lambda_1(s+t)$. 
As in~\cite{RSM} the {\bf Maslov index} is defined as 
the sum of the signatures of the crossing forms
$$
\mu(\Lambda_0,\Lambda_1) 
:= \sum_s\mathrm{sign}\,\Gamma(\Lambda_0,\Lambda_1,s)
$$
provided that the crossing forms are all nondegenerate and 
$\Lambda_0(s)$ is transverse to $\Lambda_1(s)$ for $s=0,1$.
Under this assumption the sum is finite.
The nondegeneracy condition can be achieved by a small
perturbation with fixed endpoints. 
The Maslov index is invariant under homotopies of 
paths of Lagrangian Fredholm triples with transverse 
endpoints. 

\smallskip\noindent{\bf (ii)}
The spectral flow in Lemma~\ref{le:maslov}
can be identified with the Maslov index
\begin{equation}\label{maslovclaim}
\mu_\spec\left(\left(D+P\right)_\Lambda\right)
= \mu\left(\Lambda_P,\Lambda\right),
\end{equation}
where $\Lambda_0(s):=\Lambda_{P(s)}=\{[x]\in V\,|\,D^*x+P(s)x=0\}$
and $\Lambda_1(s):=\Lambda$ for every~$s$.
The Fredholm property of the triples $(V,\Lambda_{P(s)},\Lambda)$
follows from Remark~\ref{rmk:P}~(iv).

To prove~(\ref{maslovclaim}), fix a real number $s$, choose $\Lambda_0'$ and 
$v_0'(t)$ as in (i), let $x_0(t)\in\dom D^*$ be the 
smooth path defined by $(D^*+P(s+t))x_0(t)=0$ and 
$[x_0(t)]=v+v_0'(t)\in\Lambda_0(s+t)$, and denote $x:=x_0(0)$ 
so that $[x]=v\in\Lambda_0(s)\cap\Lambda$.  Then 
\begin{align*}
\left.\frac{\rd}{\rd t}\right|_{t=0} 
\om(v,v_0'(t))
&=\left.\frac{\rd}{\rd t}\right|_{t=0}
  \om(v,v+v_0'(t)) \\
&=\left.\frac{\rd}{\rd t}\right|_{t=0}
  \biggl(\Inner{D^*x}{x_0(t)}-\Inner{x}{D^*x_0(t)}\biggr) \\
&= \left.\frac{\rd}{\rd t}\right|_{t=0}
  \Inner{\bigl(D^*+P(s+t)\bigr)x}{x_0(t)} \\
&= \left.\frac{\rd}{\rd t}\right|_{t=0}
  \Inner{\bigl(D^*+P(s+t)\bigr)x}{x}.
\end{align*}
This shows that the crossing forms 
$\Gamma\left(\left(D^*+P\right)_\Lambda,s\right)$ 
and $-\Gamma(\Lambda_0,\Lambda,s)$ 
agree under the isomorphism 
$
\ker\left(D^*+P(s)\right)_\Lambda
\to\Lambda_0(s)\cap\Lambda:x\mapsto[x].
$
\end{rmk}

\begin{lem}\label{lem:kernel}
Let $D,V,\Lambda_0,Q,P$ be as in Remark~\ref{rmk:loopmaslov}.
Denote by $\Upsilon$ the set of Lagrangian subspaces $\Lambda\subset V$ 
that are compact perturbations of $\Lambda_0^\perp$. 
For every $\Lambda\in\Upsilon$ there is a Fredholm operator
$$
\cD_\Lambda:=\p_s+D_\Lambda+P(s):
\cW_\Lambda\to\cH
$$
with 
$$
\cH:=\left\{\xi\in L^2_\loc(\R,H)\,|\,\xi(s+1)=Q^{-1}\xi(s)\right\}, 
$$
$$
\cW_\Lambda:=\bigl\{
\xi\in L^2_\loc(\R,\dom D_\Lambda)\cap W^{1,2}_\loc(\R,H)\,|\,
\xi(s+1)=Q^{-1}\xi(s)\bigr\} .
$$
The determinants $\det(\cD_\Lambda)$ for $\Lambda\in\Upsilon$
form a line bundle over $\Upsilon$.
\end{lem}

\begin{proof}
$\cD_\Lambda$ is Fredholm since it is the operator of Theorem~\ref{thm:indexS1} 
with $A(s)=D_\Lambda+P(s)$ and constant domain $W(s)=\dom D_\Lambda$.

We do not know if for any two subspaces $\Lambda,\Lambda'\in\Upsilon$ 
there is a Hilbert space isomorphism
of $Q:H\to H$ that identifies $\dom D_\Lambda$ with $\dom D_{\Lambda'}$,
as would be required for a homotopy of operator families 
in the sense of Appendix~\ref{app:spec}.
However, one can prove directly that the kernel
of $\cD_\Lambda$ depends continuously on $\Lambda$
(as a subspace of $\cH$) if $\cD_\Lambda$ is surjective.
This proves the lemma since the transverse situation can always be
achieved by finite dimensional stabilization.

To prove the continuous dependence of $\ker\cD_\Lambda$ on $\Lambda$
we will use the fact that every element $\xi\in\ker\,\cD_\Lambda$ 
is a smooth function from $\R$ to $\dom D_\Lambda$ 
(see~\cite[Theorem~3.13]{RS}) and satisfies an estimate of the form 
$\Norm{\xi}_{\cW_\Lambda}
+\Norm{\p_s\xi}_{\cW_\Lambda}\le c\Norm{\xi}_\cH$.
Two Lagrangian subspaces $\Lambda,\Lambda'\in\Upsilon$ are close
if there exists an isomorphism of $V$ close to the identity that
maps $\Lambda$ to $\Lambda'$. This extends to an isomorphism of
$\dom D_\Lambda=\Lambda\oplus\dom D$ 
and $\dom D_{\Lambda'}=\Lambda'\oplus\dom D$
(which does not necessarily extend to an isomorphism of $H$).
This isomorphism of domains followed by the orthogonal projection 
onto the kernel of $\cD_{\Lambda'}$ induces a map 
${\ker\cD_\Lambda\to\cW_{\Lambda'}}$,
which is an isomorphism for $\Lambda'$ sufficiently close to $\Lambda$. 
\end{proof}


\section{Unique continuation}\label{app:uc}

In this appendix we formulate a general unique continuation
theorem based on the Agmon--Nirenberg technique.  The method 
was also used by Donaldson--Kronheimer \cite[pp150]{DK}
and Taubes \cite{T:ucon} to prove unique continuation results 
for anti-self-dual instantons
and by Kronheimer--Mrowka~\cite{KM} and in \cite{SaSW}
for the Seiberg--Witten equations. 

Let $H$ be a Hilbert space and $A(s)$ be a family 
of (unbounded) symmetric operators on $H$ 
with domains $\dom(A(s))\subset H$.
The operators $A(s)$ are not required 
to be self-adjoint although in the main applications
they will be and, moreover, their domains will be independent of $s$.  
However, in some interesting cases these operators 
are symmetric with respect to time-dependent 
inner products. The following theorem is a special case 
of a result by Agmon and Nirenberg~\cite{AN}. 

\begin{thm}[Agmon-Nirenberg]\label{thm:agni1}
Let $H$ be a real Hilbert space and let $A(s):\dom(A(s))\to H$
be a family of symmetric linear operators.
Assume that $x:[0,T)\to H$ for $0<T\leq\infty$ is continuously 
differentiable in the weak topology
such that $x(s)\in\dom(A(s))$ and
\begin{equation}\label{eq:agni1}
    \left\|\dot x(s)+A(s)x(s)\right\|
    \le c_1(s)\left\|x(s)\right\|
\end{equation}
for every $s\in[0,T)$, where $\dot x(s):=\pd_s x(s) \in H$
denotes the time derivative of~$x$.
Assume further that the function 
$s\mapsto\inner{x(s)}{A(s)x(s)}$ is also continuously 
differentiable and satisfies
\begin{equation}\label{eq:agni2}
\frac{d}{ds}\inner{x}{Ax} - 2\inner{\dot x}{Ax}
\le 2 c_2(s)\left\|Ax\right\|\left\|x\right\|
  +  c_3(s)\left\|x\right\|^2.
\end{equation}
Here $c_1, c_2, c_3 :[0,T)\to\R$ are continuous nonnegative
functions satisfying
\begin{equation*}\label{eq:ass c}
a_0:= 2\int_0^T c_2 <\infty, \quad
b_0:= \int_0^T ( c_1^2 + c_2^2 + c_3) <\infty, \quad
c_0:= \sup c_1 <\infty .
\end{equation*}
Then the following holds.

\smallskip\noindent{\bf (i)}
If $x(0)=0$ then $x(s)=0$ for all $s\in[0,T)$.

\smallskip\noindent{\bf (ii)}
If $x(0)\ne0$ then $x(s)\ne 0$ for all 
$s\in[0,T)$ and, moreover,
$$
\left\|x(s)\right\|
\ge e^{-c s} \left\|x(0)\right\| ,\qquad
c:= c_0 + e^{a_0}\bigl(b_0 + \|x(0)\|^{-1} \| A(0) x(0) \| \bigr) .
$$
\end{thm}

\begin{proof}
The basic idea of the proof is to use the convexity of the function 
$t\mapsto\log\left\|x(t)\right\|^2$.
Assume that $x(0)\neq 0$ and define 
$$
\phi(s) := \log\,\left\|x(s)\right\|
- \int_0^s\frac{\inner{x(\sigma)}{\dot x(\sigma)+A(\sigma)x(\sigma)}}
{\left\|x(\sigma)\right\|^2}\,\rd\sigma
$$
for $0\le s < T$ wherever $x(s)\ne 0$. 
Then we prove that $\phi$ is twice continuously 
differentiable and satisfies the differential inequality
\begin{equation}\label{eq:ab}
\ddot\phi + a\left|\dot\phi\right| + b\ge 0 ,
\qquad\quad
a := 2 c_2 , \qquad b := c_1^2+c_2^2 + c_3 .
\end{equation}
Define $f(s):=\dot x(s)+A(s)x(s)$.
Then the derivative of $\phi$ is given by 
$$
\dot\phi
= \frac{\inner{x}{\dot x}}{\left\|x\right\|^2}
  - \frac{\inner{x}{f}}{\left\|x\right\|^2}
= - \frac{\inner{x}{Ax}}{\left\|x\right\|^2}.
$$
Hence 
\begin{align*}
\ddot\phi
&=
- \frac{\frac{d}{ds}\inner{Ax}{x}}{\left\|x\right\|^2}
+ \frac{2\inner{Ax}{x}\inner{\dot x}{x}}{\left\|x\right\|^4} \\
&\ge
\frac{2\inner{Ax}{Ax-f}
-2 c_2\left\|Ax\right\|\left\|x\right\|
-c_3\left\|x\right\|^2}
{\left\|x\right\|^2}
- \frac{2\inner{Ax}{x}\inner{Ax-f}{x}}{\left\|x\right\|^4}.
\end{align*}
Here the second step follows from the inequality~(\ref{eq:agni2})
and the definition of $f$.  The terms on the right hand side 
can now be organized as follows
\begin{align*}
\ddot\phi
&\ge 
\frac{2}{\left\|x\right\|^2}
\left(\left\|Ax\right\|^2 
     - \frac{\inner{Ax}{x}^2}{\left\|x\right\|^2}\right)
- \frac{2}{\left\|x\right\|^2}
   \Inner{Ax-\frac{\inner{Ax}{x}}{\left\|x\right\|^2}x}{f}  \\
&
\quad - 2 c_2 \frac{\left\|Ax\right\|}{\left\|x\right\|} - c_3.  \\
&=
\frac{2}{\left\|x\right\|^2}
\left\|Ax-\frac{\inner{Ax}{x}}{\left\|x\right\|^2}x\right\|^2
- \frac{2}{\left\|x\right\|^2}
\Inner{Ax-\frac{\inner{Ax}{x}}{\left\|x\right\|^2}x}{f}  \\
&
\quad - 2 c_2 \frac{\left\|Ax\right\|}{\left\|x\right\|} - c_3.
\end{align*}
Now abbreviate
$$
    \xi = \frac{x}{\left\|x\right\|},\qquad
    \eta = \frac{Ax}{\left\|x\right\|}.
$$
Then $\dot\phi=-\inner{\xi}{\eta}$ and
the previous inequality can be written 
in the form
\begin{align*}
\ddot\phi
&\ge
2\left\|\eta-\inner{\eta}{\xi}\xi\right\|^2
- 2\Inner{\eta-\inner{\eta}{\xi}\xi}{\frac{f}{\left\|x\right\|}} 
- 2 c_2\left\|\eta\right\| - c_3 \\
&\ge
2\left\|\eta-\inner{\eta}{\xi}\xi\right\|^2
- 2\left\|\eta-\inner{\eta}{\xi}\xi\right\|
       \frac{\left\|f\right\|}{\left\|x\right\|} 
     - 2 c_2\left\|\eta\right\| - c_3 \\
&\ge
\left\|\eta-\inner{\eta}{\xi}\xi\right\|^2
- \frac{\left\|f\right\|^2}{\left\|x\right\|^2} 
- 2 c_2\left\|\eta\right\| - c_3 \\
&\ge
\left\|\eta-\inner{\eta}{\xi}\xi\right\|^2
- c_1^2 - 2 c_2\left\|\eta\right\| - c_3.
\end{align*}
The last but one inequality uses the fact 
that $\alpha\beta\le\alpha^2/2+\beta^2/2$ 
and the last inequality uses $\|f\|\le c_1\|x\|$. 
To obtain (\ref{eq:ab}) it remains to prove that 
$$
\left\|\eta-\inner{\eta}{\xi}\xi\right\|^2
- 2 c_2\left\|\eta\right\| 
\ge - 2 c_2\left|\dot\phi\right|- c_2^2.
$$
Since $\dot\phi=-\inner{\xi}{\eta}$ this is equivalent to
$$
2 c_2\left\|\eta\right\|
\le \left\|\eta-\inner{\eta}{\xi}\xi\right\|^2
+ 2 c_2\left|\inner{\eta}{\xi}\right| + c_2^2.
$$
Now the norm squared of $\eta$ can be expressed in the form
$$
\left\|\eta\right\|^2 
  = u^2+v^2,\qquad
 u = \left\|\eta-\inner{\eta}{\xi}\xi\right\|,\qquad
 v = \left|\inner{\eta}{\xi}\right|.
$$
Hence the desired inequality has the form
$$
2 c_2\sqrt{u^2+v^2}
 \le u^2 + 2 c_2 v + c_2^2
$$
This follows from the inequalities $\sqrt{u^2+v^2}\le u+v$
and $2c_2u\le u^2+{c_2}^2$.  Thus we have proved~(\ref{eq:ab}).

Define $\alpha(s):=\int_0^s a(\sigma) \rd\sigma$. Then
$\alpha$ is nonnegative and $\dot\alpha= a$. 
Hence at each point $s\in[0,T)$ 
with $\dot\phi(s)\leq 0$ we have
$$
\frac{\rd}{\rd s} \bigl( e^{-\alpha}\dot\phi \bigr)
= e^{-\alpha} \bigl( \ddot\phi + a|\dot\phi| \bigr)
\geq - b.
$$
Integrating this inequality over maximal intervals where
$\dot\phi$ is negative we obtain
$$
e^{-\alpha(s)}\dot\phi(s)
\geq \min\{0,\dot\phi(0)\} 
- \int_0^sb(\sigma)\rd\sigma,\qquad
\text{for}\;\; 0\le s<T.
$$
This implies
$
\dot\phi(s) \ge - e^{a_0} (  b_0  + |\dot\phi(0)| ) ,
$
hence
$
\phi(s) \geq \phi(0) - e^{a_0}(  b_0  + |\dot\phi(0)| ) s ,
$
and hence, again for $0\le s<T$, 
$$
\log\|x(s)\| \geq 
\phi(s) - \int_0^s \|x\|^{-1}\|\dot x + A x\|
\geq \phi(0)-e^{a_0}(b_0+|\dot\phi(0)|)s - c_0 s.
$$
Now we can use 
$
\log\|x(0)\|=\phi(0)
$
and
$
|\dot\phi(0)|\leq\|x(0)\|^{-1}\|A(0)x(0)\|
$ 
to prove (ii):
$$
\|x(s)\| \geq
e^{\phi(0) - e^{a_0}(  b_0  + |\dot\phi(0)| ) s - c_0 s }
\geq
\|x(0)\| e^{-cs}.
$$
To prove (i) we assume by contradiction that $x(s_0)\neq 0$
for some $s_0\in(0,T)$. Then part (ii) applies to the path
$s\mapsto x(s_0-s)$ and the operator family $s\mapsto -A(s_0-s)$.
It implies $\|x(\sigma)\|\ge e^{c\sigma-cs_0} \|x(s_0)\|$ for
all $\sigma\in (0,s_0]$, so by continuity
$\|x(0)\|\ge e^{-cs_0} \|x(s_0)\|\neq 0$
in contradiction to the assumption.
\end{proof}

\subsection*{Time-dependent inner products}

There are interesting applications to operator
families $A(s)$ on a Hilbert space which are self-adjoint
with respect to a time-dependent family of inner products
which are all compatible with the standard inner product
on $H$.  Any such family of inner products can be expressed
in the form 
\begin{equation}\label{eq:Qt}
     \la x \,,\, y\ra_s = \la Q(s)x \,,\, Q(s)y \ra
\end{equation}
for some invertible bounded linear operators $Q(s):H\to H$.
Without loss of generality one can consider operators
$Q(s)$ which are self-adjoint.  Assume throughout that these 
operators satisfy the following conditions.

\begin{description}
\item[(Q1)]
  The operator $Q(s)$ is self-adjoint for every 
$s$ and there exists a constant $\delta>0$ such that 
for all $x\in H$ and $s\in[0,T)$
$$
     \delta\left\|x\right\|
     \le  \left\|Q(s)x\right\|
     \le  \delta^{-1}\left\|x\right\| .
$$
Moreover, the  map $[0,T)\to\cL(H):s\mapsto Q(s)$ is continuously 
differentiable in the weak operator topology 
and there exists a continuous function $c_Q:[0,T)\to[0,\infty)$ 
such that 
$$
\bigl\|\dot Q(s)\bigr\|_{\cL(H)} \le c_Q(s)
\quad\forall s\in[0,T),\qquad\quad
C_Q:=\int_0^T c_Q <\infty .
$$
\end{description}

\begin{thm}\label{thm:agni2}
Let $H$ be a real Hilbert space, $Q(s)\in\cL(H)$
a family of (boun\-ded) self-adjoint operators satisfying~$(Q1)$, 
and $A(s):\dom(A(s))\to H$ a family of (unbounded) linear operators 
such that $A(s)$ is symmetric with respect to the 
inner product~(\ref{eq:Qt}).
Assume that $x:[0,T)\to H$ is continuously 
differentiable in the weak topology
such that $x(s)\in\dom(A(s))$ and
\begin{equation*}\label{eq:agni1'}
    \left\|\dot x(s)+A(s)x(s)\right\|_s
    \le c_1(s)\left\|x(s)\right\|_s
\end{equation*}
for every $s\in[0,T)$.
Assume further that the function 
$s\mapsto\inner{x(s)}{A(s)x(s)}_s$ is also continuously 
differentiable and satisfies
\begin{multline*}\label{eq:agni2'}
 \frac{\rd}{\rd s}\inner{x(s)}{A(s)x(s)}_s 
   - 2\inner{\dot x(s)}{A(s)x(s)}_s  \\
    \le  2 c_2(s) \left\|A(s)x(s)\right\|_s\left\|x(s)\right\|_s
       + c_3(s) \left\|x(s)\right\|_s^2
\end{multline*}
for every $s\in[0,T)$.
Here $c_1, c_2, c_3 :[0,T)\to\R$ are continuous nonnegative
functions satisfying
\begin{equation*}\label{eq:ass cc}
\begin{split}
a_0&:= 2\int_0^T ( c_2 + \delta^{-1} c_Q ) <\infty, \\
b_0&:= \int_0^T \bigl( (c_1+\delta^{-1}c_Q)^2 
+ (c_2+\delta^{-1} c_Q)^2 + c_3 \bigr) <\infty,\\
c_0&:= \sup(c_1+\delta^{-1}c_Q) <\infty .
\end{split}
\end{equation*}
Then the following holds.

\smallskip\noindent{\bf (i)}
If $x(0)=0$ then $x(s)=0$ for all $s\in[0,T)$.

\smallskip\noindent{\bf (ii)}
If $x(0)\ne0$ then $x(s)\ne 0$ for all 
$s\in[0,T)$ and, moreover,
$$
\left\|x(s)\right\|_s
\ge e^{-c s} \left\|x(0)\right\|_0 ,\qquad
c:= c_0 + e^{a_0}\bigl(b_0 + \|x(0)\|_0^{-1} \| A(0) x(0) \|_0 \bigr) .
$$
\end{thm}

\begin{proof}  
The result reduces to Theorem~\ref{thm:agni1}.
Define 
$$
\tilde{A} := QAQ^{-1},\qquad
\tilde{x} := Qx,\qquad
\tilde{f} := \dot Qx+Qf
$$
with $\dom(\tilde{A}(s))=Q(s)\dom(A(s))$ and $f=\dot x+Ax$.
Then the operator $A(s)$ is symmetric with respect to
the inner product~(\ref{eq:Qt}) if and only
if $\tilde{A}(s)$ is symmetric with respect to the 
standard inner product.  (Moreover, one can easily
check that $A(s)$ is self-adjoint with respect to~(\ref{eq:Qt})
if and only if $\tilde{A}(s)$ is self-adjoint with respect 
to the standard inner product. However, this is not needed 
for the proof.)  It also easy to see that 
$$
\dot x+Ax=f\qquad\iff\qquad
\dot{\tilde{x}}+\tilde{A}\tilde{x}=\tilde{f}.
$$
It remains to show that under the assumptions of 
Theorem~\ref{thm:agni2} the triple $\tilde{A}$, $\tilde{x}$,
$\tilde{f}$ satisfies the requirements of Theorem~\ref{thm:agni1}.
First, note that
$$
\|\tilde f\| = \|\dot Qx + Qf\| 
\le c_Q\left\|x\right\| + \left\|f\right\|_s 
\le c_Q\delta^{-1}\left\|x\right\|_s + c_1\left\|x\right\|_s
$$
and hence $\tilde x$ satisfies~(\ref{eq:agni1})
with $c_1$ replaced by $\tilde{c}_1=c_1+c_Q/\delta$.
Secondly, the function 
$$
s\mapsto
\inner{\tilde x(s)}{\tilde{A}(s)\tilde x(s)}
= \inner{x(s)}{A(s)x(s)}_s
$$
is continuously differentiable and
a simple calculation shows that 
$$
\frac{\rd}{\rd s}\inner{\tilde x}{\tilde{A}\tilde x}
- 2\inner{\dot{\tilde{x}}}{\tilde{A}\tilde x}  
=
\frac{\rd}{\rd s}\inner{x}{Ax}_s
- 2\inner{\dot x}{Ax}_s 
- 2\inner{\dot Qx}{QAx}.
$$
Hence 
\begin{align*}
     \frac{\rd}{\rd s}\inner{\tilde x}{\tilde{A}\tilde x}
     - 2\inner{\dot{\tilde{x}}}{\tilde{A}\tilde x}  
&\le
      2 c_2\left\|x\right\|_s\left\|Ax\right\|_s
     + c_3\left\|x\right\|_s^2     
     + 2\|\dot Qx\|
        \|\tilde{A}\tilde{x}\|  \\
&\le
      2 c_2\left\|\tilde x\right\|
        \|\tilde{A}\tilde{x}\|
     + c_3\left\|\tilde x\right\|^2     
     + 2c_Q\delta^{-1}\left\|\tilde x\right\|
        \|\tilde{A}\tilde{x}\|.
\end{align*}
This shows that $\tilde x$ satisfies~(\ref{eq:agni2})
with $c_2$ and $c_3$ replaced by 
$\tilde{c}_2=c_2+c_Q/\delta$ and $\tilde{c}_3=c_3$.
Hence $\tilde x$ and $\tilde{A}$ satisfy the requirements 
of Theorem~\ref{thm:agni1} and this proves 
Theorem~\ref{thm:agni2}.
\end{proof}


\section{Holonomy perturbations}\label{app:Xf}

In this appendix we review the properties of the 
holonomy perturbations used in this paper. Throughout
this appendix $Y$ is a compact oriented $3$-manifold,
$\D\subset\C$ is the closed unit disc, 
and we identify the circle $S^1$ with $\R/\Z$.  
The elements of $S^1$ will be denoted by $\theta$ 
and those of $\D$ by $z$.
Fix a finite sequence of orientation preserving embeddings 
$\gamma_i:S^1\times\D\to Y$ for $i=1,\dots,N$ that coincide 
in a neighbourhood of $\{0\}\times\D$.  Define the holonomy maps
$$
g_i:\R\times\D\times\cA(Y)\to\rG,\qquad \rho_i:\D\times\cA(Y)\to\rG
$$
by
\begin{equation*}\label{eq:grhoi}
\p_\theta g_i + A(\p_\theta\gamma_i)g_i = 0,\qquad
g_i(0,z;A) = \one,\qquad \rho_i(z;A) := g_i(1,z;A).
\end{equation*}
and abbreviate $\rho:=(\rho_1,\dots,\rho_N):\D\times\cA(Y)\to\rG^N$.
Fix a smooth conjugation invariant function $f:\D\times\rG^N\to\R$
that vanishes near the boundary, and define the perturbation 
$h_f:\cA(Y)\to\R$ by 
$$
h_f(A):=\int_\D f(z,\rho(z;A))\,\rd^2 z.
$$
This map is smooth and its derivative has the form
\begin{equation}\label{eq:dhf}
\rd h_f(A)\alpha
= \frac\rd{\rd s}\biggr|_{s=0} \int_\D f(z,\rho(z;A+s\alpha))\,\rd^2 z
= \int_Y \winner{X_f(A)}{\alpha}
\end{equation}
for $\alpha\in\Om^1(Y,\cg)$. The map $X_f:\cA(Y)\to\Om^2(Y,\cg)$
is uniquely determined by~(\ref{eq:dhf});  it has the form
$$
X_f(A) = \sum_{i=1}^N {\gamma_i}_*\bigl(X_{f,i}(A)\rd^2 z\bigr),
$$
where $X_{f,i}(A)\in\Om^0(S^1\times\D,\cg)$ is given by
\begin{align}\label{eq:Xfi}
X_{f,i}(A) (\theta,z) &
= -g_i(\theta,z;A)\nabla_i f(z,\rho(z;A))g_i(\theta,z;A)^{-1} .
\end{align}
Here the gradient $\nabla_i f:\D\times\rG^N\to\cg$ is defined by 
$$
\la \nabla_i f(z,g),\xi\ra
:= \left.\frac\rd{\rd t}\right|_{t=0} 
f(z,g_1,\dots,g_{i-1},g_i \exp(t\xi),g_{i+1},\dots,g_N)
$$
for $z\in\D$, $g=(g_1,\dots,g_N)\in\rG^N$, and $\xi\in\cg$. 
It vanishes near the boundary of $\D$ and, 
since $f$ is conjugation invariant, it satisfies 
\begin{equation}\label{eq:hf}
\nabla_i f(z,hgh^{-1})=h\nabla_i f(z,g) h^{-1}
\end{equation}
for $h\in\rG$. If follows from~(\ref{eq:hf}) that $X_{f,i}(A):\R\times\D\to\cg$
descends to a function on $S^1\times\D$.  
If the center of $\rG$ is discrete then equation~(\ref{eq:hf}) 
implies that $\nabla_i f (z,(\one,\ldots,\one))=0$ and hence 
$X_f(0)=0$ for every $f\in\Cinf_0(\D\times\rG^N)^\rG$.
Thus, for $\rG=\SU(2)$ the trivial connection is always a critical point 
of the perturbed Chern-Simons functional $\CS_\cL+h_f$.
The next proposition summarizes the properties of $X_f$. 
We denote the space connections of class $W^{k,p}$ by
$$
\cA^{k,p}(Y):=W^{k,p}(Y,\rT^*Y\otimes\cg).
$$

\begin{proposition}\label{prop:Xf}
Let $f\in\cC^{\ell+1}_0(\D\times\rG^N)^\rG$
for some integer $\ell\ge 0$. Then the following holds 
(with uniform constants independent of $f$).

\smallskip\noindent{\bf (i)}
For every integer $\ell\ge k\ge1$ and every $p>2$ with $kp>3$, 
$X_f$ extends to a $\cC^{\ell-k}$ map from
$\cA^{k,p}(Y)$ to $W^{k,p}(Y,\Lambda^2\rT Y\otimes\cg)$,
mapping bounded sets to bounded sets.

\smallskip\noindent{\bf (ii)} 
For all $A\in\cA(Y)$, $u\in\cG(Y)$, $\xi\in\Om^0(Y,\cg)$,
and $\alpha\in\Om^1(Y,\cg)$ we have
\begin{equation*}\label{Xf1}
\rd_A (X_f(A))=0,\qquad
X_f(u^*A) = u^{-1}X_f(A)u,
\end{equation*}
\begin{equation*}\label{Xf2}
\rd X_f(A) \rd_A\xi = [X_f(A),\xi],\qquad
\rd_A(\rd X_f(A)\alpha) = [X_f(A)\wedge\alpha].
\end{equation*}

\smallskip\noindent{\bf (iii)}
For every $k\in\{0,\dots,\ell\}$ and every $p\in[1,\infty]$
there is a constant $c$ such that 
$$
\Norm{X_f(A)}_{W^{k,p}}
\le c \Norm{\nabla f}_{\cC^k}\left(1
+\sum_{j_0+\dots+j_s=k\atop s\ge0,\,j_\nu\ge1}
\Norm{A}_{W^{j_0,p}}\Norm{A}_{\cC^{j_1}}\cdots\Norm{A}_{\cC^{j_s}}
\right)
$$
for every $A\in\cA(Y)$. If $k=0$ 
then $\Norm{X_f(A)}_{L^p}\le c\Norm{\nabla f}_{L^p}$.

\smallskip\noindent{\bf (iv)}
For every $k\in\{0,\dots,\ell-1\}$ and $p\in[1,\infty]$
there is a constant $c$ such that 
\begin{equation*}\label{eq:Xf3}
\Norm{\rd X_f(A)\alpha}_{W^{k,p}(Y)}
\leq c \Norm{\nabla f}_{\cC^{k+1}} 
\bigl(1+\Norm{A}_{\cC^k}\bigr)^k\Norm{\alpha}_{W^{k,p}(Y)}
\end{equation*}
for all $A\in\cA(Y)$ and $\alpha\in\Om^1(Y,\cg)$.

\smallskip\noindent{\bf (v)}
For all $p,q,r\in[1,\infty]$ with $q^{-1}+r^{-1}=p^{-1}$
there is a constant $c$ such that 
\begin{align*}
\Norm{\rd^2 X_f(A)(\alpha,\beta)}_{L^{p}(Y)}
&\leq c \Norm{\nabla f}_{\cC^1}\Norm{\alpha}_{L^{q}(Y)}
\Norm{\beta}_{L^r(Y)},\\
\Norm{X_f(A+\alpha)-X_f(A)-\rd X_f(A)\alpha}_{L^p(Y)}
&\leq c \Norm{\nabla f}_{\cC^1}
\Norm{\alpha}_{L^q(Y)}\Norm{\alpha}_{L^{r}(Y)}
\end{align*}
for all $A\in\cA(Y)$ and $\alpha,\beta\in\Om^1(Y,\cg)$

\smallskip\noindent{\bf (vi)}
For every $p\in[1,\infty]$ there is a uniform constant $c$ 
such that 
\begin{equation}\label{eq:froyshov}
\Norm{\Nabla{A}(X_f(A))}_{L^p}\le 
c\bigl( 1+\Norm{\nabla f}_{\cC^1}\bigr)\bigl( 1+\Norm{F_A}_{L^p}\bigr)
\end{equation}
for every $A\in\cA(Y)$ 
\end{proposition}

\begin{rmk}\rm \label{rmk:Xf}
Consider a connection $\Xi=\Phi\ds+A\in\cA(I\times Y)$ 
for a compact interval $I$, given by $A:I\to\cA(Y)$ 
and $\Phi:I\to\Om^0(Y,\cg)$.
Proposition~\ref{prop:Xf} extends to the perturbation 
$X_f(\Xi):=X_f\comp A\in\Om^2(I\times Y,\cg)$
-- except for (ii), and in (i) we need to assume $kp>4$.
In particular, for every $k\ge1$ and $p>2$,
$X_f$ maps bounded sets in $\cA^{k,p}(I\times Y)$ 
to bounded sets in $W^{k,p}(I\times Y,\Lambda^2\rT Y\otimes\cg)$.
In the case $k=1$ and $kp\leq 4$ this follows from 
Proposition~\ref{prop:Xf}~(iii).
\end{rmk}

The proof requires some preparation.
We begin by considering connections on the circle.
The canonical $1$-form $d\theta\in\Om^1(S^1)$ 
allows us to identify the space $\cA(S^1)=\Om^1(S^1,\cg)$
of $\rG$-connections on $S^1$ with the space 
$\Om^0(S^1,\cg)$ of Lie algebra valued functions.
The holonomy of a connection $A=\eta\rd\theta\in\cA(S^1)$
with ${\eta:S^1\to\cg}$ is the solution $g:\R\to\rG$
of the differential equation
\begin{equation}\label{eq:holonomy}
\pd_\theta g + \eta g=0,\qquad g(0)=\one.
\end{equation}
The solutions give rise to a map
$
\hol :\R\times \Om^0(S^1,\cg) \to \rG
$
which assigns to each pair $(\theta,\eta)\in\R\times\Om^0(S^1,\cg)$
the value $\hol(\theta;\eta):=g(\theta)$ of the unique
solution of~(\ref{eq:holonomy}) at $\theta$. 
The gauge invariance of the holonomy takes the form
\begin{equation*}\label{eq:ghol}
\hol(\theta;u^{-1}\p_\theta u+u^{-1}\eta u)
= u(\theta)^{-1}\hol(\theta;\eta)u(0)
\end{equation*}
for $u:S^1\to\rG$.
One can think of $\hol$ as a map from
$\Om^0(S^1,\cg)$ to $\Cinf([0,1],\rG)$ defined by 
$\hol(\eta)(\theta):=\hol(\theta;\eta)$.  The holonomy then induces 
a map between Sobolev completions, 
for every integer $k\geq 0$ and every $p\geq 1$,
\begin{equation}\label{eq:hol}
\hol : W^{k,p}(S^1,\cg) \to W^{k+1,p}([0,1],\rG) .
\end{equation}
This map is continuously differentiable and its derivative 
at $\eta\in W^{k,p}(S^1,\cg)$ is the bounded linear operator
$
\rd\,\hol (\eta):W^{k,p}(S^1,\cg)
\to W^{k+1,p}([0,1],\hol(\eta)^*\rT\rG)
$
given by 
\begin{equation}\label{eq:dhol}
(\hol(\eta)^{-1}\rd\,\hol (\eta)\hat{\eta})(\theta) 
= - \int_0^\theta \hol(t;\eta)^{-1}\hat{\eta}(t)\hol(t;\eta) \dt 
\end{equation}
for $\hat{\eta}\in W^{k,p}(S^1,\cg)$. 
The formula~(\ref{eq:dhol}) shows, by induction, 
that the map~(\ref{eq:hol}) is smooth.  The next lemma is a 
parametrized version of this observation.

\begin{lemma}\label{le:holonomy}
Let $\Om$ be a compact Riemannian manifold.

\smallskip\noindent{\bf (i)}
For every integer $k\ge1$ and every $p>\dim\,\Om$,
composition with the holonomy induces smooth maps
\begin{align*}
\Hol:W^{k,p}(S^1\times \Om,\cg)&\to 
W^{k,p}([0,1]\times\Om,\rG) , \\
\Hol_1:W^{k,p}(S^1\times \Om,\cg)&\to W^{k,p}(\Om,\rG),
\end{align*}
given by $\Hol(\eta):=g_\eta$ 
with $g_\eta(\theta,x) := \hol(\theta,\eta(\cdot,x))$
and $\Hol_1(\eta)(x):=g_\eta(1,x)$
for $x\in\Om$ and $\theta\in[0,1]$.
These map $W^{k,p}$-bounded sets to $W^{k,p}$-bounded sets.

\smallskip\noindent{\bf (ii)}
For every integer $k\ge 1$ there is a constant $c$ such that 
$$
\Norm{\Hol(\eta)}_{W^{k,p}} + \Norm{\Hol_1(\eta)}_{W^{k,p}}
\le c 
\left(1+\sum_{j_0+\dots+j_s=k\atop s\ge0,\, j_\nu\ge1}
\Norm{\eta}_{W^{j_0,p}}\Norm{\eta}_{\cC^{j_1}}\cdots\Norm{\eta}_{\cC^{j_s}}
\right)
$$
for every $\eta\in \cC^k(S^1\times\Om,\cg)$
and every $p\in[1,\infty]$. 

\smallskip\noindent{\bf (iii)}
For every integer $k\ge0$ and every $p\in[1,\infty]$
there is a uniform constant $c$ such that,
for every $\eta\in C^k(S^1\times \Om,\cg)$, 
the derivatives
\begin{align*}
\Hol(\eta)^{-1}\rd\Hol(\eta):W^{k,p}(S^1\times \Om,\cg)
&\to W^{k,p}([0,1]\times\Om,\cg), \\
\Hol_1(\eta)^{-1}\rd\Hol_1(\eta):W^{k,p}(S^1\times \Om,\cg)
&\to W^{k,p}(\Om,\cg)
\end{align*}
are bounded linear operators with norms 
less than or equal to $c\,( 1 + \|\eta\|_{\cC^k})^k$. 
\end{lemma}

\begin{proof}
Think of $\eta$ as a map from $\Om$ to $W^{j,p}(S^1,\cg)$ and
of $\Hol(\eta)$ as a map from $\Om$ to $W^{j+1,p}([0,1],\rG)$.
Then $\Hol(\eta)$ is the composition 
$$
\Om\stackrel{\eta}{\longrightarrow}W^{j,p}(S^1,\cg)
\stackrel{\hol}{\longrightarrow}W^{j+1,p}([0,1],\rG).
$$
Since $\hol:W^{j,p}(S^1,\cg)\to W^{j+1,p}([0,1],\rG)$ is smooth
the composition induces a smooth map 
$$
\Hol:W^{\ell,p}(\Om,W^{j,p}(S^1,\cg))
\to W^{\ell,p}(\Om,W^{j+1,p}([0,1],\rG))
$$
for $\ell p>\dim\,\Om$ and any $j$; hence it defines a smooth map from
$$
W^{k,p}(S^1\times\Om,\cg)
= \bigcap_{\ell=0}^k W^{\ell,p}(\Om,W^{k-\ell,p}(S^1,\cg))
$$
to
$$
\bigcap_{\ell=1}^k W^{\ell,p}(\Om,W^{k-\ell+1,p}([0,1],\rG))
\subset W^{k,p}([0,1]\times\Om,\rG)
$$
for $k\ge 1$ and $p>\dim\,\Om$.  This proves~(i) for $\Hol$.
To prove~(i) for $\Hol_1$, take $\ell=k$ and note that
evaluation at $\theta=1$ gives a smooth map from
$W^{k,p}(\Om,W^{1,p}([0,1],\rG))$ to $W^{k,p}(\Om,\rG)$.
The boundedness of $\Hol$ and $\Hol_1$ is a consequence 
of (ii) and (iii).

To prove~(ii) we differentiate the function 
$g(\theta,x)=\hol(\theta,\eta(\cdot,x))$:
$$
g^{-1}\pd_\theta g = - g^{-1}\eta g,\qquad
(g^{-1}\pd_xg)(\theta,x) = -\int_0^\theta
g(t,x)^{-1}\pd_x\eta(t,x)g(t,x)\,\dt.
$$
Hence there are constants $c_1,c_2,c_3,\dots$ such that
\begin{equation}\label{eq:getakp}
\Norm{g}_{W^{k,p}} \le c_k
\left(1+\sum_{j_0+\dots+j_s=k\atop s\ge 0,\,j_\nu\ge1}
\Norm{\eta}_{W^{j_0,p}}\Norm{g}_{\cC^{j_1}}\cdots\Norm{g}_{\cC^{j_s}}
\right)
\end{equation}
for every smooth function $\eta:S^1\times\Om\to\cg$, 
every integer $k\ge1$, and every $p\in[1,\infty]$. 
For $p=\infty$ assertion~(ii) now follows by induction on $k$.
Inserting the resulting estimate into~(\ref{eq:getakp}) proves~(ii)
for all $p$. 
For $k=0$ assertion~(iii) follows immediately from~(\ref{eq:dhol})
with $c=1$. To prove~(iii) for $k\ge 1$ differentiate equation~(\ref{eq:dhol})
with respect to $\theta$ and $x$ and use~(ii). 
This proves the lemma.
\end{proof}

\begin{proof}[Proof of Proposition~\ref{prop:Xf}.]
The map $X_{f,i}:\cA^{k,p}(Y)\to W^{k,p}(S^1\times\D,\cg)$ can be 
expressed as composition of three maps.  The first is the product
of the $N$ maps
$$
\cA^{k,p}(Y)\to W^{k,p}(S^1\times\D,\cg):
A\mapsto \eta_j:=A(\p_\theta\gamma_j),
$$
the second is given by composition with the holonomy
$$
W^{k,p}(S^1\times\D,\cg)\to W^{k,p}([0,1]\times\D,\rG^2):
\eta_j\mapsto (g_j,\rho_j),
$$
where $g_j(\theta,z):=\hol(\eta_j(\cdot,z))(\theta)$ and
$\rho_j(\theta,z):=\hol(\eta_j(\cdot,z))(1)$, and the third 
map has the form
$$
W^{k,p}([0,1]\times\D,\rG^{2N})\to W^{k,p}(S^1\times\D,\cg):
(g_1,\rho_1,\dots,g_N,\rho_N)\mapsto \zeta ,
$$
with 
$$
\zeta:=g_i \nabla_if(\rho_1,\dots,\rho_N)g_i^{-1}
$$ 
(see equation~(\ref{eq:Xfi})). The first map is bounded linear 
(and hence smooth) for all $k$ and $p$ because composition 
with a smooth embedding at the source and multiplication 
with a smooth function define bounded linear maps 
between $W^{k,p}$-spaces.  
The second map is smooth and bounded for $k\ge 1$ and $p>2$ 
by Lemma~\ref{le:holonomy}. 
The third map is bounded and $\cC^{\ell-k}$ 
because composition with a $\cC^{k}$-map at the target defines a 
continuous map from $W^{k,p}$ to $W^{k,p}$ for all 
$kp\ge\dim(\R\times\D)$ (or $kp\ge\dim(\R^2\times\D)$ 
in the case of Remark~\ref{rmk:Xf}).
This proves~(i). 
Assertion~(ii) follows by straight forward calculations and~(iii)
follows from~(\ref{eq:Xfi}) and Lemma~\ref{le:holonomy}~(ii). 

To prove~(iv) we abbreviate $A_i:=\gamma_i^*A_i$, 
$\alpha_i:=\gamma_i^*\alpha_i$, and differentiate 
equation~(\ref{eq:Xfi}) to obtain
\begin{align*}
&\Hol(A_i)^{-1} \bigl( \rd X_{f,i}(A)\alpha \bigr) \,\Hol(A_i) \\
&\qquad\qquad\qquad
= \bigl[ \Hol(A_i)^{-1} X_{f,i}(A) \,\Hol(A_i) ,
         \Hol(A_i)^{-1} \rd\Hol(A_i)\alpha_i \bigr] \\
&\qquad\qquad\qquad\quad 
-  \sum_{j=1}^N \nabla_j\nabla_i f((\Hol_1(A_\ell))_{\ell=1,\dots N})  
\,\Hol_1(A_j)^{-1}\rd\Hol_1(A_j)\alpha_j .
\end{align*}
The estimate now follows from Lemma~\ref{le:holonomy} 
and the uniform bounds in (iii).

To prove (v) we differentiate the last equation again 
and obtain the inequality
\begin{align*}
&\bigl| \rd^2 X_f(A)(\alpha,\beta) \bigr| \\
&\leq
\bigl| \rd X_f(A)\alpha \bigr| \bigl| \rd\Hol(A)\beta \bigr|
+ \bigl| \rd X_f(A)\beta \bigr| \bigl| \rd\Hol(A)\alpha \bigr| \\
&\quad
+ |X_f(A)| 
\bigl| \rd\Hol(A)\alpha \bigr| \bigl| \rd\Hol(A)\beta \bigr|
+ |X_f(A)| 
\bigl|\rd\bigl(\Hol(A)^{-1}\rd\Hol(A)\alpha\bigr)\beta \bigr| \\
&\quad
+ |\nabla^2 f|
\bigl| \rd\Hol_1(A)\alpha \bigr| \bigl| \rd\Hol_1(A)\beta \bigr|
+ |\nabla^2 f|
\bigl|\rd\bigl(\Hol_1(A)^{-1}\rd\Hol_1(A)\alpha\bigr)\beta \bigr| 
\end{align*}
with
$$
\rd\bigl(\Hol(A)^{-1}\rd\Hol(A)\alpha\bigr) \beta
=
\int_0^\theta \bigl[ \Hol(A)^{-1} \rd\Hol(A)\alpha ,
\Hol(A)^{-1} \rd\Hol(A)\beta \bigr].
$$
A similar inequality holds for $\Hol_1$.  
The first estimate in~(v) now 
follows from the $L^q$- and $L^r$-bounds in (iv) and 
Lemma~\ref{le:holonomy} and the $L^\infty$-bounds on $X_f$
and $\nabla^2f$.  The second estimate in~(v) follows
from the first and
$$
X_f(A+\alpha)-X_f(A)-\rd X_f(A)\alpha 
= \int_0^1 \int_0^\tau \rd^2 X_f(A+t\alpha)(\alpha,\alpha) \dt \,\rd\tau .
$$
Assertion~(vi) is a result of Froyshov~\cite{FROY}. 
The proof uses the formula 
\begin{eqnarray} \label{eq:magic}
&&\p_tg(\theta,t)
+A(\p_t\gamma(\theta,t))g(\theta,t)  \nonumber \\
&&=g(\theta,t)
\biggl(\int_0^\theta 
g(s,t)^{-1}F_A(\p_\theta\gamma(s,t),\p_t\gamma(s,t))g(s,t)\,ds
\biggr)
\end{eqnarray}
for $\gamma:[0,1]^2\to Y$ and $g:[0,1]^2\to\rG$
with 
$$
\p_\theta g+A(\p_\theta\gamma)g=0,\qquad {g(0,t)=\one}.
$$
Namely, inserting a $t$-dependent parameter $z=z(t)$ into~(\ref{eq:Xfi}),
abbreviating 
$$
g(\theta,t):=g_i(\theta,z(t);A),\qquad
\gamma_i(\theta,t):=\gamma_i(\theta,z(t)),
$$ 
$$
\xi(\theta,t) 
:= X_{f,i}(\theta,z(t))  
 = -g(\theta,t)\nabla_i f(z(t),\rho(z(t);A))g(\theta,t)^{-1} ,
$$
and differentiating $\xi$ covariantly with respect to 
$\gamma_i^*A$ we find that $\Nabla{\theta}\xi=0$ and
\begin{equation}\label{eq:nablaxi}
\begin{split}
\Nabla{t}\xi \nonumber
&= \pd_t\xi + [A(\pd_t\gamma_i),\xi]  \\
&= \bigl[\bigl(\pd_tg\,g^{-1}+A(\pd_t\gamma_i)\bigr),\xi\bigr] 
- g^{-1}\bigl((\pd_1\Nabla{i}f)(z,\rho(z;A))\pd_t z\bigr)g \\
&\quad -\,g^{-1}\Biggl(
\sum_{j=1}^N (\Nabla{j}\Nabla{i}f)(z,\rho(z;A))
\rho_j(z;A)^{-1} \pd_t\rho_j(z;A )\Biggr)g. 
\end{split}
\end{equation}
Since the estimate~(\ref{eq:froyshov})
is gauge invariant and the $\gamma_j$
all coincide near $\gamma_j(0,z)=\gamma_j(1,z)$ 
we can assume that $A(\pd_t\gamma_j(1,z(t)))=0$ 
for all $j$ and $t$. Then it follows
from~(\ref{eq:magic}) that
$$
\rho_j(z(t);A)^{-1} \pd_t\rho_j(z(t);A)
=\int_0^1
g_j(s,t)^{-1}F_A(\pd_\theta\gamma_j(s,t),\pd_t\gamma_j(s,t))g_j(s,t)\,ds .
$$
So the first and third term on the right hand side 
of~(\ref{eq:nablaxi}) can be estimated by the curvature 
of $A$, and the second term is uniformly bounded. 
This proves the proposition. 
\end{proof}

In the remainder of this section we give a proof of the 
basic compactness result for solutions $\Xi\in\cA(\R\times Y)$
of the perturbed anti-self-duality equation 
\begin{equation}\label{eq:ASDp}
\bigl( F_{\Xi} + X_f(\Xi) \bigr)^+ = 0
\end{equation}
with bounded energy
$$
E_f(\Xi) = \tfrac 12 \int_{\R\times Y} \bigl| F_\Xi + X_f(\Xi) \bigr|^2 .
$$
A similar proof for somewhat different perturbations can be found
in~\cite{Kronheimer}.

\begin{thm}\label{thm:compactness}
There exists a universal constant $\hbar>0$ such that 
the following holds for every perturbation $X_f$, 
every real number $E>0$, and every $p>1$.

Let $\Xi_\nu\in\cA(\R\times Y)$ be a sequence of solutions 
of~(\ref{eq:ASDp}) with bounded energy 
$$
\sup_\nu E_f(\Xi_\nu)\leq E.
$$
Then there exists a subsequence (again denoted $(\Xi_\nu)$)
and a finite set of bubbling points 
$S=\{x_1,\ldots,x_N\}\subset\R\times{\rm int}(Y)$ with
$$
\liminf_{\nu\to\infty} \; \frac 12 
\int_{B_\delta(x_j)} \bigl| F_{\Xi_\nu} + X_f(\Xi_\nu) \bigr|^2 
\geq\hbar
\qquad \forall \delta>0 , x_j\in S .
$$
Moreover, there is a sequence of gauge transformations 
$u_\nu\in\cG((\R\times{\rm int}(Y))\setminus S)$ and a limit 
connection $\Xi_\infty\in\cA(\R\times{\rm int}(Y))$ such that
$u_\nu^*\Xi_\nu$ converges to $\Xi_\infty$ in the 
$W^{1,p}$-norm on every compact subset of
$\R\times{\rm int}(Y)\setminus S$.
The limit $\Xi_\infty$ solves (\ref{eq:ASDp}) 
and has energy 
$$
E_f(\Xi_\infty) \;\leq\;
\limsup_{\nu\to\infty} E_f(\Xi_\nu) - N\hbar .
$$
\end{thm}

\begin{rmk}\label{rmk:compacntess}\rm
If $S\subset(T_-,T_+)\times Y$ in Theorem~\ref{thm:compactness}, 
then the convergence can be improved to the 
$\cC^\infty$-topology on every compact subset of 
${(-\infty,T_-]\times{\rm int}(Y)}$ and
$[T_+,\infty)\times{\rm int}(Y)$
(in particular on $\R\times{\rm int}(Y)$ if $S=\emptyset)$.
This follows from the standard bootstrapping
techniques (e.g.~\cite{DK}, \cite{W}) and Remark~\ref{rmk:Xf}.
The crucial point is that a $W^{k,p}$-bound on
${u^\nu}^*\Xi^\nu$ implies a $W^{k,p}$-bound on
$X_{f}({u^\nu}^*\Xi^\nu)$ and thus on
$F_{{u^\nu}^*\Xi^\nu}^+$.
The appropriate gauge transformations can be interpolated  
to the ones over $(T_-,T_+)\times{\rm int}(Y)$.
\end{rmk}

\begin{proof}[Proof of Theorem~\ref{thm:compactness}.]
Without loss of generality we prove the theorem 
for a fixed constant $p>4$.
We follow the line of argument in \cite[4.4.4]{DK}.
Let ${\eps_{Uh}>0}$ and $C_{Uh}$ be the (universal) 
constants in Uhlenbeck's gauge fixing theorem 
(see~\cite{U} or \cite[Theorem~B]{W}). 
Then for each $x\in\R\times Y$, each sufficiently 
small constant $\delta>0$ with $B_\delta(x)\subset\R\times Y$, 
and each connection $\Xi\in\cA(\R\times Y)$ with energy
$$
\int_{B_\delta(x)}|F_\Xi|^2\leq\eps_{Uh}
$$ 
on the geodesic ball $B_\delta(x)$
there is a gauge transformation 
${u\in\cG(\R\times Y)}$ such that 
$$
\Norm{u^*\Xi}_{L^4(B_\delta(x))} + 
\Norm{u^*\Xi}_{W^{1,2}(B_\delta(x))}
\le C_{Uh} \Norm{F_\Xi}_{L^2(B_\delta(x))}.
$$

\noindent{\bf Step~1.}
{\it 
For every $\eps>0$ there is a finite set of bubbling points 
$S_\eps\subset\R\times\mathrm{int}(Y)$ and a subsequence,
still denoted by $\Xi_\nu$, such that the following holds.

\smallskip\noindent{\bf (a)}
If $x\in(\R\times\mathrm{int}(Y))\setminus S_\eps$ 
then there is a $\delta>0$ with
$
\sup_\nu\int_{B_\delta(x)}\Abs{F_{\Xi_\nu}}^2 \le \eps.
$

\smallskip\noindent{\bf (b)}
If $x\in S_\eps$ then
$
\inf_{\delta>0}\liminf_{\nu\to\infty} 
\int_{B_\delta(x)}\Abs{F_{\Xi_\nu}}^2\ge \eps/2.
$}

\medskip\noindent
Let $S_\eps$ be the set of points 
$x\in\R\times\mathrm{int}(Y)$ that satisfy
the inequality in~(b). Since $X_f(\Xi_\nu)$ 
is uniformly bounded we have
$$
\inf_{\delta>0}\liminf_{\nu\to\infty} 
\int_{B_\delta(x)}\Abs{F_{\Xi_\nu}+X_f(\Xi_\nu)}^2 
= \inf_{\delta>0}\liminf_{\nu\to\infty} 
\int_{B_\delta(x)}\Abs{F_{\Xi_\nu}}^2
\ge \frac{\eps}{2}
$$
for every $x\in S_\eps$ and hence the energy bound 
guarantees that $S_\eps$ contains at most $4E/\eps$ elements. 
If each point in $(\R\times\mathrm{int}(Y))\setminus S_\eps$
satisfies~(a) we are done. Otherwise there is a point
$x\in(\R\times{\rm int}(Y))\setminus S_\eps$ with
$$
\inf_{\delta>0}\sup_\nu 
\int_{B_\delta(x)} \Abs{F_{\Xi_\nu}}^2 
\ge \eps.
$$
In this case we can choose a subsequence 
(still denoted by $\Xi_\nu$) such that 
$$
\int_{B_{1/\nu}(x)}\Abs{F_{\Xi_\nu}}^2 
\ge \frac{\eps}{2}
$$
for all $\nu$. 
After passing to this subsequence we obtain a new
strictly larger set $S_\eps$.  Continue by induction.
The induction terminates when each point 
$x\in(\R\times{\rm int}(Y))\setminus S_\eps$ 
satisfies~(a).  It must terminate because in each step the 
the set $S_\eps$ contains at most $4E/\eps$ points. 

\medskip\noindent{\bf Step~2.}
{\it
We denote $q:=4p/(p+4)\in (2,4)$.
If $\eps>0$ is sufficiently small and $S=S_\eps$ is as in Step~1, 
then there exists a subsequence, still denoted by $\Xi_\nu$, 
and a sequence of gauge transformations 
${u_\nu\in\cG((\R\times Y)\setminus S)}$ 
such that $u_\nu^*\Xi_\nu$ converges to 
${\Xi_\infty\in\cA^{1,q}_{\rm loc}((\R\times {\rm int}(Y))\setminus S)}$ 
in the $W^{1,q}$-norm on every compact subset of 
$(\R\times\mathrm{int}(Y))\setminus S$.}

\medskip\noindent
There are universal constants $C_0\ge1$ and $C_1\ge1$ such that
\begin{equation}\label{eq:C01}
\Norm{\nabla\alpha}_{L^2}\leq 
C_0 \left( \Norm{\rd^+\alpha}_{L^2} 
+ \Norm{\rd^*\alpha}_{L^2}\right),\qquad
\Norm{\alpha}_{L^4}\leq C_1\Norm{\nabla\alpha}_{L^2}
\end{equation}
for ${\alpha\in\Om^1(B_{1}(0))}$ supported in the interior of the
Euclidean unit ball.  These inequalities are scale invariant, 
and for $\delta>0$ sufficiently small the metric in geodesic 
coordinates on $B_\delta(x)$ is $\cC^1$-close up to a conformal
factor to the Euclidean metric on $B_1(0)$. Hence the 
estimates~(\ref{eq:C01}) continue to hold with the same 
constants $C_0$ and $C_1$ for every compactly 
supported $1$-form on a geodesic ball $B_\delta(x)\subset\R\times Y$,
provided that $\delta>0$ is sufficiently small. 

Now fix
$0 <  \eps \le (4C_0C_1C_{Uh})^{-1}$
and choose a finite set $S=S_\eps\subset\R\times\mathrm{int}(Y)$ 
and a subsequence (still denoted by $\Xi_\nu$) 
as in Step~1.  
Since $\eps\le\eps_{Uh}$ it follows from Uhlenbeck's gauge that,
for every ${x\in(\R\times\mathrm{int}(Y))\setminus S}$,
there is a radius $\delta>0$ and a gauge transformation
$u_{\nu,x}\in\cG(B_\delta(x))$ such that
\begin{equation}\label{eq:Uh}
\Norm{u_{\nu,x}^*\Xi_\nu}_{W^{1,2}(B_\delta(x))}
\leq C_{Uh}\eps,\qquad
\rd^*(u_{\nu,x}^*\Xi_\nu)=0.
\end{equation}
By a global patching argument as in~\cite[Lemma~4.4.5]{DK} 
or~\cite[Proposition 7.6]{W}, it suffices to construct gauge 
transformations, limit connections, and establish the convergence
on every compact deformation retract $K\subset(\R\times{\rm int}(Y))\setminus S$.
We fix $K$ and find a covering by finitely many of the 
Uhlenbeck gauge neighbourhoods $B_{\delta_i}(x_i)$. 
On each of these $u_{\nu,x_i}^*\Xi_\nu$ satisfies~(\ref{eq:Uh}).
Now we fix a smooth cutoff function $h:B_{\delta_i}(x_i)\to[0,1]$ that
vanishes near the boundary. Then
\begin{align*}
& \tfrac12\Norm{h\cdot u_{\nu,x_i}^*\Xi_\nu}_{W^{2,2}} \\
&\leq C_0 \| h \nabla \rd^+ (u_{\nu,x_i}^*\Xi_\nu) \|_{L^2}
+ C \|u_{\nu,x_i}^*\Xi_\nu\|_{W^{1,2}} \\
&\leq C_0 \bigl\| h\nabla\bigl( u_{\nu,x_i}^{-1} X_f(\Xi_\nu)^+ u_{\nu,x_i}
- \tfrac 12 [ u_{\nu,x_i}^*\Xi_\nu \wedge u_{\nu,x_i}^*\Xi_\nu ]^+ 
\bigr)\bigr\|_{L^2} 
+ C \Norm{u_{\nu,x_i}^*\Xi_\nu}_{W^{1,2}} \\
&\leq 
 C_0 \Norm{u_{\nu,x_i}^*\Xi_\nu}_{L^4}
\bigl\|\nabla\bigl( h\cdot u_{\nu,x_i}^*\Xi_\nu\bigr)\bigr\|_{L^{4}} 
+ C \Norm{u_{\nu,x_i}^*\Xi_\nu}_{L^4}^2
+ C \Norm{u_{\nu,x_i}^*\Xi_\nu}_{W^{1,2}} \\
&\quad
+ C_0 \bigl( 
\Norm{u_{\nu,x_i}^{-1} 
\bigl(\nabla_{\Xi_\nu} X_f(\Xi_\nu)\bigr) u_{\nu,x_i}}_{L^2}
+ \Norm{u_{\nu,x_i}^*\Xi_\nu}_{L^2}\Norm{X_f(\Xi_\nu)}_{L^\infty}
\bigr) \\
&\leq 
 C_0 C_1C_{Uh}\eps\Norm{h\cdot u_{\nu,x_i}^*\Xi_\nu}_{W^{2,2}}
 + CC_{Uh}^2\eps^2 + CC_{Uh}\eps  \\
&\quad
+ C_0 \bigl( 
\Norm{u_{\nu,x_i}^{-1} 
\bigl(\nabla_{\Xi_\nu} X_f(\Xi_\nu)\bigr) u_{\nu,x_i}}_{L^2}
+ C_{Uh}\eps\Norm{X_f(\Xi_\nu)}_{L^\infty}
\bigr).
\end{align*}
Here all norms are in $B_{\delta_i}(x_i)$ and $C$ denotes 
a constant that only depends on $h$ and the radius $\delta_i$. 
In the first step we have used~(\ref{eq:C01}) with 
$\alpha=\p_i(h\cdot u_{\nu,x_i}^*\Xi_\nu)$,
$i=1,\dots,4$, and~(\ref{eq:Uh}).
In the last step we have used~(\ref{eq:Uh}) and the inequality
$$
\Norm{\nabla( h\cdot u_{\nu,x_i}^*\Xi_\nu)}_{L^{4}} 
\leq C_1 \Norm{h\cdot u_{\nu,x_i}^*\Xi_\nu}_{W^{2,2}}
$$
of~(\ref{eq:C01}). Since $C_0C_1C_{Uh}\eps\le1/4$ and 
$$
\Norm{\nabla_{\Xi_\nu} X_f(\Xi_\nu)}_{L^2(B_{\delta_i}(x_i))}
\leq C\bigl(1+\Norm{F_{\Xi_\nu}}_{L^2(I\times Y)}\bigr)
$$ 
for an interval $I\subset\R$ 
with $B_{\delta_i}(x_i)\subset I\times Y$
we obtain a $W^{2,2}$-bound on $u_{\nu,x_i}^*\Xi_\nu$
over a slightly smaller ball in $B_{\delta_i}(x_i)$
where $h\equiv 1$.

By Uhlenbeck's patching procedure 
(see~\cite{U} or~\cite[Chapter~7]{W})
the gauge transformations $u_{\nu,x_i}$ can then 
be interpolated to find $u_\nu\in\cG(K)$ such that 
$u_\nu^*\Xi_\nu$ is bounded in $W^{2,2}(K)$.
The compact Sobolev embedding ${W^{2,2}(K)\hookrightarrow W^{1,q}(K)}$
for $q<4$ then provides a $W^{1,q}$-convergent subsequence 
${u_\nu^*\Xi_\nu\to\Xi_\infty\in\cA^{1,q}(K)}$.

\medskip\noindent{\bf Step~3.}
{\it
We prove the theorem with $\hbar=\eps/4$ 
where $\eps$ is as in Step~2.
In particular, we remove the singularities to find
$\tilde\Xi_\infty\in\cA(\R\times{\rm int}(Y))$, a subsequence,
and gauge transformations 
$\tilde u_\nu\in\cG((\R\times{\rm int}(Y))\setminus S)$ 
such that $\tilde u_\nu^*\Xi_\nu\to\tilde\Xi_\infty$ 
in the $W^{1,p}$-norm on every compact subset
of $(\R\times {\rm int}(Y))\setminus S$.}

\medskip\noindent
Step~2 gives
$u_\nu^*\Xi_\nu\to\Xi_\infty\in
\cA^{1,q}_{\rm loc}((\R\times{\rm int}(Y))\setminus S)$
with $q>2$.
This implies $L^2$-convergence of the curvature 
on every compact subset of ${(\R\times{\rm int}(Y))\setminus S}$, 
and hence with the exhausting sequence
$K_\delta:=([-\delta^{-1},\delta^{-1}]\times Y)\setminus 
B_\delta(S\cup\R\times\pd Y)$
\begin{align*}
\int_{\R\times Y} | F_{\Xi_\infty} |^2
\;=\; \lim_{\delta\to 0}
\int_{K_\delta} | F_{\Xi_\infty} |^2 
\;=\; \lim_{\delta\to 0} \lim_{\nu\to\infty}
\int_{K_\delta} | F_{\Xi_\nu} |^2 \;\leq\; E . 
\end{align*}
Next we consider small annuli around the singularities and denote their
union, for $k\in\N$ sufficiently large, by 
$$
A_k:=B_{2^{1-k}}(S)\setminus B_{2^{-k}}(S).
$$
Then $\int_{A_k} |F_{\Xi_\infty}|^2 \to 0$ as $k\to\infty$
since the above limit exists.
For sufficiently large $k$ we can now patch Uhlenbeck gauges to obtain
a gauge transformation $u_k\in\cG(A_k)$ such that 
$\|u_k^*\Xi_\infty\|_{L^4(A_k)}\leq C \|F_{\Xi_\infty}\|_{L^2(A_k)} \to 0$.
The patching procedure does not introduce $k$-dependent constants or a flat
connection since the inequality is scale invariant and each annulus can be 
covered by two balls whose intersection
is connected and simply connected (see \cite[4.4.10]{DK}).

We extend $u_k$ to $(\R\times Y)\setminus S$ and denote 
$$
\Xi'_k:=(u_{\nu_k} u_k)^*\Xi_{\nu_k}.
$$ 
Here we pick a subsequence $\nu_k\to\infty$ such that 
$$
\bigl\|F_{u_{\nu_k}^*\Xi_{\nu_k}}\bigr\|_{L^2(A_{k_0})}^2
\le 2\Norm{F_{\Xi_\infty}}_{L^2(A_{k_0})}^2
$$ 
for all $k\geq k_0$ sufficiently large, and 
$$
\lim_{k\to\infty}\sup_{\ell\geq k}
\Norm{u_{\nu_\ell}^*\Xi_{\nu_\ell}-\Xi_\infty}_{L^4(A_{k})} = 0
$$
In particular, we have $\|\Xi'_k\|_{L^4(A_k)}\to 0$ as $k\to\infty$.
Now consider the sequence of extended connections
$$
\tilde\Xi_k:= h_k \cdot \Xi'_k \in\cA(\R\times Y), 
$$
where $h_k:\R\times Y\to[0,1]$ is a cutoff function that
vanishes on $B_{2^{-k}}(S)$, varies smoothly on $A_k$
with $|\rd h_k|\leq 2^{k+1}$, and equals to~$1$ on the complement
of $B_{2^{1-k}}(S)$.
The curvature of the extended connections is
\begin{align*}
F_{\tilde\Xi_k}
&= h_k\cdot F_{\Xi'_k} + \tfrac 12 (h_k^2-h_k) [\Xi'_k\wedge\Xi'_k] 
+ \rd h_k\wedge\Xi'_k .
\end{align*}
So for $\delta=2^{-\ell}$, and $k\geq\ell+1$
we have
\begin{align*}
&\int_{B_\delta(S)} |F_{\tilde\Xi_k}|^2 \\
&\leq
2 \int_{B_\delta(S)\setminus B_{2^{1-k}}(S)} |F_{\Xi_\infty}|^2
+ \int_{A_k} \Bigl( h_k^2|F_{\Xi'_k}|^2 + |h_k^2-h_k|^2 |\Xi'_k|^4 
+ |\rd h_k|^2\cdot |\Xi'_k|^2 \Bigr) \\
&\leq
2 \int_{B_\delta(S)\setminus B_{2^{-k}}(S)} |F_{\Xi_\infty}|^2
+ \|\Xi'_k\|_{L^4(A_k)}^4
+ 2^{2k+2} {\rm Vol}(A_k)^{\frac 12} \|\Xi'_k\|_{L^4(A_k)}^2 .
\end{align*}
The right hand side converges to 
$2 \int_{B_\delta(S)} |F_{\Xi_\infty}|^2$
as $k\to\infty$, so for sufficiently small 
$\delta=2^{-\ell}$ we have locally small energy
$\sup_k \int_{B_\delta(x)} |F_{\tilde\Xi_k}|^2\leq \eps$
at every $x\in\R\times{\rm int}(Y)$ 
for the subsequence $(\tilde\Xi_k)_{k\geq\ell}$.
(For $x\notin S$ this is true by Step~1.)

Now we can find an Uhlenbeck gauge $v_k\in\cG(B_{\delta}(S))$
such that 
\begin{equation}\label{eq:vkUh}
\rd^*(v_k^*\tilde\Xi_k)= 0,\qquad
\bigl\| v_k^*\tilde\Xi_k \bigr\|_{W^{1,2}(B_\delta(S))}\leq C_{Uh}\eps.
\end{equation}
The $W^{1,2}$-bound allows us to choose a 
$W^{1,2}$-weakly convergent subsequence 
$$
v_k^*\tilde\Xi_k\to\tilde\Xi_\infty\in\cA^{1,2}(B_\delta(S)).
$$
On the other hand, for every closed ball 
$D\subset B_\delta(S)\setminus S$
and every sufficiently large $k$ (such that $h_k|_D\equiv 1$)
the same estimate as in Step~2 provides
$W^{2,2}$-bounds on $v_k^*\tilde\Xi_k|_D$ and thus 
$W^{1,q}$-convergence $v_k^*\tilde\Xi_k\to\tilde\Xi_\infty
\in\cA^{1,q}_{\rm loc}(B_\delta(S)\setminus S)$
on every compact subset.

We can extend the gauge transformations $v_k\in\cG(B_\delta(S))$
by Uhlenbeck's patching procedure to a compact deformation retract
$S\subset K\subset\R\times{\rm int}(Y)$ 
(which is covered by $B_\delta(S)$ and finitely many balls in
$(\R\times{\rm int}(Y))\setminus S$ on which we also have an Uhlenbeck gauge
and hence $W^{2,2}$-bounds), and to $\R\times{\rm int}(Y)$ by the general 
extension procedure \cite[Proposition 7.6]{W}.
This provides a subsequence and gauge transformations 
$v_k\in\cG(\R\times{\rm int}(Y))$ such that the $v_k^*\tilde\Xi_k$ 
converge in the $W^{1,q}$-norm on every compact subset
of $(\R\times{\rm int}(Y))\setminus S$ to a limit connection
$\tilde\Xi_\infty\in\cA^{1,q}_{\rm loc}((\R\times {\rm int}(Y))\setminus S)$.
In particular, this means that
$$
\tilde u_{\nu_k}^*\Xi_{\nu_k}\to\tilde\Xi_\infty,\qquad
\tilde u_{\nu_k}:=u_{\nu_k}u_k v_k
\in\cG((\R\times{\rm int}(Y))\setminus S),
$$
because $\tilde\Xi_k=(u_{\nu_k}u_k)^*\Xi_{\nu_k}$
on compact subsets of $(\R\times{\rm int}(Y))\setminus S$.
Moreover, the limit connection extends to $S$ such that
$v_k^*\tilde\Xi_k\to\tilde\Xi_\infty\in\cA^{1,2}(B_\delta(S))$
converges $W^{1,2}$-weakly and $L^4$-weakly.

Since $\tilde\Xi_\infty$ is of class $W^{1,2}$, 
the perturbation 
$X_f(\tilde\Xi_\infty)\in{L^\infty(\R\times Y)}$
is well defined, and we claim that 
\begin{equation}\label{eq:Xflim}
 X_f(v_k^*\tilde\Xi_k) \to X_f(\tilde\Xi_\infty) ,
\qquad
\tu_{\nu_k}^{-1} X_f(\Xi_{\nu_k}) \tu_{\nu_k} \to X_f(\tilde\Xi_\infty)
\end{equation}
in the $L^p$-norm on every compact subset of $\R\times Y$.
If $S$ does not intersect the support 
$
\supp X_f:=\bigcup_{i=1}^N \R\times\im\gamma_i
$
of the perturbation 
then $\tilde\Xi_\infty|_{\supp X_f}$ is the 
$W^{1,q}_{\rm loc}$-limit of
$
v_k^*\tilde\Xi_k|_{\supp X_f}
=v_k^*\Xi_k'|_{\supp X_f}
=\tu_{\nu_k}^* \Xi_{\nu_k}|_{\supp X_f}
$
and the claim follows directly from Remark~\ref{rmk:Xf}
and the Sobolev embedding ${W^{1,q}\hookrightarrow L^p}$ on compact
subsets of $\R\times Y$.
If $S$ does intersect the set $\supp X_f$ at some points
${(s_j,\gamma_{i_j}(\theta_j,z_j))_{j=1,\dots,n}\subset S}$, 
then we have
$$
X_f(v_k^*\tilde\Xi_k)
=v_k^{-1} X_f(h_k\Xi_k') v_k
=v_k^{-1} X_f(\Xi_k') v_k
=\tu_{\nu_k}^{-1} X_f(\Xi_{\nu_k}) \tu_{\nu_k}
$$
only on the complement of a solid cylinder neighbourhood $Z_k$
of the loops $(s_j,\gamma_{i_j}(S^1,z_j))\subset\R\times Y$.
More precisely, $Z_k\subset\R\times{\rm int}(Y)$ is given by
the union of all loops $(s,\gamma_i(S^1,z))$
that intersect the support of $1-h_k$.
It thus is a union of solid cylinders whose width is 
of order $2^{1-k}$.
If we fix the cylinder neighbourhood $Z_{k_0}$,
then the previous argument still applies for $k\geq k_0$
to give $L^p$-convergence on the complement of $Z_{k_0}$.
The remaining $Z_{k_0}$ has volume of
order $2^{3-3k_0}$, and the perturbations 
$X_f(\tilde\Xi_\infty)$, $X_f(v_k^*\tilde\Xi_k)$, 
and $X_f(\Xi_{\nu_k})$
are all 
uniformly bounded by Proposition~\ref{prop:Xf}~(iii) 
(with $k=0$). 
So we see that
$\|X_f(v_k^*\tilde\Xi_k) - X_f(\tilde\Xi_\infty)\|_{L^p(Z_{k_0})}$
and
$\|\tu_{\nu_k}^{-1} X_f(\Xi_{\nu_k}) \tu_{\nu_k} 
- X_f(\tilde\Xi_\infty)\|_{L^p(Z_{k_0})}$ 
also converge to zero as we let $k\geq k_0\to\infty$.
This proves~(\ref{eq:Xflim}).

A first consequence is that the limit connection satisfies
\begin{equation}\label{eq:asdlim}
\bigl( F_{\tilde\Xi_\infty} + X_f(\tilde\Xi_\infty) \bigr)^+ = 0
\end{equation}
because this is the local weak $L^2$-limit of
$\bigl(F_{v_k^*\tilde\Xi_k}+ X_f(v_k^*\tilde\Xi_k)\bigr)^+$
and

\begin{align*}
& \bigl\| \bigl(F_{v_k^*\tilde\Xi_k} 
+ X_f(v_k^*\tilde\Xi_k)\bigr)^+\bigr\|_{L^2(\R\times Y)} \\
&= \bigl\| \bigl(F_{h_k\Xi_k'} + X_f(h_k\Xi_k') \bigr)^+
-  (u_{\nu_k}u_k)^{-1}\bigl(F_{\Xi_{\nu_k}} 
+ X_f(\Xi_{\nu_k}) \bigr)^+ (u_{\nu_k}u_k)
\bigr\|_{L^2(\R\times Y)} \\
&\leq 
\bigl\| F_{h_k\Xi_k'} -  F_{\Xi_k'} \bigr\|_{L^2(A_k)} 
+ \bigl\| X_f(h_k\Xi_k') 
- (u_{\nu_k}u_k)^{-1} X_f(\Xi_{\nu_k})  (u_{\nu_k}u_k)
\bigr\|_{L^2(Z_k)} ,
\end{align*}
which converges to zero by similar estimates as before.
Another consequence is the energy identity:
We have 
$$
F_{v_k^*\tilde\Xi_k}+X_f(v_k^*\tilde\Xi_k)\to
F_{\tilde\Xi_\infty}+X_f(\tilde\Xi_\infty)
$$
in the $L^2$-norm on every compact subset 
of $(\R\times{\rm int}(Y))\setminus S$.
So, exhausting $\R\times Y$ with 
$$
K_\delta:=\bigl([-\delta^{-1},\delta^{-1}]\times Y \bigr) 
\setminus B_\delta(S\cup\R\times\pd Y),
$$
we have
\begin{align*}
&E_f(\tilde\Xi_\infty) \\
&= \lim_{\delta\to 0} \;\frac 12 \int_{K_\delta} 
\bigl| F_{\tilde\Xi_\infty} + X_f(\tilde\Xi_\infty) \bigr|^2 
\;=\;
 \lim_{\delta\to 0}\lim_{k\to\infty} \;\frac 12 \int_{K_\delta} 
\bigl| F_{\tilde\Xi_k} + X_f(\tilde\Xi_k) \bigr|^2 \\
&\leq \lim_{\delta\to 0} \lim_{k\to\infty} \;\frac 12 \biggl(
\int_{K_\delta} \bigl| F_{\Xi_{\nu_k}}+X_f(\Xi_{\nu_k}) \bigr|^2
+ \int_{K_\delta} \bigl| F_{h_k\Xi_k'} -  F_{\Xi_k'} \bigr|^2 \\
&\qquad\qquad\qquad\quad
+ \int_{K_\delta} \bigl| X_f(h_k\Xi_k') 
- (u_{\nu_k}u_k)^{-1} X_f(\Xi_{\nu_k})  (u_{\nu_k}u_k)\bigr|^2
 \biggr) \\
&= \lim_{\delta\to 0} \lim_{k\to\infty} \frac 12 \biggl(
\int_{[-\delta^{-1},\delta^{-1}]\times Y} 
\bigl| F_{\Xi_{\nu_k}}+X_f(\Xi_{\nu_k}) \bigr|^2 
- \int_{B_\delta(S)} \bigl| F_{\Xi_{\nu_k}}+X_f(\Xi_{\nu_k}) \bigr|^2 
\biggr) \\
&\leq \limsup_{\nu\to\infty} E_f(\Xi_\nu) - N\hbar .
\end{align*}
Here $\hbar:=\eps/4$ with $\eps>0$ as in Step~2.

It follows from~(\ref{eq:vkUh}) and~(\ref{eq:asdlim}) that
$$
\rd^*\tilde\Xi_\infty = 0,\qquad
\bigl\|\tilde\Xi_\infty\bigr\|_{L^4(B_\delta(S))}\leq C_{Uh}\eps,\qquad
F_{\tilde\Xi_\infty}^+\in L^\infty(\R\times Y).
$$
This implies $\tilde\Xi_\infty\in\cA^{1,3}(B_{\delta/2}(S))$ 
by a standard argument as in~\cite[Proposition~4.4.13]{DK}, 
using the estimate
$$
\bigl\|\A\bigr\|_{W^{1,3}}
\le C\left(
\bigl\|\rd^*\A\bigr\|_{L^3}
+ \bigl\|F_\A^+\bigr\|_{L^3}
+ \bigl\|\A\bigr\|_{L^4}\bigl\|\A\bigr\|_{W^{1,3}}
\right)
$$
for compactly supported $\A\in\cA(B_\delta(S))$.  Hence we have 
$\tilde\Xi_\infty\in\cA^{1,3}_{\rm loc}(\R\times{\rm int}(Y))$.
Now the standard regularity theory for anti-self-dual connections
(e.g.\ \cite[Chapter~9]{W}) together with Remark~\ref{rmk:Xf},
for control of the perturbation, provides another gauge transformation 
that makes $\tilde\Xi_\infty$ smooth and does not affect the convergence.

It remains to strengthen the convergence
$$
\tilde u_{\nu_k}^*\Xi_{\nu^k}\to\tilde\Xi_\infty
$$
on $(\R\times{\rm int}(Y))\setminus S$ to the 
$W^{1,p}_{\rm loc}$-topology.
Again, it suffices to construct the required subsequence
and gauge transformations on a compact deformation retract 
$K\subset(\R\times{\rm int}(Y))\setminus S$.
We pick a compact submanifold 
$$
M\subset(\R\times{\rm int}(Y))\setminus S
$$
such that $K\subset{\rm int}(M)$ and apply the local slice theorem 
(e.g.\ \cite[Theorem~8.1]{W}) to find gauge transformations 
$u_{\nu_k}\in\cG(M)$ such that
$$
\rd_{\tilde\Xi_\infty}^*(u_{\nu_k}^*\Xi_{\nu_k} 
- \tilde\Xi_\infty) = 0,\qquad
\lim_{k\to\infty}\bigl\|u_{\nu_k}^*\Xi_{\nu_k} 
- \tilde\Xi_\infty\bigr\|_{W^{1,q}(M)}=0.
$$
Since $\tu_{\nu_k}^*\Xi_{\nu_k}|_M$ has the same $W^{1,q}$-limit,
the gauge transformations $\tu_{\nu_k}^{-1}u_{\nu_k}\in\cG(M)$
converge, for a further subsequence, in the weak $W^{2,q}(M)$-topology 
to an element $u_\infty$ of the isotropy subgroup of $\tilde\Xi_\infty$.
We can make sure that this limit is in fact $\one$, by modifying
$u_{\nu_k}$ to $u_{\nu_k}u_\infty^{-1}$ in the local slice gauge.
With this we have
$$
\lim_{k\to\infty}
\bigl\|u_{\nu_k}^*\Xi_{\nu_k}
-\tilde\Xi_\infty\bigr\|_{L^p(M)}=0,\qquad
\lim_{k\to\infty}
\bigl\|\rd^*(u_{\nu_k}^*\Xi_{\nu_k}-\tilde\Xi_\infty)\bigr\|_{L^p(M)}=0,
$$
so we can use the elliptic estimate for $\rd^+\oplus\rd^*$ on $M$.
For that purpose fix a cutoff function $h:M\to[0,1]$ with $h|_K\equiv 1$
and $h\equiv 0$ near $\pd M$.
Then
\begin{align*}
\bigl\| \rd^+\bigl(
&h (u_{\nu_k}^*\Xi_{\nu_k}-\tilde\Xi_\infty)\bigr) \bigr\|_{p} \\
&\leq
C_h \bigl\| u_{\nu_k}^*\Xi_{\nu_k} - \tilde\Xi_\infty \bigr\|_{p}
+ 
\bigl\| X_f(\tilde\Xi_\infty) 
- u_{\nu_k}^{-1} X_f(\Xi_{\nu_k}) u_{\nu_k} \bigr\|_{p} \\
&\quad
+ \bigl\| h[\tilde\Xi_\infty\wedge \tilde\Xi_\infty]^+ 
- h[u_{\nu_k}^*\Xi_{\nu_k}\wedge u_{\nu_k}^*\Xi_{\nu_k}]^+ 
\bigr\|_{p}.
\end{align*}
Here the constant $C_h:=\|\nabla h\|_\infty$ is finite, so the first
term converges to zero as $k\to\infty$.
The second term also converges to zero due to (\ref{eq:Xflim}) 
and the $\cC^0$-convergence $\tu_{\nu_k}^{-1}u_{\nu_k}\to\one$.
Finally, the third term can be bounded by the constant
${\bigl( 2 \|\tilde\Xi_\infty\|_{L^\infty}
+ C_S\| h( u_{\nu_k}^*\Xi_{\nu_k} - \tilde\Xi_\infty ) \|_{W^{1,p}} \bigr)
\| u_{\nu_k}^*\Xi_{\nu_k} - \tilde\Xi_\infty \|_{L^p}}$
with a constant $C_S$ from the Sobolev embedding
$W^{1,p}(M)\hookrightarrow\cC^0(M)$.
Now apply the elliptic estimate for $\rd^+\oplus\rd^*$ to the 
compactly supported $1$-form
$\eta_k:=h( u_{\nu_k}^*\Xi_{\nu_k} - \tilde\Xi_\infty )$ 
to obtain
$$
\|\eta_k\|_{W^{1,p}} \leq 
C \bigl( 1 + \| \eta_k \|_{W^{1,p}} \bigr) 
\bigl\|u_{\nu_k}^*\Xi_{\nu_k} - \tilde\Xi_\infty \bigr\|_p
+ C \bigl\| X_f(\tilde\Xi_\infty) 
- u_{\nu_k}^{-1} X_f(\Xi_{\nu_k}) u_{\nu_k} \bigr\|_p
$$
with a finite constant $C$.
Since $\|u_{\nu_k}^*\Xi_{\nu_k} - \tilde\Xi_\infty\|_{L^p(M)}\to 0$
this can be rearranged to prove that
$$
\bigl\|u_{\nu_k}^*\Xi_{\nu_k} - \tilde\Xi_\infty\bigr\|_{W^{1,p}(K)}
\leq \Norm{\eta_k}_{W^{1,p}(M)}\to 0.
$$
This finishes the proof of Step~3 and the theorem.
\end{proof}


\section{The Lagrangian and its tangent bundle}\label{app:Lag}
 
For any compact manifold $X$, any integer $k\geq 0$, and any $p>1$ 
we denote the space of $W^{k,p}$-connections by
$$
\cA^{k,p}(X) := W^{k,p}(X,\rT^*X\otimes\cg).
$$
If $(k+1)p>{\rm dim}\, X$ then the gauge group
$$
\cG^{k+1,p}(X) := W^{k+1,p}(X,\rG) 
$$
acts smoothly on $\cA^{k,p}(X)$. For $p=\infty$ we denote by
$\cA^{k,\infty}(X)$ the space of $\cC^k$-connections;
similarly for $\cG^{k,\infty}(X)$.

Let $Y$ be a compact oriented Riemannian $3$-manifold with boundary 
$\pd Y=\Sigma$ and $\cL\subset\cA(\Sigma)$ be a gauge invariant
Lagrangian submanifold 
(in the sense of (L1) of the introduction)
such that $\cL/\cG_z(\Sigma)$ is compact.
For $(k+1)p>2$ the $W^{k,p}$-closure of 
$\cL$ is a Banach submanifold of $\cA^{k,p}(\Sigma)$, which we denote
by $\cL^{k,p}$. (This follows from the Sobolev embedding
$W^{k,p}(\Sigma)\hookrightarrow L^q(\Sigma)$ with $q>2$
and the fact that the $L^q$-Banach submanifold coordinates 
in~\cite[Lemma~4.3]{W Cauchy} restrict to $W^{k,p}$-coordinates.)
Again, we denote by $\cL^{k,\infty}$ the $\cC^k$-completion.
Denote
\begin{align*}
 \cA^{k,p}(Y,\cL) &:= \{ A\in\cA^{k,p}(Y) \st A|_{\pd Y} \in \cL^{0,q} \} .
\end{align*}
This is a Banach submanifold of $\cA^{k,p}(Y)$ for $(k+1)p>3$ since
the restriction map $\cA^{k,p}(Y)\hookrightarrow \cA^{0,q}(\Sigma)$ 
with $q>2$ is smooth and transverse to~$\cL$.
Theorem~\ref{thm:expmap} will provide a gauge equivariant 
exponential map for $\cA^{1,p}(Y,\cL)$, 
from which we construct an exponential map for $\cA^{1,p}(\R\times Y,\cL;B_-,B_+)$
in Corollary~\ref{cor:expmap4}.

Moreover, consider the vector bundle $\cE\to\cA(Y,\cL)$ with fibre
$$
\cE_A:=\Om^1_A(Y,\cg) 
=\bigl\{\alpha \in \Om^1(Y,\cg) \st *\alpha|_{\pd Y} = 0,\, 
\alpha|_{\pd Y} \in \rT_A\cL \bigr\} .
$$
In Theorem~\ref{thm:Q} below we construct local trivializations of $\cE$.
In a preliminary step we construct local trivializations of the tangent 
bundle of $\cL$.
Note that these trivializations extend to the fibrewise $L^2$-closure of the 
tangent bundle although it is not known whether the $L^2$-closure of $\cL$ 
is smooth.

\begin{thm} \label{thm:P}
For every $A_0\in\cL$ there exists a neighbourhood
$\cU\subset\cL$ of $A_0$ (open in the $\cC^0$-topology)
and a family of bijective linear operators 
$$
P_A: \Om^1(\Sigma,\cg) \to \Om^1(\Sigma,\cg),
$$
parametrized by $A\in\cU$, such that the following holds.
\begin{description}
\item[(i)]
$P_{A_0}=\one$.
\item[(ii)]
For every $A\in\cU$ and every $\alpha\in\Om^1(\Sigma,\cg)$ 
we have 
$$
P_A \alpha \in\rT_A\cL 
\qquad\iff \qquad \alpha \in \rT_{A_0}\cL.
$$
\item[(iii)]
For every integer $k\ge 0$ and every $p>1$ the operator $P_A$ extends 
to a Banach space isomorphism from $W^{k,p}(\Sigma,\rT^*\Sigma\otimes\cg)$ 
to itself; this extended operator depends smoothly on $A\in\cL^{k,\infty}$
with respect to the operator norm on $P_A$.
\item[(iv)]
For every integer $k\ge 0$, every $p>1$, every $\lambda\in[0,1]$,
and every $A\in\cU^{k,\infty}$ the operator $\lambda\one + (1-\lambda)P_A$ 
extends to a Banach space isomorphism from 
$W^{k,p}(\Sigma,\rT^*\Sigma\otimes\cg)$ to itself.
Here $\cU^{k,\infty}$ denotes the interior 
of the closure of $\cU$ in $\cL^{k,\infty}$.
\end{description}
\end{thm}

\begin{proof}
Choose a $3$-dimensional subspace $E\subset\Om^0(\Sigma,\cg)$ such that 
the restriction of $\rd_{A_0}:\Om^0(\Sigma,\cg)\to\Om^1(\Sigma,\cg)$
to $E^\perp$ (the $L^2$-orthogonal complement of $E$) is injective.  
Then there is a constant $C$ such that 
$\|\xi\|_{W^{1,2}} \leq C \|\rd_{A_0}\xi\|_{L^2}$
for all $\xi\in E^\perp$.
This estimate continues to hold for each $A\in\cL$ 
that is sufficiently close to $A_0$ in the $\cC^0$-norm.
Hence there is a $\cC^0$-open neighbourhood $\cU\subset\cL$ of $A_0$ 
such that $\rd_{A}: E^\perp \to\Om^1(\Sigma,\cg)$ is injective
for every $A\in\cU$.  
Define 
$$
H^1_{A,E}:=\left\{\alpha\in\Om^1(\Sigma,\cg)\st
*\rd_{A}\alpha\in E,\,
\rd_{A}^*\alpha\in E\right\}.
$$
Then, for every $A\in\cU$, there is a generalized Hodge decomposition
\begin{equation}\label{eq:hodge}
\Om^1(\Sigma,\cg) = H^1_{A,E} 
\oplus \rd_{A}(E^\perp) \oplus *\rd_{A}(E^\perp).
\end{equation}
The three summands in~(\ref{eq:hodge}) are orthogonal 
to each other and the generalized Hodge decomposition
extends to each Sobolev completion 
$\cA^{k,p}(\Sigma)$ in the usual fashion.
This uses the fact that the operator 
$$
\laplace_{A,E}:=\rd_{A}^*\rd_{A}:
\Om^0(\Sigma,\cg)\supset E^\perp\to\Om^0(\Sigma,\cg)/E
$$
extends to an isomorphism from $W^{k+2,p}$ to $W^{k,p}$ 
(with $p>1$) for every ${A\in\cU}$. 
(The operators $\laplace_{A,E}$ are all injective and 
compact perturbations of the isomorphism $\laplace_{A_0,E}$.)
The standard Hodge decomposition corresponds 
to the case $E=\ker\,\rd_{A}$. 
The reason for our construction with $E$ independent of $A$
is the need for a Hodge decomposition
which depends smoothly on $A$.

The Lagrangian submanifold $\cL$ gives rise to another 
$L^2$-orthogonal decomposition, 
$\Om^1(\Sigma,\cg)=\rT_A\cL \oplus * \rT_A\cL$,
see \cite[Lemma~4.2]{W Cauchy}.
Since $\rd_{A}(E^\perp)\subset\rT_{A}\cL$ 
and $*\rd_{A}(E^\perp)$ is perpendicular 
to $\rT_{A}\cL$ it follows from~(\ref{eq:hodge}) that 
we have $\rT_{A}\cL = \Lambda_A \oplus \rd_{A}(E^\perp)$,
where 
$$
\Lambda_A:=H^1_{A,E}\cap \rT_{A}\cL
$$ 
is a Lagrangian subspace of $H^1_{A,E}$.
Hence there is a refined Hodge decomposition
\begin{equation}\label{eq:hodge L}
\Om^1(\Sigma,\cg) = \Lambda_A \oplus *\Lambda_A
\oplus \rd_{A}(E^\perp) \oplus *\rd_{A}(E^\perp).
\end{equation}
For $A\in\cU$ we define a bijective linear operator 
$P_A:\Om^1(\Sigma,\cg)\to\Om^1(\Sigma,\cg)$ by
$$
P_A(\alpha_0 + *\beta_0 +\rd_{A_0}\xi + *\rd_{A_0}\eta)
:=\Pi_A\alpha_0 + *\Pi_A\beta_0 +  \rd_{A}\xi + *\rd_{A}\eta
$$
for $\alpha_0,\beta_0\in \Lambda_{A_0}$ and
 $\xi,\eta\in E^\perp\subset\Om^0(\Sigma,\cg)$,
where
$$
\Pi_A:\Om^1(\Sigma,\cg)\to\Lambda_A
$$
denotes the $L^2$-orthogonal projection.
(Shrink $\cU$, if necessary, so that the restriction of $\Pi_A$ to 
$\Lambda_{A_0}$ is a vector space isomorphism for every $A\in\cU$.)
Note that $P_{A_0}=\Id$ and $P_A\alpha\in\rT_A\cL$ 
iff $\alpha\in\rT_{A_0}\cL$.  We claim that each operator 
$P_A$ extends to a Banach space automorphism of 
$\rT_A\cA^{k,p}(\Sigma)=W^{k,p}(\Sigma,\rT^*\Sigma\otimes\cg)$ 
for all $k$ and $p$, and this automorphism depends smoothly on 
$A\in\cL^{k,\infty}$.  To prove this we write $P_A$ as the composition 
of three linear operators.  The first is the Banach space 
isomorphism 
$$
W^{k,p}(\Sigma,\rT^*\Sigma\otimes\cg) 
\to \Lambda_{A_0} \times \Lambda_{A_0} 
\times W^{k+1,p}_E(\Sigma,\cg)\times W^{k+1,p}_E(\Sigma,\cg)
$$
induced by the Hodge decomposition for $A_0$.
Here $W^{k+1,p}_E(\Sigma,\cg)$ denotes the $L^2$-orthogonal complement of $E$
in $W^{k+1,p}(\Sigma,\cg)$.
The second operator is the restriction of $\Pi_A$ on the factors 
$\Lambda_{A_0}$ and is the identity on the factors 
$W^{k+1,p}_E(\Sigma,\cg)$.  We think of the target space of this 
second operator as the product
$$
W^{k,p}(\Sigma,\rT^*\Sigma\otimes\cg) 
\times W^{k,p}(\Sigma,\rT^*\Sigma\otimes\cg)
\times W^{k+1,p}_E(\Sigma,\cg)\times W^{k+1,p}_E(\Sigma,\cg).
$$
The third operator maps this product to
to $W^{k,p}(\Sigma,\rT^*\Sigma\otimes\cg)$ via
$$
(\alpha,\beta,\xi,\eta)\mapsto \alpha + * \beta + \rd_A\xi + *\rd_A\eta.
$$
The first operator is independent of $A$ and the third 
depends smoothly on ${A\in\cL^{k,\infty}}$. By the Hodge 
decomposition for $A$ it restricts to an isomorphism from
$\Lambda_A \times \Lambda_A \times 
W^{k+1,p}_E(\Sigma,\cg)\times W^{k+1,p}_E(\Sigma,\cg)$
to $W^{k,p}(\Sigma,\rT^*\Sigma\otimes\cg)$. 
It remains to prove that the map
$$
\cU^{k,\infty}\to\Hom(\Lambda_{A_0},W^{k,p}(\Sigma,\rT^*\Sigma\otimes\cg)):
A\mapsto\Pi_A
$$
is smooth.  To see this we 
write $\Pi_A$ as the composition of two projections
$$
\Pi_A = \Pi_{H^1_{A,E}}\circ \Pi_{\rT_A\cL}|_{\Lambda_{A_0}}.
$$
Here 
$
\Pi_{H^1_{A,E}} : 
W^{k,p}(\Sigma,\rT^*\Sigma\otimes\cg) 
\to W^{k,p}(\Sigma,\rT^*\Sigma\otimes\cg)
$
denotes the $L^2$-orthogonal projection onto $H^1_{A,E}$
given by
$$
\Pi_{H^1_{A,E}}\alpha := \alpha 
- \rd_A \laplace_{A,E}^{-1} (\rd_A^*\alpha) 
+ *\rd_A \laplace_{A,E}^{-1} (*\rd_A\alpha) .
$$
It depends smoothly on $A\in\cL^{k,\infty}\cap\cU$ since the same holds for 
the operator $\laplace_{A,E}:W^{k+1,p}_E(\Sigma,\cg)\to W^{k-1,p}(\Sigma,\cg)/E$
and its inverse. The operator
$$
\Pi_{\rT_A\cL}: 
W^{k,p}(\Sigma,\rT^*\Sigma\otimes\cg) \to W^{k,p}(\Sigma,\rT^*\Sigma\otimes\cg) 
$$
denotes the $L^2$-orthogonal projection onto $\rT_A \cL^{k,p}$.
For $(k+1)p>2$ we know that $\cL^{k,p}\subset\cA^{k,p}(\Sigma)$ 
is a Banach submanifold, so $\Pi_{\rT_A\cL}$ depends smoothly 
on $A\in\cL^{k,p}$, and this proves that 
$\Pi_A$ depends smoothly on $A\in\cL^{k,\infty}$.
In the case $(k+1)p\le2$, i.e.\ $k=0$, $p\leq 2$, we have 
${\cA^{0,3}(\Sigma)\subset\cA^{k,p}(\Sigma)}$.
The $L^p$- and the $L^3$-norm are 
equivalent on the finite dimensional space
$\Lambda_{A_0}\subset\Om^1(\Sigma,\cg)$.
Hence $\Pi_A$ is the composition of the projection
${\Pi_{\rT_A\cL} : L^3(\Sigma,\rT^*\Sigma\otimes\cG)
 \to L^3(\Sigma,\rT^*\Sigma\otimes\cG)}$,
restricted to $\Lambda_{A_0}$, the inclusion
$L^3(\Sigma,\rT^*\Sigma\otimes\cG)\hookrightarrow 
L^p(\Sigma,\rT^*\Sigma\otimes\cG)$,
and the projection $\Pi_{H^1_{A,E}}:L^p(\Sigma,\rT^*\Sigma\otimes\cG)
\to L^p(\Sigma,\rT^*\Sigma\otimes\cG)$.
All of these depend smoothly on $A\in\cL^{0,\infty}$.

To prove (iv) shrink $\cU$ such that $\|\one-P_A\|_{\cL(L^2)}\leq1/2$
for all $A\in\cU$. Then $\lambda\one + (1-\lambda)P_A$ is invertible on $L^2$
for every $\lambda\in[0,1]$ and every $A\in\cU^{0,\infty}$.
Invertibility on $W^{k,p}$ for $A\in\cU^{k,\infty}$ now follows from elliptic
regularity for the Laplace operator.
This proves the theorem.
\end{proof}

\begin{thm} \label{thm:Q}
For every $A_0\in\cA(Y,\cL)$ there is a neighbourhood
${\cU\subset\cA(Y,\cL)}$ of $A_0$ (open in the $\cC^0$-topology)
and a family of bijective linear operators 
$$
Q_A: \Om^1(Y,\cg) \to \Om^1(Y,\cg),
$$
parametrized by $A\in\cU$, such that the following holds.
\begin{description}
\item[(i)]
$Q_{A_0}=\one$.
\item[(ii)]
For every $A\in\cU$ and every $\alpha\in\Om^1(Y,\cg)$ 
we have 
$$
Q_A \alpha \in\Om^1_A(Y,\cg) 
\qquad\iff \qquad \alpha \in \Om^1_{A_0}(Y,\cg).
$$
Moreover, $*(Q_A\alpha)|_{\pd Y}=*\alpha|_{\pd Y}$ and
$(Q_A \alpha)|_{\pd Y}=0$ iff $\alpha|_{\pd Y}=0$.
\item[(iii)]
For every integer $k\ge 0$ and every $p>1$ the operator $Q_A$ extends 
to a Banach space isomorphism from $W^{k,p}(Y,\rT^*Y\otimes\cg)$ to itself; 
this extended operator depends smoothly on $A\in\cA^{k,\infty}(Y,\cL)$
with respect to the operator norm on $Q_A$.
\end{description}
\end{thm}

\begin{proof}
Choose geodesic normal coordinates to identify a neighbourhood 
of $\pd Y$ with the product $(-\eps,0]\times\Sigma$ via an orientation 
preserving embedding
$$
\iota:(-\eps,0]\times\Sigma\to Y.
$$
For a connection $A\in\cA(Y,\cL)$ and a $1$-form $\alpha\in\Om^1(Y,\cg)$
we write the pullbacks under $\iota$ in the form
\begin{equation}\label{eq:alpha}
\iota^*A =: B(t) +\Psi(t)\,\dt,\qquad
\iota^*\alpha =: \beta(t)+\psi(t)\,\dt.
\end{equation}
Then  $B(0)=A|_\Sigma\in\cL$.  
Choose a neighbourhood $\cU_0\subset\cL$ of $B_0:=A_0|_\Sigma$
(open in the $\cC^0$-topology) and an operator family 
$P_B:\Om^1(\Sigma,\cg)\to\Om^1(\Sigma,\cg)$,
parametrized by $B\in\cU_0$, which satisfies the requirements of
Theorem~\ref{thm:P}. Then we have $P_{B_0}=\one$.
Now
$$
\cU:=\left\{A\in\cA(Y,\cL)\st A|_\Sigma\in\cU_0\right\}
$$
is a $\cC^0$-open neighbourhood of $A_0$.
For $A\in\cU$ we define the bijective linear operator 
${Q_A:\Om^1(Y,\cg)\to\Om^1(Y,\cg)}$ 
by 
$$
\iota^*(Q_A\alpha)
:=h(t)\beta(t) + (1-h(t)) P_{A|_\Sigma}\beta(t)+ \psi(t)\,\dt
$$
for $\iota^*\alpha$ of the form~(\ref{eq:alpha}),
and by $Q_A\alpha :=\alpha$ outside of the image of $\iota$.  
Here $h:(-\eps,0]\to[0,1]$ is a smooth cutoff function that 
vanishes near $0$ and equals to $1$ near $-\eps$.
The operator family $\{Q_A\}_{A\in\cU}$ satisfies conditions~(i)-(iii).
\end{proof} 
 
The construction of exponential maps will be based on the following.

\begin{lem}\label{lem:lagg}
Fix a constant $p>2$. There is an open neighbourhood 
$$
\cU^{0,p}\subset L^{p}(\Sigma,\rT^*\Sigma\otimes\cg)
$$ 
of zero and a smooth map
$$
\cL^{0,p}\times\cU^{0,p}\to \cA^{0,p}(\Sigma):
(A,\alpha) \mapsto \Theta_A(\alpha)
$$
satisfying the following conditions:
\begin{description}
\item[(i)]
For every $A\in\cL^{0,p}$ the map $\Theta_A:\cU^{0,p}\to\cA^{0,p}(\Sigma)$
is a diffeomorphism from $\cU^{0,p}$ onto an $L^p$-open neighbourhood 
of $A$ in $\cA^{0,p}(\Sigma)$ such that $\Theta_A(0)=A$ and
$D\Theta_A(0)=\Id$. In particular, there is a uniform constant $C$ such that
\begin{align*}
\| \Theta_A(\alpha) - \Theta_A(\alpha') \|_{L^p} &\leq C \|\alpha-\alpha'\|_{L^p} \\
\| D\Theta_A(\alpha)\beta - D\Theta_A(\alpha')\beta \|_{L^p} &\leq C \|\alpha-\alpha'\|_{L^p} \|\beta\|_{L^p}
\end{align*}
for all $A \in \cL^{0,p}$, $\alpha\in\cU^{0,p}$, $\beta\in L^p(\Sigma,\rT^*\Sigma\otimes\cg)$.
\item[(ii)]
$\Theta$ is gauge equivariant in the sense that
for $u\in\cG^{1,p}(\Sigma)$
$$
\Theta_{u^*A}(u^{-1}\alpha u) = u^*\Theta_A(\alpha) .
$$
\item[(iii)]
For every $A\in\cL^{0,p}$
$$
\Theta_A (\rT_A\cL^{0,p}\cap\cU^{0,p}) 
= \cL^{0,p}\cap\Theta_A(\cU^{0,p}).
$$
\item[(iv)]
For every integer $k\ge 1$ and every $A\in\cL^{k,p}$ the restriction of $\Theta_A$
to the intersection $\cU^{k,p}:=\cU^{0,p}\cap W^{k,p}$
is a diffeomorphism onto its (open) image in $\cA^{k,p}(\Sigma)$.
It depends smoothly on $A\in\cL^{k,p}$ and satisfies 
\begin{align*}
\| \Theta_A(\alpha) - \Theta_A(\alpha') \|_{W^{1,p}} &\leq C (1 + \|A\|_{L^\infty} ) \|\alpha-\alpha'\|_{W^{1,p}} , \\
\| D\Theta_A(\alpha)\beta - D\Theta_A(\alpha')\beta \|_{W^{1,p}} &\leq C(1 + \|A\|_{L^\infty} ) 
\|\alpha-\alpha'\|_{W^{1,p}} \|\beta\|_{W^{1,p}}
\end{align*}
for all $A\in\cL^{1,p}$, $\alpha,\alpha'\in\cU^{1,p}$, and $\beta\in W^{1,p}(\Sigma,\rT^*\Sigma\otimes\cg)$
with a uniform constant~$C$.
\item[(v)]
The restriction of $\Theta$ to an open neighbourhood of the 
zero section in the subbundle $*\rT\cL^{0,p}\subset\cL^{0,p}\times\cU^{0,p}$
is a diffeomorphism onto an open neighbourhood 
$\cW^{0,p}\subset\cA^{0,p}(\Sigma)$ of $\cL^{0,p}$.
The composition of its inverse with the projection onto $\cL^{0,p}$
$$
\pi : \cW^{0,p} \to \cL^{0,p} 
$$
is gauge equivariant and maps $\cW^{k,p}:=\cW^{0,p}\cap W^{k,p}$
to $\cL^{k,p}$ for every $k$.
\end{description}
\end{lem}

\begin{proof}
Since $\cL^{0,p}/\cG^{1,p}(\Sigma)$ is compact it suffices to provide the construction
for smooth $A\in\cL$. The smooth extension to $\cL^{0,p}$ is then provided by the equivariance (ii).
For every smooth connection $A\in\cL$ we have an $L^2$-orthogonal
direct sum decomposition from \cite[Lemma~4.2]{W Cauchy},
\begin{equation}\label{eq:TL}
L^{p}(\Sigma,\rT^*\Sigma\otimes\cg) 
= \rT_A\cL^{0,p} \oplus * \rT_A\cL^{0,p}.
\end{equation}
Moreover, $\rT_A\cL^{0,p}=L_A\oplus\rd_A W^{1,p}(\Sigma,\cg)$,
where $L_A:=\rT_A\cL^{0,p}\cap h^1_A\subset\Om^1(\Sigma,\cg)$ is
the intersection of $\rT_A\cL$ with the harmonic (and thus smooth)
$1$-forms 
$$
{h^1_A:=\ker\rd_A\cap\ker\rd_A^*\subset\Om^1(\Sigma,\cg)}.
$$
We denote the $L^2$--orthogonal projection in (\ref{eq:TL}) by
$$
\pi_A:L^{p}(\Sigma,\rT^*\Sigma\otimes\cg)\to \rT_A\cL^{0,p} .
$$
It smoothly depends on $A\in\cL$, is gauge equivariant
$\pi_{u^*A}(u^{-1}\alpha u)=u^{-1}\pi_A(\alpha) u$, and satisfies 
${\rd_A^*\comp\pi_A = \rd_A^*}$ because ${\im\rd_A\subset\rT_A\cL^{0,p}}$.
By standard Hodge theory, this projection restricts to a 
bounded linear operator from the subspace
$W^{k,p}(\Sigma,\rT^*\Sigma\otimes\cg)$
to $\rT_A\cL^{k,p}=L_A\oplus\rd_A W^{k+1,p}(\Sigma,\cg)$ 
for every integer $k\ge1$. 
For each $A\in\cL$ the map 
$$
\cL^{0,p}\to\rT_A\cL^{0,p}:B\mapsto\pi_A(B-A)
$$
is smooth and its differential at $B=A$ is the identity.
Hence it restricts to a diffeomorphism from an 
$L^p$--open neighbourhood of $A$ onto an open
set 
$$
\cV_A^{0,p}\subset \rT_A\cL^{0,p}. 
$$
We denote its inverse by 
$$
\psi_A:\cV_A^{0,p}\to\cL^{0,p}.
$$
It follows immediately from the definition that $\psi$ is smooth 
and gauge equivariant in the sense that
$$
\psi_{u^*A}(u^{-1}\alpha u) = u^*\psi_A(\alpha)
$$
for all $A\in\cL$, $u\in\cG(\Sigma)$ and $\alpha\in\cV^{0,p}_A$.
Its differential at $0$ is the identity, $D\psi_A(0)={\rm Id}$, hence on 
a small ball $\{\|\alpha\|_{L^p}\leq\delta\}\supset\cV^{0,p}_A$, we 
can bound the $L^p$-operator norm $\|D\psi_A(\alpha)\|\leq 2$, and thus
obtain a linear estimate for all $\alpha,\alpha'\in\cV^{0,p}_A$
$$
\| \psi_A(\alpha) - \psi_A(\alpha') \|_{L^p} \leq 
\int_0^1 \| D\psi_A(t\alpha+(1-t)\alpha')\| \|\alpha-\alpha'\|_{L^p}
\leq 2 \|\alpha-\alpha'\|_{L^p} .
$$
Similarly, since $D\psi_A$ is continuously differentiable, we obtain 
for all $\alpha,\alpha'$ in (the possibly smaller) $\cV^{0,p}_A$ and
all $\beta\in \rT_A\cL^{0,p}$
$$
\| D\psi_A(\alpha)\beta - D\psi_A(\alpha')\beta \|_{L^p} \leq C \|\alpha-\alpha'\|_{L^p} \|\beta\|_{L^p}
$$
with a uniform constant $C$.
(In fact, $C$ is also independent of $A\in\cL$ since the estimates are gauge invariant 
and $\cL/\cG(\Sigma)$ is compact).
In particular, we have
$$
\| \psi_A(\alpha) - A \|_{L^p} \leq 2 \|\alpha\|_{L^p} ,\qquad 
\| D\psi_A(\alpha)\beta - \beta \|_{L^p} \leq C \|\alpha\|_{L^p} \|\beta\|_{L^p} .
$$
Moreover, $\psi_A$ maps the intersection $\cV_A^{k,p}:=\cV_A^{0,p}\cap W^{k,p}$
to $W^{k,p}$-regular points in $\cL^{k,p}$ because $F_{\psi_A(\alpha)}=0$ and 
$$
\rd_A^*(\psi_A(\alpha)-A)
\;=\;\rd_A^*\bigl(\pi_A(\psi_A(\alpha)-A)\bigr)
\;=\;\rd_A^*\alpha
\;\in\; W^{k-1,p}(\Sigma,\cg).
$$
In fact, we obtain an estimate for all $A\in\cL^{1,p}$, $\alpha,\alpha'\in\cV_A^{1,p}$ 
(denoting all uniform constants by $C$)
\begin{align*}
& \| \psi_A(\alpha) - \psi_A(\alpha') \|_{W^{1,p}}  \\
&\leq C \bigl( 
 \bigl\| \rd \bigl(  \psi_A(\alpha) - \psi_A(\alpha') \bigr) \bigr\|_{p} 
+ \bigl\| \rd^* \bigl(  \psi_A(\alpha) - \psi_A(\alpha') \bigr) \bigr\|_{p} 
+ \bigl\|  \psi_A(\alpha) - \psi_A(\alpha') \bigr\|_{p} \bigr) \\
&\leq
C \bigl( \| \psi_A(\alpha) - \psi_A(\alpha') \bigr\|_{p} + \|\psi_A(\alpha')-A\|_{p} \bigr)
 \| \psi_A(\alpha) - \psi_A(\alpha') \bigr\|_{\infty} \\
&\quad   
+ C\bigl\| \rd_A^* (\alpha - \alpha') \bigr\|_{p} 
+ C\|A\|_{\infty} \|\psi_A(\alpha) - \psi_A(\alpha') \|_{p}   
+ C\| \alpha -\alpha' \|_{p}  \\
&\leq C (1+ \|A\|_{\infty}) \| \alpha - \alpha' \|_{W^{1,p}}
+ C \bigr( \| \alpha - \alpha' \|_{p} + \| \alpha' \|_{p} \bigl)  
\| \psi_A(\alpha) - \psi_A(\alpha') \bigr\|_{W^{1,p}} .
\end{align*}
If we choose $\cV_A^{0,p}$ sufficiently small, 
then the second term can be absorbed into the left hand side, which proves 
$$
\| \psi_A(\alpha) - \psi_A(\alpha') \|_{W^{1,p}} 
\leq C (1+ \|A\|_{L^\infty}) \| \alpha - \alpha' \|_{W^{1,p}} 
\qquad \forall A\in\cL^{1,p} , \alpha,\alpha'\in\cV_A^{1,p} .
$$
Note that this estimate does not simply follow from smoothness of $\psi_A$ since
$\cV_A^{1,p}$ is not even bounded in the $W^{1,p}$-norm.
Similarly, we obtain uniform estimates for the linearization $D\psi_A$ of $\psi_A$ 
using the identities $\rd_{\psi_A(\alpha)}( D\psi_A(\alpha)\beta ) = 0 = \rd_A\beta$
and $\rd_A^*(D\psi_A(\alpha)\beta) = \rd_A^*\beta$,
\begin{align*}
& \| D\psi_A(\alpha)\beta -  D\psi_A(\alpha')\beta \|_{W^{1,p}}  \\
&\leq C \bigl( 
 \bigl\| \rd \bigl(  D\psi_A(\alpha)\beta - D\psi_A(\alpha')\beta \bigr) \bigr\|_{p} 
+ \bigl\| \rd^* \bigl(  D\psi_A(\alpha)\beta - D\psi_A(\alpha')\beta \bigr) \bigr\|_{p} \\
&\qquad
+ \bigl\|  D\psi_A(\alpha)\beta - D\psi_A(\alpha')\beta \bigr\|_{p} \bigr) \\
&\leq C( 1+\|A\|_{\infty} ) \| (D\psi_A(\alpha) - D\psi_A(\alpha'))\beta \|_{p}  
+ \| \psi_A(\alpha)-\psi_A(\alpha') \|_{\infty} \|D\psi_A(\alpha)\beta \|_{p}  \\
&\qquad
+ C\| \psi_A(\alpha')-A \|_{p} \bigl\| D\psi_A(\alpha)\beta - D\psi_A(\alpha')\beta \bigr\|_{\infty} \\
&\leq C ( 1+\|A\|_{\infty}) \| \alpha- \alpha'\|_{p} \|\beta \|_{p} 
+ C\|\alpha-\alpha' \|_{W^{1,p}} (1+ \| \alpha' \|_{p}) \|\beta \|_{p}  \\
&\quad
+ C\| \alpha \|_{L^p} \bigl\| D\psi_A(\alpha)\beta - D\psi_A(\alpha')\beta \bigr\|_{W^{1,p}} .
\end{align*}
For $\cV_A^{0,p}$ sufficiently small, this can be rearranged to
$$
\| D\psi_A(\alpha)\beta -  D\psi_A(\alpha')\beta \|_{W^{1,p}} 
\leq C ( 1+\|A\|_{L^\infty}) \| \alpha- \alpha'\|_{W^{1,p}} \|\beta \|_{L^p} .
$$
Now choose an open neighbourhood 
$\cU^{0,p}\subset L^p(\Sigma,\rT^*\Sigma\otimes\cg)$ of $0$
such that 
$$
\pi_A(\cU^{0,p})\subset\cV^{0,p}_A
$$
for every $A\in\cL$.  Then the map
$\Theta_A:\cU^{0,p}\to\cA^{0,p}(\Sigma)$ defined by
$$
\Theta_A(\alpha) := \psi_A(\pi_A(\alpha))+\alpha-\pi_A(\alpha)
$$
has the required properties. 
The estimates for $\Theta_A$ follow from the linearity of $\pi_A$ 
and the linear estimates for $\psi_A$.
To check (v) note that the differential of $\Theta|_{*\rT\cL^{0,p}}$ at $(A,0)$
is the isomorphism 
${\rT_A\cL^{0,p}\times*\rT_A\cL^{0,p}\to L^p(\Sigma,\rT^*\Sigma\otimes\cg)}$,
${(\eta,\beta)\mapsto \eta + \beta}$.
So the restriction of $\Theta$ to $*\rT\cL^{0,p}$ is a local diffeomorphism
near the zero section. 
To see that it is globally injective we assume by contradiction that
${\Theta_{A_i}(\alpha_i)=\Theta_{B_i}(\beta_i)}$ for some 
$A_i,B_i\in\cL^{0,p}$ and some $\alpha_i,\beta_i\in*\rT_{A_i}\cL^{0,p}$ 
with ${\|\alpha_i\|_{L^p}+ \|\beta_i\|_{L^p}\to 0}$. 
Since $\Theta$ is equivariant and $\cL^{0,p}/\cG^{1,p}(\Sigma)$ is compact,
we can assume w.l.o.g.\ $A_i\to A_\infty$ and $u_i^*B_i\to \tilde B_\infty$ 
in the $\cC^\infty$-topology for some $u_i\in\cG^{1,p}(\Sigma)$.
Then $\Theta_{A_i}(\alpha_i)\to A_\infty$ and
${u_i^*\Theta_{A_i}(\alpha_i) 
= \Theta_{u_i^*B_i}(u_i^{-1}\beta_i u_i)\to\tilde B_\infty}$,
so we can find a convergent subsequence $u_i\to u_\infty\in\cG(\Sigma)$.
Consequently $B_i\to u_\infty^{-1\;*}\tilde B_\infty=A_\infty$ 
has the same limit as $A_i$, in contradiction 
to the local injectivity of $\Theta|_{*\rT\cL^{0,p}}$.
\end{proof}

\begin{thm} \label{thm:expmap}
Fix a constant $p>2$ and a compact subset $N\subset\cA^{1,p}(Y,\cL)$. 
Then there is an open neighbourhood 
$\cU\subset \rT\cA^{1,p}(Y,\cL)$
of the zero section over $N$ and a smooth map
$$
\cU \to \cA^{1,p}(Y,\cL) : (A,\alpha)\mapsto E_A(\alpha)
$$
satisfying the following conditions:
\begin{description}
\item[(i)]
For every $A\in\cA^{1,p}(Y,\cL)$ the map 
$E_A:\cU\cap\rT_A\cA^{1,p}(Y,\cL)\to\cA^{1,p}(Y,\cL)$
is a diffeomorphism from a neighbourhood of $0$ onto a 
neighbourhood of $A$ such that $E_A(0)=A$ and
$\rd E_A(0)=\Id$.
\item[(ii)]
$E$ is gauge equivariant in the sense that
for $u\in\cG^{1,p}(Y)$
$$
E_{u^*A}(u^{-1}\alpha u) = u^*E_A(\alpha) .
$$
\end{description}
\end{thm}  
\begin{proof} 
Our construction will be based on the two maps 
from Lemma~\ref{lem:lagg},
$$
{\Theta:\cL^{0,p}\times\cU^{0,p}\to\cA^{0,p}(\Sigma)},\qquad
\pi:\cW^{0,p}\to\cL^{0,p} .
$$ 
We start by fixing a tubular neighbourhood 
$\tau:(-1,0]\times\Sigma\hookrightarrow Y$ of
the boundary $\pd Y\cong\{0\}\times\Sigma$
such that $\tau^*A|_{\{t\}\times\Sigma}\in\cW^{0,p}$
for all $A\in N$ and $t\in(-1,0]$.
This is possible since 
$\tau^*N\subset W^{1,p}((-1,0]\times\Sigma)
\subset\cC^0((-1,0],\cA^{0,p}(\Sigma))$
is compact.

On the complement of the image of $\tau$ we define 
$E_A(\alpha):=A+\alpha$.
On the image of $\tau$ write 
$\tau^*A=B(t)+\Psi(t)\dt$ and $\tau^*\alpha=\beta(t)+\psi(t)\dt$,
where $\beta(t)\in\cU^{0,p}$ can be ensured by the choice of
neighbourhood $\cU\ni\alpha$ of the zero section.
With this we can
define $\tau^*E_A(\alpha) := \tB + (\Psi+\psi)\dt $ by
$$
\tB(t) 
:= B(t) + \rho(t)\bigl( \Theta_{\pi(B(t))}(\beta(t)) - \pi(B(t)) \bigr)
+ (1-\rho(t)) \beta(t) ,
$$
where $\rho:(-1,0]\to[0,1]$ is a smooth cutoff function satisfying
$\rho\equiv 1$ near $0$ and $\rho\equiv 0$ near $-1$.
The claimed properties of $E$ now simply follow from the properties
of $\Theta$ and $\pi$ in Lemma~\ref{lem:lagg}.
\end{proof}

\begin{cor} \label{cor:expmap4}
Let $B_-,B_+\in\cA(Y,\cL)$ and $\Xi=A+\Phi\ds\in\cA(\R\times Y,\cL;B_-,B_+)$.
Fix $p>2$, then there is an open neighbourhood 
$\tilde{\cU}\subset \rT_\Xi\cA^{1,p}(\R\times Y,\cL;B_-,B_+)$ of zero such that
$$
\tilde E : \tilde{\cU} \to \cA^{1,p}(\R\times Y,\cL;B_-,B_+) ,\quad
\tilde E(\alpha + \phi\ds) := E_{A}(\alpha) + (\Phi+\phi)\ds
$$
defines a continuously differentiable homeomorphism onto a neighbourhood of~$\Xi$.
\end{cor}
\begin{proof}
Here we follow the construction of the exponential map of Theorem~\ref{thm:expmap} 
over the compact subset $N:=\{A(s)|s\in\R\}\cup\{B_-,B_+\}\subset\cA(Y,\cL)$. 
We fix the tubular neighbourhood $\tau:(-1,0]\times\Sigma\hookrightarrow Y$ of the boundary
such that $\tau^*A(s)=B(s,t)+\Psi(s,t)\dt$ with $B(s,t)\in\cW^{0,p}(\Sigma)$ 
for all $(s,t)\in \R\times(-1,0]$.
For $\alpha+\phi\ds\in \rT_\Xi\cA^{1,p}(\R\times Y,\cL;B_-,B_+)$ with
$\|\alpha+\phi\ds\|_{W^{1,p}(\R\times Y)}$ sufficiently small
the Sobolev embedding 
$W^{1,p}(\R\times(-1,0]\times\Sigma)\hookrightarrow \cC^0(\R\times(-1,0],L^p(\Sigma))$ 
ensures that $\tau^*\alpha = \beta(s,t) + \psi(s,t)\dt$ with $\beta(s,t)\in\cU^{0,p}$ 
for all $(s,t)\in\R\times(-1,0]$.

Thus we have $\tilde E(\alpha+\phi\ds)=A+\alpha + (\Phi+\phi)\ds$
on $\R\times (Y\setminus\im\tau)$ and
$\tau^*\tilde E(\alpha+\phi\ds) = \tB + (\Psi+\psi)\dt + (\Phi+\phi)\ds$
on $\R\times(-1,0]\times\Sigma$ with
$$
\tB(s,t) = B(s,t) + \rho(t)\bigl( \Theta_{\pi(B(s,t))}(\beta(s,t)) - \pi(B(s,t)) \bigr)
+ (1-\rho(t)) \beta(s,t)  .
$$
That $\tilde E$ is a bijection to a neighbourhood of $\Xi$ follows directly from Theorem~\ref{thm:expmap}.
For a restriction to a compact subset of $\R\times Y$ the smoothness of $\tilde E$ 
follows directly from the smoothness of the 3-dimensional exponential map.
To see that the 4-dimensional exponential map also is continuously differentiable with respect to
the $W^{1,p}(\R\times Y)$-norm on the noncompact domain, it suffices to drop linear terms
and the cutoff function $\rho$ and check that $\beta\mapsto\Theta_{\pi(B)}(\beta)-\pi(B)$
defines a $\cC^1$-map
$W^{1,p}(\R\times(-1,0],L^p(\Sigma))\to W^{1,p}(\R\times(-1,0],L^p(\Sigma))$
and also induces a $\cC^1$-map
$L^p(\R\times(-1,0],W^{1,p}(\Sigma))\to L^p(\R\times(-1,0],W^{1,p}(\Sigma))$.
This follows from the linear bounds for $\Theta$ and $\pi$ in Lemma~\ref{lem:lagg}, as follows.
For all $\beta,\beta'\in W^{1,p}(\R\times(-1,0]\times\Sigma,\rT^*\Sigma\otimes\cg)$ we have
\begin{align*}
\bigl\| \Theta_{\pi(B(s,t))}(\beta(s,t)) -  \Theta_{\pi(B(s,t))}(\beta'(s,t)) \bigr\|_{L^p(\Sigma)}  
&\leq C \| \beta(s,t) - \beta'(s,t) \|_{L^p(\Sigma)}, \\
\bigl\| \Theta_{\pi(B(s,t))}(\beta(s,t)) - \Theta_{\pi(B(s,t))}(\beta'(s,t))  \bigr\|_{W^{1,p}(\Sigma)}  
&\leq C \| \beta(s,t) - \beta'(s,t) \|_{W^{1,p}(\Sigma)} .
\end{align*}
For the $(s,t)$-derivatives we use the smoothness 
of $\Theta$ in the $L^p$-norm to obtain uniform
continuity for the derivative by $A$ in the 
$L^p$-operator norm, i.e.\
$
\|D_1\Theta(A,\alpha)-D_1\Theta(A,\alpha')\| 
\leq C\|\alpha-\alpha'\|_{L^p(\Sigma)}
$ 
for all sufficiently small 
$\alpha,\alpha'\in L^p(\Sigma,\rT^*\Sigma\otimes\cg)$. 
Since $\|\beta(s,t)\|_{L^p(\Sigma)}\to 0$ for $s\to\pm\infty$
this applies for all $t\in(-1,0]$ and $|s|$ sufficiently large, 
so that
\begin{align*}
&\bigl\| \partial_s \bigl( \Theta_{\pi(B(s,t))}(\beta(s,t)) 
- \Theta_{\pi(B(s,t))}(\beta'(s,t)) \bigr) \bigr\|_{L^p(\Sigma)}   \\
&\leq \bigl\| D\Theta_{\pi(B(s,t))}(\beta) (\partial_s\beta(s,t) - \partial_s\beta'(s,t)) \bigr\|_{L^p(\Sigma)} \\
&\quad
+\bigl\|  \bigl( D_1\Theta(\pi(B(s,t)),\beta(s,t)) - D_1\Theta(\pi(B(s,t)),\beta'(s,t)) \bigr) \partial_s \pi(B(s,t)) \bigr\|_{L^p(\Sigma)}  \\
&\leq C \bigl( \|\partial_s\beta(s,t) - \partial_s\beta'(s,t) \|_{L^p(\Sigma)} 
+\|\beta(s,t)-\beta'(s,t)\|_{L^p(\Sigma)} \| \partial_s B(s,t) \|_{L^p(\Sigma)} \bigr).
\end{align*}
(The same holds for $\partial_t(\ldots)$.) 
Integrating these estimates over $(s,t)\in(-1,0]\times\R$ 
proves $W^{1,p}$-continuity of $\tilde E(\alpha+\phi\ds)$.
To check continuity of the differential 
we use the analogous estimates for $D\Theta$, in particular
we use uniform continuity for the second derivatives 
of $\Theta$ (which again hold for 
$\|\beta(s,t)\|_{L^p(\Sigma)}$ sufficiently small, 
i.e.\ $|s|$ sufficiently large) to obtain
\begin{align*}
&\bigl\| \partial_s \bigl( D\Theta_{\pi(B(s,t))}(\beta(s,t)) 
- D\Theta_{\pi(B(s,t))}(\beta'(s,t)) \bigr) \gamma(s,t) \bigr\|_{L^p(\Sigma)}   \\
&\leq 
\bigl\| \bigl( D\Theta_{\pi(B)}(\beta) - D\Theta_{\pi(B)}(\beta') \bigr) \partial_s\gamma \bigr\|_{L^p(\Sigma)}  
+ \bigl\| D^2\Theta_{\pi(B)}(\beta) (\partial_s\beta - \partial_s\beta',\gamma) \bigr\|_{L^p(\Sigma)} \\
&\quad
+\bigl\|  \bigl( D_1 D_2\Theta(\pi(B),\beta) - D_1 D_2\Theta(\pi(B),\beta') \bigr) (\partial_s \pi(B),\gamma) \bigr\|_{L^p(\Sigma)}  \\
&\leq C \bigl( \|\beta - \beta' \|_{L^p(\Sigma)}  \| \partial_s \gamma \|_{L^p(\Sigma)} 
+ \|\partial_s\beta - \partial_s\beta' \|_{L^p(\Sigma)}  \| \gamma \|_{L^p(\Sigma)}  \\
&\qquad 
+  \|\beta-\beta' \|_{L^p(\Sigma)} \| \partial_s B \|_{L^p(\Sigma)} 
\| \gamma \|_{L^p(\Sigma)} \bigr) .
\end{align*}
Integration then proves the continuity of $D\tilde E$ 
in $W^{1,p}(\R\times Y)$.  (Strictly speaking, 
we can only integrate the above estimate 
over the complement of a compact interval in $\R$. 
However, the same estimate holds on the compact part due to
the smoothness of $\Theta$.)
\begin{align*}
&\bigl\| \partial_s \bigl( D\Theta_{\pi(B)}(\beta) 
- D\Theta_{\pi(B)}(\beta') \bigr) \gamma \bigr\|_{L^p(\R\times(-1,0]\times\Sigma)}   \\
&\leq C \|\beta - \beta' \|_{L^\infty(\R\times(-1,0]\times\Sigma)}  
\| \partial_s \gamma\|_{L^p(\R\times(-1,0]\times\Sigma)} \\
&\quad 
+C \|\partial_s\beta - \partial_s\beta' \|_{L^p(\R\times(-1,0]\times\Sigma)} 
\| \gamma \|_{L^\infty(\R\times(-1,0]\times\Sigma)} \\
&\quad
+ C \|\beta - \beta' \|_{L^\infty(\R\times(-1,0]\times\Sigma)}  
\| \partial_s B \|_{L^p(\R\times(-1,0]\times\Sigma)} \| \gamma \|_{L^\infty(\R\times(-1,0]\times\Sigma)} \\
&\leq C 
 \|\beta - \beta' \|_{W^{1,p}(\R\times(-1,0]\times\Sigma)}  
\bigl(1 + \| \partial_s B \|_{L^p(\R\times(-1,0]\times\Sigma)} \bigr) 
\| \gamma \|_{W^{1,p}(\R\times(-1,0]\times\Sigma)} .
\end{align*}
\end{proof}  
 
 
\bibliographystyle{alpha}

\end{document}